\documentclass[onecolumn,11pt]{article}
\usepackage[top=1in, bottom=1in, left=1in, right=1in]{geometry}
\usepackage{amsmath,amsfonts,amscd,amssymb,amsthm}
\usepackage{fancyhdr}
\usepackage{enumitem}
\usepackage[percent]{overpic}
\usepackage{tikz}
\usetikzlibrary{shapes.misc}
\usepackage{mathrsfs}
\usepackage{ifthen,tabularx,graphicx,multirow}
\usepackage{xargs}
\usepackage[colorinlistoftodos,prependcaption,textsize=tiny]{todonotes}
\usepackage{tcolorbox}
\usepackage{enumerate}
\usepackage{comment}
\usepackage{algorithm}
\usepackage{booktabs}
\usepackage{varwidth}
\usepackage{graphicx}
\usepackage{epstopdf}
\usepackage{mathtools}
\usepackage{setspace}
\usepackage{palatino} 
\usepackage{natbib}
\usepackage{bm}
\usepackage[bottom,flushmargin,hang,multiple]{footmisc}

\usepackage[pdfencoding=auto, psdextra]{hyperref}
\usepackage{algpseudocode}
\usepackage[capitalise]{cleveref}
\usepackage{lscape}
\usepackage{afterpage}
\usepackage{multirow}
\newcommand{\minitab}[2][l]{\begin{tabular}{#1}#2\end{tabular}}

\makeatletter
\renewcommand\paragraph{\@startsection{paragraph}{4}{\z@}%
            {-2.5ex\@plus -1ex \@minus -.25ex}%
            {1.25ex \@plus .25ex}%
            {\normalfont\normalsize\bfseries}}
\makeatother

\makeatletter
\def\hlinewd#1{%
\noalign{\ifnum0=`}\fi\hrule \@height #1 \futurelet
\reserved@a\@xhline}
\makeatother

\tikzset{cross/.style={cross out, draw=blue, minimum size=2*(#1-\pgflinewidth), inner sep=0pt, outer sep=0pt},
	cross/.default={1pt}}

\newcommand\blfootnote[1]{%
  \begingroup
  \renewcommand\thefootnote{}\footnote{#1}%
  \addtocounter{footnote}{-1}%
  \endgroup
}

\newcommand{\Av}{\mathbf{A}}
\newcommand{\Bv}{\mathbf{B}}
\newcommand{\Cv}{\mathbf{C}}
\newcommand{\Dv}{\mathbf{D}}
\newcommand{\Ev}{\mathbf{E}}
\newcommand{\Fv}{\mathbf{F}}
\newcommand{\Gv}{\mathbf{G}}

\newcommand{\Iv}{\mathbf{I}}

\newcommand{\Kv}{\mathbf{K}}
\newcommand{\Lv}{\mathbf{L}}

\newcommand{\Nv}{\mathbf{N}}

\newcommand{\Qv}{\mathbf{Q}}
\newcommand{\Rv}{\mathbf{R}}
\newcommand{\Sv}{\mathbf{S}}

\newcommand{\Uv}{\mathbf{U}}
\newcommand{\Vv}{\mathbf{V}}
\newcommand{\Wv}{\mathbf{W}}
\newcommand{\Xv}{\mathbf{X}}
\newcommand{\Yv}{\mathbf{Y}}
\newcommand{\Zv}{\mathbf{Z}}
\newcommand{\Psiv}{\mathbf{\Psi}}

\newcommand{\bv}{\mathbf{b}}
\newcommand{\cv}{\mathbf{c}}

\newcommand{\gv}{\mathbf{g}}
\newcommand{\hv}{\mathbf{h}}

\newcommand{\uv}{\mathbf{u}}
\newcommand{\vv}{\mathbf{v}}
\newcommand{\wv}{\mathbf{w}}
\newcommand{\xv}{\mathbf{x}}
\newcommand{\yv}{\mathbf{y}}
\newcommand{\zv}{\mathbf{z}}

\numberwithin{equation}{section}

\numberwithin{theorem}{section}
\numberwithin{lemma}{section}
\numberwithin{corollary}{section}
\numberwithin{proposition}{section}
\numberwithin{definition}{section}
\numberwithin{example}{section}
\makeatletter 
\newcounter{algorithmicH}
\let\oldalgorithmic\algorithmic
\renewcommand{\algorithmic}{
  \stepcounter{algorithmicH}
  \oldalgorithmic}
\renewcommand{\theHALG@line}{ALG@line.\thealgorithmicH.\arabic{ALG@line}}
\makeatother

\title{{\fontsize{16}{16}\selectfont \textbf{The Multiverse of Dynamic Mode Decomposition Algorithms}}\vspace{-.15in}}

\author{\normalsize{Matthew J. Colbrook$^{1*}$}\\
\footnotesize{$^1$DAMTP, University of Cambridge, Cambridge, CB3 0WA, United Kingdom}}
\date{}
\begin{document}
\maketitle

\blfootnote{$^*$ Corresponding author (m.colbrook@damtp.cam.ac.uk).}
\vspace{-.2in}
\begin{abstract}
Dynamic Mode Decomposition (DMD) is a popular data-driven analysis technique used to decompose complex, nonlinear systems into a set of modes, revealing underlying patterns and dynamics through spectral analysis. This review presents a comprehensive and pedagogical examination of DMD, emphasizing the role of Koopman operators in transforming complex nonlinear dynamics into a linear framework. A distinctive feature of this review is its focus on the relationship between DMD and the spectral properties of Koopman operators, with particular emphasis on the theory and practice of DMD algorithms for spectral computations. We explore the diverse ``multiverse'' of DMD methods, categorized into three main areas: linear regression-based methods, Galerkin approximations, and structure-preserving techniques. Each category is studied for its unique contributions and challenges, providing a detailed overview of significant algorithms and their applications as outlined in Table 1. We include a MATLAB package with examples and applications to enhance the practical understanding of these methods. This review serves as both a practical guide and a theoretical reference for various DMD methods, accessible to both experts and newcomers, and enabling readers to delve into their areas of interest in the expansive field of DMD.
\end{abstract}
\noindent\emph{Keywords --} dynamical systems, Koopman operator, data-driven discovery,
dynamic mode decomposition, spectral theory, spectral computations
{\small\selectfont{}\tableofcontents}

\section{Introduction}

Dynamical systems provide a powerful framework for modeling the evolution of various scientific and engineering systems over time. They are crucial for understanding complex phenomena ranging from weather patterns and population growth to stock market fluctuations. We consider discrete-time dynamical systems represented as:
\begin{equation}
\label{eq:DynamicalSystem} 
\xv_{n+1} = \Fv(\xv_n), \qquad n= 0,1,2,\ldots, 
\end{equation}
where $\xv\in\Omega$ denotes the state of the system, and $\Omega\subseteq\mathbb{R}^d$ is the statespace. The function $\Fv:\Omega \rightarrow \Omega$ governs the system's evolution. The classical approach to analyzing such systems, tracing back over a century to the seminal work of \cite{poincare1899methodes}, is geometric. It involves local analysis of fixed points, periodic orbits, and stable or unstable manifolds. While Poincaré's framework has significantly advanced our understanding of dynamical systems, it faces two main challenges in modern applications:
\begin{itemize}
	\item \textbf{Global understanding of nonlinear dynamics:} Unlike linear systems, there is no comprehensive mathematical framework for nonlinear systems. The principle of linear superposition is not applicable in this context. Local models can predict long-term dynamics near fixed points and attracting manifolds but have limited predictive power for other initial conditions. Consequently, the global understanding of nonlinear dynamics in state space is predominantly qualitative.
\item \textbf{Incomplete knowledge of evolution:} Many systems cannot be analytically described due to their complexity or our incomplete understanding. Typically, our knowledge is limited to discrete-time snapshots of the system, i.e., a finite dataset
$$
\left\{\xv^{(m)},\yv^{(m)}\right\}_{m=1}^M\quad \text{such that}\quad \yv^{(m)}=\Fv(\xv^{(m)}),\quad m=1,\ldots,M.
$$
We concisely write this data in the form of snapshot matrices
\begin{equation}
\label{snapshot_data}
\Xv=\begin{pmatrix} 
\xv^{(1)} & \xv^{(2)} & \cdots & \xv^{(M)}
\end{pmatrix}\in\mathbb{R}^{d\times M},\quad \Yv=\begin{pmatrix} 
\yv^{(1)} & \yv^{(2)} & \cdots & \yv^{(M)}
\end{pmatrix}\in\mathbb{R}^{d\times M}.
\end{equation}
Advances in measurement technologies have significantly enhanced our ability to collect detailed multimodal and multi-fidelity snapshot data. Data could be collected from one long trajectory or multiple shorter trajectories. It can come from experimental observations or numerical simulations. The question becomes how to use this data to meaningfully study the dynamical system in \eqref{eq:DynamicalSystem}.
\end{itemize}
The advent of big data \citep{hey2009fourth}, coupled with strides in modern statistical learning \citep{friedman2017elements} and machine learning \citep{mohri2018foundations}, has heralded a new era of data-driven algorithms to address these issues. This review will focus on one of the most prominent of these algorithms, \textit{Dynamic Mode Decomposition} (DMD), closely connected with \textit{Koopman operators}.

\textit{Koopman Operators} --- 
In 1931, Koopman introduced his operator-theoretic approach to dynamical systems, initially to describe Hamiltonian systems \citep{koopman1931hamiltonian}. This theory was further expanded by \cite{koopman1932dynamical} to include systems with continuous spectra. Koopman operators offer a powerful alternative to the classical geometric view of dynamical systems by addressing the fundamental issue of \textit{nonlinearity}. We lift a nonlinear system \eqref{eq:DynamicalSystem} into an infinite-dimensional space of observable functions $g:\Omega\rightarrow\mathbb{C}$ using a Koopman operator $\mathcal{K}$:
$$
[\mathcal{K}g](\xv) = g(\Fv(\xv)), \quad \text{so that}\quad [\mathcal{K}g](\xv_n)=g(\xv_{n+1}).
$$
Through this approach, the evolution dynamics become linear, enabling the use of generic solution techniques based on spectral decompositions. Initially, the primary application of Koopman operators was in ergodic theory \citep{eisner2015operator}, notably playing a pivotal role in proving the ergodic theorem by von Neumann \citep{neumann1932proof} and Birkhoff \citep{birkhoff1931proof, birkhoff1932recent}. More recently, they have been extensively used in data-driven methods for studying dynamical systems.

\textit{Dynamic Mode Decomposition} --- A significant objective of modern Koopman operator theory is to identify a coordinate transformation under which even strongly nonlinear dynamics may be approximated by a linear system. This coordinate system is related to the spectrum of the Koopman operator. DMD was initially developed by \cite{schmid2009dynamic,schmid2010dynamic} (see also \citep{schmid2008dynamic}) in the context of fluid dynamics. \cite{mezic2005spectral} introduced the Koopman mode decomposition, providing a theoretical basis for \cite{rowley2009spectral} to connect DMD with Koopman operators. This connection validated DMD's application in nonlinear systems and offered a powerful yet straightforward, data-driven approach for approximating Koopman operators. The fusion of contemporary Koopman theory with an efficient numerical algorithm has led to significant advancements and a surge in research. DMD is now the central algorithm for computational approximations of Koopman operators with applications in various fields beyond fluid mechanics, such as neuroscience, disease modeling, robotics, video processing, power grids, financial markets, and plasma physics. The simplicity and effectiveness of DMD have led to numerous innovations, giving rise to a diverse array of DMD methods, playfully described here as a ``multiverse'', aimed at addressing specific challenges.

\textit{This Review} ---  We provide a comprehensive tour of this ``multiverse'' of DMD methods. Our primary focus is on the interplay between DMD, the spectral properties of Koopman operators, and their numerical computations. At the time of writing, these methods can be broadly categorized into three main areas:
\begin{itemize}
	\item DMD methods based on linear regression;
	\item DMD methods utilizing Galerkin approximations;
	\item DMD methods aimed at preserving structures or symmetries of \eqref{eq:DynamicalSystem}.
\end{itemize}
These distinctions are not rigid, and some methods encompass multiple flavors. This review navigates these key areas and variants, summarized in \cref{tab:summary}, where we also highlight the unique challenges each algorithm addresses (see also \cref{sec:challenges_of_DMD}). We provide detailed summaries and examples of these algorithms in action. Accompanying this review is a MATLAB package:
\begin{center}
\textcolor[rgb]{0,0,1}{https://github.com/MColbrook/DMD-Multiverse}
\end{center}
featuring user-friendly implementations and examples from the paper, most of which are new. We aim for readers to utilize this paper as a practical manual for various DMD methods. Although extensive, the review is structured modularly, enabling readers to selectively engage with DMD versions and topics that interest them most.

\afterpage{%
    \clearpage
    \thispagestyle{empty}
\begin{landscape}
\begin{table}
\renewcommand{\arraystretch}{1.2}
\centering\scriptsize
\begin{tabular}{|c|l|l|c|}
\hline
\textbf{DMD Method} & \multicolumn{1}{c|}{\textbf{Challenges Overcome}} & \multicolumn{1}{c|}{\textbf{Key Insight/Development}} & \textbf{Key Reference(s)} \\
\hlinewd{2pt}
\multirow{2}{*}{Forward-Backward DMD} & \multirow{2}{*}{Sensor noise bias.} & Take geometric mean of forward and& \multirow{2}{*}{\cite{dawson2016characterizing}}\\
&& backward propagators for the data. & \\
\hline
\multirow{2}{*}{Total Least-Squares DMD} & \multirow{2}{*}{Sensor noise bias.} & Replace least-squares problem with & \cite{hemati2017biasing}\\
& & total least-squares problem. & \cite{dawson2016characterizing}
\\
\hline
\multirow{3}{*}{\minitab[c]{Optimized DMD\\Bagging Optimized DMD}} & \multirow{3}{*}{\hspace{-2mm}\minitab[l]{Sensor noise bias.\\Optimal collective processing of snapshots.}} & Exponential fitting problem,& \cite{chen2012variants}\\
&  &solve using variable projection method.  & \cite{askham2018variable}\\
& & Statistical bagging sampling strategy. & \cite{sashidhar2022bagging}\\
\hline
\multirow{3}{*}{Compressed Sensing}&\multirow{3}{*}{\hspace{-2mm}\minitab[l]{Computational efficiency.\\Temporal or spatial undersampling.}} & Unitary invariance of DMD extended to& \cite{tu2014spectral}\\
& &settings of compressed sensing & \cite{brunton2016compressed}\\
&&(e.g., RIP, sparsity-promoting regularizers).& \cite{erichson2019compressed}\\
\hline
\multirow{2}{*}{Randomized DMD} & Computational efficiency. &  Sketch data matrix for computations & \multirow{2}{*}{\cite{erichson2019randomized}}\\
& Memory usage.  &  in reduced-dimensional space. &\\
\hline
Multiresolution DMD & Multiscale dynamics. & Filtered decomposition across scales. & \cite{kutz2016multiresolution}\\
\hline
\multirow{2}{*}{DMD with Control} & \multirow{2}{*}{Separation of unforced dynamics and actuation.} & Generalized regression for globally &\multirow{2}{*}{\cite{proctor2016dynamic}}\\
&& linear control framework. & \\
\hlinewd{2pt}
\multirow{2}{*}{Extended DMD} & \multirow{2}{*}{Nonlinear observables.} & Arbitrary (nonlinear) dictionaries, &\multirow{2}{*}{\cite{williams2015data}} \\
& & recasting of DMD as a Galerkin method. & \\
\hline
\multirow{2}{*}{Hankel DMD} & Delay-embedding for ergodic systems. & Connection with Krylov subspace methods&\multirow{2}{*}{\cite{arbabi2017ergodic}}\\
& Convergence under invariant subspace assumption. & and Birkhoff's ergodic theorem.& \\
\hline
HAVOK & Lack of closed linear models for chaotic systems. & Delay-embedding with chaos as forcing. &  \cite{brunton2017chaos}\\
\hline
\multirow{3}{*}{Residual DMD} & 
Inf.-dim. projection errors, verification (general systems). & \multirow{3}{*}{\hspace{-2mm}\minitab[l]{Append EDMD with additional matrix\\(available from the snapshot data)\\ to compute infinite-dimensional residuals.}}& \multirow{3}{*}{\hspace{-2mm}\minitab[c]{\cite{colbrook2021rigorous}\\\cite{colbrook2023residualJFM}}}\\
&Computation of Koopman spectra (general systems). &&\\
&Spectral measures (measure-preserving systems).&&\\
\hlinewd{2pt}
\multirow{2}{*}{Physics-Informed DMD} & Preserving structure of dynamical systems. & Restrict the least-squares optimization & \multirow{2}{*}{\cite{baddoo2023physics}}\\
& Numerous instances given in general framework. & to lie on a matrix manifold. &\\
\hline
\multirow{3}{*}{\minitab[c]{Measure-Preserving\\Extended DMD}} & Measure-preserving discretizations of system. & Alter EDMD to be measure-preserving& \multirow{3}{*}{\cite{colbrook2023mpedmd}}\\
 & Convergence to Koopman spectral properties & with respect to a learned inner product. & \\
& (including continuous spectra/spectral measures). & (via a polar decomposition of EDMD) & \\
\hline
\multirow{4}{*}{\minitab[c]{Compactification\\Methods}} & Cts.-time generator of measure-preserving system.
& \multirow{4}{*}{\hspace{-2mm}\minitab[l]{Compactification of generator or its\\ resolvent using  kernel integral operators,\\dictionary of kernel eigenvectors.}}& 
\multirow{4}{*}{\hspace{-2mm}\minitab[c]{\cite{das2021reproducing}\\
\cite{valva2023consistent}}}\\
 & Convergence to Koopman spectral properties  & & \\
& (including continuous spectra/spectral measures).&  &\\
&Conditioning of dictionary. & &\\
\hline
\end{tabular}
\caption{Executive summary of the DMD methods discussed in detail in this review. Numerous others are also discussed. The bold horizontal lines separate the different flavors of regression (top), Galerkin (middle), and structure-preserving (bottom). The fundamental DMD algorithm, exact DMD, is given in \cref{alg:DMD_vanilla}.\\\textbf{NB:} For measure-preserving systems, discretizations that preserve the measure are crucial for convergence, recovering the correct dynamical behavior, stability, robustness to noise, and improved qualitative and long-time behavior.}
\label{tab:summary}
\end{table}
\normalsize
\end{landscape}
\clearpage
}

Differentiating itself from prior reviews, this review specifically focuses on the theory and practice of DMD algorithms for computing spectral properties of Koopman operators, complementing other reviews on the subject. \cite{mezic2013analysis} and the more recent review of \cite{schmid2022dynamic} (see also \citep{taira2017modal,taira2020modal}) focus on developments associated with applications in fluid dynamics. While the initial applications of Koopman and DMD techniques were in fluid problems, their utility has been demonstrated in a broader range of fields. We also briefly explore the applications of DMD in control theory, and readers seeking further exploration in this area are encouraged to read \cite{otto2021koopman}. An excellent early review of ``Applied Koopmanism'' is presented by \cite{budivsic2012applied}. More recently, \cite{brunton2021modern} have provided a broad overview of Koopman operators and their applications, with connections to other fields.

This review is organized into several sections. In \cref{sec:basic_DMD}, we introduce the basic DMD algorithm, offering a concise introduction to Koopman operators and their spectral properties, before delving into the fundamental DMD algorithm and its two key interpretations: regression and projection. We discuss three canonical examples, followed by an examination of the goals and challenges of DMD. \cref{sec:regression_variants} focuses on variants from the regression viewpoint, including noise reduction, compression, randomized linear algebra, multiscale dynamics, and control. The connection with Koopman operators is further explored in \cref{sec:Galerkin_perspective}, where we discuss nonlinear observables, time-delay embedding methods, and methods for controlling the infinite-dimensional projection error of DMD. In \cref{sec:structure_preserving_methods}, we review recent methods that preserve the structures of dynamical systems. These methods often lead to greater noise resistance, improved generalization, and reduced data demands for training. We conclude in \cref{sec:conclusion} by discussing further connections and open problems. I hope the reader enjoys this tour of the DMD multiverse as much as I have enjoyed writing it!

Due to the sheer breadth and thousands of papers written on DMD, it is impossible for this review to cover every version of DMD in great detail. Significant DMD algorithms that are not discussed in their own sections are still discussed in some detail. If the reader searches this paper, they will find dozens of DMD algorithms. I have included all significant references I am aware of, but many others may not have been included. I apologize for that in advance and encourage all readers to inform the author about results that deserve more discussion.

\section{The Basics of DMD}
\label{sec:basic_DMD}

To understand the DMD ``multiverse'', we must first study the basic DMD algorithm. We begin with Koopman operator theory, the theoretical underpinning of DMD, before moving on to the fundamental DMD algorithm and two important viewpoints. We then provide three canonical examples and discuss the goals and challenges of DMD.

\subsection{The underlying theory: Koopman operators and spectra}
\label{sec:Koopman_operators}

In this section, we recall the definition of Koopman operators and equip the reader with a crash course on their relevant spectral properties. At its core, DMD is an algorithm that uses the snapshot data in \eqref{snapshot_data} to approximate the spectral properties of Koopman operators.

\begin{figure}
\centering
\includegraphics[width=0.9\textwidth,trim={0mm 0mm 0mm 0mm},clip]{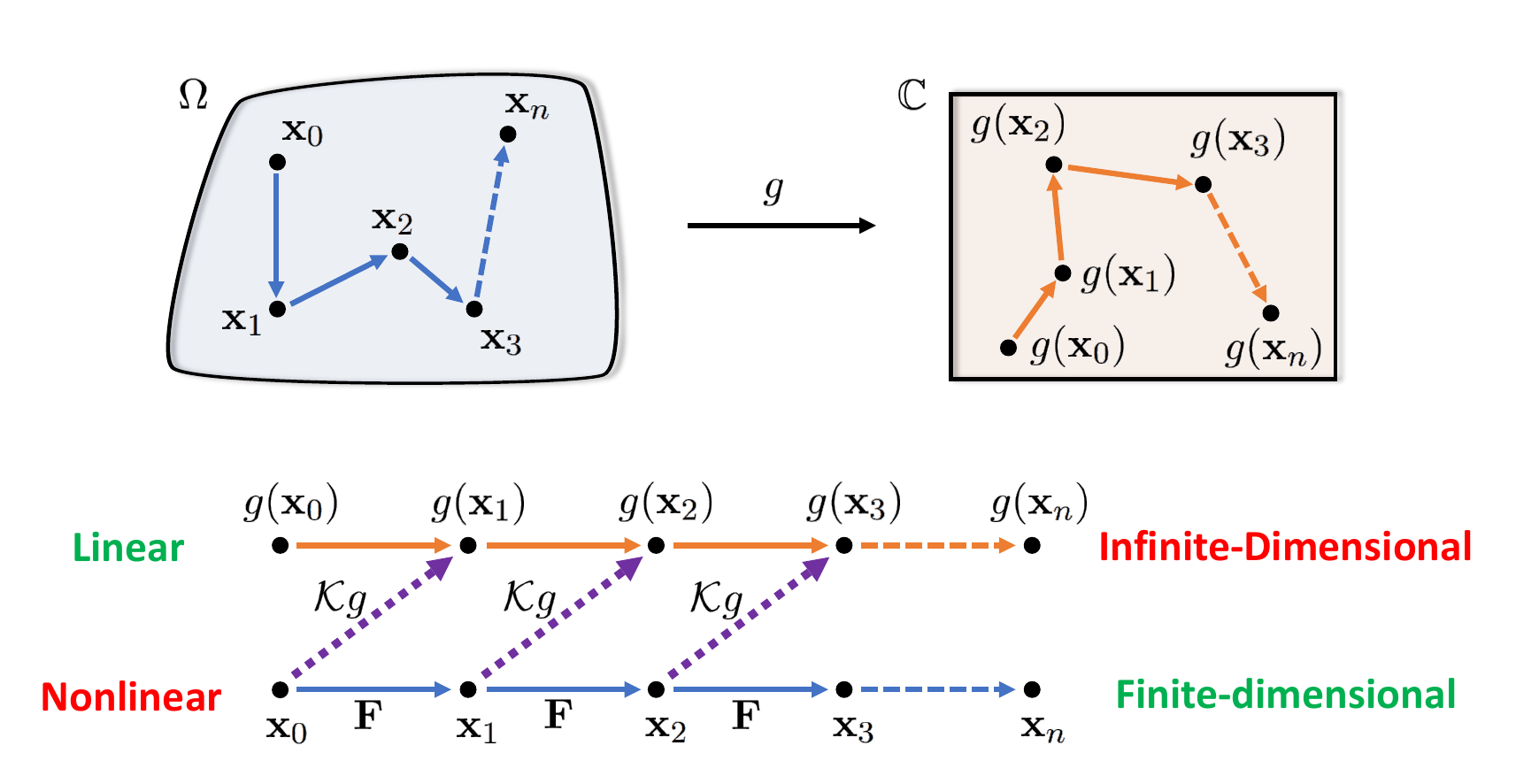}
\caption{Summary of the idea of Koopman operators. By lifting to a space of observables, we trade a nonlinear finite-dimensional system for a linear infinite-dimensional system.}
\label{schematic}
\end{figure}

\subsubsection{What is a Koopman operator?}
To define a Koopman operator, we begin with a space $\mathcal{F}$ of functions $g: \Omega \rightarrow \mathbb{C}$, where $\Omega$ is the state space of our dynamical system. The functions $g$, referred to as \textit{observables}, serve as tools for indirectly measuring the state of the system described in \eqref{eq:DynamicalSystem}. Specifically, $g(\xv_n)$ indirectly measures the state $\xv_n$. Koopman operators enable us to capture the time evolution of these observables through a linear operator framework. For a suitable domain $\mathcal{D}(\mathcal{K}) \subset \mathcal{F}$, we define the Koopman operator via the composition formula:
\begin{equation} 
[\mathcal{K}g](\xv) = [g\circ \Fv](\xv)=g(\Fv(\xv)), \qquad g\in \mathcal{D}(\mathcal{K}).
\label{eq:KoopmanOperator} 
\end{equation}
In this context, $[\mathcal{K}g](\xv_n)= g(\Fv(\xv_n))=g(\xv_{n+1})$ represents the measurement of the state one time step ahead of $g(\mathbf{x}_n)$. This process effectively captures the dynamic progression of the system. The overarching concept is summarized in \cref{schematic}.

The key property of the Koopman operator $\mathcal{K}$ is its \textit{linearity}. This linearity holds irrespective of whether the system's dynamics, as represented in \eqref{eq:DynamicalSystem}, are linear or nonlinear. Consequently, the spectral properties of $\mathcal{K}$ become a powerful tool in analyzing the dynamical system's behavior. To study spectra, we assume that $\mathcal{F}$ is a Banach space.\footnote{A Banach space is a normed vector space that is complete, i.e., every Cauchy sequence converges. Thus, a Banach space has no `holes' in it. We have deliberately kept the background functional analysis to a minimum in this review.} For the spectrum of $\mathcal{K}$ to be meaningful and nontrivial, we assume that its domain, $\mathcal{D}(\mathcal{K})$, is dense in $\mathcal{F}$ and that $\mathcal{K}$ itself is a closed operator.\footnote{An operator being 'closed' means that its graph $
\left\{(g,\mathcal{K}g):g\in\mathcal{D}(\mathcal{K})\right\}$ is a closed subset within the product space $\mathcal{F}\times\mathcal{F}$.} If these conditions are not met, the spectrum would encompass the entirety of $\mathbb{C}$. It is crucial to recognize that the Koopman operator is not uniquely defined by the dynamical system in \eqref{eq:DynamicalSystem}; rather, it is fundamentally dependent on the choice of the space of observables $\mathcal{F}$. In this review, we focus on cases where $\mathcal{F}$ is defined as the following Hilbert space:
$$
\mathcal{F}=L^2(\Omega,\omega)\quad \text{with inner product}\quad \langle g_1,g_2 \rangle=\int_{\Omega} g_1(\xv)\overline{g_2(\xv)}\ \mathrm{d}\omega(\xv)\quad \text{and norm}\quad \|g\|=\sqrt{\langle g,g \rangle},
$$
for some positive measure $\omega$.\footnote{We do not assume that this measure is invariant. For Hamiltonian systems, a common choice of $\omega$ is the standard Lebesgue measure, for which the Koopman operator is unitary on $L^2(\Omega,\omega)$. For other systems, we can select $\omega$ according to the region where we wish to study the dynamics, such as a Gaussian measure. In many applications, $\omega$ corresponds to an unknown ergodic measure on an attractor.} In going from a pointwise definition in \eqref{eq:KoopmanOperator} to the space $L^2(\Omega,\omega)$, a little care is needed since $L^2(\Omega,\omega)$ consists of equivalence classes of functions. We assume that the map $\Fv$ is nonsingular with respect to $\omega$, meaning that
$$
\omega(E)=0\quad \text{implies that}\quad \omega(\{\xv:\Fv(\xv)\in E\})=0.
$$
This ensures that the Koopman operator is well-defined since $g_1(\xv)=g_2(\xv)$ for $\omega$-almost every $\xv$ implies that $g_1(\Fv(\xv))=g_2(\Fv(\xv))$ for $\omega$-almost every $\xv$. The above Hilbert space setting is standard in most of the Koopman literature for two reasons. First, it is a reasonable assumption for many dynamical systems, particularly if we study the dynamics on an attractor. Second, working with operators in a Hilbert space is much easier computationally than in a Banach space. For Koopman operators on more general spaces, see \citep{mezic2020spectrum}. Practical algorithms for Koopman operators on more general Banach spaces remain a largely open problem (see \cref{sec:open_problems}). 

Since $\mathcal{K}$ acts on an \textit{infinite-dimensional} function space, we have exchanged the nonlinearity in \eqref{eq:DynamicalSystem} for an infinite-dimensional linear system. This means that the spectral properties of $\mathcal{K}$ can be significantly more complex than those of a finite matrix, making them more challenging to compute. While this might seem disheartening, as we will explore in \cref{sec:Galerkin_perspective}, methodologies exist that enable the analysis of infinite-dimensional spectral properties through a series of finite-dimensional approximations.

\subsubsection{Crash course on spectral properties of Koopman operators}
\label{sec:spectra_crash_course}

We will now review the relevant spectral properties of $\mathcal{K}$. Readers primarily interested in applying DMD algorithms will still find the dynamical interpretations of these properties insightful. The sole assumption made throughout this paper is that $\mathcal{K}$ is a closed and densely defined operator. Specifically, unless stated otherwise, we do not presuppose that $\mathcal{K}$ possesses a nontrivial\footnote{If the measure $\omega$ is finite, then the constant function $g(\xv)=1$ is a trivial eigenfunction with eigenvalue $1$. It is deemed trivial because the dynamics of a constant observable lack informative content.} finite-dimensional invariant subspace, nor do we assume it has an eigenvector basis. These two assumptions are often implicitly (and sometimes wrongly) assumed in DMD papers and can lead to confusion if care is not taken.

\paragraph{Koopman spectra}

If $g\in L^2(\Omega,\omega)$ is an \textit{eigenfunction} of $\mathcal{K}$ with \textit{eigenvalue} $\lambda$, then $g$ exhibits perfect coherence\footnote{In the setting of dynamical systems, coherent sets or structures are subsets of the phase space where elements (e.g., particles, agents, etc.) exhibit similar behavior over some time interval. This behavior remains relatively consistent despite potential perturbations or the chaotic nature of the system.} with
\begin{equation}
\label{eq:perfectly_coherent}
g(\xv_n)=[\mathcal{K}^ng](\xv_0)=\lambda^n g(\xv_0)\quad \forall n\in\mathbb{N}.
\end{equation}
The oscillation and decay/growth of the observable $g$ are dictated by the complex argument and absolute value of the eigenvalue $\lambda$, respectively. In infinite dimensions, the appropriate generalization of the set of eigenvalues of $\mathcal{K}$ is the \textit{spectrum}, denoted by $\mathrm{Sp}(\mathcal{K})$, and defined as
$$
\mathrm{Sp}(\mathcal{K})=\left\{z\in\mathbb{C} :(\mathcal{K} - zI)^{-1}\text{ does not exist as a bounded operator}\right\}\subset\mathbb{C}.
$$
Here, $I$ denotes the identity operator. The spectrum $\mathrm{Sp}(\mathcal{K})$ includes the set of eigenvalues of $\mathcal{K}$, but in general, $\mathrm{Sp}(\mathcal{K})$ contains points that are not eigenvalues. This is because there are more ways for $(\mathcal{K} - zI)^{-1}$ to not exist in infinite dimensions than in finite dimensions. For example, we may have continuous spectra. The standard Lorenz system on the Lorenz attractor gives rise to a Koopman operator that has no nontrivial eigenvalues, yet the spectrum is the whole unit circle!

In general, we cannot numerically approximate an eigenfunction perfectly. Moreover, the operator 
$\mathcal{K}$ may not have any nontrivial eigenfunctions, for instance, if the system is weakly mixing. Instead, the so-called \textit{approximate point spectrum} is the following subset of $\mathrm{Sp}(\mathcal{K})$:
$$
\mathrm{Sp}_{\mathrm{ap}}(\mathcal{K})=\left\{\lambda\in\mathbb{C}:\exists\{g_n\}_{n\in\mathbb{N}}\subset L^2(\Omega,\omega)\text{ such that }\|g_n\|=1,\lim_{n\rightarrow\infty}\|(\mathcal{K}-\lambda I)g_n\|=0\right\}\subset\mathbb{C}.
$$
An observable $g$ with $\|g\|=1$ and $\|(\mathcal{K}-\lambda I)g\|\leq\epsilon$ for $\lambda\in\mathbb{C}$ is known as $\epsilon$-pseudoeigenfunction. Such observables are important for the dynamical system \eqref{eq:DynamicalSystem} since
$$
\|\mathcal{K}^ng-\lambda^n g\|\lesssim \mathcal{O}(n\epsilon)\quad \forall n\in\mathbb{N}.
$$
In other words, $\lambda$ describes an approximate coherent oscillation and decay/growth of the observable $g$ with time. The pseudoeigenfunctions and $\mathrm{Sp}_{\mathrm{ap}}(\mathcal{K})$ encode information about the underlying dynamical system \citep{mezicAMS}. For example, the level sets of certain eigenfunctions determine the invariant manifolds \citep{mezic2015applications} and isostables \citep{mauroy2013isostables}, and the global stability of equilibria \citep{mauroy2016global} and ergodic partitions \citep{budivsic2012applied,mezic1999method} can be characterized by pseudoeigenfunctions and $\mathrm{Sp}_{\mathrm{ap}}(\mathcal{K})$.

\paragraph{Koopman pseudospectra}

Approximate point spectra and pseudoeigenfunctions are related to the notion of \textit{pseudospectra} \citep{trefethen2005spectra}. For a finite matrix $A\in\mathbb{C}^{n\times n}$ and $\epsilon>0$, the $\epsilon$-pseudospectrum of $A$ is the set\footnote{Some authors use a strict inequality in the definition of $\epsilon$-pseudospectra. We prefer the given definition since then the pseudospectrum is a \textit{closed} subset of $\mathbb{C}$.}
$$
\mathrm{Sp}_\epsilon(A) = \left\{\lambda\in\mathbb{C}: \|(A-\lambda I)^{-1}\|\geq1/\epsilon \right\}= \bigcup\limits_{B\in\mathbb{C}^{n\times n},\|B\|\leq\epsilon} \mathrm{Sp}(A + B).
$$
The $\epsilon$-pseudospectra of $A$ are regions in the complex plane enclosing the eigenvalues of $A$. These regions tell us how far an $\epsilon$-sized perturbation can perturb an eigenvalue. Pseudospectra of Koopman operators must be defined with some care because $\mathcal{K}$ may be an unbounded operator and hence the resolvent norm $\|(\mathcal{K}-\lambda I)^{-1}\|$ can be constant on open subsets of $\mathbb{C}\backslash\mathrm{Sp}(\mathcal{K})$ \citep{shargorodsky2008level}. We define the $\epsilon$-pseudospectrum of $\mathcal{K}$ as \cite[Prop. 4.15]{roch1996c}:\footnote{While \citep[Prop. 4.15]{roch1996c} considers bounded operators, it can be adjusted to cover unbounded operators \citep[Thm. 4.3]{trefethen2005spectra}.}
\begin{equation}
\label{eq_def_pseudospectra}
\mathrm{Sp}_\epsilon(\mathcal{K})=\mathrm{Cl}\left(\{\lambda\in\mathbb{C}:\|(\mathcal{K}-\lambda I)^{-1}\| > 1/\epsilon\}\right)=\mathrm{Cl}\left(\bigcup_{\|\mathcal{B}\|< \epsilon}\mathrm{Sp}(\mathcal{K}+\mathcal{B})\right),
\end{equation}
where $\mathrm{Cl}$ denotes the closure of a set. To see the connection with $\mathrm{Sp}_{\mathrm{ap}}(\mathcal{K})$, note that if $\|(\mathcal{K}-\lambda I)g\|\leq\epsilon$ for an observable $g$ with $\|g\|=1$, then $\|(\mathcal{K}-\lambda I)^{-1}\|\geq 1/\epsilon$. We care about pseudospectra for several reasons, but two stand out as the most important:
\begin{itemize}
	\item Pseudospectra allow us to determine which regions of computed spectra are accurate and trustworthy. This could be in terms of the numerical stability, but also pseudospectra aid in detecting so-called \textit{spectral pollution} (see \cref{duffing2,duffing4}). These are spurious eigenvalues arising from discretization that are unrelated to the underlying Koopman operator. The term spectral pollution refers to the accumulation of these spurious eigenvalues at points outside the spectrum of $\mathcal{K}$ as the discretization size increases \citep{lewin2010spectral}. This occurs even when $\mathcal{K}$ is a normal operator (see \cref{duffing3}). It is essential to realize that spectral pollution leads to spurious modes that are not linked to stability issues but are a consequence of discretizing the infinite-dimensional operator $\mathcal{K}$ to a finite matrix.
	\item If the Koopman operator is nonnormal, the system's transient behavior can differ significantly from the asymptotic behavior captured by $\mathrm{Sp}(\mathcal{K})$. Pseudospectra can be employed to detect and quantify transients not represented by the spectrum \citep{trefethen1993hydrodynamic}\citep[Section IV]{trefethen2005spectra}. 
\end{itemize}

Pseudospectra also provide a means of computing spectra since $\lim_{\epsilon\downarrow 0}\mathrm{Sp}_\epsilon(\mathcal{K})=\mathrm{Sp}(\mathcal{K})$. This convergence occurs in the so-called Attouch--Wets metric space \citep{beer1993topologies}, which roughly means that we obtain uniform convergence on any compact region of $\mathbb{C}$. This observation goes beyond Koopman operators and has been behind some recent breakthroughs in the computation of spectra in infinite dimensions \citep{SCI_ref,colbrookthesis,colbrook2022computation,colbrook2022foundations,colbrook2019compute}.

\paragraph{Koopman mode decompositions and spectral theorems beyond eigenvalues}

One of the most useful features of Koopman operators is the \textit{Koopman Mode Decomposition} (KMD) \citep{mezic2005spectral}. The KMD expresses the state $\xv$ or an observable $g(\xv)$ as a linear combination of dominant coherent structures. It can be considered a diagonalization of the Koopman operator. As a result, the KMD is invaluable for tasks such as dimensionality and model reduction. It generalizes the space-time separation of variables typically achieved through the Fourier transform or singular value decomposition (SVD). It is crucial to realize that an exact KMD is rigorously justified only if $\mathcal{K}$ possesses some form of spectral theorem, which extends the concept of diagonalization to infinite dimensions. Nevertheless, obtaining an \textit{approximate} KMD is still possible even without a spectral theorem.

For example, suppose that the system \eqref{eq:DynamicalSystem} is \textit{measure-preserving} with respect to the positive measure $\omega$. This means that $\omega(E)=\omega(\{\xv:\Fv(\xv)\in E\})$ for any measurable set $E\subset\Omega$. In other words, the dynamical system preserves a volume. Measure-preserving systems encompass many systems of interest such as Hamiltonian flows \citep{arnold1989mathematical}, geodesic flows \citep{dubrovin2012modern}, Bernoulli schemes \citep{shields1973theory}, physical systems in equilibrium \citep{hill1986introduction}, and ergodic systems \citep{walters2000introduction}. Furthermore, many dynamical systems either admit invariant measures \citep{kryloff1937theorie} or exhibit measure-preserving post-transient behavior \citep{mezic2005spectral}. In fact, if $\Omega$ is a compact metric space and $\Fv$ is continuous, then there is an invariant measure \cite[Prop. 8.1]{mane1987ergodic}.\footnote{Of course, whether or not this is useful or whether our chosen $\omega$ is invariant is another matter.} For a measure-preserving system, the Koopman operator $\mathcal{K}$ is an isometry, i.e., $\|\mathcal{K}g\|=\|g\|$ for all observables $g\in \mathcal{D}(\mathcal{K})=L^2(\Omega,\omega)$. For simplicity, we further assume that $\mathcal{K}$ is unitary, implying that it is normal (it commutes with its adjoint).\footnote{A Koopman operator that is an isometry need not be unitary, e.g., the Koopman operator associated with the tent map. However, an isometry can always be extended to a unitary operator, and the spectral measures associated with forward-time dynamics are independent of the chosen extension \citep{colbrook2023mpedmd}.}

Under these conditions, the spectral theorem \citep[Thm. X.4.11]{conway2019course} allows us to diagonalize the Koopman operator $\mathcal{K}$. There is a \textit{projection-valued measure} $\mathcal{E}$ supported on $\mathrm{Sp}(\mathcal{K})$. For readers unfamiliar with the spectral theorem, \cite{halmos1963does} provides an excellent and readable introduction. In our example, $\mathcal{K}$ is unitary, which implies that $\mathrm{Sp}(\mathcal{K})$ lies within the unit circle $\mathbb{T}$. The measure $\mathcal{E}$ associates an orthogonal projector with each Borel measurable subset of $\mathbb{T}$. For such a subset $S\subset\mathbb{T}$, $\mathcal{E}(S)$ is a projection onto the spectral elements of $\mathcal{K}$ inside $S$. For any observable $g\in L^2(\Omega,\omega)$,
$$
g=\left(\int_\mathbb{T} \ \mathrm{d}\mathcal{E}(\lambda)\right)g \qquad\text{and}\qquad \mathcal{K}g=\left(\int_\mathbb{T} \lambda\ \mathrm{d}\mathcal{E}(\lambda)\right)g.
$$
The essence of this formula is the decomposition of $g$ according to the spectral content of $\mathcal{K}$. The projection-valued measure $\mathcal{E}$ simultaneously decomposes the space $L^2(\Omega,\omega)$ and diagonalizes the Koopman operator. For example, we have
\begin{equation}
\label{KMD_spec_meas_forward}
g(\xv_n)= [\mathcal{K}^ng](\xv_0)= \left[\left(\int_\mathbb{T} \lambda^n\ \mathrm{d}\mathcal{E}(\lambda)\right)g\right](\xv_0).
\end{equation}
This directly extends \eqref{eq:perfectly_coherent}. The spectral theorem can be perceived as offering a custom Fourier-type transform specifically for the operator $\mathcal{K}$ that extracts coherent features. Of particular interest are \textit{scalar-valued} spectral measures. Given a normalized observable $g\in L^2(\Omega,\omega)$ with $\|g\| = 1$, the scalar-valued spectral measure of $\mathcal{K}$ with respect to $g$ is a probability measure defined as 
$$
\mu_g(S)=\langle \mathcal{E}(S)g,g \rangle.
$$
These measures can be further refined using Lebesgue's decomposition into a pure point part, supported on the eigenvalues of $\mathcal{K}$, and a continuous part. The continuous part can further be decomposed into an absolutely continuous part with a density function and a singular continuous part. The moments of the measure $\mu_g$ are the correlations
$$
\langle \mathcal{K}^n g,g\rangle=\int_\mathbb{T} \lambda^n \ \mathrm{d}\mu_g(\lambda),\quad n\in\mathbb{Z}.
$$
For example, if our system corresponds to the dynamics on an attractor, these statistical properties allow comparison of complex dynamics \citep{mezic2004comparison}. More generally, the spectral measure of $\mathcal{K}$ with respect to $g\in L^2(\Omega,\omega)$ is a signature for the forward-time dynamics of \eqref{eq:DynamicalSystem}.

Going one step further, $\mathcal{E}$ leads to a decomposition of $L^2(\Omega,\omega)$ into parts associated with quasiperiodic evolution and weak-mixing dynamics. Namely, we have the following orthogonal decomposition into two $\mathcal{K}$-invariant subspaces \citep{halmos2017lectures}
$$
L^2(\Omega,\omega)=\mathcal{H}_{\mathrm{pp}}\oplus \mathcal{H}_{\mathrm{c}}.
$$
Here, the subspace $\mathcal{H}_{\mathrm{pp}}$ consists of the closure of the linear span of eigenvectors and admits an orthonormal basis of eigenvectors $\{\phi_j\}$ of $\mathcal{K}$ with eigenvalues $\{\lambda_j\}$. This means that we can write
\begin{equation}
\label{pp_space_obs}
\mathcal{K}^ng=\sum_j  \lambda_j^n\langle g, \phi_j\rangle \phi_j\quad \forall g \in \mathcal{H}_{\mathrm{pp}},n\in\mathbb{N}.
\end{equation}
The spectrum of $\mathcal{K}\restriction_{\mathcal{H}_{\mathrm{pp}}}$ need not be a discrete subset of $\mathbb{T}$. For example, an ergodic rotation on the circle has eigenvalues that densely fill $\mathbb{T}$. In contrast to \eqref{pp_space_obs}, observables in the continuous part $\mathcal{H}_{\mathrm{c}}$ exhibit a decay of correlations that is typical of chaotic systems. Namely, for any $\epsilon>0$ \citep[p.45]{katznelson2004introduction},
$$
\lim_{n\rightarrow\infty} \frac{1}{n} \sum_{j=1}^n \left|\langle \mathcal{K}^jg,f \rangle\right|^\epsilon=0\quad \forall g\in \mathcal{H}_{\mathrm{c}}, f\in L^2(\Omega,\omega).
$$
This result says that $|\langle \mathcal{K}^jg,f \rangle|$ converges to zero in density, that is, for any $\delta>0$, the proportion in all sufficiently large intervals of integers $j$ such that $\left|\langle \mathcal{K}^jg,f \rangle\right|>\delta$ is arbitrarily small.

The above dichotomy is an example of how the decomposition of $\mathcal{E}$ into atomic and continuous parts often characterizes a dynamical system. For example, suppose that $\Fv$ is measure-preserving and bijective, and $\omega$ is a probability measure. Then, the dynamical system is \citep{halmos2017lectures}
\begin{itemize}
	\item \textbf{Ergodic} if and only if $\lambda=1$ is a simple eigenvalue of $\mathcal{K}$,
	\item \textbf{Weakly mixing} if and only if $\lambda=1$ is a simple eigenvalue of $\mathcal{K}$ and there are no other eigenvalues,
	\item \textbf{Mixing} if $\lambda=1$ is a simple eigenvalue of $\mathcal{K}$, and $\mathcal{K}$ has absolutely continuous spectrum on $\mathrm{span}\{1\}^\perp$.
\end{itemize}
Different spectral types find interpretations across various applications, including fluid mechanics \citep{mezic2013analysis}, anomalous transport \citep{zaslavsky2002chaos}, and the analysis of invariants/exponents related to trajectories \citep{kantz2004nonlinear}. The approximation of $\mathcal{E}$ is critical in many applications. For example, the approximate spectral projections provide reduced-order models \citep{mezic2004comparison,mezic2005spectral}. 

\subsection{The fundamental DMD algorithm}

With the definition of a Koopman operator in hand, we can now present the fundamental DMD algorithm and two interpretations. The first interpretation of DMD is as a linear regression. The second is as a projection method. Both interpretations are instrumental, and understanding their interplay is often key to unlocking the power of DMD.

\subsubsection{The linear regression interpretation}
\label{sec_DMD_regression_view}

The simplest and historically first interpretation of DMD is as a linear regression. Given the snapshot matrices $\Xv,\Yv\in\mathbb{C}^{d\times M}$ in \eqref{snapshot_data}, we seek a matrix $\Kv_{\mathrm{DMD}}$ such that $\Yv\approx \Kv_{\mathrm{DMD}}\Xv.$ We can think of this as constructing a \textit{linear and approximate} dynamical system. To find a suitable matrix $\Kv_{\mathrm{DMD}}$, we consider the minimization problem
\begin{equation}
\label{DMD_opt_vanilla}
\min_{\Kv_{\mathrm{DMD}}\in\mathbb{C}^{d\times d}} \left\|\Yv-\Kv_{\mathrm{DMD}}\Xv\right\|_{\mathrm{F}},
\end{equation}
where $\|\cdot\|_{\mathrm{F}}$ denotes the Frobenius norm. Similar optimization problems will be at the heart of the various DMD-type algorithms we consider in this review. A solution to the problem in \eqref{DMD_opt_vanilla} is
$$
\Kv_{\mathrm{DMD}}=\Yv \Xv^{\dagger}\in\mathbb{C}^{d\times d},
$$
where $\dagger$ denotes the Moore--Penrose pseudoinverse. Often, the matrices $\Xv$ and $\Yv$ are tall and skinny, meaning that $d\gg M$. In this scenario, we typically first project onto a low-dimensional subspace to reconstruct the leading nonzero eigenvalues and eigenvectors of the matrix $\Kv_{\mathrm{DMD}}$ without explicitly computing it. The standard DMD algorithm does this using an SVD and is summarized in \cref{alg:DMD_vanilla}, where we have assumed that the projected DMD matrix is diagonalizable.\footnote{We make this assumption about various matrices throughout. Mathematically, a Jordan decomposition may be substituted for an eigendecomposition, and the modes corresponding to a single Jordan block can be considered as interacting modes. However, computing a Jordan block should be avoided. A stable alternative is a Schur decomposition that provides an orthogonal set of interacting modes (in sharp contrast with what is typically considered a DMD mode) or a block-diagonal Schur decomposition with nearly confluent eigenvalues grouped together.} \cref{alg:DMD_vanilla} is known as \textit{exact DMD} \citep{tu2014dynamic} and often the modes are further scaled by $\Lambda^{-1}$. There are several remarks about this algorithm that are worth mentioning:

\begin{algorithm}[t]
\textbf{Input:} Snapshot data $\Xv\in\mathbb{C}^{d\times M}$ and $\Yv\in\mathbb{C}^{d\times M}$, rank $r\in\mathbb{N}$. \\
\vspace{-4mm}
\begin{algorithmic}[1]
\State Compute a truncated SVD of the data matrix $
\Xv \approx \Uv \mathbf{\Sigma}\Vv^*$, $\Uv\in\mathbb{C}^{d\times r}$, $\mathbf{\Sigma}\in\mathbb{R}^{r\times r}$,$\Vv\in\mathbb{C}^{M\times r}.$ The columns of $\Uv$ and $\Vv$ are orthonormal and $\mathbf{\Sigma}$ is diagonal.
\State Compute the compression $
\tilde{\Kv}_{\mathrm{DMD}}=\Uv^*\Yv\Vv\mathbf{\Sigma}^{-1}\in\mathbb{C}^{r\times r}.$
\State Compute the eigendecomposition $
\tilde{\Kv}_{\mathrm{DMD}}\Wv=\Wv\mathbf{\Lambda}$.

\noindent{}The columns of $\Wv$ are eigenvectors and $\mathbf{\Lambda}$ is a diagonal matrix of eigenvalues.
\State Compute the modes $
\mathbf{\Phi}=\Yv\Vv\mathbf{\Sigma}^{-1}\Wv.$
\end{algorithmic} \textbf{Output:} The eigenvalues $\mathbf{\Lambda}$ and modes $\mathbf{\Phi}\in\mathbb{C}^{d\times r}$.
\caption{The exact DMD algorithm \citep{tu2014dynamic}, which has become the workhorse DMD algorithm.}
\label{alg:DMD_vanilla}
\end{algorithm}

\begin{itemize}
	\item The rank $r$ is usually chosen based on the decay of singular values of $\Xv$. If low-dimensional structure is present in the data \citep{udell2019big}, the singular values decrease rapidly, and small $r$ captures the dominant modes. Moreover, the lowest energy modes may be corrupted by noise, and low-dimensional projection is a form of spectral filtering which has the positive effect of dampening the influence of noise \citep{hansen2006deblurring}.\footnote{Low-energy modes can be important though, for example, in optimized control \citep{rowley2005model,rowley2006linear}.} The question of how best to truncate is difficult to answer and is often performed heuristically. If the measurement error is additive white noise, there are algorithmic choices \citep{gavish2014optimal}. In the context of Koopman operators, $r$ is equivalent to the size of the space spanned by basis functions, and a good choice depends on the chosen observables. For example, we shall see below that \cref{alg:DMD_vanilla} corresponds to a linear set of basis functions, which may not capture the relevant nonlinear dynamics. Hence, a larger $r$ may be suitable for other basis choices. Often, the choice of $r$ is modest, meaning that randomized methods \citep{halko2011finding} for computing the SVD can significantly reduce the computational cost. We will explore this and other compression methods in \cref{sec:compression_random}.
	\item We can interpret the algorithm as constructing a linear model of the dynamical system on projected coordinates $\tilde{\xv}=\Uv^*\xv$. Namely, $\tilde{\xv}_{n+1}\approx \tilde{\Kv}_{\mathrm{DMD}} \tilde{\xv}_{n}$. The left singular vectors $\Uv$ are known as proper orthogonal decomposition (POD) modes \citep{berkooz1993proper}.
	\item If the SVD is exact, so that $\Xv = \Uv \mathbf{\Sigma}\Vv^*$, then
	$$
	\Kv_{\mathrm{DMD}}=\Yv \Vv \mathbf{\Sigma}^{-1}\Uv^*.
	$$
	Using this relation, we have
	$$
	\Kv_{\mathrm{DMD}}[\Yv\Vv\mathbf{\Sigma}^{-1}\Wv] =\Yv \Vv \mathbf{\Sigma}^{-1}\underbrace{\Uv^*\Yv\Vv\mathbf{\Sigma}^{-1}}_{\tilde{\Kv}_{\mathrm{DMD}}}\Wv= [\Yv\Vv\mathbf{\Sigma}^{-1}\Wv] \mathbf{\Lambda},
	$$
	and hence \cref{alg:DMD_vanilla} computes exact eigenvalues and eigenvectors of $\Kv_{\mathrm{DMD}}$. Moreover, one can show that this process identifies all of the nonzero eigenvalues of $\Kv_{\mathrm{DMD}}$ \citep[Thm. 1]{tu2014dynamic}. It is common to call $\Yv\Vv\mathbf{\Sigma}^{-1}\Wv$ \textit{exact modes} and $\Uv\Wv$ \textit{projected modes}.
	\item Originally, DMD was developed in connection with Krylov subspaces and the Arnoldi algorithm. In this version, it is assumed that data is gathered along a single trajectory. The SVD version, on the other hand, is capable of handling more general trajectory data. Strategies for using this flexibility to reduce computational cost and average snapshot data noise are given in \citep{tu2014dynamic}. The SVD version is also more numerically stable. Drma{\v{c}} has carefully analyzed the stability of DMD \citep{drmac2018data,drmac2019data,drmavc2020dynamic,drmac2020least}.
	\item Centering the data before applying DMD can be helpful if the mean-subtracted data have linearly dependent columns, especially if the dynamics are perturbations about an equilibrium \citep{hirsh2020centering}. This is equivalent to including an affine term in the linear regression. However, computing the DMD of centered data can be restrictive and have undesirable consequences \citep{chen2012variants}.
\end{itemize}

The above interpretation of DMD is simple and intuitive. However, using DMD to analyze nonlinear dynamics globally seems dubious, as there is an underlying assumption of approximately linear dynamics in \eqref{DMD_opt_vanilla}. Nevertheless, we shall now demonstrate that DMD can be interpreted as an approximation to Koopman spectral analysis. This provides a solid theoretical foundation for applying DMD in analyzing nonlinear dynamics, which will be further elaborated upon in \cref{sec:Galerkin_perspective}.

\subsubsection{The Galerkin interpretation: Connecting to Koopman operators}
\label{sec:galerkin_interp}

The connection between \cref{alg:DMD_vanilla} and Koopman operators is revealed once we interpret DMD as a Galerkin method. Consider the two correlation matrices
$$
\Gv=\frac{1}{M}\overline{\Uv^*\Xv} (\Uv^*\Xv)^\top,\quad \Av = \frac{1}{M}\overline{\Uv^*\Xv}(\Uv^*\Yv)^\top.
$$
We can think of the $j$th row of the POD matrix $\Uv^*\Xv$ as an affine function $u_j$ on the statespace $\Omega$ evaluated at the snapshot data:
$$
u_j(\xv)=[\Uv_{:,j}]^*\xv,\quad u_j(\xv^{(m)})=[\Uv^*\Xv]_{jm}.
$$
It follows that $\Gv$ can be interpreted as a Gram matrix with respect to the positive semi-definite Hermitian form induced by the probability measure with equal point masses at the $\{\xv^{(m)}\}$. Namely,
$$
\Gv_{jk}=\frac{1}{M}\sum_{m=1}^M \overline{u_j(\xv^{(m)})}u_k(\xv^{(m)})=\int_{\Omega} \overline{u_j(\xv)}u_k(\xv) \ \mathrm{d} \omega_M(\xv),\quad \omega_M=\frac{1}{M}\sum_{m=1}^M\delta_{\xv^{(m)}}.
$$
Writing $\langle \cdot,\cdot\rangle_{M}$ for the form induced by $\omega_M$, we can argue similarly for $\Av$ and succinctly write
$$
\Gv_{jk}=\langle u_k, u_j \rangle_M,\quad \Av_{jk}=\frac{1}{M}\sum_{m=1}^M \overline{u_j(\xv^{(m)})}u_k(\yv^{(m)})=\int_{\Omega} \overline{u_j(\xv)}u_k(\Fv(\xv)) \ \mathrm{d} \omega_M(\xv)=\langle \mathcal{K}u_k, u_j \rangle_{M}.
$$
Assuming that the matrix $\Uv^*\Xv$ is of rank $r$, and using $\Uv^*\Xv=\mathbf{\Sigma}\Vv^*$, we can write
$$
\tilde{\Kv}_{\mathrm{DMD}}^\top=\mathbf{\Sigma}^{-1}\Vv^\top\Yv^\top\overline{\Uv}=(\Xv^\top\overline{\Uv})^\dagger \Yv^\top\overline{\Uv}= \Gv^{-1}\Av.
$$
The matrix $\Gv^{-1}\Av$ is an approximation of the action of $\mathcal{K}$ on the subspace spanned by the functions $\{u_j\}_{j=1}^r$. Namely, if $g$ is an observable that can be expressed as the linear combination
$$
g(\xv)=\sum_{j=1}^r u_j(\xv)\gv_j,\quad\text{for some}\quad \gv\in\mathbb{C}^r,
$$
then
$$
[\mathcal{K}g](\xv)\approx \sum_{j=1}^r u_j(\xv)(\Gv^{-1}\Av\gv)_j=\sum_{j=1}^r u_j(\xv)(\tilde{\Kv}_{\mathrm{DMD}}^\top\gv)_j.
$$
In other words, $\tilde{\Kv}_{\mathrm{DMD}}^\top$ is a matrix that approximates the action of the Koopman operator on expansion coefficients. More precisely, it is a Galerkin method corresponding to $\mathcal{K}$ and the form $\langle \cdot , \cdot\rangle_M$. This connection is explored more deeply in \cref{sec:EDMD}.

If $\langle \cdot , \cdot\rangle_M$ converges to $\langle \cdot , \cdot\rangle$ in the large data limit $M\rightarrow\infty$, then DMD can be considered to be a numerical approximation to Koopman spectral analysis. The terms \textit{DMD mode} and \textit{Koopman mode} are often used interchangeably in the literature. It is important to note that the Koopman modes and eigenfunctions are distinct mathematical objects, requiring different approaches for approximation. The right eigenvectors of $\tilde{\Kv}_{\mathrm{DMD}}$ give rise to time-invariant directions in the state space $\xv$, whereas the left-eigenvectors give rise to Koopman eigenfunctions, which are similarly time-invariant directions in the space of observables.

\subsubsection{The Koopman mode decomposition}\label{sec:KMD_linear}

We can now connect DMD with the spectral expansions discussed in \cref{sec:spectra_crash_course}. First, we approximate an initial condition $\xv_0$ in the eigenvector coordinates via
$$
\xv_0\approx\mathbf{\Phi}\bv,\quad \bv=\mathbf{\Phi}^\dagger \xv_0.
$$
This is not the only choice, but it is the simplest. The KMD then provides an approximation of the dynamics by
\begin{equation}
\label{linear_KMD}
\xv_n\approx {\Kv}_{\mathrm{DMD}}^n\xv_0\approx  {\Kv}_{\mathrm{DMD}}^n \mathbf{\Phi}\bv= \mathbf{\Phi}\mathbf{\Lambda}^n\bv,
\end{equation}
which echoes \eqref{KMD_spec_meas_forward}. Since zero eigenvalues do not contribute to the dynamics, this decomposition further justifies the compression in \cref{alg:DMD_vanilla}. The KMD has also been related to other decompositions in various situations, particularly those that have arisen in the fluid dynamics community \citep{taira2017modal,taira2020modal}. These include POD \citep{towne2018spectral}, optimal mode decomposition \citep{wynn2013optimal}, and resolvent analysis \citep{sharma2016correspondence,herrmann2021data}. Under suitable conditions, the KMD converges as we increase the dimension of the projected Koopman operator (see \cref{sec:EDMD_convergence} and \cref{sec:mpEDMD_convergence}).

\subsection{Three canonical examples}

Having grasped the notion of Koopman operators and the basic DMD algorithm, it is time for some examples. As a warm-up for the reader, we consider three canonical well-studied examples of \cref{alg:DMD_vanilla}, each with a unitary Koopman operator:
\begin{itemize}
	\item The flow past a cylinder wake at $Re=100$ with a state space dimension $d=160,000$ that corresponds to the number of spatial measurement points in the flow. The associated Koopman operator has a pure point spectrum consisting of powers of a fundamental eigenvalue.
	\item The Lorenz system on the Lorenz attractor with a state space dimension $d=3$. The associated Koopman operator possesses no eigenvalues, except for the simple eigenvalue $\lambda=1$ whose eigenfunction is the constant function. The rest of the spectrum is continuous.
	\item The Duffing oscillator with a state space dimension $d=2$. The associated Koopman operator possesses no eigenvalues, except for $\lambda=1$, whose eigenspace now corresponds to the conserved Hamiltonian energy of the system and indicator functions associated with invariant sets of positive area.
\end{itemize}
Despite its large ambient space dimension, the first example is the easiest to address using \cref{alg:DMD_vanilla}. This is because the cylinder wake exhibits an attracting limit cycle, and the Koopman operator has a basic spectrum. The other two examples demonstrate three difficulties of DMD: noise, projection error (which can lead to spurious modes and missing parts of the spectrum), and continuous spectra. These and further challenges are discussed in \cref{sec:challenges_of_DMD}.

\subsubsection{Flow past a cylinder wake}
\label{sec:example_cylinder}

We first consider the classic DMD example of low Reynolds number flow past a circular cylinder. Due to its simplicity and relevance in engineering, this is one of the most studied examples in modal-analysis techniques \citep[Table 3]{rowley2017model}\citep{chen2012variants,taira2020modal}. $Re=100$ is chosen so that it is larger than the critical Reynolds number at which the flow undergoes a supercritical Hopf bifurcation, resulting in laminar vortex shedding \citep{jackson1987finite,zebib1987stability}. This limit cycle is stable and is representative of the three-dimensional flow \citep{noack1994global,noack2003hierarchy}. The Koopman operator of the post-transient flow has a pure point spectrum with a lattice structure on the unit circle \citep{bagheri2013koopman}.

\begin{figure}
\centering
\raisebox{-0.5\height}{\includegraphics[width=0.28\textwidth,trim={0mm 0mm 0mm 0mm},clip]{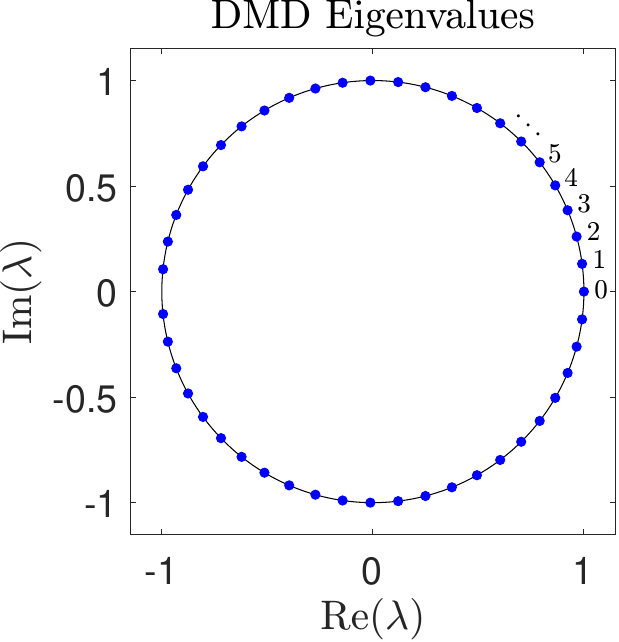}}
\raisebox{-0.5\height}{\includegraphics[width=0.32\textwidth,trim={0mm 0mm 0mm 0mm},clip]{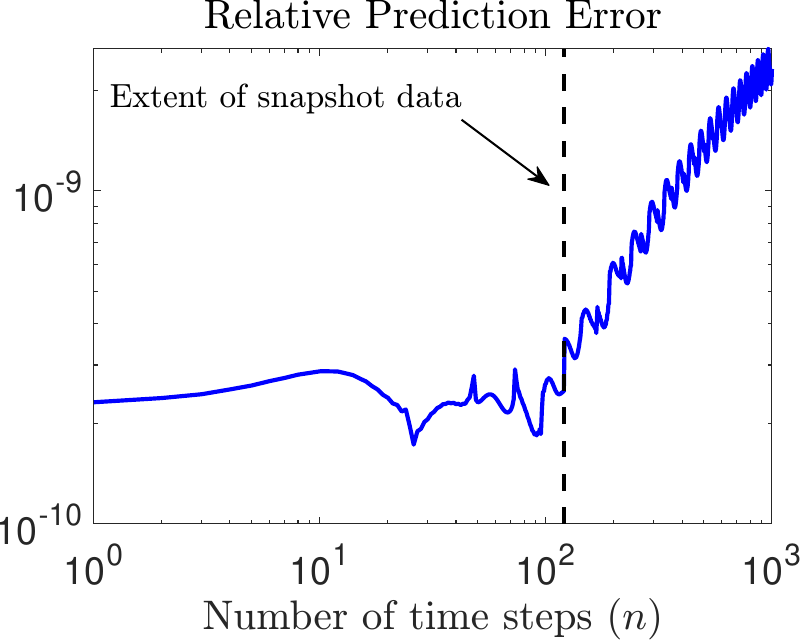}}
\raisebox{-0.5\height}{\includegraphics[width=0.38\textwidth,trim={15mm 0mm 12mm 0mm},clip]{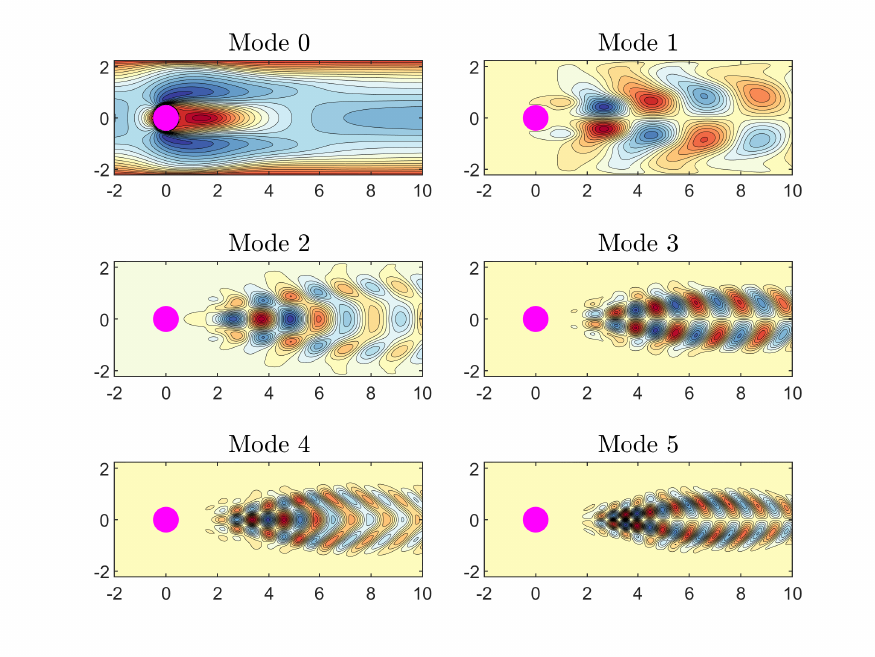}}
\caption{Output of \cref{alg:DMD_vanilla} for the flow past a cylinder wake. Left: The DMD eigenvalues. Middle: The total relative prediction error of \eqref{linear_KMD}. Right: Example modes of the horizontal component of the velocity field. The magenta disc corresponds to the cylinder. The zeroth mode corresponds to the time-averaged flow.}
\label{wake1}
\end{figure}

To collect snapshot data, we numerically compute the velocity field of a flow around a circular cylinder of diameter $D=1$ using an incompressible, two-dimensional lattice-Boltzmann solver \citep{jozsa2016validation,szHoke2017performance}. The temporal resolution of the flow is chosen so that approximately 24 snapshots of the flow field correspond to the period of vortex shedding. The computational domain size is $18D$ in length and $5D$ in height, with a $800\times 200$ grid resolution. The cylinder is positioned $2D$ downstream of the inlet at the mid-height of the domain. The cylinder side walls are defined as bounce-back and no-slip walls, and a parabolic velocity profile is given at the inlet of the domain. The outlet is defined as a non-reflecting outflow. After simulations converge to steady-state vortex shedding, we collect $M=120$ snapshots for the DMD algorithm and a further $880$ snapshots to test the prediction of the KMD. One should think of this as training data and test data, respectively. Letting $\Vv_x(t)$ denote the vectorized horizontal velocity field at time $t$, our snapshot matrices have the form
$$
\Xv=\begin{pmatrix} 
\Vv_x(0) & \Vv_x(\Delta t) &\cdots & \Vv_x(119 \Delta t)
\end{pmatrix}, \quad \Yv=\begin{pmatrix} 
\Vv_x(\Delta t) & \Vv_x(2\Delta t) &\cdots & \Vv_x(120 \Delta t)
\end{pmatrix}.
$$
We use a rank of $r=47$ to recover the trivial mode corresponding to $\lambda=1$ and $24$ conjugate pairs of modes up to the timescale of vortex shedding. The eigenvalues come in conjugate pairs due to processing real-valued data $\Xv$ and $\Yv$.

\cref{wake1} shows the output of \cref{alg:DMD_vanilla}. In the left panel, we see the lattice structure of the DMD modes correctly identified by \cref{alg:DMD_vanilla}. In the middle plot, we show the predictive error of \eqref{linear_KMD}. The relative error is computed by taking the $2$-norm of the error in the velocity field $\Vv_x$ and normalizing it by the $2$-norm of the mean-subtracted flow at each time step. Due to the periodic nature of the flow, there is excellent agreement between the KMD and flow, with slow algebraic growth of the error beyond the snapshot data time window. The right panel of \cref{wake1} shows the real part of some of the Koopman modes for the horizontal velocity field.

\subsubsection{Lorenz system}
\label{sec:intro_lorenz}

The Lorenz (63) system \citep{lorenz1963deterministic} is the following three coupled ordinary differential equations:
$$
\dot{x}=10\left(y-x\right),\quad\dot{y}=x\left(28-z\right)-y,\quad \dot{z}=xy-8z/3.
$$
We consider the dynamics of $\xv=(x,y,z)$ on the Lorenz attractor. The system is chaotic and strongly mixing \citep{luzzatto2005lorenz}. It follows that $\lambda=1$ is the only eigenvalue of $\mathcal{K}$, corresponding to a constant eigenfunction, and that this eigenvalue is simple. We consider a discrete-time dynamical system by sampling with a time-step $\Delta t=0.001$. We use time-delay embedding, which is a popular method for DMD-type algorithms \citep{susuki2015prony,arbabi2017ergodic,brunton2017chaos,das2019delay,kamb2020time,pan2020structure} and corresponds to building a Krylov subspace. This technique is justified through Takens' embedding theorem \citep{takens2006detecting}, which says that under certain technical conditions, delay embedding a signal coordinate of the system can reconstruct the attractor of the original system up to a diffeomorphism. In this example, we augment $\xv$ by $N-1$ further time-delays of length $\Delta t'=0.2$ and consider $M=5\times 10^5$ snapshots along a single trajectory. Specifically, our snapshot matrices have the form
\begin{align*}
\Xv&=\begin{pmatrix} 
\xv(0) & \xv(\Delta t) &  \cdots & \xv( (M{-}1)\Delta t)\\ 
\xv(\Delta t') & \xv(\Delta t'{+} \Delta t)  &\cdots & \xv(\Delta t'{+}(M{-}1)\Delta t)\\
\vdots & \vdots & \vdots & \vdots  \\
\xv((N{-}1)\Delta t') & \xv((N{-}1)\Delta t'{+} \Delta t) & \cdots & \xv((N{-}1)\Delta t'{+}(M{-}1)\Delta t)
\end{pmatrix}\in\mathbb{R}^{3N\times M},\\
\Yv&=\begin{pmatrix} 
\xv(\Delta t) & \xv(2\Delta t) &  \cdots & \xv( M\Delta t)\\ 
\xv(\Delta t'{+}\Delta t) & \xv(\Delta t'{+} 2\Delta t) & \cdots & \xv(\Delta t'{+}M\Delta t)\\
\vdots & \vdots & \vdots & \vdots  \\
\xv((N{-}1)\Delta t'{+}\Delta t) & \xv((N{-}1)\Delta t'{+} 2\Delta t) &\cdots & \xv((N{-}1)\Delta t'{+}M\Delta t)
\end{pmatrix}\in\mathbb{R}^{3N\times M}.
\end{align*}
We use the \texttt{ode45} command in MATLAB to collect the data after an initial burn-in time to ensure that the initial point $\xv(0)$ is (approximately) on the Lorenz attractor. The system is chaotic, so we cannot hope to integrate for long periods accurately numerically. However, convergence is still obtained in the large data limit $M\rightarrow\infty$ due to an effect known as shadowing.

\begin{figure}
\centering
\raisebox{-0.5\height}{\includegraphics[width=0.31\textwidth,trim={0mm 0mm 0mm 0mm},clip]{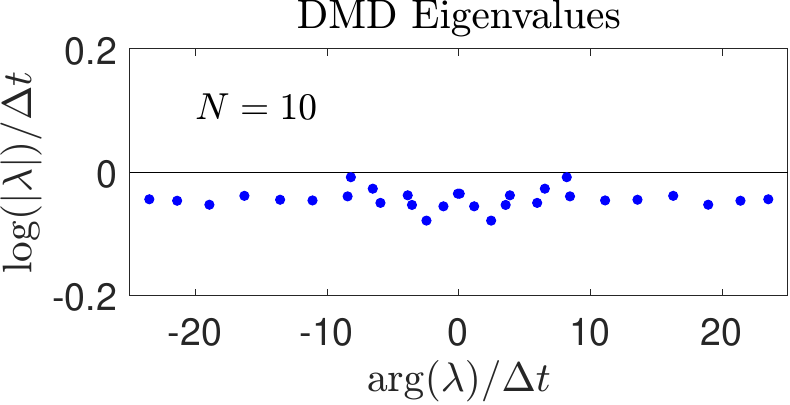}}
\hfill
\raisebox{-0.5\height}{\includegraphics[width=0.31\textwidth,trim={0mm 0mm 0mm 0mm},clip]{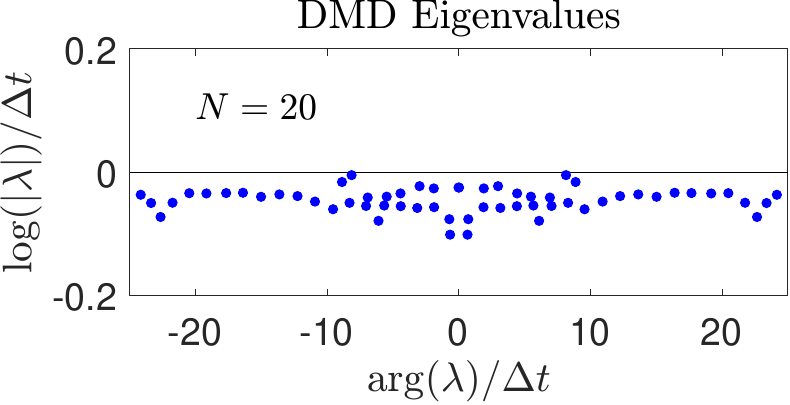}}
\hfill
\raisebox{-0.5\height}{\includegraphics[width=0.31\textwidth,trim={0mm 0mm 0mm 0mm},clip]{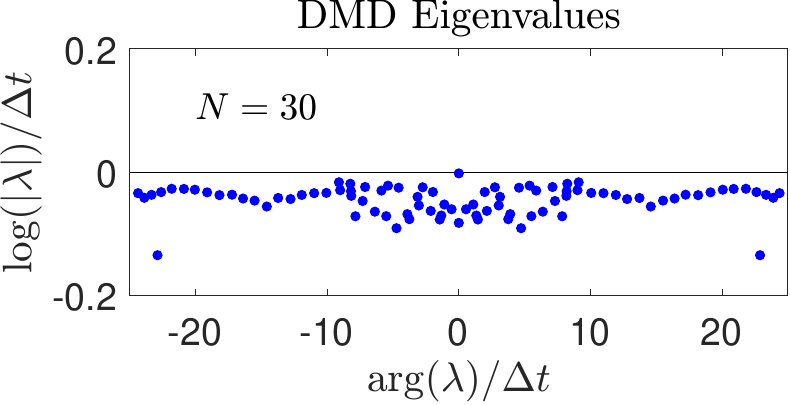}}\\\vspace{5mm}
\raisebox{-0.5\height}{\includegraphics[width=0.31\textwidth,trim={0mm 0mm 0mm 0mm},clip]{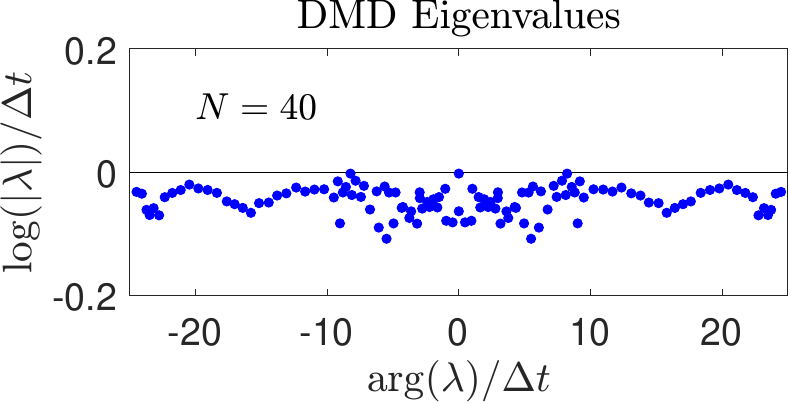}}
\hfill
\raisebox{-0.5\height}{\includegraphics[width=0.31\textwidth,trim={0mm 0mm 0mm 0mm},clip]{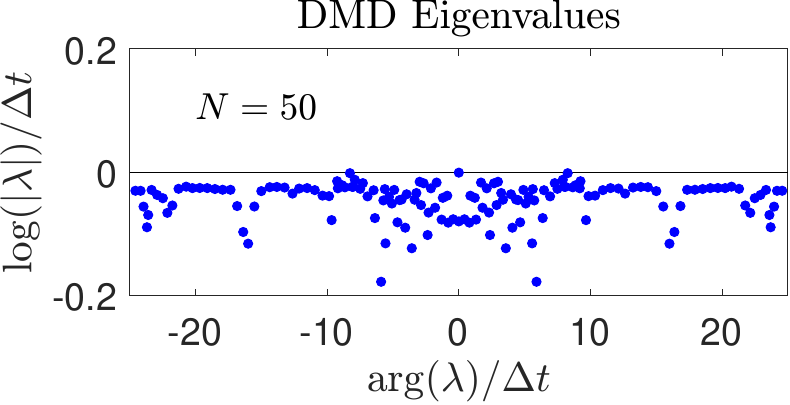}}
\hfill
\raisebox{-0.5\height}{\includegraphics[width=0.31\textwidth,trim={0mm 0mm 0mm 0mm},clip]{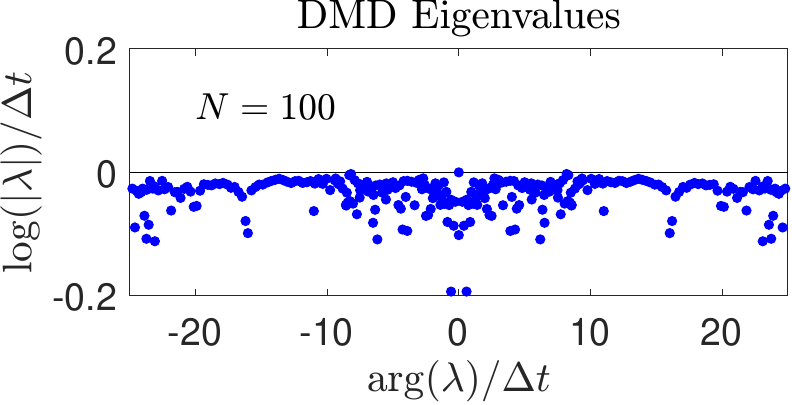}}
\caption{DMD eigenvalues for the Lorenz system and different choices of $N$ (the number of eigenvalues is $3N$). The logarithm of the eigenvalues is plotted to align with the continuous-time system. As $N$ increases, the eigenvalues cluster to approximate the continuous spectrum.}
\label{lorenz1}
\end{figure}

In addition to a discrete time-step $\Delta t$, we can consider Koopman operators associated with continuous-time dynamical systems. The continuous-time infinitesimal generator is defined by
\begin{equation}
\label{eq:gen_def}
\mathcal{L} g = \lim_{\Delta t\downarrow 0}\frac{\mathcal{K}_{\Delta t}g-g}{\Delta t},
\end{equation}
where $\mathcal{K}_{\Delta t}$ is the Koopman operator corresponding to a time-step $\Delta t$. The generator satisfies
$$
\mathcal{K}_{\Delta t}=\exp(\Delta t\mathcal{L}),
$$
which can be made precise through the theory of semigroups \citep{pazy2012semigroups}. Hence, in this example, we consider the following time-scaled logarithms of the eigenvalues:
$$
{\log(\lambda)}/{\Delta t}= {\log(|\lambda|)}/{\Delta t} +i {\mathrm{arg}(\lambda)}/{\Delta t}.
$$
\cref{lorenz1} shows the DMD eigenvalues for various choices of $N$. The horizontal line corresponds to a portion of the spectrum of the Koopman operator. The DMD eigenvalues fall below this line, corresponding to a dampening effect in the dynamics encapsulated by DMD. For these choices of parameters, this error is due to the finite amount of trajectory data and noise in the data matrices from the numerical solver. In \cref{sec:noise_examples}, we shall see that this effect can be reduced using DMD methods designed to be robust to noise and by increasing $M$. For other parameter choices, DMD approximations for the Lorenz system (and other systems) can also suffer from projection errors. Another interesting feature of \cref{lorenz1} is the clustering of the DMD eigenvalues with increasing $N$ as they attempt to approximate the continuous spectrum. Recall from above that the Koopman operator for this example has no eigenvalues except the trivial eigenvalue $\lambda=1$. In \cref{sec:mpEDMD}, we shall see that \textit{Measure-Preserving Extended DMD} \citep{colbrook2023mpedmd} can deal with continuous spectra (see also the discussion in \cref{sec:methods_for_spec_meas} for further methods).

\begin{figure}
\centering
\raisebox{-0.5\height}{\includegraphics[width=1\textwidth,trim={0mm 0mm 0mm 0mm},clip]{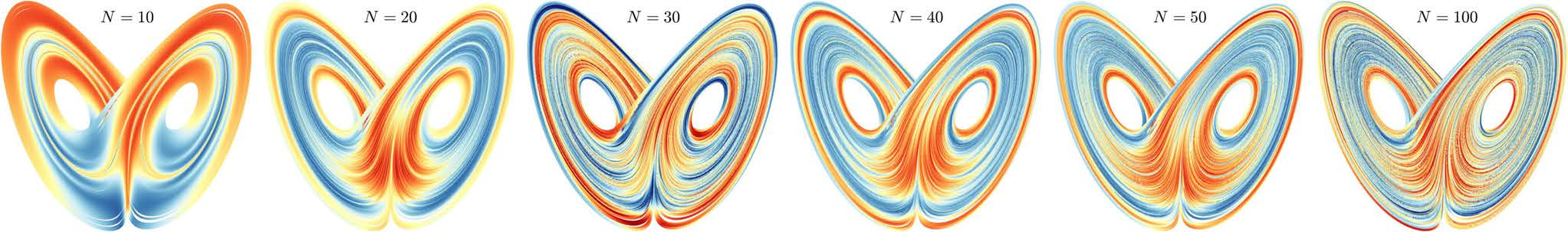}}
\caption{DMD eigenfunctions corresponding to $\log(\lambda)/(i\Delta t)\approx \pm 8.2\ \mathrm{rad}/\mathrm{s}$. These are similar to the singularity in spectral measures detected in \citep{korda2020data}.}
\label{lorenz2}
\end{figure}

\begin{figure}
\centering
\raisebox{-0.5\height}{\includegraphics[width=1\textwidth,trim={0mm 0mm 0mm 0mm},clip]{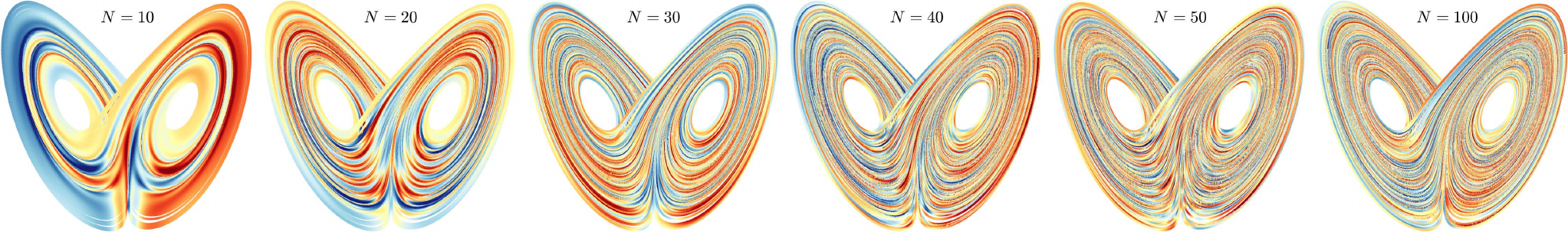}}
\caption{DMD eigenfunctions corresponding to $\log(\lambda)/(i\Delta t)\approx 6\ \mathrm{rad}/\mathrm{s}$. We see increasing oscillations as $N$ gets larger and the eigenfunctions resemble unstable periodic orbits, see also \citep{colbrook2021rigorous}.}
\label{lorenz3}
\end{figure}

We next show DMD eigenfunctions but associated with the matrix $\Xv^\dagger \Yv$. The discussion in \cref{sec:galerkin_interp} shows that these correspond to pseudoeigenfunctions of $\mathcal{K}$. Note that these pseudoeigenfunctions do not approximate true eigenfunctions - since eigenfunctions do not exist for this system! For visualization over the attractor, we plot function values along the trajectory of snapshot data. In \cref{lorenz1} there are DMD eigenvalues close to the horizontal line with $\log(\lambda)/(i\Delta t)\approx \pm 8.2\ \mathrm{rad}/\mathrm{s}$. These correspond to an apparent singularity in the spectral measure detected in \citep{korda2020data}. \cref{lorenz2} shows the corresponding pseudoeigenfunctions. These bear a striking resemblance to the local spectral projections in \citep[Figure 13]{korda2020data}, which the authors attributed to an almost-periodic motion of the $z$ component during the time that the state resides in either of the two lobes of the Lorenz attractor. In \cref{lorenz3}, we plot the pseudoeigenfunctions corresponding to the DMD eigenvalue with $\log(\lambda)/(i\Delta t)$ closest to $6\ \mathrm{rad}/\mathrm{s}$. In this case, we see increasing oscillations as $N$ gets larger, and the pseudoeigenfunctions resemble unstable periodic orbits, which in a sense, form a backbone of the attractor \citep{eckmann1985ergodic,tufillaro1992experimental}. For further examples of these kinds of pseudoeigenfunctions, see \citep{colbrook2021rigorous}.

\subsubsection{Duffing oscillator}
\label{sec:duffing_spectral_pollution}

We now consider the Hamiltonian system
$$
\dot{x}=y,\quad\dot{y}=x-x^3,
$$
known as the (undamped nonlinear) Duffing oscillator with state $\xv=(x,y)\in\Omega=\mathbb{R}^2$. This dynamical system has three fixed points at $\xv=(0,0)$ (a saddle), and $\xv=(\pm1,0)$ (centers). The Hamiltonian for this system is $H=y^2-x^2/2+x^4/2$. We consider the corresponding discrete-time dynamical system by sampling with a time-step $\Delta t=0.25$. Instead of using the statespace $\xv$ or time-delay embedding to form our snapshot matrices, we consider an example of \textit{Extended DMD} \citep{williams2015data}, discussed in more detail in \cref{sec:EDMD}. Specifically, we consider $10^3$ random points sampled uniformly in $[-2,2]^2$, and then the trajectory of these points for $50$ times steps. This leads to $M=5\times 10^4$ snapshots $\{\xv^{(m)},\yv^{(m)}\}_{m=1}^M$. We then partition these into $N$ clusters using k-means, and use these as centers $\cv_j$ for $N$ radial basis functions of the form
$$
\psi_j(\xv)=\exp(-\gamma\|\xv-\cv_j\|),
$$
where $\gamma$ is the squared reciprocal of the average $\ell^2$-norm of the snapshot data after it is shifted to mean zero. Our snapshot matrices are then given by
\begin{align*}
\Xv&=\begin{pmatrix}
\psi_1(\xv^{(1)}) & \psi_1(\xv^{(2)}) &  \cdots & \psi_1(\xv^{(M)})\\ 
\psi_2(\xv^{(1)}) & \psi_2(\xv^{(2)}) &  \cdots & \psi_2(\xv^{(M)})\\
\vdots & \vdots & \vdots & \vdots  \\
\psi_N(\xv^{(1)}) & \psi_N(\xv^{(2)}) &  \cdots & \psi_N(\xv^{(M)})
\end{pmatrix},\quad \Yv&=\begin{pmatrix}
\psi_1(\yv^{(1)}) & \psi_1(\yv^{(2)}) &  \cdots & \psi_1(\yv^{(M)})\\ 
\psi_2(\yv^{(1)}) & \psi_2(\yv^{(2)}) &  \cdots & \psi_2(\yv^{(M)})\\
\vdots & \vdots & \vdots & \vdots  \\
\psi_N(\yv^{(1)}) & \psi_N(\yv^{(2)}) &  \cdots & \psi_N(\yv^{(M)})
\end{pmatrix}.
\end{align*}

\begin{figure}
\raisebox{-0.5\height}{\includegraphics[width=1\textwidth,trim={0mm 0mm 0mm 0mm},clip]{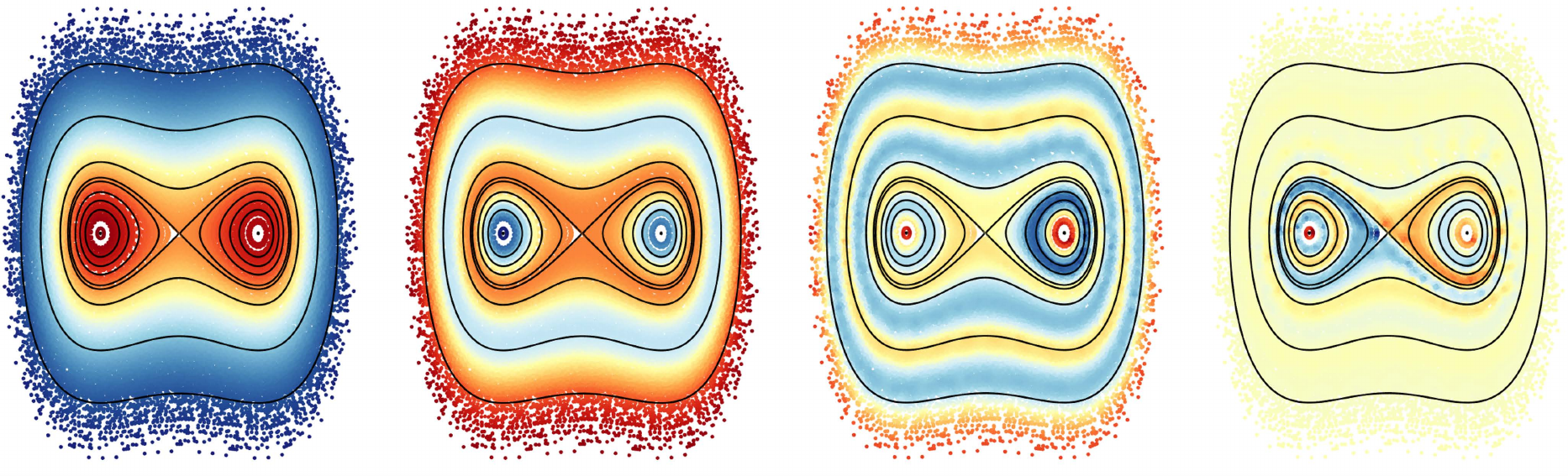}}
\caption{Some examples of computed eigenfunctions of the Duffing oscillator corresponding to $\lambda=1$. The black lines show trajectory orbits and correspond to level sets of the eigenfunctions, which are invariant in time.}
\label{duffing1}
\end{figure}

\begin{figure}
\centering
\raisebox{-0.5\height}{\includegraphics[width=0.32\textwidth,trim={0mm 0mm 0mm 0mm},clip]{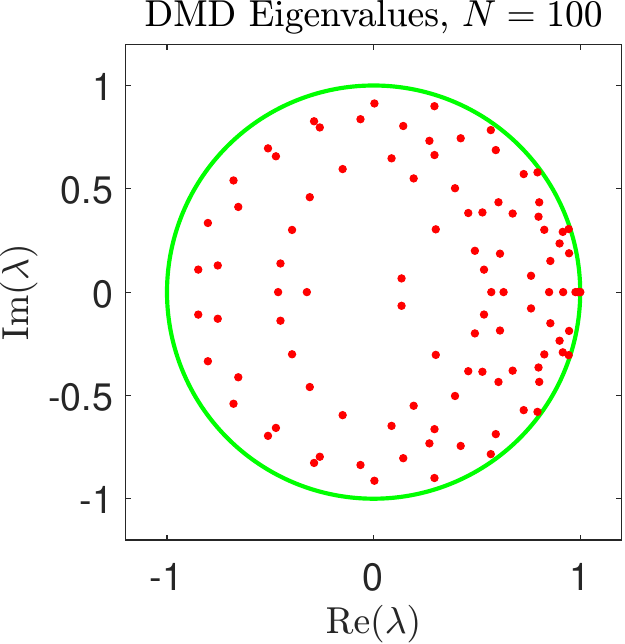}}
\raisebox{-0.5\height}{\includegraphics[width=0.32\textwidth,trim={0mm 0mm 0mm 0mm},clip]{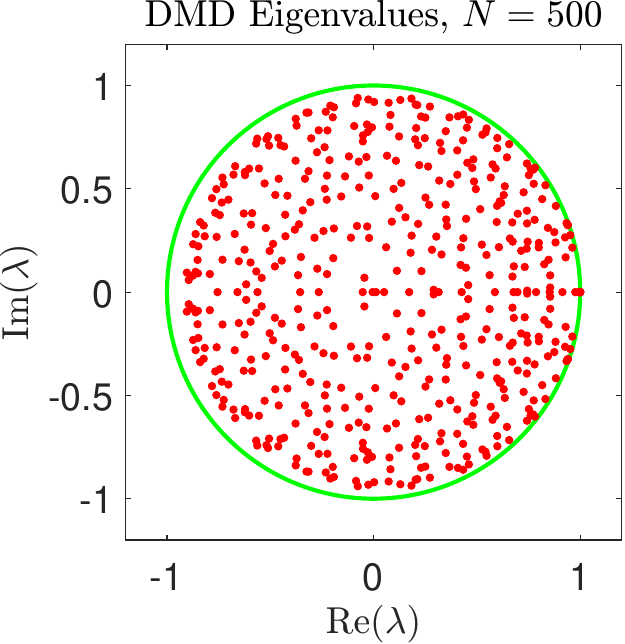}}
\raisebox{-0.5\height}{\includegraphics[width=0.32\textwidth,trim={0mm 0mm 0mm 0mm},clip]{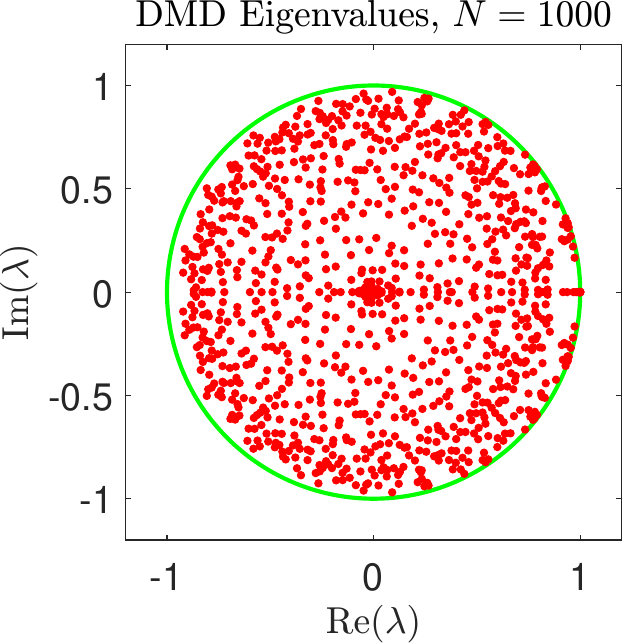}}
\caption{DMD eigenvalues (red dots) computed for various choices of $N$. The spectrum is the unit circle (green); hence, most eigenvalues are spurious. This occurs because a projection error occurs when the Koopman operator $\mathcal{K}$ is approximated by a finite DMD matrix.}
\label{duffing2}
\end{figure}

\cref{duffing1} shows some of the Koopman eigenfunctions corresponding to $\lambda=1$ and computed using $N=1000$. For visualization, the values of the functions are plotted at the data points. We see that the level sets of the eigenfunctions correspond to trajectories, as expected. However, moving away from $\lambda=1$ in the spectral plane becomes more challenging. \cref{duffing2} shows the DMD eigenvalues for various choices of $N$, along with the unit circle, which is the spectrum of the Koopman operator. Most of the DMD eigenvalues are spurious and correspond to spectral pollution. This occurs because we have approximated the infinite-dimensional Koopman operator $\mathcal{K}$, by a finite matrix. These errors persist, even as we increase $N$.

\begin{figure}
\centering
\raisebox{-0.5\height}{\includegraphics[width=0.32\textwidth,trim={0mm 0mm 0mm 0mm},clip]{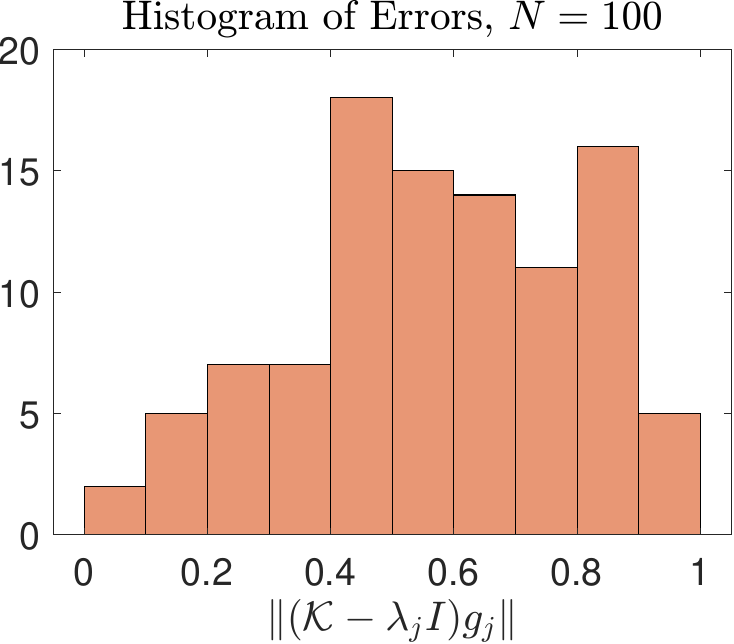}}
\raisebox{-0.5\height}{\includegraphics[width=0.32\textwidth,trim={0mm 0mm 0mm 0mm},clip]{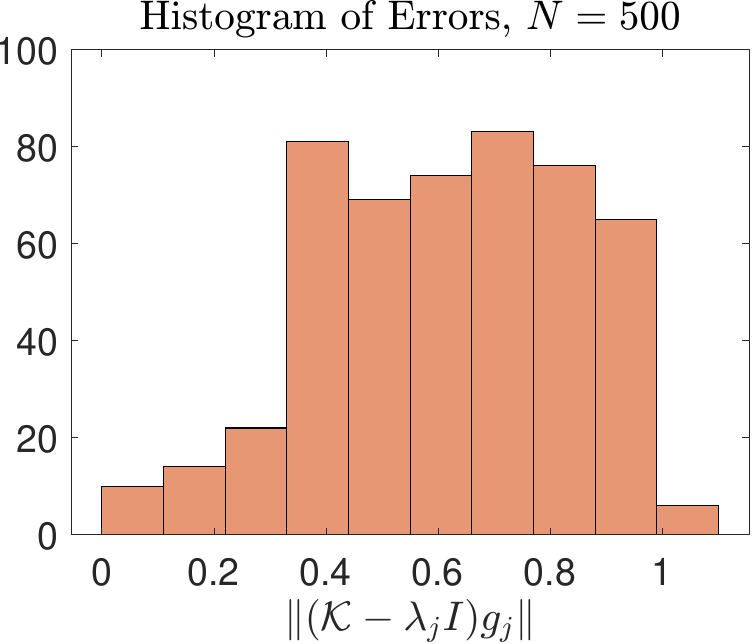}}
\raisebox{-0.5\height}{\includegraphics[width=0.32\textwidth,trim={0mm 0mm 0mm 0mm},clip]{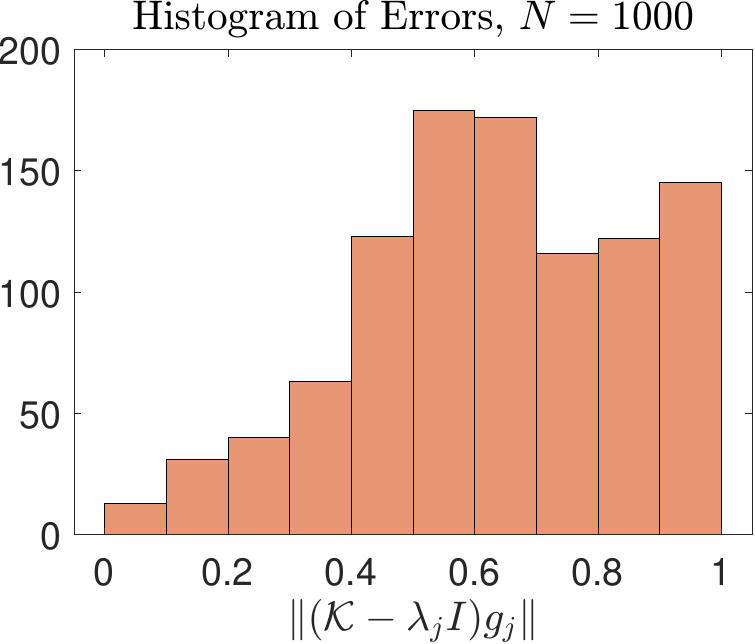}}
\caption{Histograms of the errors of the DMD eigenpairs. The errors are computed using ResDMD in \cref{alg:ResDMD1} and show the persistence of heavy spectral pollution as $N$ increases.}
\label{duffing3}
\end{figure}

\begin{figure}
\centering
\raisebox{-0.5\height}{\includegraphics[width=0.32\textwidth,trim={0mm 0mm 0mm 0mm},clip]{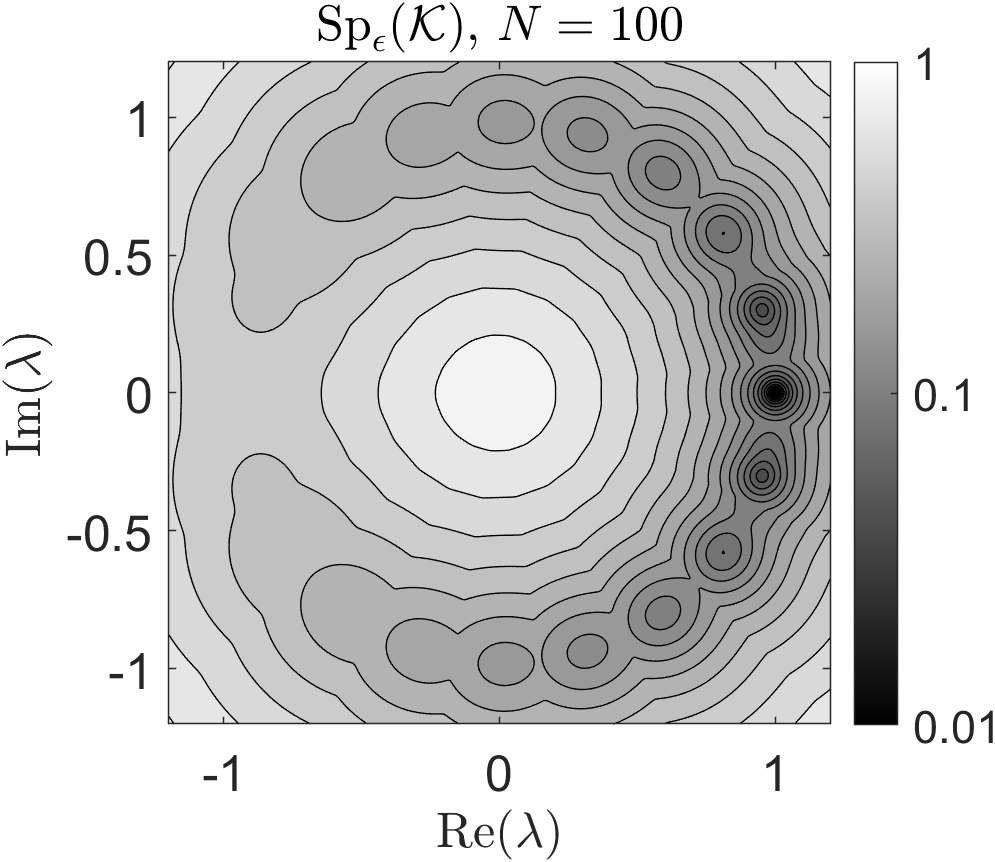}}
\raisebox{-0.5\height}{\includegraphics[width=0.32\textwidth,trim={0mm 0mm 0mm 0mm},clip]{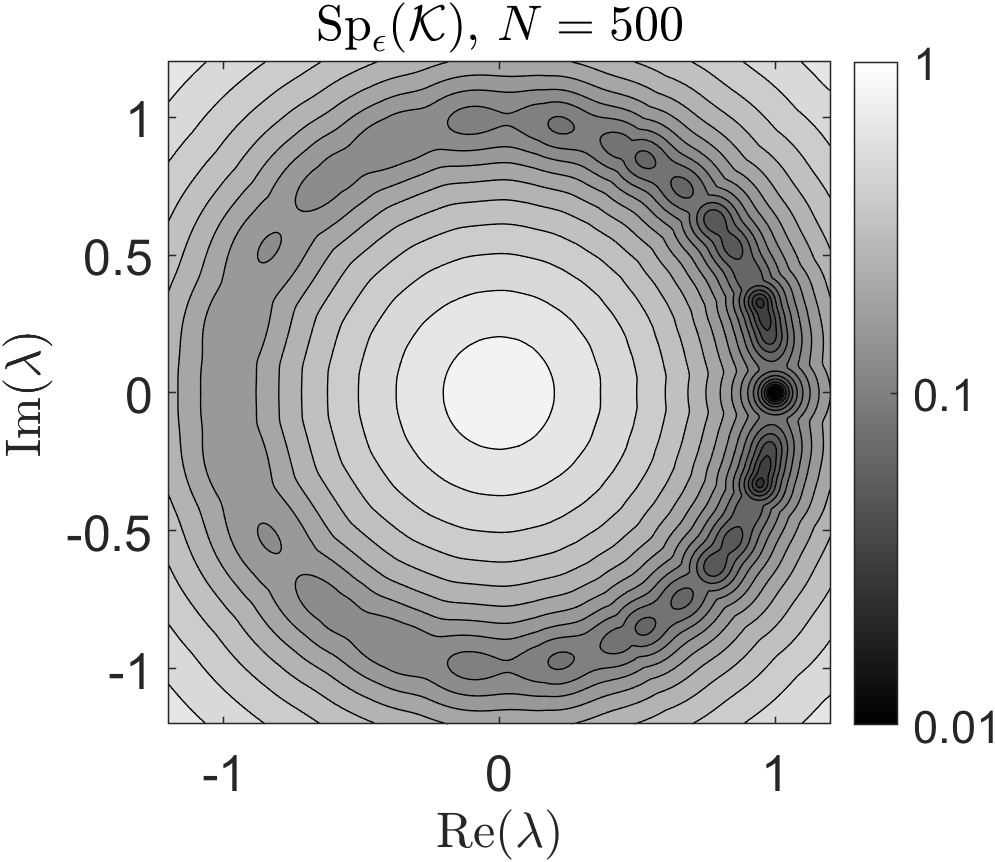}}
\raisebox{-0.5\height}{\includegraphics[width=0.32\textwidth,trim={0mm 0mm 0mm 0mm},clip]{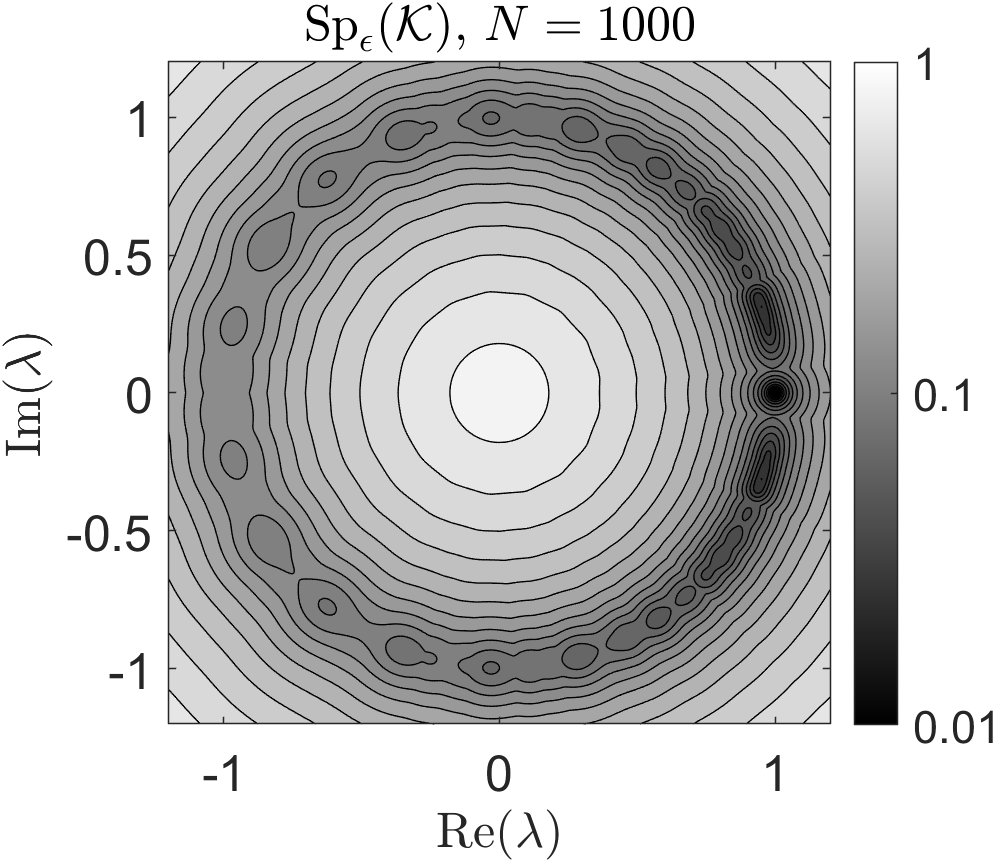}}
\caption{Pseudospectra (see \eqref{eq_def_pseudospectra}) computed using ResDMD (see \cref{alg:ResDMD2}) and visualized by plotting several contour plots of $\epsilon$ on a logarithmic scale. The pseudospectra demonstrate the heavy spectral pollution present in \cref{duffing2}. Note also that these pseudospectra are computed using the same snapshot data and dictionary used for \cref{duffing1,duffing2,duffing3}. As $N\rightarrow\infty$, the algorithm converges to the pseudospectra.}
\label{duffing4}
\end{figure}

To measure the errors, we can use \textit{Residual DMD} (ResDMD) \citep{colbrook2021rigorous,colbrook2023residualJFM} to compute the error $\|(\mathcal{K}-\lambda_jI)g_j\|$ associated with a DMD eigenfunction $g_j$ and eigenvalue $\lambda_j$. In other words, we can compute the projection error of DMD. This is detailed in \cref{sec:ResDMD}. \cref{duffing3} shows the histograms of these projection errors. Only a small proportion of reliable DMD eigenvalues persist as $N$ increases. Finally, we can also use ResDMD to compute pseudospectra. \cref{alg:ResDMD2} converges to the pseudospectrum as $N\rightarrow\infty$. The output is shown in \cref{duffing4}, where we visualize pseudospectra by plotting several contour plots of $\epsilon$ on a logarithmic scale. For this example, the pseudospectra are the annuli $\mathrm{Sp}_\epsilon(\mathcal{K})=\{\lambda\in\mathbb{C}:||\lambda|-1|\leq \epsilon\}$. The projection errors are computed using the same snapshot data and dictionary used for \cref{duffing1,duffing2,duffing3}. ResDMD allows us to compute and minimize projection errors directly in infinite dimensions to avoid spectral pollution and spurious modes.

In summary, DMD can suffer from closure issues (projection errors) associated with approximating the infinite-dimensional Koopman operator by a finite-dimensional matrix. This phenomenon is well-known \citep{brunton2016koopman,kaiser2017data}. Nevertheless, by computing projection errors, we can avoid difficulties (such as spurious modes) associated with the infinite-dimensional nature of the Koopman operator.

\subsection{The goals and challenges of DMD}
\label{sec:challenges_of_DMD}

The core goal of DMD is to apply linear algebra and spectral techniques to the analysis, prediction, and control of nonlinear dynamical systems. However, DMD  often faces several challenges \citep{kutz2016dynamic}, many of which are discussed in \cref{tab:summary}. These challenges have been a driving force for the many versions of the DMD algorithm that have appeared.

For example, the KMD in \eqref{linear_KMD} highlights the potential usefulness of DMD in forecasting. In instances where DMD is applied to noise-free data, such as in generating reduced-order models from high-fidelity numerical simulations \citep{kutz2016multiresolution,alla2017nonlinear,lu2020prediction}, DMD proves effective for both reconstruction and accurate forecasting of the solutions. However, practitioners familiar with DMD's performance in noisy conditions recognize its shortcomings; the algorithm often fails to forecast and reconstruct even the time series it was trained on. In particular, the prediction error in \cref{wake1} is somewhat misleading of what a user might expect in the general case. Even after over a decade, the application of DMD for forecasting or reconstructing time-series data remains limited, typically restricted to high-quality, low-noise scenarios. In \cref{sec:noise_DMD_robustness}, we will focus on methods that mitigate the effect of noise in snapshot data. Many of the structure-preserving methods we discuss in \cref{sec:structure_preserving_methods} have an in-built robustness to noise.

Generally speaking, the error of DMD and its approximate KMD can be split into three types:
\begin{itemize}
	\item The \textit{projection error} is due to projecting/truncating the Koopman operator onto a finite-dimensional space of observables. This is linked to the issue of closure and lack of (or lack of knowledge of) non-trivial finite-dimensional Koopman invariant subspaces.
	\item The \textit{estimation error} is due to estimating the matrices that represent the projected Koopman operator from a finite set of potentially noisy trajectory data.
	\item \textit{Numerical errors} (e.g., roundoff, stability, further compression, etc.) incurred when processing the finite DMD matrix.
\end{itemize}
In particular, \cite{wu2021challenges} highlight the issues of robustness to noise and closure/projection errors as the two fundamental challenges for DMD methods. In \cref{sec:Galerkin_perspective} we will consider methods that directly connect DMD with Koopman operators through the Galerkin perspective. Ways of controlling and measuring the projection error are discussed in \cref{sec:ResDMD}.

DMD's primary value has been as a diagnostic tool, and the interpretability of DMD modes and frequencies is crucial to this role. Most DMD papers focus on analyzing DMD modes and eigenvalues. This emphasis shapes much of this review. The KMD approximated by DMD modes and eigenvalues facilitates dimensionality reduction and model simplification, analogous to classical methods like the Fourier transform or SVD \citep{brunton2022data}. There are numerous software packages for DMD methods, including \textcolor[rgb]{0,0,1}{https://github.com/dynamicslab/pykoopman}, \textcolor[rgb]{0,0,1}{https://github.com/mathLab/PyDMD}, and \textcolor[rgb]{0,0,1}{https://github.com/decargroup/pykoop}. Moreover, there are numerous repositories connected to the papers cited below. Drmač has implemented the DMD algorithm and extensions in LAPACK \citep{drmac2022lapack,drmac2022lapack2}.

\section{Variants from the Regression Perspective}
\label{sec:regression_variants}

This section gives the reader a flavor of DMD variants from the regression perspective.\footnote{As we shall see, there is less convergence theory (e.g., in the large data limit $M\rightarrow\infty$ or as the number of observables increases) for DMD methods based on this viewpoint than for those based on the Galerkin viewpoint in \cref{sec:Galerkin_perspective}. This is due to a looser connection with Koopman operators.} We focus on four key aspects that have proved influential over the last decade or so:
\begin{itemize}
	\item Noise reduction;
	\item Compression and randomized linear algebra;
	\item Multiscale dynamics; and
	\item Control.
\end{itemize}
The methods we discuss are only some of the variants - it is impossible to do justice to the breadth of techniques! Notable omissions include the following.
\textit{Bayesian DMD} \citep{takeishi2017bayesian} transfers the Bayesian formulation into DMD.
\textit{Higher order DMD} \citep{le2017higher} applies time-delay embedding to build a larger state space after projecting onto POD modes.
\textit{Parametric DMD} \citep{huhn2023parametric} performs DMD independently per parameter realization and interpolates the resulting Koopman operators. See also \citep{andreuzzi2023dynamic}.
\textit{Refined Rayleigh--Ritz data driven modal decomposition} \citep{drmac2018data} produces refined Ritz pairs of finite DMD matrix $\Kv_{\mathrm{DMD}}$.
\textit{Spatio-temporal Koopman decomposition} \citep{clainche2018spatio} approximates spatio-temporal data as a linear combination of (possibly growing or decaying exponentially) standing or traveling waves.
\cite{klus2018tensor,klus_towards} develop tensor-based DMD methods for computing eigenfunctions of the Koopman operator. For example, \textit{tensor-based DMD} exploits low-rank tensor decompositions of the data matrices to improve efficiency and memory use. There are extensions of this approach based on reproducing kernel Hilbert spaces (RKHSs) \citep{fujii2019dynamic} and EDMD  \citep{nuske2021tensor}. Recent work has also explored connections between DMD and tensor factorizations \citep{redman2021koopman}.

\subsection{Increasing robustness to noise}
\label{sec:noise_DMD_robustness}

A challenge of DMD is that the computed eigenvalues are biased in the presence of sensor noise. Noise typically dampens the eigenvalues, meaning that for discrete-time systems, the absolute values of the eigenvalues are decreased, and the Koopman modes become distorted. For studies of this effect in various physical systems, see \citep{duke2012error,bagheri2014effects,pan2015accuracy}. \cite{dawson2016characterizing} provide an exceptionally clear discussion of this topic. The bias occurs because standard algorithms treat the data ``snapshot to snapshot'' rather than as a whole and favor one direction (forward in time). Several variants of DMD aim to address this bias. In addition to the methods presented below, other techniques include utilizing Kalman filters \citep{nonomura2018dynamic,nonomura2019extended,jiang2022correcting}, adapting DMD to online data \citep{hemati2014dynamic,hemati2016improving}, robust principal component analysis \citep{scherl2020robust}, and using a second set of noisy observables that meet some independence requirements \citep{wanner2022robust}. Moreover, the structure-preserving methods we discuss in \cref{sec:structure_preserving_methods} often have an inbuilt robustness to noise. Finally, in \cref{sec:stochastic_systems}, we discuss the stochastic Koopman operator, which can handle both system and sensor noise.

\subsubsection{The problem of noise}
We can understand the bias often encountered in DMD as follows. Assume that the snapshots come with additive sensor noise that affects only our measurements of a given system and does not interact with the true dynamics. This means that we have access to noisy data matrices
$$
\Xv_s=\Xv +\Nv_X,\quad \Yv_s=\Yv +\Nv_Y,
$$
where $\Nv_X$ and $\Nv_Y$ are random matrices representing sensor noise, and $\Xv$ and $\Yv$ are the noise-free snapshots. We then represent the data in a truncated POD mode basis so that
$$
\tilde{\Xv}_s=\tilde{\Xv} +\tilde{\Nv}_X,\quad \tilde{\Yv}_s=\tilde{\Yv} +\tilde{\Nv}_Y,
$$
and assume that a subset of POD modes has been selected so that $\tilde{\Xv}_s\tilde{\Xv}_s^*$ is invertible. Assuming the noise is sufficiently small, the DMD matrix can be expanded as
\begin{align*}
\tilde{\Kv}_{\mathrm{DMD}}&=\tilde{\Yv}_s\tilde{\Xv}_s^\dagger=\tilde{\Yv}_s\tilde{\Xv}_s^*(\tilde{\Xv}_s\tilde{\Xv}_s^*)^{-1}\\
&=(\tilde{\Yv} +\tilde{\Nv}_Y)(\tilde{\Xv} +\tilde{\Nv}_X)^*\left[(\tilde{\Xv} +\tilde{\Nv}_X)(\tilde{\Xv} +\tilde{\Nv}_X)^*\right]^{-1}\\
& = (\tilde{\Yv} +\tilde{\Nv}_Y)(\tilde{\Xv} +\tilde{\Nv}_X)^*(\tilde{\Xv}\tilde{\Xv}^*)^{-1}\left[\Iv-( \tilde{\Nv}_X\tilde{\Xv}^*+\tilde{\Xv}\tilde{\Nv}_X^*+\tilde{\Nv}_X\tilde{\Nv}_X^*)(\tilde{\Xv}\tilde{\Xv}^*)^{-1}+\cdots\right].
\end{align*}
\cite{dawson2016characterizing} discard high-order terms in the expectation of this expansion to arrive at
\begin{equation}
\label{DMD_bias}
\mathbb{E}(\tilde{\Kv}_{\mathrm{DMD}})\approx\tilde{\Yv}\tilde{\Xv}^{-1}(\Iv - \mathbb{E}( \tilde{\Nv}_X\tilde{\Nv}_X^* ) (\tilde{\Xv}\tilde{\Xv}^*)^{-1} ).
\end{equation}
This indicates that DMD has an inherent bias due to sensor noise, causing a dampening effect. Interestingly, this bias depends only on $\tilde{\Nv}_X$ and not $\tilde{\Nv}_Y$. The reason is that the least squares problem in \eqref{DMD_opt_vanilla} is optimal only when assuming that all of the noise is in $\tilde{\Yv}$, but not in $\tilde{\Xv}$. Another way of seeing this is that the expression $\tilde{\Yv}\tilde{\Xv}^{-1}$ is linear in $\tilde{\Yv}$, but not in $\tilde{\Xv}$, which is why perturbations to $\tilde{\Xv}$ do not have to propagate through the equation in an unbiased manner.

If the noise structure is known, DMD can be adjusted using a method called \textit{noise-corrected DMD} (ncDMD) \citep{dawson2016characterizing}. However, it is preferable to have methods that correct for noise without requiring explicit knowledge of its structure. We will now outline three popular DMD variants that address this bias without specific assumptions about the noise. The first two can be executed directly using SVDs. The final method requires an iterative method for solving an optimization problem and is more expensive yet more robust.

\subsubsection{Forward-Backward Dynamic Mode Decomposition (fbDMD)}

\begin{algorithm}[t]
\textbf{Input:} Snapshot data $\Xv\in\mathbb{C}^{d\times M}$ and $\Yv\in\mathbb{C}^{d\times M}$, rank $r\in\mathbb{N}$. \\
\vspace{-4mm}
\begin{algorithmic}[1]
\State Compute a truncated SVD $\Xv \approx \Uv \mathbf{\Sigma}\Vv^*, \Uv\in\mathbb{C}^{d\times r},\mathbf{\Sigma}\in\mathbb{R}^{r\times r},\Vv\in\mathbb{C}^{M\times r}.$
\State Compute the projected data matrices $\tilde{\Xv}=\Uv^*\Xv$, $\tilde{\Yv}=\Uv^*\Yv$ and their economized SVDs $\tilde{\Xv}=\Uv_X \mathbf{\Sigma}_X\Vv_X^*$, $\tilde{\Yv}=\Uv_Y \mathbf{\Sigma}_Y\Vv_Y^*$.
\State Compute the forward and backward matrices $\tilde{\Kv}_{f}={\mathbf U}_{ X}^* {\tilde{\Yv}} {\mathbf V}_{ X} {\boldsymbol \Sigma}_{ X}^{-1},
\tilde{\Kv}_{b}={\mathbf U}_{ Y}^* {\tilde{\Xv}} {\mathbf V}_{ Y} {\boldsymbol \Sigma}_{ Y}^{-1}.
$
\State Compute the matrices $\Sv_f = \tilde{\Yv}\Vv_X\mathbf{\Sigma}_X^{-1}, \Sv_b = \tilde{\Xv}\Vv_Y\mathbf{\Sigma}_Y^{-1}$, and ${\Kv}_{f}=\Sv_f\tilde{\Kv}_{f}\Sv_f^\dagger,{\Kv}_{b}=\Sv_b\tilde{\Kv}_{b}\Sv_b^\dagger.$
\State Compute the DMD matrix $\tilde{\Kv} = \left ({\Kv}_{f} {\Kv}_{b}^{-1} \right)^{1/2}$ and its eigendecomposition $\tilde{\Kv}\Wv=\Wv\mathbf{\Lambda}.$
\State Compute the modes $\mathbf{\Phi}=\Yv\Vv\mathbf{\Sigma}^{-1}\Wv.$
\end{algorithmic} \textbf{Output:} The eigenvalues $\mathbf{\Lambda}$ and modes $\mathbf{\Phi}\in\mathbb{C}^{d\times r}$.
\caption{Forward-backward DMD \citep{dawson2016characterizing}.}
\label{alg:fbDMD}
\end{algorithm}

\textit{Forward-Backward DMD} (fbDMD) can be considered a correction to the unidirectional bias of \cref{alg:DMD_vanilla} \citep{dawson2016characterizing}. Let ${\mathbf X}={\mathbf U}_{ X}{\boldsymbol \Sigma}_{ X} {\mathbf V}^*_{ X}$ and ${\mathbf Y}={\mathbf U}_{ Y} {\boldsymbol \Sigma}_{ Y} {\mathbf V}^*_{ Y} $ be truncated SVDs of the matrices ${\mathbf X}$ and ${\mathbf Y}$, respectively. We define
$$
\tilde{\Kv}_{f}={\mathbf U}_{ X}^* {\mathbf Y} {\mathbf V}_{ X} {\boldsymbol \Sigma}_{ X}^{-1},\quad 
\tilde{\Kv}_{b}={\mathbf U}_{ Y}^* {\mathbf X} {\mathbf V}_{ Y} {\boldsymbol \Sigma}_{ Y}^{-1},
$$
which represent forward and backward propagators for the data, analogous to \cref{alg:DMD_vanilla}. Assuming the system's dynamics are invertible and $\tilde{\Kv}_{b}$ is also invertible, the matrix
$$
\tilde{\Kv} = \left (\tilde{\Kv}_{f} \tilde{\Kv}_{b}^{-1} \right)^{1/2}
$$
provides a debiased estimate of the forward propagator. The method is presented in \cref{alg:fbDMD}. Nonetheless, caution is required due to the nonuniqueness of the matrix square root \citep{higham2008functions}. \cite{dawson2016characterizing} suggest selecting the square root that is closest to $\tilde{\Kv}_{f}$ in norm, although this can be computationally costly. A more economical alternative involves measuring closeness in the computed eigencoordinates. Sometimes, the nonuniqueness can be avoided. For instance, if the samples are snapshots from a continuous system whose signal has a bandwidth of $\lambda_B$ and the time-step satisfies $\Delta t < \pi/(2\lambda_B)$, then the discrete eigenvalues expected to be recovered will have a positive real part, which resolves the ambiguity mentioned previously. The square root issue is further analyzed in \cite[Section 5.4]{drmac2018data}. Finally, \cite{askham2018variable} recommend first projecting onto $r$ POD modes before applying fbDMD, an alteration that has demonstrated superior performance in practice. For a variational problem involving forward and backward dynamics, see \textit{Consistent DMD }\citep{azencot2019consistent}.

\subsubsection{Total Least-Squares Dynamic Mode Decomposition (tlsDMD)}
\label{sec:tlsDMD}

\begin{algorithm}[t]
\textbf{Input:} Snapshot data $\Xv\in\mathbb{C}^{d\times M}$ and $\Yv\in\mathbb{C}^{d\times M}$, rank $r\in\mathbb{N}$. \\
\vspace{-4mm}
\begin{algorithmic}[1]
\State Compute a truncated SVD $\Xv \approx \Uv \mathbf{\Sigma}\Vv^*, \Uv\in\mathbb{C}^{d\times r},\mathbf{\Sigma}\in\mathbb{R}^{r\times r},\Vv\in\mathbb{C}^{M\times r}.$
\State Compute the projected data matrices $\tilde{\Xv}=\Uv^*\Xv$, $\tilde{\Yv}=\Uv^*\Yv$.
\State Form the matrix
$
\Zv = \begin{pmatrix} \tilde{\Xv}^\top & \tilde{\Yv} ^\top
  \end{pmatrix}^\top
$
and compute its reduced SVD $\Zv=\Uv_{ Z} \mathbf{\Sigma}_{ Z} \Vv_{ Z}^* $.
\State Set $\Uv_1=\Uv_{ Z}(1 :  r,1 :  r)$, $\Uv_2=\Uv_{ Z}(r+1 :  2r,1 :  r)$.
\State Compute the DMD matrix $\tilde{\Kv} =\Uv_2 \Uv_1^{-1}$ and its eigendecomposition $\tilde{\Kv}\Wv=\Wv\mathbf{\Lambda}.$
\State Compute the modes $\mathbf{\Phi}=\Yv\Vv\mathbf{\Sigma}^{-1}\Wv.$ 
\end{algorithmic} \textbf{Output:} The eigenvalues $\mathbf{\Lambda}$ and modes $\mathbf{\Phi}\in\mathbb{C}^{d\times r}$.
\caption{Total least-squares DMD \citep{dawson2016characterizing,hemati2017biasing}.}
\label{alg:tlsDMD}
\end{algorithm}

\textit{Total Least-Squares DMD} (tlsDMD) addresses the asymmetric treatment of noise in $\Xv$ and $\Yv$ by \cref{alg:DMD_vanilla}. The least-squares problem in \eqref{DMD_opt_vanilla} can be formulated as
$$
\min_{\Kv} \|\Ev_Y\|_\mathrm{F}\quad \text{such that}\quad \Yv+ \Ev_Y=\Kv\Xv.
$$
Considering the reverse time direction, as in fbDMD, leads to the problem
$$
\min_{\Kv} \|\Ev_X\|_\mathrm{F}\quad \text{such that}\quad\Yv=\Kv(\Xv + \Ev_X).
$$
While fbDMD accounts for both directions of error, a more direct approach utilizes the total least-squares problem \citep{van1991total}:
$$
\min_{\Kv}\left\|\begin{pmatrix}\Ev_X\\
\Ev_Y\end{pmatrix}\right\|_\mathrm{F}\quad \text{such that}\quad \Yv+\Ev_Y=\Kv(\Xv + \Ev_X).
$$
This problem can be solved via an SVD, and we follow the version presented by \cite{dawson2016characterizing}, which is similar in spirit to that of \cite{hemati2017biasing}. First, we project ${\mathbf X}$ and ${\mathbf Y}$ onto $r < M/2$ POD modes to obtain $\tilde{{\mathbf X}}$ and $\tilde{{\mathbf Y}}$. We then define
$$
\Zv = \begin{pmatrix} \tilde{\Xv} \\ \tilde{\Yv} 
  \end{pmatrix}
$$
and compute its reduced SVD $\Zv=\Uv_{ Z} \mathbf{\Sigma}_{ Z} \Vv_{ Z}^* $. The matrix $\tilde{\Kv} = \Uv_{ Z}(r+1 :  2r,1 :  r)\Uv_{ Z}(1 :  r,1 :  r)^{-1}$ then provides a debiased estimate of the forward propagator. The method is summarized in \cref{alg:tlsDMD}.

\subsubsection{Optimized Dynamic Mode Decomposition (optDMD)}

\textit{Optimized DMD} (optDMD) is a variation of DMD that processes all data snapshots collectively \citep{chen2012variants}. This approach reduces much of the bias associated with exact DMD. Nonetheless, it necessitates solving a nonlinear optimization problem, initially thought to hinder its practical application. However, \cite{askham2018variable} demonstrated that an approximate solution to the optimization problem can be efficiently computed using the variable projection method \citep{golub1973differentiation}. In this framework, DMD is reformulated as an exponential data fitting problem \citep{pereyra2010exponential}, which brings an additional advantage: the data snapshots do not have to be equidistant in time. For further DMD methodologies tailored for data with irregular time intervals, see \citep{tu2014spectral,gueniat2015dynamic,leroux2016dynamic}.

Initially, we project onto $r$ POD modes to construct the data matrix $\tilde{\Xv}=[\zv_0 \hspace{1mm}\zv_1\cdots \zv_{M}]$, which corresponds to the projected data at times $t_0,t_1,\ldots, t_{M}$. Depending on the data structure, the projected matrix $\Yv$ may also be incorporated into this matrix. We posit that the data represents the solution to a linear system of differential equations, expressed as
$$
  {\mathbf z}(t) \approx {\mathbf S} e^{{\boldsymbol \Lambda} t} {\mathbf S}^\dagger {\mathbf z}_0 \; ,
$$
where ${\mathbf S} \in \mathbb{C}^{r\times r}$ and ${\boldsymbol \Lambda} \in \mathbb{C}^{r\times r}$. This representation can be reformulated to
$$
\tilde{\Xv}^\top \approx {\boldsymbol \Phi}({\boldsymbol \alpha}) {\mathbf B},\quad { \Bv}_{i,j} = { \Sv}_{j,i} \left ({\mathbf S}^\dagger {\mathbf z}_0 \right)_i,
$$
where ${\boldsymbol \Phi}({\boldsymbol \alpha}) \in \mathbb{C}^{(m+1)\times r}$, whose elements are ${ \Phi}({\boldsymbol \alpha})_{i,j} = \mathrm{exp} ({\boldsymbol \alpha}_j t_i)$. From this, we arrive at an exponential fitting problem:
$$
  \min_{\boldsymbol\alpha\in\mathbb{C}^r,\Bv\in\mathbb{C}^{r\times r}} \left\| \tilde{\Xv}^\top - {\boldsymbol \Phi}({\boldsymbol \alpha}) {\mathbf B} \right\|_{\mathrm{F}}.
$$
The optimized DMD eigenvalues are determined by $\lambda_j = {\boldsymbol \alpha}_j$. This optimization problem is solved using the variable projection method, which exploits the specific structure of the exponential data fitting problem to eliminate many of the variables from the optimization process. A summary is provided in \cref{alg:optDMD}, with practical details given in \citep{askham2018variable}, including strategies for selecting the initial guess (e.g., employing an alternate DMD algorithm).

While this nonlinear, nonconvex optimization problem is not guaranteed to be solved globally, and the method may be computationally intensive due to its iterative nature, optDMD has been proven to yield significant enhancements over traditional DMD approaches. Moreover, optDMD's efficacy can be further heightened by employing Breiman's statistical bagging sampling strategy \citep{breiman2017classification}, which assembles a collection of models to reduce model variance, mitigate overfitting, and facilitate uncertainty quantification. This augmented method is referred to as \textit{bagging optimized DMD} (bopDMD) \citep{sashidhar2022bagging}.

\begin{algorithm}[t]
\textbf{Input:} Snapshot data $\Xv\in\mathbb{C}^{d\times (M+1)}$, rank $r\in\mathbb{N}$, initial guess for $\boldsymbol \alpha$. \\
\vspace{-4mm}
\begin{algorithmic}[1]
\State Compute a truncated SVD $\Xv \approx \Uv \mathbf{\Sigma}\Vv^*, \Uv\in\mathbb{C}^{d\times r},\mathbf{\Sigma}\in\mathbb{R}^{r\times r},\Vv\in\mathbb{C}^{(M+1)\times r}.$
\State Compute the projected data matrix $\tilde{\Xv}=\Uv^*\Xv$.
\State Solve the problem
$$
  \min_{\boldsymbol\alpha\in\mathbb{C}^r,\Bv\in\mathbb{C}^{r\times r}} \left\| \tilde{\Xv}^\top - {\boldsymbol \Phi}({\boldsymbol \alpha}) {\mathbf B} \right\|_{\mathrm{F}}
$$
using a variable projection algorithm.
\State Set $\lambda_j = {\boldsymbol \alpha}_j$ and $\mathbf{\Phi}(:,i) = ({\|\Uv\Bv^\top( : ,i) \|_{\ell^2}})^{-1}
      \Uv\Bv^{\top}( : ,i)$.
\end{algorithmic} \textbf{Output:} The eigenvalues $\mathbf{\Lambda}$ and modes $\mathbf{\Phi}\in\mathbb{C}^{d\times r}$.
\caption{Optimized DMD, algorithmic details are given in \citep{askham2018variable}.}
\label{alg:optDMD}
\end{algorithm}

\subsubsection{Examples}
\label{sec:noise_examples}

\begin{figure}
\centering
\raisebox{-0.5\height}{\includegraphics[width=0.45\textwidth,trim={0mm 0mm 0mm 0mm},clip]{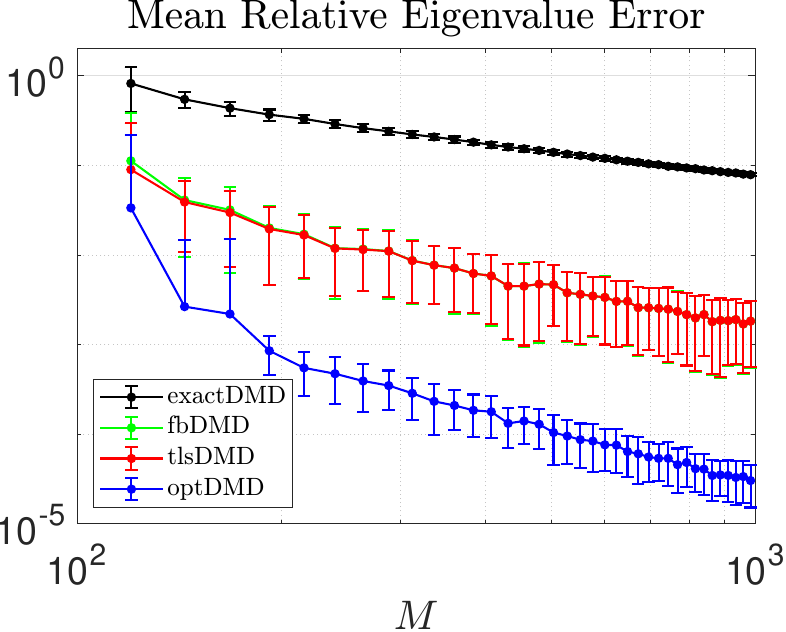}}\hfill
\raisebox{-0.5\height}{\includegraphics[width=0.45\textwidth,trim={0mm 0mm 0mm 0mm},clip]{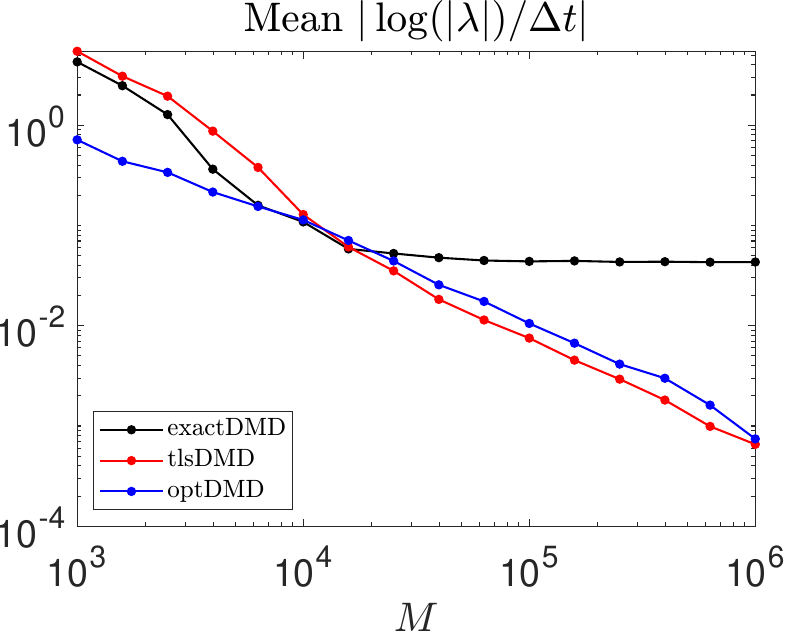}}
\caption{Left: Mean error (error bars correspond to the standard deviation across noise realizations) in the first 11 eigenvalues of the cylinder example. Right: Mean value of $|\log(|\lambda|)/\Delta t|$ for the DMD eigenvalues for the Lorenz system. We have not shown the results for fbDMD since they are almost identical to tlsDMD.}
\label{fig_denoised}
\end{figure}

For simplicity, we focus on the error associated with the approximated eigenvalues. Other error metrics related to the modes or the accuracy of the decomposition in fitting the data or forecasting are also frequently considered in the literature, often yielding similar results.

\paragraph{Noisy cylinder wake}

We revisit the example of flow past a cylinder from \cref{sec:example_cylinder}. We center and normalize the data grid-wise before adding 40\% Gaussian random noise to the measurements. \cref{fig_denoised} (left) shows the mean relative $\ell^2$ error of the first 11 eigenvalues (see \cref{wake1}), averaged over 100 realizations of random noise. The errors are calculated by comparison with eigenvalues computed from noise-free snapshots that have converged in terms of both the size of the truncated SVD and the number of snapshots. The error bars represent one standard deviation from the mean. All methods exhibit a decreasing error as $M$ increases, which is largely attributable to the truncation in the SVD used in DMD. As often noted in the literature, the fbDMD and tlsDMD methods perform comparably. However, optDMD demonstrates a significantly smaller error.

\paragraph{Lorenz system}

We revisit the Lorenz system example from \cref{sec:intro_lorenz}. Since the spectrum is continuous (apart from the trivial eigenvalue $\lambda=1$), measuring the error of individual DMD eigenvalues is meaningless unless methods such as the residual in \cref{sec:ResDMD} are used. However, $|\log(|\lambda|)/\Delta t|$ vanishes on the spectrum of the Koopman operator. Therefore, we select $N=10$ and compute the mean value of $|\log(|\lambda|)/\Delta t|$ over the DMD eigenvalues. \cref{fig_denoised} (right) shows the results, averaged over 50 randomly selected initial conditions on the attractor for the initial value $\xv(0)$. For exact DMD, this error metric plateaus as $M$ increases. Generally, the eigenvalues computed using DMD with delay embedding are damped and lie strictly within the unit disk \citep[Corollary 2]{korda2020data}. Consequently, their logarithms are in the left-half plane, corresponding to positions below the horizontal line in \cref{lorenz1}. Conversely, the eigenvalues computed by tlsDMD and optDMD approach the unit disk with increasing $M$ and exhibit greater robustness to noise in the measurements.

\subsection{Compression and randomized linear algebra}
\label{sec:compression_random}

With ever-increasing volumes of measurement data from simulations and experiments, modal extraction algorithms such as DMD can become prohibitively expensive, particularly for online or real-time analysis. Dynamics often evolve on low-dimensional attractors, indicating sparsity in a suitable coordinate system or an intrinsic low rankness. However, the SVD used in \cref{alg:DMD_vanilla} scales with the dimension of the measurements, not with the intrinsic dimension of the data. This section explores two principles aimed at mitigating this computational cost:
\begin{itemize}
	\item \textbf{Compressed sensing} \citep{donoho2006compressed,candes2006robust} facilitates the reconstruction of sparse signals from a limited number of measurements, allowing for undersampling below traditional Shannon--Nyquist limits \citep{nyquist1928certain,shannon1948mathematical}. Applying compressed sensing to DMD can substantially improve computational efficiency, particularly during the SVD step of the algorithm. Acquiring high-resolution, time-resolved measurements can be challenging. Nevertheless, temporally and spatially sparse signals may be sampled less frequently than traditionally expected, which is crucial if data acquisition is costly.
	\item \textbf{Randomized numerical linear algebra} \citep{martinsson2020randomized}  offers a way to solve certain linear algebra problems much faster than classical methods. The randomized SVD is a fast and straightforward technique for computing an approximate low-rank SVD \citep{halko2011finding}. It is robust and amenable to parallelization and can benefit from GPU architectures. When coupled with randomized SVD, DMD scales with the intrinsic rank of the data matrices rather than the measurement dimension. The approximation error is manageable through oversampling and power iterations, providing a balance between computational speed and accuracy. Moreover, it can accommodate large datasets that exceed the capacity of fast memory by using a blocked matrix approach.
\end{itemize}

Beyond the methods detailed here, several DMD variants are based on related principles. For instance, due to the non-orthogonality of DMD modes, choosing an appropriate low-rank representation can be difficult \citep{kou2017improved}. \textit{Sparsity-Promoting DMD} \citep{jovanovic2014sparsity} aims to strike a balance between accuracy and the number of modes by identifying a sparse subset of modes. Other techniques for selecting dominant modes include ranking each DMD mode's importance by time integration \citep{kou2017improved} or by assessing the time-averaged modal energy contribution \citep{tissot2014model}. Furthermore, one can apply DMD recursively to achieve orthogonality \citep{noack2016recursive}, a method termed \textit{Recursive DMD}, which blends the principles of POD and DMD. Additionally, regularization terms can be imposed to encourage sparsity in the Koopman matrix \citep{sinha2019computation}.

Finally, it is crucial to recognize that the usefulness of the methods in this section presumes the dynamics are evolving on a low-dimensional subspace characterized by a quickly decaying singular value spectrum. While not a fundamental limitation of DMD, this is a common underlying assumption which may not hold for all dynamical systems. \cite{erichson2019randomized} provide a turbulent flow example that demonstrates the limits of the approaches in this section when this assumption does not hold.

\subsubsection{Compressed Sensing meets DMD (cDMD and csDMD)}

A full description of the extensive field of compressed sensing is beyond the scope of this review. We outline the key points to understand its interplay with DMD. The reader is encouraged to consult the excellent textbooks \citep{foucart2013invitation,adcock2021compressive} for a comprehensive understanding or \citep{candes2008introduction} for a concise introductory tutorial. Compressed sensing is founded on two central principles: \textit{sparsity}, which pertains to the signals of interest, and \textit{incoherence}, which relates to the sensing methodology.

Consider a signal $\xv\in\mathbb{C}^d$ that is approximately sparse in some basis $\Bv\in\mathbb{C}^{d\times d}$, meaning that $\xv=\Bv\zv$, where the vector $\zv$ can be well approximated by a sparse vector. Many natural signals, such as images and audio, are approximately sparse in specific bases like the Fourier or wavelet bases. When we transform an image using Fourier or wavelet transformations, most coefficients are small and can be disregarded while still retaining the quality of the image. We assume that we have access to measurements:
$$
\xv_c=\Cv \xv = \Cv \Bv \zv,
$$
where $\Cv\in \mathbb{C}^{p\times d}$ is a measurement matrix with $p<d$. Compressed sensing theory implies that, under suitable conditions, we can recover an accurate approximation of $\zv$ (and hence $\xv$) from the subsampled measurements $\xv_c$. For example, consider the $\ell^1$-minimization problem
\begin{equation}
\label{eq:CS_min}
\min \|\zv\|_{\ell^1}\quad\text{subject to}\quad \xv_c=\Cv \Bv \zv.
\end{equation}
Specifically, the measurement matrix $\Cv$ must be incoherent with respect to the sparse basis $\Bv$, meaning that the rows of $\Cv$ are uncorrelated with the columns of $\Bv$. If the matrix $\Cv\Bv$ satisfies a restricted isometry property (RIP):\footnote{There are no known large matrices with bounded restricted isometry constants since computing these constants is NP-hard and hard to approximate. Typically, one builds random matrices so that the RIP holds with overwhelming probability. For example, Bernoulli and Gaussian random measurement matrices satisfy the RIP for a generic basis $\Bv$ with high probability \citep{candes2006near}.}
$$
(1-\delta_k)\|\zv\|_{\ell^2}^2\leq \|\Cv\Bv\zv\|_{\ell^2}^2\leq (1+\delta_k)\|\zv\|_{\ell^2}^2\quad\text{for $k$-sparse vectors $\zv$},
$$
then we can prove results about how close solutions of \eqref{eq:CS_min} are to the true $\zv$, how issues such as only approximately numerically solving \eqref{eq:CS_min} affect the solution, robustness to noise, and so forth. Beyond the above $\ell^1$-minimization problem, many successful optimization problems and algorithms approximate their solutions in compressed sensing.

\begin{algorithm}[t]
\textbf{Input:} Snapshot data $\Xv\in\mathbb{C}^{d\times M}$ and $\Yv\in\mathbb{C}^{d\times M}$, rank $r\in\mathbb{N}$, and measurement matrix $\Cv\in\mathbb{C}^{p\times d}$. \\
\vspace{-4mm}
\begin{algorithmic}[1]
\State Compress $\Xv$ and $\Yv$ to $\Xv_c=\Cv\Xv$ and $\Yv_c=\Cv\Yv$.
\State Apply \cref{alg:DMD_vanilla} with input $\Xv_c$ and $\Yv_c$ and outputs $\mathbf{\Lambda}_c$, $\Wv_c$, $\Vv_c$ and $\mathbf{\Sigma}_c$.
\State Reconstruct full-state modes via $\mathbf{\Phi}=\Yv\Vv_c\mathbf{\Sigma}^{-1}_c\Wv_c$.
\end{algorithmic} \textbf{Output:} The eigenvalues $\mathbf{\Lambda}_c$ and DMD modes $\mathbf{\Phi}\in\mathbb{C}^{d\times r}$.
\caption{Compressed DMD \citep{brunton2016compressed}, suitable when given access to the full snapshots.}
\label{alg:compressed_DMD}
\end{algorithm}

\begin{algorithm}[t]
\textbf{Input:} Compressed snapshot data $\Xv_c\in\mathbb{C}^{p\times M}$ and $\Yv_c\in\mathbb{C}^{p\times M}$, measurement matrix $\Cv\in\mathbb{C}^{p\times d}$, and basis $\Bv\in\mathbb{C}^{d\times d}$. \\
\vspace{-4mm}
\begin{algorithmic}[1]
\State Apply \cref{alg:DMD_vanilla} with input $\Xv_c$ and $\Yv_c$ and outputs $\mathbf{\Lambda}_c$ and $\mathbf{\Phi}_c$.
\State Apply $\ell^1$-minimization \eqref{eq:CS_min} columnwise to reconstruct modes $\mathbf{\Phi}_s\in\mathbb{C}^{d\times r}$.
\State Recover full-state modes via $\mathbf{\Phi}=\Bv\mathbf{\Phi}_s$.
\end{algorithmic} \textbf{Output:} The eigenvalues $\mathbf{\Lambda}_c$ and DMD modes $\mathbf{\Phi}\in\mathbb{C}^{d\times r}$.
\caption{Compressed sensing DMD \citep{brunton2016compressed}, suitable when given access to only compressed data. The $\ell^1$-minimization can be replaced with a plethora of similar minimization problems from the compressed sensing literature.}
\label{alg:cs_DMD}
\end{algorithm}

\cite{tu2014spectral} combine \textit{temporal} compressed sensing with ideas from DMD to recover POD modes. For the remainder of this section, we focus instead on \textit{spatial} compressed sensing, following the methods of \cite{brunton2016compressed}. In essence, \cite{brunton2016compressed} demonstrated that the unitary invariance of the DMD algorithm can be extended to approximate invariance under transformations satisfying a RIP, provided that the data is sparse in a basis that is incoherent with respect to the measurements. For compressed data matrices
$$
\Xv_c=\Cv\Xv,\quad \Yv_c=\Cv \Yv,
$$
there are essentially two approaches, depending on whether one has access to the matrix $\Yv$ or not. \cref{alg:compressed_DMD} illustrates \textit{compressed DMD} (cDMD) (see also \citep{erichson2019compressed}), where one performs the standard DMD algorithm on the compressed data matrices, and then reconstructs the full-state modes using $\Yv$. If access to $\Yv$ is not available, we can use an optimization problem such as \eqref{eq:CS_min} to recover the modes in the sparse basis $\mathbf{\Phi}_s$, and then reconstruct the full-state modes. This approach, known as \textit{compressed sensing DMD} (csDMD), is outlined in \cref{alg:cs_DMD}.

\subsubsection{Randomized Dynamic Mode Decomposition (rDMD)}

Early uses of DMD with randomized SVD include \citep{erichson2016randomized}, who utilized it to expedite DMD applications in video background subtraction, and \citep{bistrian2017randomized}, who applied it as a component of a reduced-order model for two-dimensional fluid flows. Although this method is reliable and robust to noise, it only accelerates the computation of the SVD, with subsequent computational steps in the DMD algorithm remaining costly. Instead, \cite{erichson2019randomized} developed a \textit{randomized DMD} (rDMD) algorithm. This algorithm relies on sketching the range of $\Xv$ and executing the entire DMD process in a reduced-dimensional space, ultimately recovering the DMD of the original system at the end.

\begin{algorithm}[t]
\textbf{Input:} Snapshot data $\Xv\in\mathbb{C}^{d\times M}$, target rank $r\in\mathbb{N}$, oversampling factor $p\in\mathbb{N}$, and power iteration factor $q\in\mathbb{N}\cup\{0\}$. \\
\vspace{-4mm}
\begin{algorithmic}[1]
\State Generate a random Gaussian matrix $\mathbf{\Omega}\in\mathbb{R}^{M\times (r+p)}$ and form the matrix $\Zv=\Xv \mathbf{\Omega}$.
\For{$j=1,\ldots, q$}
        \State $[\Qv,\sim]=\texttt{qr}(\Zv,\text{'econ'})$
				\State $[\Cv,\sim] = \texttt{qr}(\Xv^*\Qv,\text{'econ'})$
				\State $\Zv =\Xv\Cv$
\EndFor
\State $[\Qv,\sim]=\texttt{qr}(\Zv,\text{'econ'})$.
\end{algorithmic} \textbf{Output:} Range matrix $\Qv\in\mathbb{C}^{d\times(r+p)}$.
\caption{Randomized range finder. Other choices of random test matrices may be employed for computational efficiency. The QR algorithm is written using MATLAB notation.}
\label{alg:range_finder}
\end{algorithm}

\begin{algorithm}[t]
\textbf{Input:} Snapshot data $\Xv\in\mathbb{C}^{d\times M}$ and $\Yv\in\mathbb{C}^{d\times M}$, target rank $r\in\mathbb{N}$, oversampling factor $p\in\mathbb{N}$, and power iteration factor $q\in\mathbb{N}\cup\{0\}$. \\
\vspace{-4mm}
\begin{algorithmic}[1]
\State Run \cref{alg:range_finder} to generate the matrix $\Qv$.
\State Compress $\Xv$ and $\Yv$ to $\Xv_c=\Qv^*\Xv$ and $\Yv_c=\Qv^*\Yv$.
\State Apply \cref{alg:DMD_vanilla} with input $\Xv_c$ and $\Yv_c$ and outputs $\mathbf{\Lambda}_c$ and $\mathbf{\Phi}_c$.
\State Reconstruct full-state modes via $\mathbf{\Phi}=\Qv\mathbf{\Phi}_c$.
\end{algorithmic} \textbf{Output:} The eigenvalues $\mathbf{\Lambda}_c$ and DMD modes $\mathbf{\Phi}\in\mathbb{C}^{d\times r}$.
\caption{Randomized DMD \citep{erichson2019randomized}.}
\label{alg:rDMD}
\end{algorithm}

The idea is to use randomness as a computational strategy to find a smaller representation, known as a \textit{sketch}. This smaller matrix sketch can be used to compute an approximate low-rank factorization for the high-dimensional data matrix. rDMD utilizes the off-the-shelf probabilistic framework proposed in the seminal work of \cite{halko2011finding}. Given a target rank $r$, the aim is to compute a near-optimal basis $\Qv\in\mathbb{C}^{d\times r}$ for the input matrix $\Xv$ such that $\Xv\approx \Qv\Qv^*\Xv$. A test matrix $\mathbf{\Omega}\in\mathbb{R}^{M\times r}$ is drawn from a normal Gaussian distribution to sample the range of $\Xv$ via
$$
\Zv = \Xv \mathbf{\Omega}.
$$
To mitigate the $\mathcal{O}(dMr)$ cost of dense matrix multiplication, more sophisticated random test matrices, such as the subsampled randomized Hadamard transform, can also be used, leading to a complexity of 
$\mathcal{O}(dM\log(r))$. The orthonormal basis $\Qv$ is then obtained via QR decomposition of $\Zv$. In practice, we slightly oversample the desired rank $r$ by a constant factor (typically, ten suffices). A second strategy to improve performance involves power iterations \citep{rokhlin2010randomized, gu2015subspace}. Particularly, a slowly decaying singular value spectrum of the input matrix can significantly affect the quality of the approximated basis matrix $\Qv$. Power iterations are employed to preprocess the input matrix to promote a more rapidly decaying spectrum. The sampling matrix obtained is
$$\Zv=(\Xv\Xv^*)^q\Xv\mathbf{\Omega},$$
and as few as $q=2$ power iterations can considerably improve the approximation quality, even when the singular values of the input matrix decay slowly. This procedure is outlined in \cref{alg:range_finder}, and we direct the reader to \citep[Section 11]{martinsson2020randomized} for probabilistic error bounds. With the matrix $\Qv$ in hand, we can perform DMD on the lower-dimensional space, as summarized in \cref{alg:rDMD}. One can also simultaneously sketch the range and corange of $\Xv$. This method, known as \textsl{sketchy DMD}, was proposed by \cite{ahmed2022dynamic}.

\subsubsection{Examples}
\label{sec:examples_rDMD}

\begin{figure}
\centering
\raisebox{-0.5\height}{\includegraphics[width=0.32\textwidth,trim={0mm 0mm 0mm 0mm},clip]{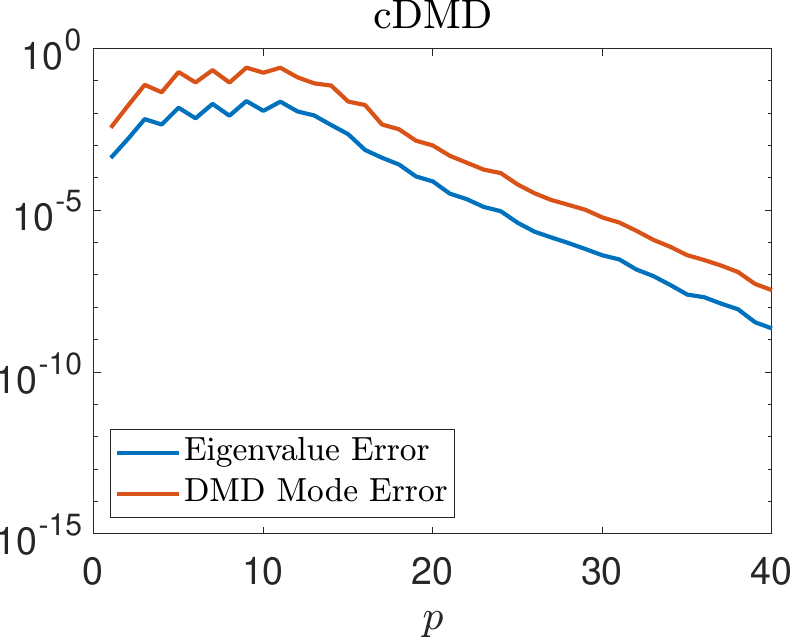}}\hfill
\raisebox{-0.5\height}{\includegraphics[width=0.32\textwidth,trim={0mm 0mm 0mm 0mm},clip]{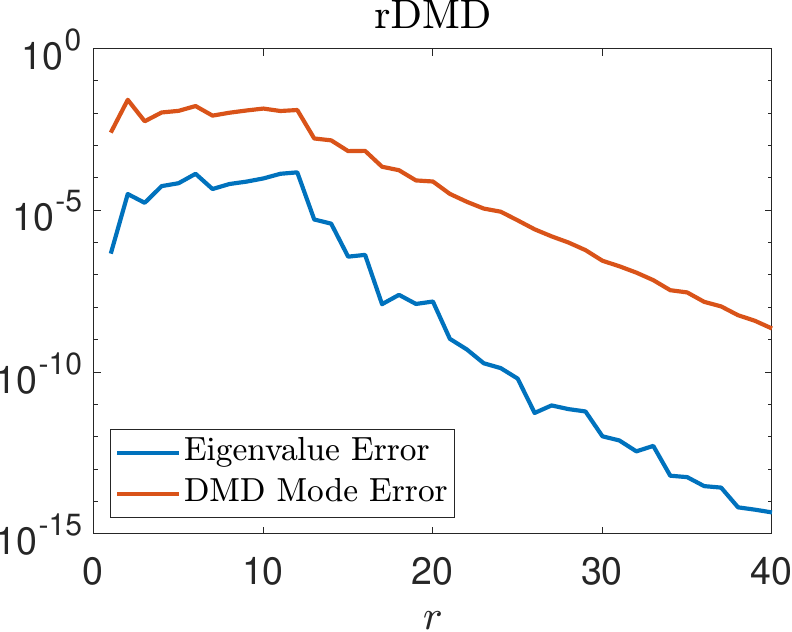}}\hfill
\raisebox{-0.5\height}{\includegraphics[width=0.32\textwidth,trim={0mm 0mm 0mm 0mm},clip]{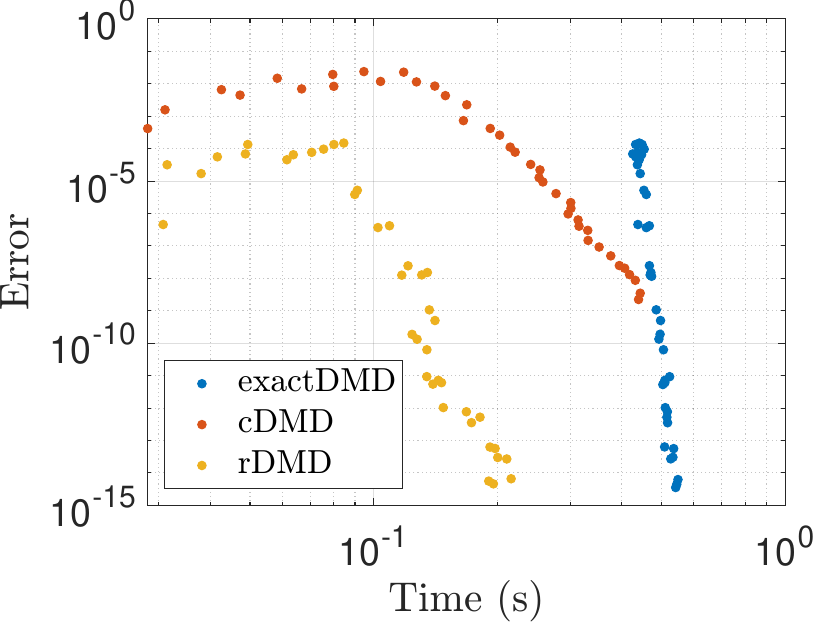}}
\caption{Left: Errors of cDMD (\cref{alg:compressed_DMD}) for the first 6 eigenvalues and DMD modes. Middle:  Errors of rDMD (\cref{alg:rDMD}) for the first six eigenvalues and DMD modes. Right: Computational times vs eigenvalue error.}
\label{compressed_wake1}
\end{figure}

\begin{figure}
\centering
\raisebox{-0.5\height}{\includegraphics[width=1\textwidth,trim={0mm 0mm 0mm 0mm},clip]{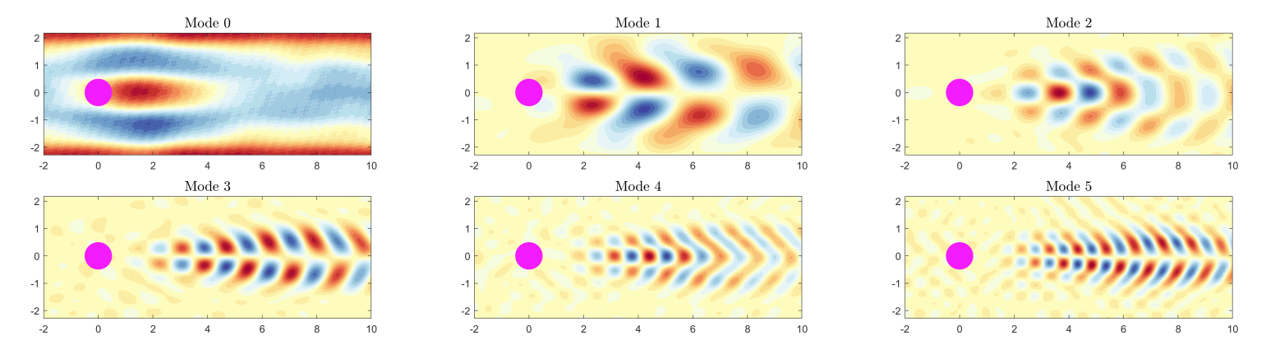}}
\caption{DMD modes constructed using csDMD (\cref{alg:cs_DMD}) with a $99.5\%$ compression in the dimension of the measurements.}
\label{compressed_wake2}
\end{figure}

\paragraph{Cylinder wake}

We return to the cylinder wake discussed in \cref{sec:example_cylinder} as an illustrative example. For cDMD and csDMD, we use the (inverse) two-dimensional discrete Fourier transform as our basis $\mathbf{B}$ and Gaussian random measurements $\mathbf{C}$. \cref{compressed_wake1} shows the results for cDMD and csDMD when recovering the first six modes plotted in \cref{wake1}. In the left and middle panels, we have shown the mean eigenvalue error and the DMD mode error (computed as a subspace angle) averaged over 20 random realizations. The errors quickly become negligible for $p,r=\mathcal{O}(10)$. In the right panel, we have displayed the mean execution times on a standard laptop (without GPU) against the eigenvalue error. Even for this simple example with rapidly decreasing singular values, cDMD performs better than exact DMD, while rDMD is the clear winner.

\cref{compressed_wake2} shows the modes computed by csDMD with $p=800$, which corresponds to a 99.5\% compression in the dimension of the snapshot matrices. We employ the CoSaMP algorithm \citep{needell2009cosamp} to perform the $\ell^1$-minimization step. When only compressed measurements are available, it is still possible to reconstruct full-state modes using compressed sensing. However, this typically requires more measurements and computational resources than cDMD or rDMD.

\paragraph{Sea surface data}
We now consider high-resolution sea surface temperature (SST) data. SST data are widely studied in climate science for climate monitoring and prediction \citep{reynolds2007daily,reynolds2002improved,smith2005global}, and measurements are constructed by combining infrared satellite data with observations provided by ships and buoys. The data are available from the National Oceanic and Atmospheric Administration at \textcolor[rgb]{0,0,1}{https://www.esrl.noaa.gov/psd/} for the years $1981$ to $2023$, with a grid resolution of $0.25^\circ$. Omitting data over land results in $d = 691,150$ spatial grid points. The following experiments, similar in spirit to \citep{erichson2019randomized}, were performed using a system with Intel(R) Xeon(R) Gold 6126 CPU at 2.60GHz (48 cores) and 767GiB system memory.

\begin{figure}
\centering
\raisebox{-0.5\height}{\includegraphics[width=1\textwidth,trim={0mm 0mm 0mm 0mm},clip]{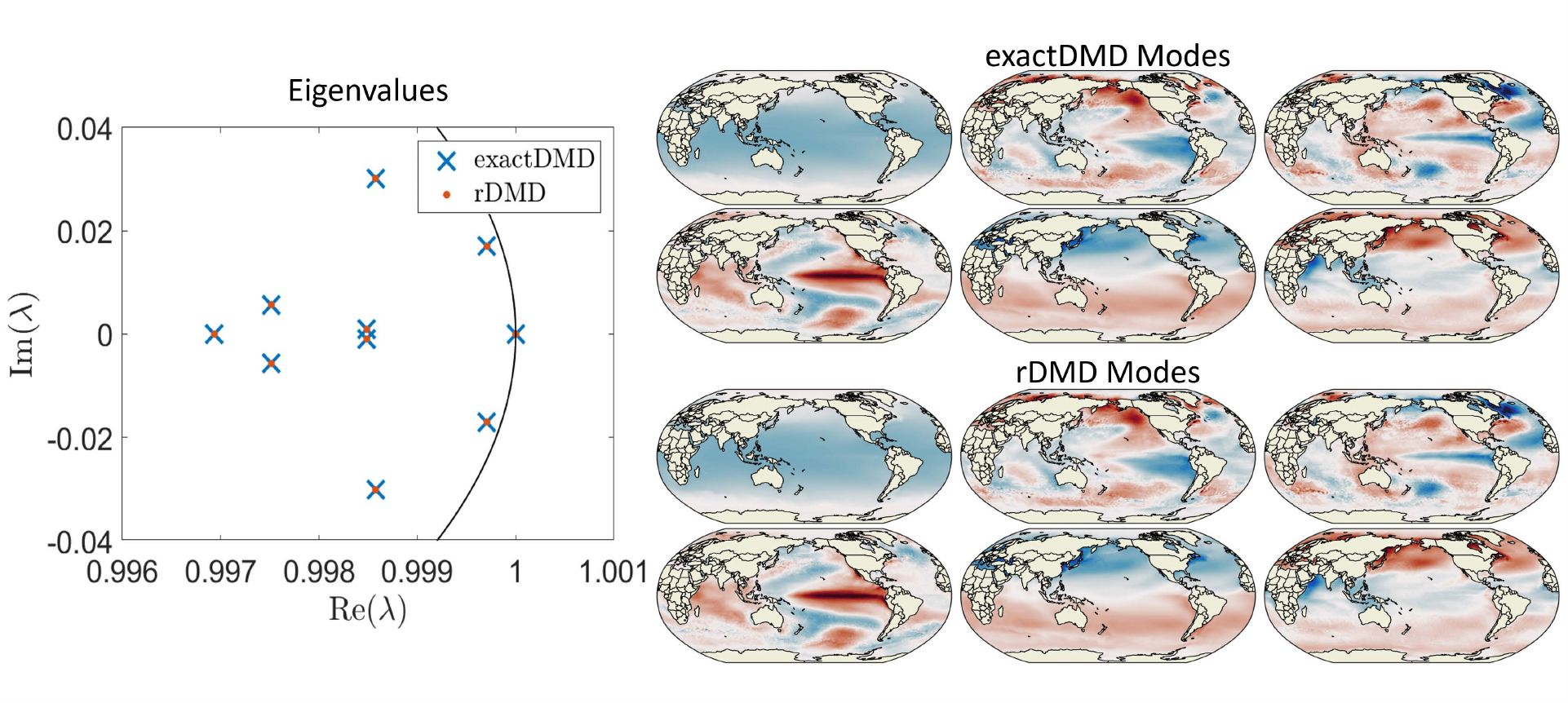}}
\caption{Eigenvalues and dynamic modes of the SST data set computed using exact DMD and rDMD. The eigenvalue plot shows the unit circle as a black line.}
\label{fig_SAT_modes}
\end{figure}

\begin{figure}
\centering
\raisebox{-0.5\height}{\includegraphics[width=0.32\textwidth,trim={0mm 0mm 0mm 0mm},clip]{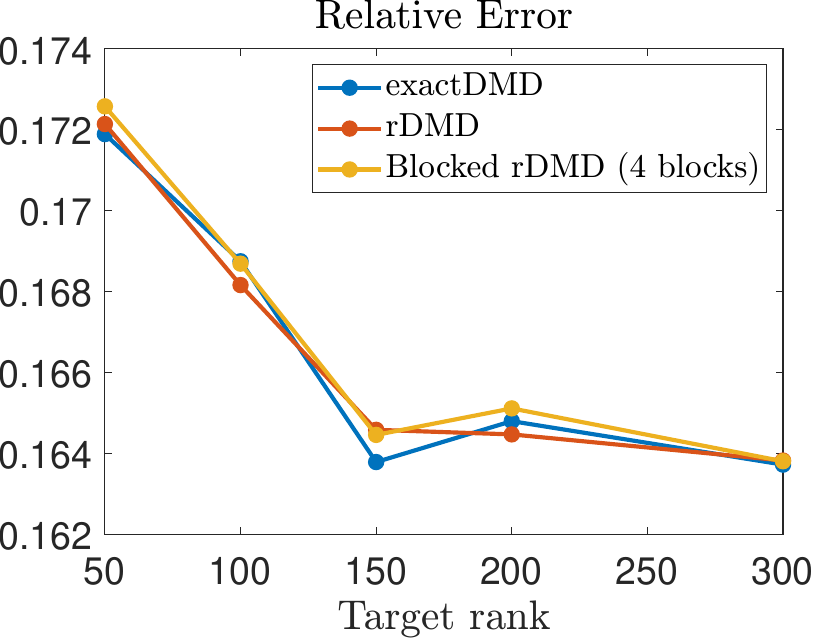}}\hfill
\raisebox{-0.5\height}{\includegraphics[width=0.32\textwidth,trim={0mm 0mm 0mm 0mm},clip]{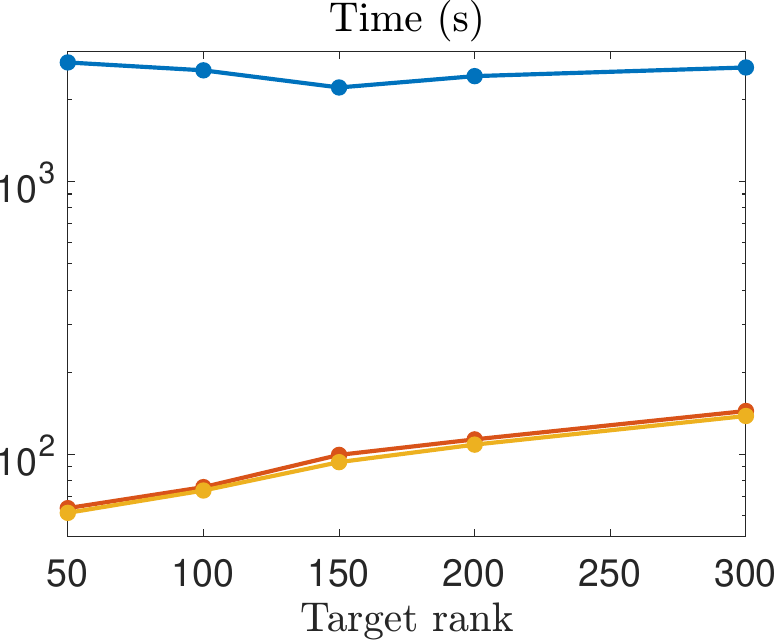}}\hfill
\raisebox{-0.5\height}{\includegraphics[width=0.32\textwidth,trim={0mm 0mm 0mm 0mm},clip]{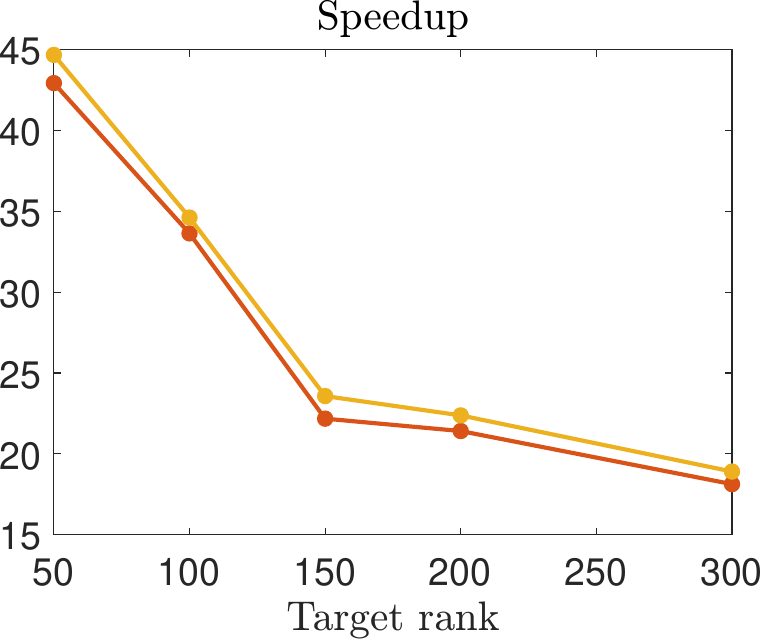}}
\caption{Left: Accuracy of exact DMD, rDMD, and batched rDMD for the SST data set. Middle: the computational ties of each method. Right: the speedup of rDMD and batched rDMD compared to exact DMD. All plots show the mean of 10 independent runs.}
\label{fig_SAT_times}
\end{figure}

We first consider a temporal resolution of one day and a data matrix $\Xv \in  \mathbb{R}^{691150\times 15097}$. \cref{fig_SAT_modes} shows the eigenvalues and dynamic modes computed using exact DMD and rDMD (with $r=10$), demonstrating the accuracy of rDMD. The bottom left mode is reminiscent of an El Niño mode generated from the El Niño-Southern Oscillation (ENSO). El Niño is the warm phase of the ENSO cycle. It is associated with a band of warm ocean water that develops in the central and east-central equatorial Pacific, including off the Pacific coast of South America (see also \cref{fig_mrDMD}).

Next, we compare the accuracy and computational times for a temporal resolution of one week and a data matrix $\Xv \in  \mathbb{R}^{691150\times 2156}$ of weekly averages. \cref{fig_SAT_times} shows the relative error in the Frobenius norm of the reconstructed data matrix and time taken for exact DMD, rDMD, and blocked rDMD (using four blocks). We observe substantial gains in computational time when using rDMD while maintaining an accuracy similar to the full deterministic exact DMD.

\subsection{Multiresolution Dynamic Mode Decomposition (mrDMD)}
\label{sec:mrDMD}

Multiscale systems are widespread across various scientific disciplines. Modeling the interactions between microscale and macroscale phenomena, which may differ by orders of magnitude either spatially or temporally, poses a considerable challenge. Wavelet-based methods and windowed Fourier transforms are well-suited for multiresolution analysis (MRA), as they systematically isolate temporal or spatial features through recursive refinement when sampling from the targeted data \citep{daubechies1992ten}. Typically, MRA is employed separately in either space or time, but it is seldom applied to both simultaneously.

\textit{Multiresolution DMD} (mrDMD) \citep{kutz2016multiresolution} integrates DMD with core principles from wavelet theory and MRA. It adjusts the sampling window of the data collection process in line with wavelet theory, filtering information across various scales. The process is iteratively refined through progressively shorter snapshot sampling windows, leading to the recursive extraction of DMD modes from slow to rapidly changing timescales. The benefits of this approach include enhanced prediction of the near-future state of the system, which is vital for control; effective management of transient phenomena; and improved handling of moving (translating/rotating) structures within the data. The latter two points underscore significant challenges inherent in standard DMD methods. The mrDMD algorithm has led to practical applications such as determining optimal sensor placement \citep{manohar2019optimized}.

\subsubsection{The algorithm}

When using mrDMD, it is typical to work with $M$ such that a full-rank approximation with $r=M$ in \cref{alg:DMD_vanilla} is feasible and such that high- and low-frequency content is present. We assume that data is collected along a single time trajectory with time-step $\Delta t$ and express the eigenvalues in terms of their time-scaled logarithms $\eta = {\log(\lambda)}/{\Delta t}$. In the first pass, mrDMD separates the approximation in \eqref{linear_KMD} into slow modes and fast modes:
\begin{equation}
\label{linear_KMD_mrDMD}
\xv(t)\approx \underbrace{\sum_{|\eta_k|\leq \tau} \pmb{\phi}_k^{(1)}\exp(\eta_kt)\bv(k)}_{\text{slow modes}}+\underbrace{\sum_{|\eta_k|> \tau} \pmb{\phi}_k^{(1)}\exp(\eta_kt)\bv(k)}_{\text{fast modes}},
\end{equation}
where $\pmb{\phi}_k^{(1)}=\mathbf{\Phi}(:,k)$, and the superscript $(1)$ indicates the level. The first sum in the expression \eqref{linear_KMD_mrDMD} represents the slow-mode dynamics, whereas the second sum is everything else. How to choose the slow modes is important in the practical implementation. In the original mrDMD paper, it is suggested to set the threshold $\tau$ to select eigenvalues whose temporal behavior allows for at most one wavelength to fit within the sampling window.

The fast modes in \eqref{linear_KMD_mrDMD} can be collected into a data matrix $\Xv_{M/2}$, where we let $m_1$ denote the number of slow modes in \eqref{linear_KMD_mrDMD}. The matrix $\Xv_{M/2}$ is now split into two matrices, where the first matrix contains the first $M/2$ snapshots, and the second matrix contains the remaining $M/2$ snapshots. The process is now repeated, where $m_2$ slow modes are collected at the second level and computed separately in the first and second intervals of snapshots. This process is repeated to obtain the decomposition
$$
\xv(t)\approx \sum_{k=1}^{m_1} b_k^{(1)}\pmb{\phi}_k^{(1)}\exp(\eta_k^{(1)}t)+ \sum_{k=1}^{m_2} b_k^{(2)}\pmb{\phi}_k^{(2)}\exp(\eta_k^{(2)}t)+\sum_{k=1}^{m_3} b_k^{(3)}\pmb{\phi}_k^{(3)}\exp(\eta_k^{(3)}t)+\cdots,
$$
where the $\pmb{\phi}_k^{(\ell)}$ and $\eta_k^{(\ell)}$ are the DMD modes and eigenvalues at the $\ell$th level of the decomposition, the $b_k^{(\ell)}$ are the initial projections of the data onto the time interval of interest, and the $m_\ell$ are the number of slow modes retained at each level. The idea is that different spatiotemporal DMD modes are used to represent key multiresolution features. Thus, no single set of modes dominates the SVD and potentially marginalizes features at other time scales.

We can make the mrDMD more precise, letting $L$ denote the number of levels of the decomposition. The solution is a sum with $\ell$ indexing the level, $j=1,\ldots,2^{(\ell-1)}$ indexing the time bins $[t_j^{(\ell)},t_{j+1}^{(\ell)}]$ in each level and $k=1,\ldots,m_\ell$ indexing the modes extracted at each level. To simplify the sum, define the following indicator function
$$
f_{\ell,j} (t) =
\begin{cases}
1,\quad &\text{if }t\in[t_j^{(\ell)},t_{j+1}^{(\ell)}]\\
0,\quad&\text{otherwise}
\end{cases},\quad\text{so that}\quad \xv_{\mathrm{mrDMD}}(t)=\sum_{\ell=1}^L\sum_{j=1}^{2^{\ell-1}}\sum_{k=1}^{m_\ell}f_{\ell,j} (t)b_k^{(\ell,j)}\pmb{\phi}_k^{(\ell,j)}\exp(\eta_k^{(\ell,j)}t).
$$
In particular, each mode is represented in its respective time bin and level. Alternatively, this solution can be interpreted as yielding the least-squares fit to the dynamics within a given time bin at each level of the decomposition.

Numerous innovations enhance the practical implementation of mrDMD. Since only slow modes matter within a window, we can limit sampling to a fixed number of points per window, reducing the data matrix size for more manageable SVD computations. The sampling window locations are flexible, and smoothing their edges can prevent the Gibbs phenomenon due to abrupt data cutoffs. One can also employ wavelet functions, like Haar, Daubechies, or Mexican Hat, for the sifting function $f_{\ell,j} (t)$. Finally, overlapping windows prevent data loss during sampling. A sliding window approach can robustly track data features, enabling pattern correlation across windows akin to a Gabor transform-generated spectrogram.

\subsubsection{Example}

\begin{figure}
\centering
\raisebox{-0.5\height}{\includegraphics[width=0.9\textwidth,trim={0mm 0mm 0mm 0mm},clip]{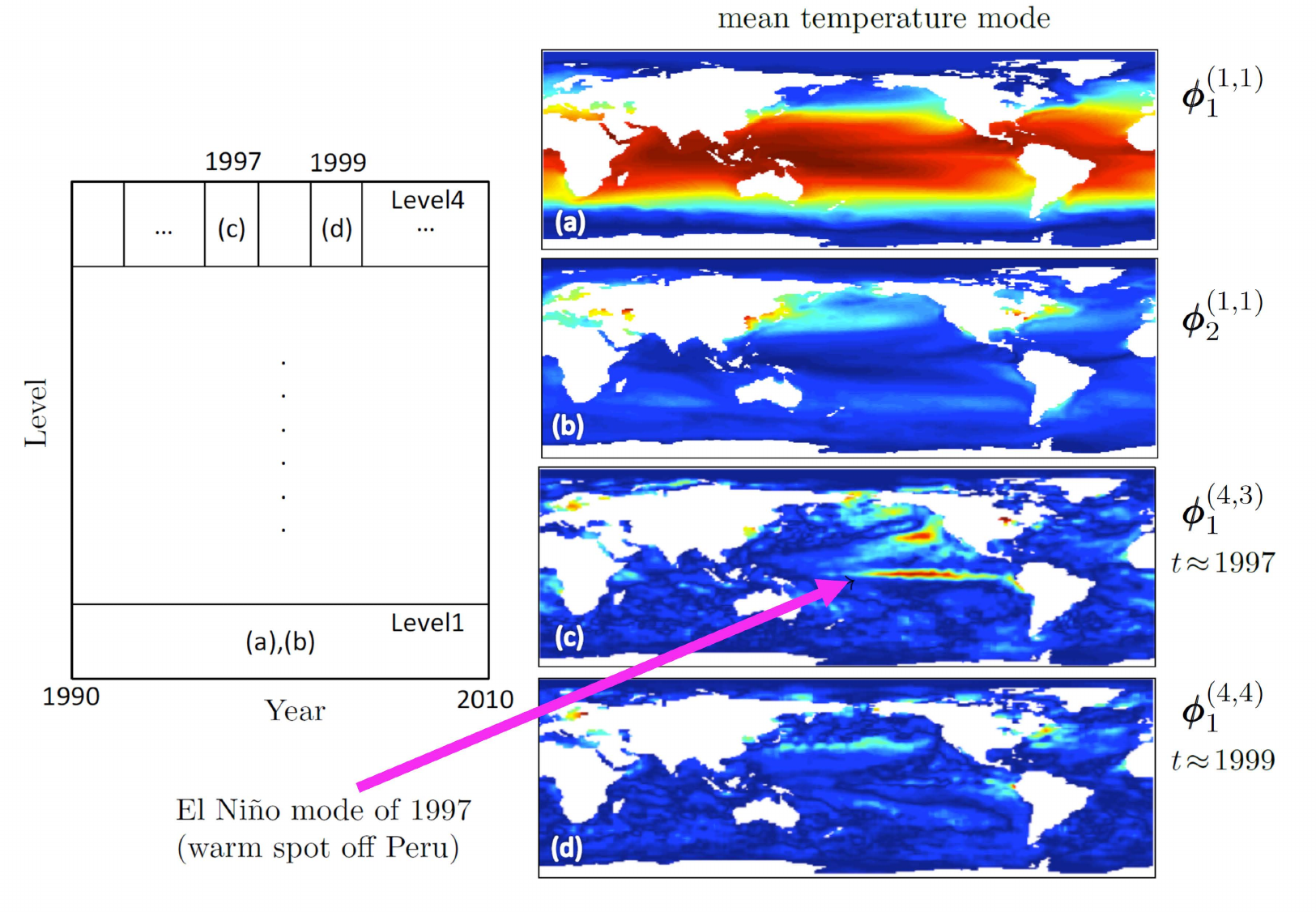}}
\caption{Application of mrDMD on sea surface temperature data from 1990 to 2010. The left panel illustrates the process for a 4-level decomposition. At each level, the slowest modes are extracted. Mode (c) clearly shows the El Niño mode of interest that develops in the central and east-central equatorial Pacific. The El Niño mode was absent in 1999, as is clear from mode (d). Reproduced with permission from \citep{kutz2016multiresolution}, copyright © 2016 Society for Industrial and Applied Mathematics, all rights reserved.}
\label{fig_mrDMD}
\end{figure}

The following example is from \citep{kutz2016multiresolution}. We consider the global sea surface temperature (see \cref{sec:examples_rDMD}) with data that spans 20 years from 1990 to 2010. \cref{fig_mrDMD} illustrates the outcomes of employing a 4-level mrDMD decomposition. At the first level, mrDMD identifies two modes: the mean ocean temperature, denoted as $\pmb{\phi}_1^{(1,1)}$, and an annual cycle, represented by $\pmb{\phi}_2^{(1,1)}$. Intriguingly, at the fourth level, the approximate zero mode of the sampling window uncovers noteworthy phenomena; specifically, it isolates the 1997 El Niño event. In contrast, when the same sampling window is applied to the year 1999, the El Niño mode is absent, aligning with the recognized oceanic patterns of that year. These insights would not have been obtainable using traditional DMD without pre-selecting the correct sampling windows. Moreover, even if such a step were taken, the slow modes identified at the first level, as shown in \cref{fig_mrDMD} (a) and (b), would pollute the data at the level of investigation. For additional techniques identifying approximate eigenfunctions that provide a rectified representation of the ENSO and function as (approximate) semi-conjugacies or factor maps with circle rotations, see \citep{froyland2021spectral}.

\subsubsection{Nonautonomous systems}

The use of mrDMD to detect transient behavior is tantalizing! Most DMD methods are designed for autonomous dynamical systems, where the function $\Fv$ on the right-hand side of \eqref{eq:DynamicalSystem} has no time dependence. However, there has been some recent initial work on nonautonomous systems, and we expect this area to grow significantly over the next few years. \cite{mezic2016koopman} were the first to extend the Koopman operator framework to nonautonomous dynamical systems, applying the methodology to linear-periodic and quasi-periodic nonautonomous systems. \cite{giannakis2019data} developed a strategy inspired by time-changed dynamical systems that involves rescaling the generator; this can be applied to a class of time-changed mixing systems. Further development of this approach, using delay-coordinate maps for recovering the dynamical system on tori with multiple time scales, is presented in \citep{das2019delay}. The extraction of spatiotemporal patterns using an extension of approximation techniques developed in \citep{giannakis2019data} on the space-time manifold defined as a skew-product structure is considered in \citep{giannakis2020extraction,giannakis2019spatiotemporal}. For the online computation of windowed DMD using rank-one updates, see \textit{Online DMD} \citep{zhang2019online}. \cite{macesic2018koopman} provides an error analysis for DMD with moving stencils. For extensions to actuated systems, see \citep{williams2016extending,bai2020dynamic}. \cite{redman2023koopman} develop an episodic memory approach that saves spectral objects associated with temporally local approximations of the Koopman operator, and utilizes this information to make new predictions. Nonautonomous systems have also been studied using transfer operators, which are the dual of Koopman operators (see the discussion in \cref{sec:transfer_ops}), and space-time manifolds \citep{froyland2023detecting}.

\subsection{Control}

One of the most successful applications of the Koopman operator framework lies in control \citep{mauroy2020koopman,otto2021koopman}, with demonstrated successes in various challenging applications. These include fluid dynamics \citep{arbabi2018data,peitz2020feedback}, robotics \citep{abraham2017model,bruder2019modeling,mamakoukas2019local,haggerty2023control}, power grids \citep{korda2018power,netto2018robust}, biology \citep{hasnain2020steady}, and chemical processes \citep{narasingam2019koopman}. The key point is that Koopman operators represent nonlinear dynamics within a globally linear framework. This approach leads to tractable convex optimization problems and circumvents theoretical and computational limitations associated with nonlinearity. Moreover, it is amenable to data-driven, model-free approaches  \citep{proctor2016dynamic,williams2016extending,korda2018linear,proctor2018generalizing,surana2016koopman,kaiser2017data,kaiser2018discovering,peitz2019koopman,abraham2019active}. The resulting models reveal insights into global stability properties \citep{sootla2016properties,mauroy2016global}, observability/controllability \citep{vaidya2007observability,goswami2017global,yeung2018koopman}, and sensor/actuator placement \citep{sinha2016operator,sharma2019transfer} for the underlying nonlinear systems.

Koopman operator theory was first extended to actuated systems by \cite{mezic2004comparison}, with stochastic forcing interpreted as actuation. \cite{proctor2016dynamic} developed the first control schemes based on DMD. A significant strength of DMD is the ability to describe complex and high-dimensional dynamical systems with a few dominant modes. Reducing the system's dimensionality enables faster and lower-latency prediction and estimation, leading to high-performance, robust controllers.

\subsubsection{Dynamic Mode Decomposition with Control (DMDc)}

We will focus on the \textit{DMD with control} (DMDc) algorithm \citep{proctor2016dynamic}. DMDc extends DMD to disambiguate between unforced dynamics and the effect of actuation. The DMD regression of \cref{sec_DMD_regression_view} is generalized to
$$
\xv_{n+1}=\Fv(\xv_{n},\uv_n)\approx \Av \xv_{n} +\Bv \uv_n,
$$
where $\uv_n\in\mathbb{C}^q$ is a vector of control inputs for each time-step. Here $\Av\in\mathbb{C}^{d\times d}$ and $\Bv\in\mathbb{C}^{d\times q}$ are unknown matrices. Snapshot triplets of the form $\{\xv^{(m)},\yv^{(m)},\uv^{(m)}\}_{m=1}^M$ are collected, where we assume that
$$
\yv^{(m)}\approx \Fv(\xv^{(m)},\uv^{(m)}),\quad m=1,\ldots, M.
$$
The control portion of the snapshots is arranged into the matrix $\mathbf{\Upsilon}=\begin{pmatrix} 
\uv^{(1)} \, \uv^{(2)} \, \cdots \, \uv^{(M)}
\end{pmatrix}\in\mathbb{C}^{q\times M}.$ The optimization problem in \eqref{DMD_opt_vanilla} is replaced by
$$
\min_{(\Av \,\,\, \Bv)} \left\|\Yv-(\Av \,\,\, \Bv)\mathbf{\Omega}\right\|_{\mathrm{F}}^2,\quad \text{where}\quad
\mathbf{\Omega}=\begin{pmatrix} 
\Xv \\ \mathbf{\Upsilon}
\end{pmatrix}\in\mathbb{C}^{(d+q)\times M}.
$$
A solution is given as $(\Av \,\,\, \Bv) = \Yv\mathbf{\Omega}^{\dagger}$. In practice, we seek a reduced-order model by performing a truncated SVD on both the input and output space. The full algorithm is summarized in \cref{alg:DMD_control} and is an extension of \cref{alg:DMD_vanilla}. DMDc has been used with Model-Predictive Control (MPC) for enhanced control of nonlinear systems in \citep{korda2018linear,kaiser2018sparse}, with the DMDc method performing surprisingly well, even for strongly nonlinear systems. Extensions are discussed in \cref{sec:extensions_for_control}.

\begin{algorithm}[t]
\textbf{Input:} Snapshot data $\Xv\in\mathbb{C}^{d\times M}$, $\Yv\in\mathbb{C}^{d\times M}$ and $\mathbf{\Upsilon}\in\mathbb{C}^{q\times M}$, ranks $r,p\in\mathbb{N}$. \\
\vspace{-4mm}
\begin{algorithmic}[1]
\State Compute a truncated SVD of the input matrix $\begin{pmatrix} 
\Xv \\ \mathbf{\Upsilon}
\end{pmatrix} \approx \tilde{\Uv} \tilde{\mathbf{\Sigma}}\tilde{\Vv}^*,$ $\tilde{\Uv}\in\mathbb{C}^{(d+q)\times p}$, $\tilde{\mathbf{\Sigma}}\in\mathbb{R}^{p\times p}$, $\tilde{\Vv}\in\mathbb{C}^{M\times p}.$ Break up the matrix $\tilde{\Uv}$ into $\tilde{\Uv}^*=[\tilde{\Uv}_1^* \, \tilde{\Uv}_2^*]$ where $\tilde{\Uv}_1\in\mathbb{C}^{d\times p}$ and $\tilde{\Uv}_2\in\mathbb{C}^{q\times p}$.
\State Compute a truncated SVD of $
\Yv \approx \widehat{\Uv} \widehat{\mathbf{\Sigma}}\widehat{\Vv}^*$, $\widehat{\Uv}\in\mathbb{C}^{d\times r}$, $\widehat{\mathbf{\Sigma}}\in\mathbb{R}^{r\times r}$, $\widehat{\Vv}\in\mathbb{C}^{M\times r}.$
\State Compute the compressions $
\tilde{\Av}= \widehat{\Uv}^*\Yv\tilde{\Vv}\tilde{\mathbf{\Sigma}}^{-1}\tilde{\Uv}_1^*\widehat{\Uv}\in\mathbb{C}^{r\times r}$ and $\tilde{\Bv}= \widehat{\Uv}^*\Yv\tilde{\Vv}\tilde{\mathbf{\Sigma}}^{-1}\tilde{\Uv}_2^*\in\mathbb{C}^{r\times q}.$
\State Compute the eigendecomposition $
\tilde{\Av}\Wv=\Wv\mathbf{\Lambda}$.
\noindent{}The columns of $\Wv$ are eigenvectors and $\mathbf{\Lambda}$ is a diagonal matrix of eigenvalues.
\State Compute the modes $
\mathbf{\Phi}=\Yv\tilde{\Vv}\tilde{\mathbf{\Sigma}}^{-1}\tilde{\Uv}_1^*\widehat{\Uv}\Wv.$
\end{algorithmic} \textbf{Output:} The eigenvalues $\mathbf{\Lambda}$ and modes $\mathbf{\Phi}\in\mathbb{C}^{d\times r}$.
\caption{DMD with control \citep{proctor2016dynamic}.}
\label{alg:DMD_control}
\end{algorithm}

\subsubsection{Example}

We illustrate DMDc for system identification on a high-dimensional, linear system with spectral sparsity following \citep[Section 4.3]{proctor2016dynamic}. We consider a two-dimensional torus discretized by a $128\times 128$ equispaced grid such that $\xv\in\mathbb{R}^{128\times 128}\cong \mathbb{R}^{16384}$. Taking the two-dimensional discrete Fourier transform of $\xv$, we obtain $\hat{\xv}$. The system evolves according to
$$
\hat{\xv}_{n+1}=\hat{\Av}\hat{\xv}_n+\hat{\Bv}\hat{\uv}_n.
$$
Here, $\hat{\Av}$ is a diagonal matrix with five non-zero entries representing the modes, each with a randomly chosen frequency and small damping. The random input signal, $\hat{\uv}$, is one-dimensional and directly influences the sparse modes, resulting in a localized negative control input when transformed back to the spatial domain. This back-transformation yields our dynamical system in physical space. The system is constructed by sampling a continuous-time system at time steps of $\Delta t=0.01$. We collect $M=400$ snapshots of the data for our analysis. Further details of this system can be found in \citep{brunton2016compressed}.

\begin{figure}
\centering
\raisebox{-0.5\height}{\includegraphics[width=0.5\textwidth,trim={0mm 0mm 0mm 0mm},clip]{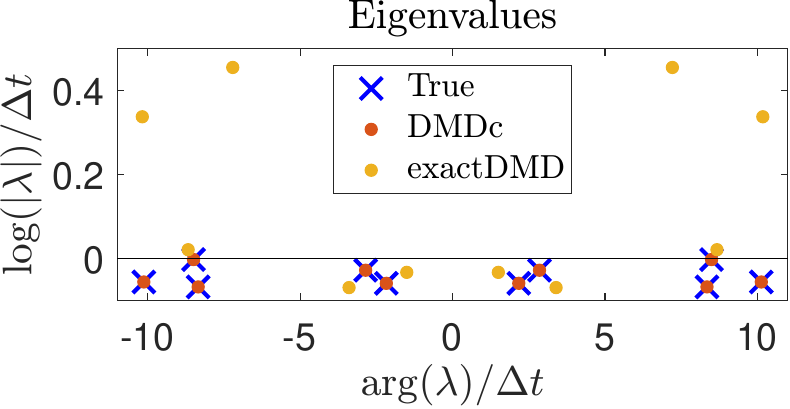}}
\caption{True eigenvalues of the torus example and those computed by DMDc and exact DMD. The logarithm of the eigenvalues are plotted to align with the continuous-time system.}
\label{fig_DMDc1}
\end{figure}

\begin{figure}
\centering
\raisebox{-0.5\height}{\includegraphics[width=0.9\textwidth,trim={0mm 0mm 0mm 0mm},clip]{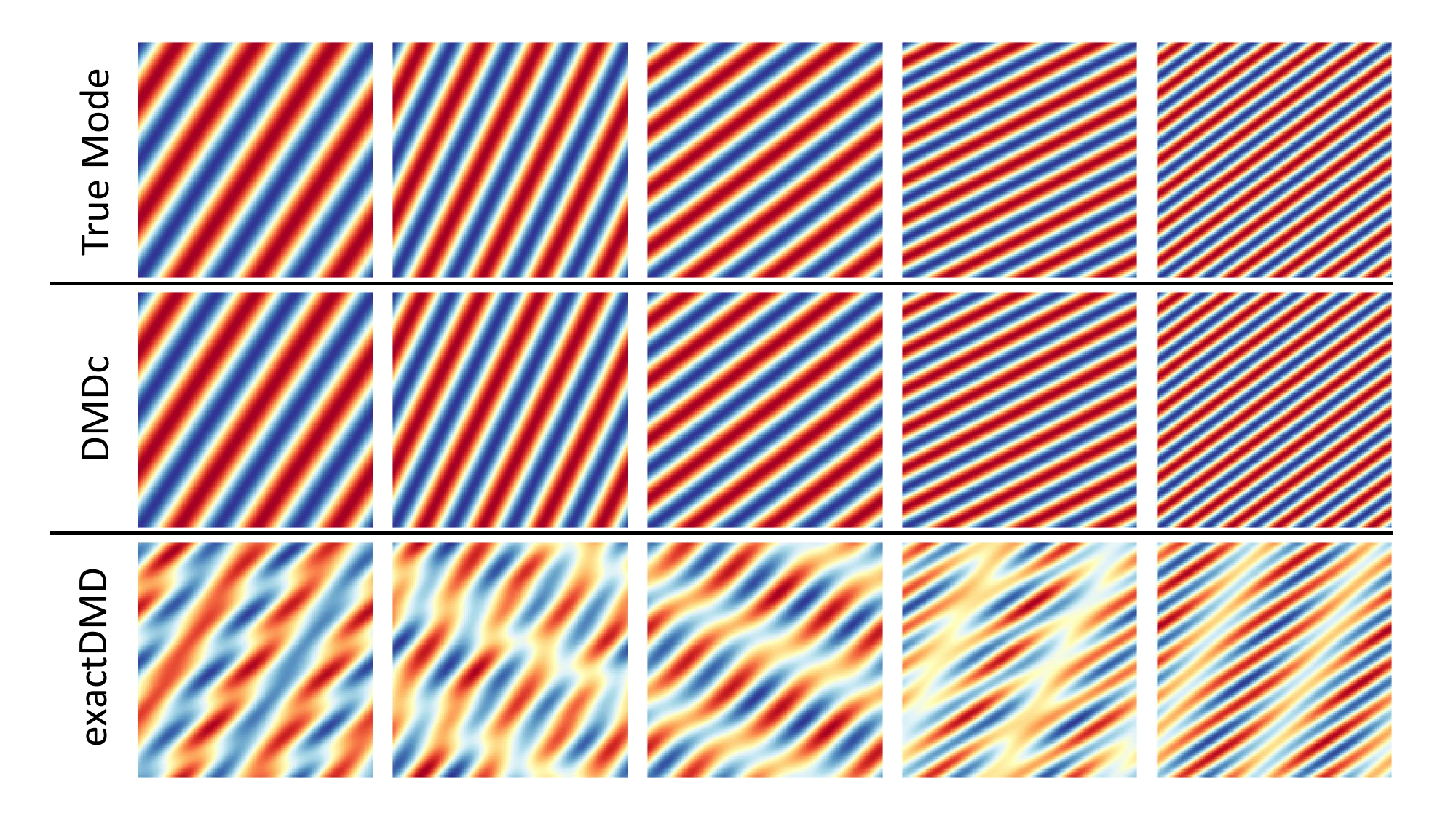}}
\caption{The true DMD modes for the torus example, alongside those computed by DMDc and exact DMD. The modes obtained from DMDc are accurate to machine precision, whereas those computed using exact DMD are significantly distorted.}
\label{fig_DMDc2}
\end{figure}

\cref{fig_DMDc1} displays the true eigenvalues alongside those computed by DMDc and exact DMD. Ten eigenvalues are present in conjugate pairs due to processing real-valued data. DMDc demonstrates greater accuracy than exact DMD, which inaccurately estimates some eigenvalues and generates unstable modes. The true DMD modes for this system appear at the top \cref{fig_DMDc2}. The DMDc modes in the middle row align almost perfectly with the true modes. The subspace angle between the true modes and the DMDc-computed modes is on the order of machine precision. In contrast, the modes produced by exact DMD show significant distortion.

\subsubsection{Extensions and connection with Koopman operators}
\label{sec:extensions_for_control}

Koopman theory has been used in combination with the
Linear Quadratic Regulator (LQR) \citep{brunton2016koopman,mamakoukas2019local,mamakoukas2021derivative},
state-dependent LQR \citep{kaiser2017data},
and MPC \citep{korda2018linear,kaiser2018sparse}.
Other noteworthy directions include
optimal control for switching control problems \citep{peitz2019koopman,peitz2020feedback},
Lyapunov-based stabilization \citep{huang2018feedback,huang2020data},
eigenstructure assignment \citep{hemati2017dynamic},
and active learning \citep{abraham2019active}.
Additionally, deep learning architectures have been employed to represent the nonlinear observables in combination with MPC \citep{li2019learning}, see also \citep{liu2018decomposition,han2020deep}, and \citep{peitz2019koopman,peitz2020data,klus2020data} for parametrized models.

Koopman theory is closely related to Carleman linearization \citep{carleman1932application}, which embeds finite-dimensional dynamics into infinite-dimensional linear systems using a polynomial basis. Carleman linearization has been used for decades to obtain truncated linear (and bilinear) state estimators \citep{krener1974linearization,brockett1976volterra} and to examine stability, observability, and controllability of the underlying nonlinear system \citep{loparo1978estimating}.

The DMDc framework may be extended to nonlinear observables using EDMD (see \cref{sec:Galerkin_perspective}), an approach called eDMDc \citep{williams2016extending}. \cite{korda2018linear} integrated eDMDc into MPC. Here, the Koopman operator is characterized as an autonomous operator on the extended state vector $(\mathbf{x}^\top, \mathbf{u}^\top)^\top$, with observables that may be nonlinear functions of both the state and the input. In practical applications, simplifications are employed to ensure the control problem remains convex \citep{korda2018linear,proctor2018generalizing}. This method has been applied for control in the coordinates of Koopman eigenfunctions \citep{kaiser2017data,kaiser2018discovering,folkestad2020extended} and in interpolated Koopman models \citep{peitz2018controlling,peitz2020data}. Convergence can be established under the assumption of an infinite amount of data and an infinite number of basis functions. Koopman Lyapunov-based MPC guarantees closed-loop stability and controller feasibility \citep{narasingam2019koopman,son2020handling}. However, general guarantees regarding the optimality, stability, and robustness of the controlled dynamical system are still limited.

The Koopman operator's eigenfunctions (or approximate eigenfunctions) are a natural choice of observables due to their simple temporal behavior. It is crucial to validate computed eigenfunctions to ensure that their evolution is consistent with the predictions of their associated eigenvalues, particularly for prediction tasks. They have been used for observer design within the Koopman canonical transform \citep{surana2016koopman,surana2016linear} and within the Koopman reduced-order nonlinear identification and control framework \citep{kaiser2017data}, which both typically yield a global bilinear representation of the underlying system. Subsequent research has focused on directly identifying Koopman eigenfunctions \citep{korda2020optimal,pan2021sparsity} and approximate invariant subspaces \citep{haseli2023modeling}.

The efficacy of Koopman-based MPC is currently at odds with the difficulties of approximating the Koopman operator and its spectra. Only a limited number of systems with a known Koopman-invariant subspace and verifiable eigenfunctions exist for model analysis and evaluation. Furthermore, the linearity of Koopman eigenfunctions is seldom validated. Nevertheless, Koopman-based MPC demonstrates remarkable resilience with models of marginal predictive ability. Despite notable successes, understanding how well the Koopman operator is actually approximated and producing error bounds remains largely incomplete.

\section{Variants from the Galerkin Perspective}
\label{sec:Galerkin_perspective}

We now explore variants of DMD from the Galerkin (or projection) perspective, building on the connection established in \cref{sec:galerkin_interp}. This approach particularly focuses on addressing the infinite-dimensional nature of Koopman operators. Given that a Koopman operator transforms a finite-dimensional nonlinear system into an infinite-dimensional linear one, a significant part of this section will address nonlinear observables. We will concentrate on three methods designed to tackle these challenges:
\begin{itemize}
	\item \textbf{Extended DMD:} This represents a fundamental extension of DMD that treats it as a Galerkin method.\footnote{Though once nonlinear observables have been chosen, one can also apply the regression interpretation of \cref{sec:regression_variants}.} In particular, it introduces nonlinear observables to form a dictionary, which generates a subspace within $L^2(\Omega,\omega)$. Adopting the Galerkin perspective enables the application of numerical tools for addressing infinite-dimensional spectral problems. However, the well-studied challenges of infinite-dimensional spectral computations are significant. Generally, EDMD will not converge to the spectral properties of the Koopman operator, either theoretically or practically (see \cref{sec:EDMD_convergence} and common pitfalls in \cref{sec:ResDMD}).
	\item \textbf{Time-delay Embedding:} This technique is commonly used to construct a dictionary of observables for EDMD and generates a Krylov subspace. Our focus will be on two methods: Hankel-DMD, which is a widely used technique suitable for ergodic systems that have a low-dimensional attractor, and HAVOK (Hankel Alternative View Of Koopman) analysis, which produces a linear model using the leading delay coordinates and includes forcing terms represented by low-energy delay coordinates.
	\item \textbf{Residual DMD:} This algorithm computes verified spectral properties of Koopman operators via an infinite-dimensional residual corresponding to the projection error of (E)DMD. This residual is computed from the snapshot data by augmenting EDMD with an additional matrix. This leads to the computation of spectra and pseudospectra without spectral pollution (general systems) and can be used to compute spectral measures (measure-preserving systems). Since the algorithms have error control, ResDMD allows aposteri verification of spectral quantities, Koopman mode decompositions, and learned dictionaries.
\end{itemize}

\subsection{Nonlinear observables: Extended Dynamic Mode Decomposition (EDMD)}
\label{sec:EDMD}

The standard DMD algorithm can accurately characterize periodic and quasi-periodic behaviors in nonlinear systems. However, DMD models based on linear observables generally fail to capture truly nonlinear phenomena. To address this limitation, \cite{williams2015data} introduced \textit{Extended DMD} (EDMD), which also elucidated the interpretation of DMD as a Galerkin method. Specifically, they demonstrated that EDMD converges to the numerical approximation obtained by a Galerkin method in the limit of large data sets. Prior research in a similar vein includes \citep{tu2014dynamic}. Moreover, the connection between EDMD and the earlier variational approach of conformation dynamics \citep{noe2013variational,nuske2014variational} from molecular dynamics is explored in \citep{wu2017variational,klus2018data}.

\subsubsection{The algorithm}
\label{sec:EDMD_algorithm}

Following the discussion of Koopman operators in \cref{sec:Koopman_operators}, the objective of EDMD is to approximate the Koopman operator with a matrix. For the sake of simplicity, the initial formulation of EDMD assumes that the columns of the snapshot matrix $\mathbf{X}$ are independently sampled from the distribution $\omega$. In our discussion, we extend EDMD to accommodate any given snapshot matrices
$$
\Xv=
\begin{pmatrix} 
\xv^{(1)} & \xv^{(1)} &\cdots & \xv^{(M)}
\end{pmatrix}
\quad \text{and} \quad
\Yv=
\begin{pmatrix} 
\yv^{(1)} & \yv^{(1)} &\cdots & \yv^{(M)}
\end{pmatrix},
$$
and consider the $\xv^{(m)}$ as \textit{quadrature nodes} used for integration with respect to $\omega$. This adaptability permits the application of various quadrature weights tailored to the specific scenario. It will be shown that EDMD generalizes the setup of \cref{sec:galerkin_interp}.

One first chooses a dictionary $\{\psi_1,\ldots,\psi_{N}\}$, i.e., a list of observables, in the space $L^2(\Omega,\omega)$. These observables form a finite-dimensional subspace $V_N=\mathrm{span}\{\psi_1,\ldots,\psi_{N}\}$. EDMD selects a matrix $\Kv\in\mathbb{C}^{N\times N}$ that approximates the action of $\mathcal{K}$ confined to this subspace. We desire that
$$
{[\mathcal{K}\psi_j](\xv) = \psi_j(\Fv(\xv)) \approx \sum_{i=1}^{N} \Kv_{ij} \psi_i(\xv)},\quad 1\leq j\leq {N}.
$$
Define the vector-valued feature map
$$
\Omega\ni\xv\mapsto\Psiv(\xv)=\begin{bmatrix}\psi_1(\xv) & \cdots& \psi_{{N}}(\xv) \end{bmatrix}\in\mathbb{C}^{1\times {N}}.
$$
Any $g\in V_{N}$ can be written as $g(\xv)=\sum_{j=1}^{N}\psi_j(\xv)\gv_j=\Psiv(\xv)\,\gv$ for some vector $\gv\in\mathbb{C}^{N}$. Hence
$$
	[\mathcal{K}g](\xv)=\Psiv(\Fv(\xv))\,\gv=\Psiv(\xv)(\Kv\,\gv)+\underbrace{\left(\sum_{j=1}^{N}\psi_j(\Fv(\xv))\gv_j-\Psiv(\xv)(\Kv\,\gv)\right)}_{=:R(\gv,\xv)}.
$$
Typically, $V_{N}$ is not an invariant subspace of $\mathcal{K}$. Hence, there is no choice of $\Kv$ that makes $R(\gv,\xv)$ zero for all $g\in V_N$ and $\omega$-almost every $\xv\in\Omega$. Instead, it is natural to select $\Kv$ as a solution of
\begin{equation}
	\min_{\Kv\in\mathbb{C}^{N\times N}} \left\{\int_\Omega \max_{\gv\in\mathbb{C}^{N},\|\Cv\gv\|_{\ell^2}=1}|R(\gv,\xv)|^2\ \mathrm{d}\omega(\xv)=\int_\Omega \left\|\Psiv(\Fv(\xv))\Cv^{-1} - \Psiv(\xv)\Kv\Cv^{-1}\right\|^2_{\ell^2}\ \mathrm{d}\omega(\xv)\right\}.
	\label{eq:ContinuousLeastSquaresProblem}
\end{equation}
Here, $\|\cdot\|_{\ell^2}$ denotes the standard Euclidean norm of a vector, and $\Cv$ is a positive self-adjoint matrix that controls the size of $g=\Psiv\gv$. One should think of this $\Cv$ as choosing an appropriate norm. This is important since not all norms on an infinite-dimensional vector space are equivalent (for an example in DMD analysis of fluid flow, see \citep[Figure 7]{colbrook2023mpedmd}).

In practical, data-driven contexts, it is not possible to directly evaluate the integral in \eqref{eq:ContinuousLeastSquaresProblem}. Instead, we approximate it via a quadrature rule with nodes $\{\xv^{(m)}\}_{m=1}^{M}$ and weights $\{w_m\}_{m=1}^{M}$.  For notational convenience, let $\Dv=\mathrm{diag}(w_1,\ldots,w_{M})$ and
\begin{equation}
	\begin{split}
		\Psiv_X=\begin{pmatrix}
			\Psiv(\xv^{(1)}) \\
			\vdots              \\
			\Psiv(\xv^{(M)})
		\end{pmatrix}\in\mathbb{C}^{M\times N},\quad
		\Psiv_Y=\begin{pmatrix}
			\Psiv(\yv^{(1)}) \\
			\vdots              \\
			\Psiv(\yv^{(M)})
		\end{pmatrix}\in\mathbb{C}^{M\times N}.
		\label{psidef}
	\end{split}
\end{equation}
The discretized version of \eqref{eq:ContinuousLeastSquaresProblem} is the following weighted least-squares problem:
\begin{equation}
	\label{EDMD_opt_prob2}
	\min_{\Kv\in\mathbb{C}^{N\times N}}\left\{\sum_{m=1}^{M} w_m\left\|\Psiv(\yv^{(m)})\Cv^{-1}-\Psiv(\xv^{(m)})\Kv\Cv^{-1}\right\|^2_{\ell^2}=\left\|\Dv^{1/2}\Psiv_Y\Cv^{-1}-\Dv^{1/2}\Psiv_X \Kv\Cv^{-1}\right\|_{\mathrm{F}}^2\right\},
\end{equation}
where we remind the reader that $\|\cdot\|_{\mathrm{F}}$ denotes the Frobenius norm. By reducing the size of the dictionary if necessary, we may assume without loss of generality that $\Dv^{1/2}\Psiv_X$ has rank $N$. For example, we can do this in DMD by projecting onto POD modes. Regularization through a truncated singular value decomposition may also be considered. A solution to \eqref{EDMD_opt_prob2} is
$$
	\Kv= (\Dv^{1/2}\Psiv_X)^\dagger \Dv^{1/2}\Psiv_Y=(\Psiv_X^*\Dv\Psiv_X)^\dagger\Psiv_X^*\Dv\Psiv_Y,
$$
where `$\dagger$' denotes the pseudoinverse. The second equality follows since $\Dv^{1/2}\Psiv_X$ has linearly independent columns. Note that this solution is independent of the matrix $\Cv$. However, a suitable choice of $\Cv$ is vital once we add constraints to the optimization problem in \eqref{eq:ContinuousLeastSquaresProblem}, see \cref{sec:mpedmd_alg}. As observed in \cref{sec:galerkin_interp}, if the quadrature weights are equal and $\Psiv=\begin{bmatrix}u_1 & \cdots& u_{r} \end{bmatrix}$ constitutes an appropriate linear dictionary, then $\Kv$ is the transpose of the DMD matrix. Conceptually, DMD can be regarded as a particular instance of EDMD employing a set of linear basis functions.

We now generalize \cref{sec:galerkin_interp} by defining the two correlation matrices
\begin{equation}
\label{eq_EDMD_corr_matrices}
\Gv=\Psiv_X^*\Dv\Psiv_X=\sum_{m=1}^{M} w_m \Psiv(\xv^{(m)})^*\Psiv(\xv^{(m)}),\quad
\Av=\Psiv_X^*\Dv\Psiv_Y=\sum_{m=1}^{M} w_m \Psiv(\xv^{(m)})^*\Psiv(\yv^{(m)}).
\end{equation}
If we consider the discrete measure $\omega_M=\sum_{m=1}^Mw_m\delta_{\xv^{(m)}}$, then
$$
\Gv_{jk}=\int_{\Omega} \overline{\psi_j(\xv)}\psi_k(\xv) \ \mathrm{d} \omega_M(\xv),\quad \Av_{jk}=\int_{\Omega} \overline{\psi_j(\xv)}\psi_k(\Fv(\xv)) \ \mathrm{d} \omega_M(\xv).
$$
If the quadrature converges, then
\begin{equation}
	\label{quad_convergence}
	\lim_{M\rightarrow\infty}\Gv_{jk} = \langle \psi_k,\psi_j \rangle\quad \text{ and }\quad \lim_{M\rightarrow\infty}\Av_{jk} = \langle \mathcal{K}\psi_k,\psi_j \rangle,
\end{equation}
where $\langle \cdot,\cdot \rangle$ is the inner product associated with $L^2(\Omega,\omega)$. Hence, in the large data limit, $\Kv=\Gv^\dagger\Av$ approaches a matrix representation of $\mathcal{P}_{V_{N}}\mathcal{K}\mathcal{P}_{V_{N}}^*$, where $\mathcal{P}_{V_{N}}$ denotes the orthogonal projection onto $V_{N}$. In essence, EDMD is a Galerkin method. The EDMD eigenvalues thus approach the spectrum of $\mathcal{P}_{V_{N}}\mathcal{K}\mathcal{P}_{V_{N}}^*$, and EDMD is an example of the so-called \textit{finite section method} \citep{bottcher1983finite}\citep[Section 4]{mezic2020numerical}. Since the finite section method can suffer from spectral pollution (spurious modes), spectral pollution is a concern for EDMD \citep{williams2015data}. We saw an explicit example in \cref{sec:duffing_spectral_pollution}. See also \citep[Example 2]{mezic2020numerical} for the worked example $\Fv(\xv)=\xv^2$ on the circle.

\begin{algorithm}[t]
\textbf{Input:} Snapshot data $\Xv\in\mathbb{C}^{d\times M}$ and $\Yv\in\mathbb{C}^{d\times M}$, quadrature weights $\{w_m\}_{m=1}^{M}$, and a dictionary of functions $\{\psi_j\}_{j=1}^{N}$.\\
\vspace{-4mm}
\begin{algorithmic}[1]
\State Compute the matrices $\Psiv_X$ and $\Psiv_Y$ defined in \eqref{psidef} and $\Dv=\mathrm{diag}(w_1,\ldots,w_{M})$.
\State Compute the EDMD matrix $\Kv= (\Dv^{1/2}\Psiv_X)^\dagger \Dv^{1/2}\Psiv_Y\in\mathbb{C}^{N\times N}$.
\State Compute the eigendecomposition $
\Kv\Vv=\Vv\mathbf{\Lambda}$.

\noindent{}The columns of $\Vv$ are eigenvector coefficients and $\mathbf{\Lambda}$ is a diagonal matrix of eigenvalues.
\end{algorithmic} \textbf{Output:} The eigenvalues $\mathbf{\Lambda}$ and eigenvector coefficients $\Vv\in\mathbb{C}^{N\times N}$.
\caption{The EDMD algorithm \citep{williams2015data}.}
\label{alg:EDMD}
\end{algorithm}

\cref{alg:EDMD} summarizes the procedure for computing eigenvalues and eigenvectors. We can also use EDMD to compute Koopman modes. Given an observable $g=\Psiv\gv\in V_N$, we may expand $g$ in terms of the eigenvectors of $\Kv$ as
\begin{equation}
\label{KMD_exact}
g=\Psiv\gv = \Psiv \Vv\left[\Vv^{-1}\gv\right],
\end{equation}
where $\Vv$ is the matrix of eigenvectors of $\Kv$ with eigenvalues $\{\lambda_j\}_{j=1}^N$. Similarly, for general $g\in L^2(\Omega,\omega)\backslash V_N$, we obtain an approximate expansion
\begin{equation}
\label{KMD_approx}
g\approx \Psiv \Vv\left[\Vv^{-1}(\Dv^{1/2}\Psiv_X)^{\dagger}\Dv^{1/2}\left(g(\xv^{(1)}),\ldots,g(\xv^{(M)})\right)^\top\right].
\end{equation}
This expansion is called the KMD of $g$.\footnote{Unfortunately, there are numerous meanings of the term KMD in the literature. There is the KMD of \cite{mezic2005spectral}, which we discussed in \cref{sec:spectra_crash_course} and is based on the spectral theorem for unitary Koopman operators. There is also the (typically approximate) KMD produced by DMD and EDMD.} With an abuse of notation, if $g\notin V_N$, we set
$$
\gv=(\Dv^{1/2}\Psiv_X)^{\dagger}\Dv^{1/2}\left(g(\xv^{(1)}),\ldots,g(\xv^{(M)})\right)^\top.
$$
As $M\rightarrow\infty$, assuming that the quadrature rule underlying EDMD converges, the approximation $g\approx \Psiv \gv$ converges to the projected observable $\mathcal{P}_{V_N}g$. As a particular case, we can vectorize and obtain
$$
\xv\approx \Psiv(\xv) \Vv\left[\Vv^{-1}(\Dv^{1/2}\Psiv_X)^{\dagger}\Dv^{1/2}\left(\xv^{(1)},\ldots,\xv^{(M)}\right)^\top\right].
$$
The $j$th row of the matrix in square brackets is known as the $j$th Koopman mode, which we denote as $\pmb{\xi}_j\in\mathbb{C}^{1\times d}$. Generalizing \eqref{linear_KMD}, the KMD provides an approximation of the dynamics by
$$
\xv_n\approx \Psiv(\xv_0)\Kv^n\Vv\begin{pmatrix}
\pmb{\xi}_1\\
\vdots \\
\pmb{\xi}_N
\end{pmatrix}= \Psiv(\xv_0)\Vv\mathbf{\Lambda}^n\begin{pmatrix}
\pmb{\xi}_1\\
\vdots \\
\pmb{\xi}_N
\end{pmatrix}=\sum_{j=1}^N  \Psiv(\xv_0)\Vv(:,j)\lambda_j^n\pmb{\xi}_j.
$$
Similarly, for general $g$, we obtain
$$
g(\xv_n)\approx \Psiv(\xv_0)\Kv^n\Vv\Vv^{-1}\gv= \Psiv(\xv_0)\Vv\mathbf{\lambda}^n\Vv^{-1}\gv=\sum_{j=1}^N  \Psiv(\xv_0)\Vv(:,j)\lambda_j^n[\Vv^{-1}\gv]_j,
$$
which includes the triple of Koopman eigenvectors, eigenvalues, and modes.

\subsubsection{Choices of dictionary}

We have already met two examples of EDMD in this review: the Lorenz system discussed in \cref{sec:intro_lorenz}, where a dictionary was constructed from delay embedding, and the Duffing oscillator discussed in \cref{sec:duffing_spectral_pollution}, where we utilized a dictionary of radial basis functions. The selection of the dictionary significantly affects the efficacy of EDMD. In their original formulation, \cite{williams2015data} suggest various dictionary choices, such as polynomials, Fourier modes, spectral elements, and radial basis functions. Subsequent extensions have primarily focused on addressing the challenges of large state-space dimensions and mitigating the curse of dimensionality.

\textit{Kernelized EDMD} \citep{williams2015kernel} (developed in parallel in \citep{kawahara2016dynamic}) uses the kernel trick \citep{scholkopf2001kernel} to perform EDMD with a choice of dictionary determined implicitly by a choice of a kernel function. This approach can help circumvent the curse of dimensionality and can be very effective when the state-space dimension $d$ is large. Numerous papers have been written on the approximation of Koopman operators in a RKHS \citep{klus2018kernel,fujii2019dynamic,degennaro2019scalable,alexander2020operator,das2020koopman,klus2020eigendecompositions,klus2020kernel,mezic2020spectrum,burov2021kernel,baddoo2021kernel,kostic2022learning,khosravi2023representer,philipp2023error}. This also includes methods for continuous-time dynamical systems \citep{das2021reproducing,rosenfeld2022dynamic}. A challenge associated with RKHS techniques is that a general RKHS does not exhibit invariance under the action of the Koopman operator. This situation renders the selection of a reproducing kernel a delicate task. Ideally, one should choose the kernel so that the Koopman operator on the RKHS is not only densely defined but also closable. Finding such a kernel is generally non-trivial, as indicated in \citep{ikeda2022koopman}.

\textit{Kernel analog forecasting} (KAF) \citep{zhao2016analog} is a kernel method used for nonparametric statistical forecasting of dynamically generated time series data. Under measure-preserving and ergodic dynamics, KAF consistently approximates the conditional expectation of observables that are acted upon by the Koopman operator of the dynamical system and are conditioned on the observed data at forecast initialization \citep{alexander2020operator}. KAF yields optimal predictions in the sense of minimal root mean square error with respect to the invariant measure in the asymptotic limit of large data. This connection facilitates the analysis of generalization error and uncertainty quantification. KAF has been used with streaming kernel regression \citep{giannakis2021learning} and for multiscale systems \citep{burov2021kernel}.

\textit{Diffusion forecasting} \citep{berry2015nonparametric} uses the diffusion maps algorithm \citep{coifman2006diffusion} to construct a data-driven basis. Leveraging spectral convergence results for kernel integral operators \citep{garcia2020error,von2008consistency}, this approach produces a well-conditioned and consistent approximation as both the amount of training data and the number of basis functions increase. \cite{giannakis2015spatiotemporal,giannakis2019data} use the diffusion forecasting technique in a framework that approximates the generator $\mathcal{L}$ of measure-preserving ergodic flows on manifolds by an advection-diffusion operator $\mathcal{L}_\tau = \mathcal{L} - \tau \Delta$, where $\tau$ is a regularization parameter, and $\Delta$ is a Laplace-type diffusion operator. A Galerkin method was developed for the eigenvalue problem of $\mathcal{L}_\tau$, which was observed to perform efficiently for systems with a pure point spectrum, such as ergodic rotations on tori. The most straightforward case for analyzing the spectral properties of diffusion-regularized generators arises when the regularizing operator $\Delta$ commutes with $\mathcal{L}$. \cite{das2019delay,giannakis2019data} demonstrated that a commuting operator $\Delta$ can be derived from the infinite-delay limit of a family of kernel integral operators constructed using time-delay embedding.

Another prevalent method involves training \textit{neural networks} as a suitable dictionary to construct Koopman forecasts, as demonstrated in several studies \citep{li2017extended,takeishi2017learning,wehmeyer2018time,yeung2019learning,azencot2020forecasting,eivazi2021recurrent,li2021deep,alford2022deep}. This approach is typically implemented in two ways: by identifying a few key latent variables or by lifting to a higher-dimensional input space. Variational autoencoders (VAMPnets) have been employed for stochastic dynamical systems such as in molecular dynamics \citep{mardt2018vampnets,wehmeyer2018time}, wherein the mapping back to the physical configuration space from the latent variables is probabilistic. The integration of Koopman analysis with graph convolutional neural networks has been explored to learn the dynamics of atoms within materials \citep{xie2019graph}. \cite{lusch2018deep} employ an auxiliary network to parameterize the continuously varying spectral parameter, enabling a network structure that offers both parsimony and interpretability. A notable challenge when incorporating EDMD with neural networks is the trade-off between representing data accurately and the potential for overfitting, particularly with limited data. To address this issue, \cite{otto2019linearly} proposed an architecture that combines an autoencoder with linear recurrent dynamics in the encoded space. Beyond employing neural networks for learning Koopman embeddings, Koopman theory has also been applied to understand the behavior of neural networks themselves \citep{manojlovic2020applications,dogra2020optimizing} and algorithms more broadly \citep{dietrich2020koopman,redman2022algorithmic}.

\subsubsection{Convergence theory}
\label{sec:EDMD_convergence}

We now outline the convergence theory for EDMD. Unfortunately, the type of convergence (in the strong operator topology) is too weak to ensure the convergence of spectral properties. To take into account the snapshot data and dictionary, we let $\Kv_{N,M}$ denote the EDMD matrix. When considering the convergence of EDMD and related methods, there are two limits of interest:
\begin{itemize}
	\item The large-data limit which corresponds to $M\rightarrow\infty$;
	\item The large-subspace limit which corresponds to $N\rightarrow\infty$.
\end{itemize}
To compute the spectral properties of $\mathcal{K}$, a double limit
$$
\lim_{N\rightarrow\infty}\lim_{M\rightarrow\infty}\Kv_{N,M}
$$
must be considered. Generally, these limits do not commute. More broadly, the use of successive limits is a common occurrence in spectral problems and other areas of scientific computation and cannot be overcome regardless of the choice of algorithm \citep{colbrookthesis,colbrook2022computation,colbrook2022foundations,SCI_ref}.

We saw above that if the quadrature rule converges, i.e., \eqref{quad_convergence} holds, then $\lim_{M\rightarrow\infty}\Kv_{N,M}=\Kv_N$ is a Galerkin matrix. There are essentially three options for the quadrature rule:
\begin{itemize}
	\item \textbf{Random sampling:} We may draw $\xv^{(m)}$ at random according to a probability measure that is absolutely continuous with respect to $\omega$, and select the quadrature weights according to the corresponding Radon--Nikodym derivative. This was essentially observed in \citep{williams2015data}. Convergence holds with probability one \citep[Section 3.4]{2158-2491_2016_1_51} provided that $\omega$ is not supported on a zero level set that is a linear combination of the dictionary \citep[Section 4]{korda2018convergence}. The convergence rate is typically $\mathcal{O}(M^{-1/2})$ \citep{caflisch1998monte}, but is a practical approach if the state-space dimension is large. One could also consider quasi-Monte Carlo integration, which can achieve a faster rate of $\mathcal{O}(M^{-1})$ (up to logarithmic factors) under suitable conditions \citep{caflisch1998monte}.
	\item \textbf{Ergodic sampling:} If the system is ergodic, then we can replace the strong law of large numbers with Birkhoff's Ergodic theorem \citep{birkhoff1931proof}:
	\begin{equation}
	\label{eq:birkhoff}
\lim_{n\rightarrow\infty} \frac{1}{n} \sum_{j=0}^{n-1} [\mathcal{K}^jg](\xv_0)=\lim_{n\rightarrow\infty} \frac{1}{n} \sum_{j=0}^{n-1} g(\xv_j)=\int_{\Omega} g(\xv)\ \mathrm{d} \omega(\xv)\quad \forall g\in L^1(\Omega,\omega).
	\end{equation}
We may select $\xv^{(m)}=\xv_{m-1}$ from a single trajectory starting at $\omega$-almost any initial condition $\xv_0$ and $w_m=1/M$. Often, the measures are `physical,' meaning that the set of initial points with convergence has a positive Lebesgue measure.\footnote{There is also the notion of SRB measure, which often coincides. For a survey of these measures and their definitions, see \citep{young2002srb}.} For example, taking $g=[\mathcal{K}\psi_k]\cdot \overline{\psi_j}$ in \eqref{eq:birkhoff}, we obtain
$$
\lim_{M\rightarrow\infty}\underbrace{\frac{1}{M}\sum_{n=0}^{M-1} \psi_k(\xv_{n+1})\overline{\psi_j(\xv_n)}}_{=\Av_{jk}}=\langle \mathcal{K}\psi_k,\psi_j \rangle.
$$
Convergence in this scenario is analyzed in \citep{arbabi2017ergodic,korda2018convergence}. The convergence rate in $M$ is problem dependent \citep{kachurovskii1996rate,mezic2002ergodic}. For periodic and quasi-periodic attractors, the error of approximating the inner products is generally $\mathcal{O}(M^{-1})$. For strongly mixing systems, the rate of convergence slows down to $\mathcal{O}(M^{-1/2})$. However, convergence rates cannot be established for the
general class of ergodic systems. For convergence rates of von Neumann's ergodic theorem in the context of Koopman operators, see \citep{aloisio2022spectral}.
	\item \textbf{High-order quadrature:} If the statespace dimension $d$ is not too large and $\Omega$ is sufficiently simple, it can be effective to choose $\{(\xv^{(m)},w_{m})\}$ according to a high-order quadrature rule. Even changing the weights $\{w_m\}$ for a fixed set of sample points $\{\xv^{(m)}\}$ can lead to a considerable acceleration of the convergence \citep{colbrook2021rigorous}.
\end{itemize}
\cite{colbrook2023beyond} provide concentration bounds on the error of the finite $M$ EDMD matrix. \cite{mollenhauer2022kernel} provide a rigorous analysis of kernel autocovariance operators, including nonasymptotic error bounds under classical ergodic and mixing assumptions. \cite{nuske2023finite} presented the first rigorously derived probabilistic bounds on the finite-data approximation error for the truncated Koopman generator of stochastic differential equations (SDEs) and nonlinear control systems. Two settings were analyzed: independent and identically distributed sampling and ergodic sampling, where it was assumed that the Koopman semigroup is exponentially stable for the latter. \cite{lu2020prediction} provide bounds for parabolic PDEs.

Suppose the quadrature rule converges and we have passed to the limit $M\rightarrow\infty$. \cite{korda2018convergence} show that under a natural density assumption of $V_N$ as $N\rightarrow\infty$, $\Kv_N$ converges strongly to $\mathcal{K}$ for bounded Koopman operators. This means that
\begin{equation}
\label{SOT_Koopman}
\lim_{N\rightarrow\infty}\|\mathcal{K}g-\Psiv \Kv_N \mathcal{P}_{V_{N}}g \|_{L^2(\Omega,\omega)}=0\quad \forall g\in L^2(\Omega,\omega),
\end{equation}
where, with an abuse of notation, $\mathcal{P}_{V_{N}}g$ denotes the vector of coefficients of $\mathcal{P}_{V_{N}}g$. It is straightforward to drop the assumption that $\mathcal{K}$ is bounded by making natural assumptions on the dictionary and considering $g$ in the domain of $\mathcal{K}$ \citep{colbrook2021rigorous}. Unfortunately, strong convergence is insufficient to ensure that the spectral properties of $\Kv_N$ converge to that of $\mathcal{K}$ - \cite{mezic2020numerical} provides an explicit example. We also saw an example of this effect in \cref{sec:duffing_spectral_pollution}. In \cref{sec:ResDMD}, we will show how Residual DMD provides convergence and error control in the final limit $N\rightarrow\infty$.

\subsubsection{Infinitesimal generators}

Several methods have also been proposed for continuous-time systems and approximating the Koopman infinitesimal generator defined in \eqref{eq:gen_def}. For example, \textit{generator EDMD} (gEDMD) \citep{klus2020data} uses time derivatives of the dictionary to extend EDMD to compute the generator, see also \citep{klus2020kernel,rosenfeld2022dynamic}. Other methods include computing the matrix logarithm of the Koopman operator \citep{mauroy2019koopman,drmavc2021identification}, approximating the Koopman operator family, and using finite-differences to compute the Lie derivative of the Koopman operator \citep{giannakis2021delay,sechi2021estimation}. Finally, \cite{das2019delay,giannakis2019data,giannakis2020extraction} approach the problem of approximating both the Koopman and its generator as a manifold-learning problem on a space-time manifold. This challenge was successfully addressed for ergodic dynamical systems, such as those evolving on a chaotic attractor.

\subsubsection{Example}

\begin{figure}
\centering
\raisebox{-0.5\height}{\includegraphics[width=0.32\textwidth,trim={0mm 0mm 0mm 0mm},clip]{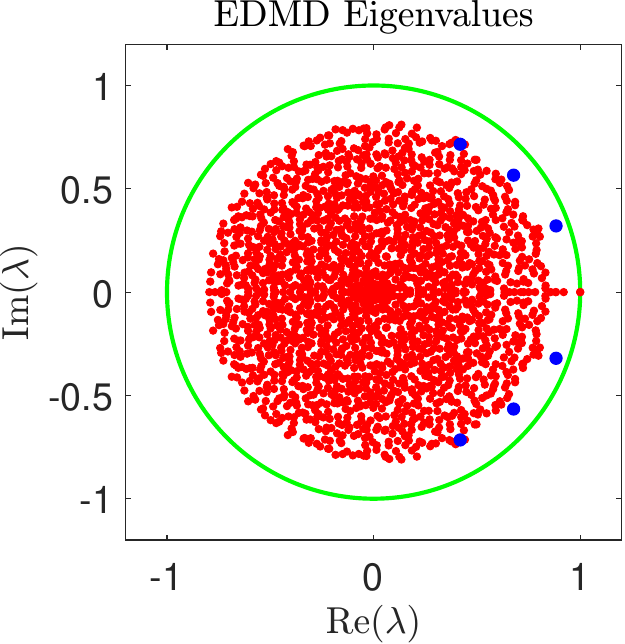}}\hfill
\raisebox{-0.5\height}{\includegraphics[width=0.30\textwidth,trim={0mm 0mm 0mm 0mm},clip]{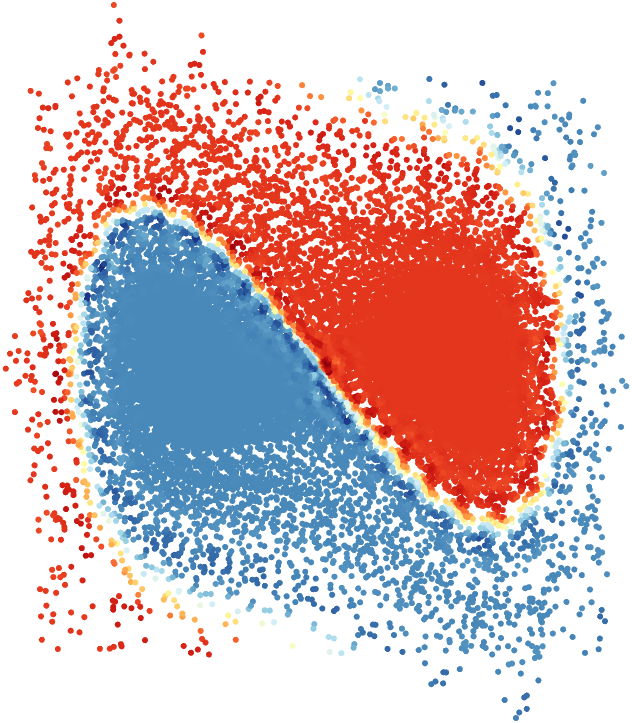}}\hfill
\raisebox{-0.5\height}{\includegraphics[width=0.32\textwidth,trim={0mm 0mm 0mm 0mm},clip]{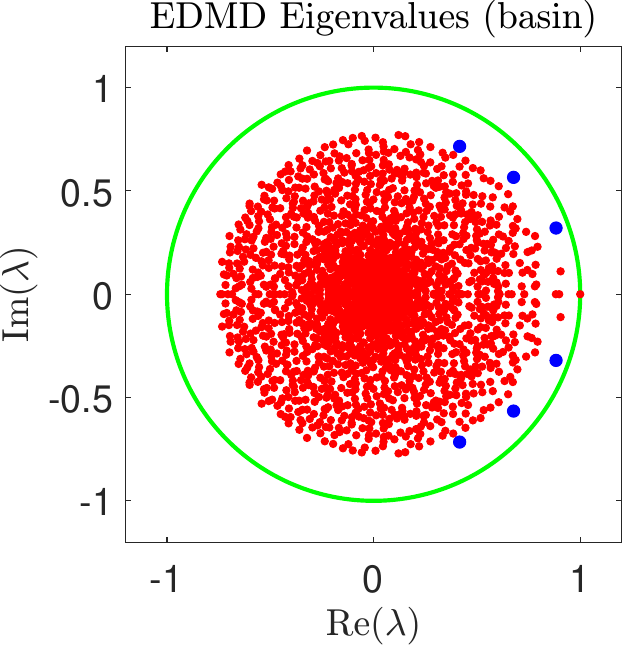}}
\caption{Left: EDMD eigenvalues computed over the entire state space. Middle: Eigenfunction that parametrizes the basins of attraction. Right: EDMD eigenvalues obtained after restricting the process to one basin of attraction. The powers of the dominant damped mode are shown in blue.}
\label{fig_EDMD1}
\end{figure}

As an example of EDMD, we revisit the Duffing oscillator from \cref{sec:duffing_spectral_pollution} but follow the experiment of \cite{williams2015data} closely. Namely, we consider the damped system:
$$
\dot{x}=y,\quad\dot{y}=-0.5y+x-x^3.
$$
In this regime, there are two stable spirals at $(\pm1,0)$ and a saddle at the origin. Almost every initial condition, except those on the stable manifold of the saddle, is drawn to one of the spirals. We collect trajectory data and form the dictionary of observables $\{\psi_j\}_{j=1}^N$ in the same manner as before. \cref{fig_EDMD1} (left) shows the eigenvalues computed using EDMD with $N=2000$. The system is now damped, and there is a lattice structure of dominant but damped modes, $\{\lambda_1^n,\overline{\lambda_1}^n:n\in\mathbb{N},\lambda_1\approx 0.8831 + 0.3203i\}$ shown in blue. The lattice structure can be understood as follows: if $g$ and $f$ are eigenfunctions of $\mathcal{K}$ corresponding to eigenvalues $\lambda$ and $\mu$, respectively, and if the product $fg$ is within the function space that forms the domain of $\mathcal{K}$, then
$$
[\mathcal{K}(fg)](\xv)=f(\Fv(\xv))g(\Fv(\xv))=[\mathcal{K}f](\xv)[\mathcal{K}g](\xv)=\lambda \mu f(\xv)g(\xv).
$$
Namely, further eigenvalues and eigenfunctions can be constructed by taking products.

\begin{figure}
\centering
\raisebox{-0.5\height}{\includegraphics[width=1\textwidth,trim={0mm 0mm 0mm 0mm},clip]{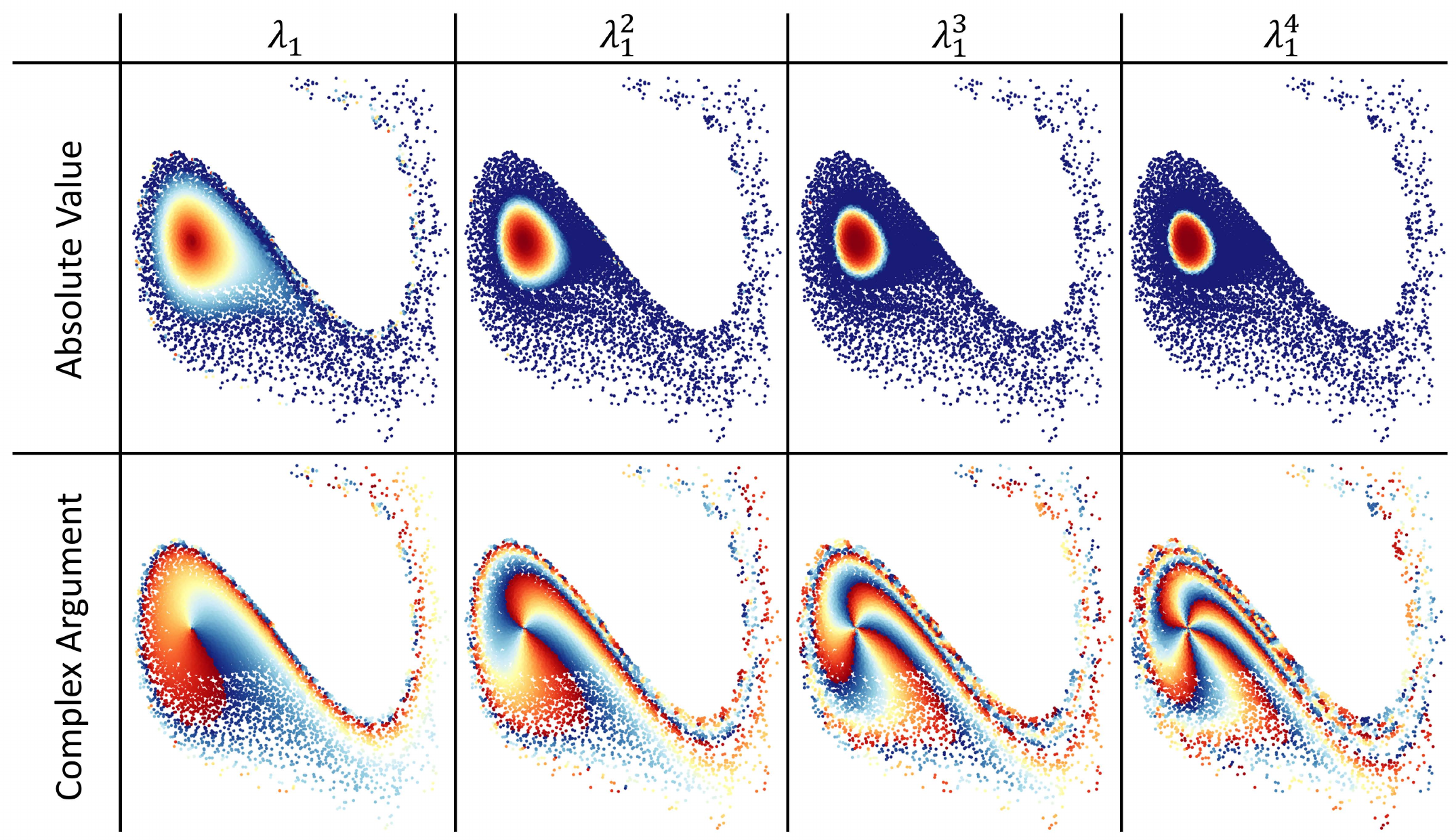}}
\caption{Top row: Absolute value of eigenfunction. Bottom row: Complex argument of eigenfunction. Left to write: Eigenfunctions corresponding to powers $\lambda_1$, $\lambda_1^2$, $\lambda_1^3$ and $\lambda_1^4$ of the fundamental eigenvalue $\lambda_1$.}
\label{fig_EDMD2}
\end{figure}

For this system, the eigenspace corresponding to $\lambda=1$ is spanned by the constant function and the indicator function of the invariant set corresponding to the two basins of attraction. This is illustrated in the middle of \cref{fig_EDMD1}. Utilizing the level sets of this eigenfunction, we limit the data to the basin of $(-1,0)$ and rerun the process to compute a new dictionary. The resulting EDMD eigenvalues are displayed on the right side of \cref{fig_EDMD1}, where the eigenvalue $\lambda_1^3$ is now more distinctly observable. In \cref{fig_EDMD2}, we plot the eigenfunctions corresponding to the powers $\lambda_1$, $\lambda_1^2$, $\lambda_1^3$, and $\lambda_1^4$ of the fundamental eigenvalue. Note that the eigenfunctions are successive powers of one another. Furthermore, the amplitude and phase of a Koopman eigenfunction are analogous to an 'action–angle' parametrization of the basin of attraction. The level sets of the absolute values of the eigenfunctions are the so-called isostables, while the level sets of the arguments of the eigenfunctions are termed isochrons \citep{mauroy2013isostables}.

\subsection{Time-delay embedding}

In many applications, only partial observations of the system are available, leading to hidden or latent variables. Additionally, the explicit construction of a robust nonlinear dictionary can be challenging, particularly when the system evolves on a low-dimensional attractor that may be unknown or fractal. Nevertheless, it is often feasible to utilize time-delayed measurements of the system to construct an augmented state vector. This approach yields an intrinsic coordinate system that is hoped to form an approximate invariant subspace. This technique was discussed in \cref{sec:intro_lorenz}, where it was justified by Takens' embedding theorem \citep{takens2006detecting}. \cite{mezic2004comparison} established the connection between delay embeddings and the Koopman operator via a statistical Takens' embedding theorem.

Employing the same time step for both the delay interval and the frequency of measurements results in a data matrix with a Hankel structure. Hankel matrices have been used in system identification for decades, as seen in the eigensystem realization algorithm \citep{juang1985eigensystem} and singular spectrum analysis \citep{broomhead1989time}. Although these early algorithms were initially developed for linear systems, they have frequently been applied to weakly nonlinear systems as well. The practice of computing DMD on a Hankel matrix was introduced by \cite{tu2014dynamic} and subsequently utilized in the field of neuroscience \citep{brunton2016extracting}. In this section, we focus on two prevalent methods: Hankel-DMD, which is essentially EDMD applied to a dictionary created from time-delay embedding, and the Hankel Alternative View of Koopman (HAVOK) framework, which enhances the DMD model by incorporating a forcing term.

\subsubsection{Hankel Dynamic Mode Decomposition (Hankel-DMD)}
\label{sec:Hankel_DMD}

\textit{Hankel-DMD}, introduced by \cite{arbabi2017ergodic} and closely related to the Prony approximation of the KMD \citep{susuki2015prony}, represents a specialized instance of EDMD where the dictionary is constructed through time-delay embedding. This approach is particularly effective for ergodic systems that exhibit low-dimensional attractors. We saw a slightly generalized variant of this algorithm in \cref{sec:intro_lorenz}, where we employed distinct time steps for sampling trajectories and the lengths of the time delays. Typically, Hankel-DMD utilizes the same time steps for delay embedding and trajectory data collection.

Suppose the map $\Fv$ in \eqref{eq:DynamicalSystem} is ergodic. We can construct a dictionary by starting with an observable $g$ and forming the \textit{Krylov subspace}
$$
V_N=\mathrm{span}\{g,\mathcal{K}g,\mathcal{K}^2g,\ldots,\mathcal{K}^{N-1}g\}.
$$
Given a single trajectory of the observable, $\{g(\xv_0),g(\xv_1),\ldots,g(\xv_{M+N-1})\}$, the matrices $\Psiv_X$ and $\Psiv_Y$ in \eqref{psidef} are given explicitly by the Hankel matrices
\begin{equation}
\label{hankel_data_matrices}
\Psiv_X=\begin{pmatrix} 
g(\xv_0) & g(\xv_1) &  \cdots & g(\xv_{N-1})\\
g(\xv_1) & g(\xv_2) &  \cdots & g(\xv_{N})\\
\vdots & \vdots & \vdots & \vdots  \\
g(\xv_{M-1}) & g(\xv_M) &  \cdots & g(\xv_{M+N-2})
\end{pmatrix},\quad
\Psiv_Y=\begin{pmatrix} 
g(\xv_1) & g(\xv_2) &  \cdots & g(\xv_{N})\\
g(\xv_2) & g(\xv_3) &  \cdots & g(\xv_{N+1})\\
\vdots & \vdots & \vdots & \vdots  \\
g(\xv_{M}) & g(\xv_{M+1}) &  \cdots & g(\xv_{M+N-1})
\end{pmatrix}.
\end{equation}
Applying Birkhoff's ergodic theorem \eqref{eq:birkhoff}, we obtain the convergence specified in \eqref{quad_convergence}. A common simplifying assumption in Hankel-DMD is the existence of a finite-dimensional $\mathcal{K}$-invariant subspace $V$ of $L^2(\Omega,\omega)$ generated by $g$. $\mathcal{K}$-invariance means that $\mathcal{K}V\subset V$ and allows us to study a portion of the spectral properties of $\mathcal{K}$ by restricting to the finite-dimensional subspace $V$. Suppose such a subspace exists and has dimension $k$, then $V_k=V$. We can identify this invariant subspace as $M\rightarrow\infty$ by selecting $N=k$ and employing the aforementioned dictionary. This is proven in \citep{arbabi2017ergodic} and is derived from the ergodic theorem in conjunction with the quadrature interpretation of EDMD. These findings also apply when constructing a Krylov subspace from multiple initial observables $g_1,\ldots,g_p$. Nonetheless, the existence of such a subspace is not guaranteed, and even if it is, the dimension $k$ is typically unknown. Practically, one postulates an \textit{approximate invariant subspace} and truncates to $r \leq N$ modes for the basis by executing an SVD.

\begin{figure}
\centering
\raisebox{-1\height}{\includegraphics[width=0.32\textwidth,trim={0mm 0mm 0mm 0mm},clip]{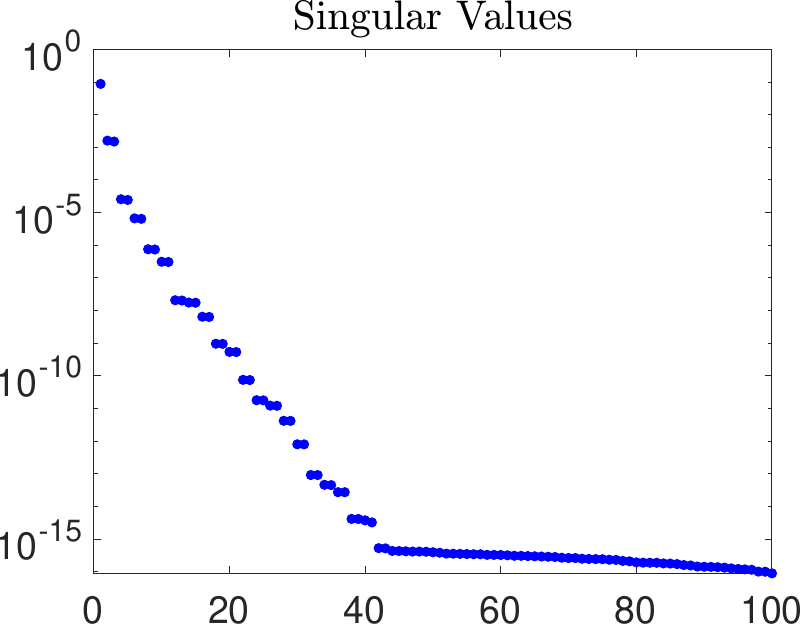}}
\raisebox{-1\height}{\includegraphics[width=0.28\textwidth,trim={0mm 0mm 0mm 0mm},clip]{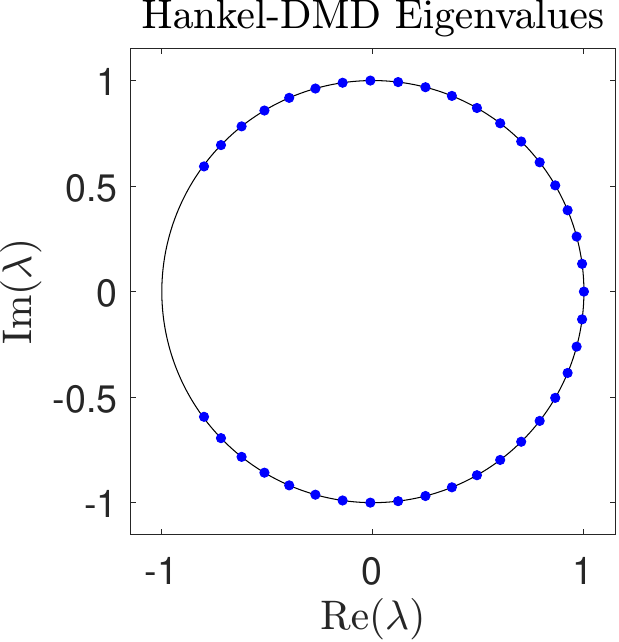}}
\raisebox{-1\height}{\includegraphics[width=0.36\textwidth,trim={0mm 0mm 0mm 0mm},clip]{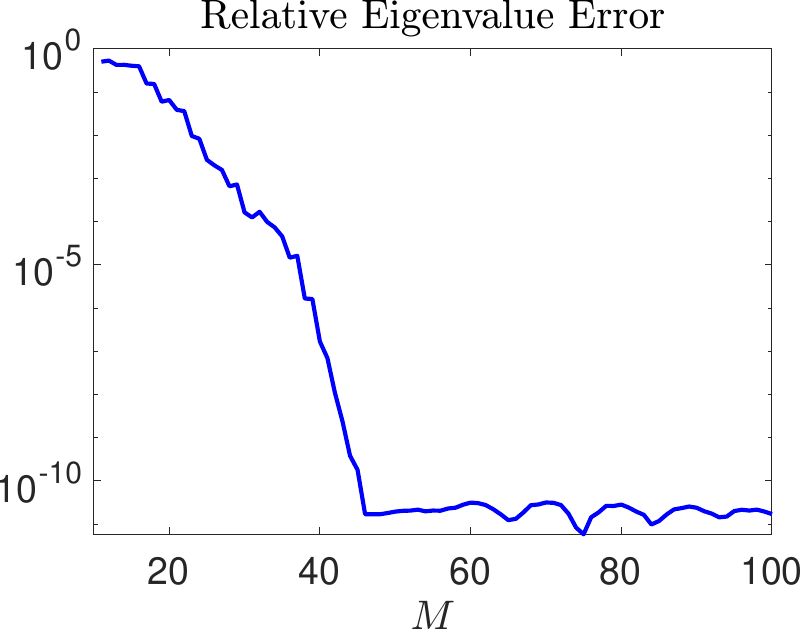}}
\caption{Left: Singular values of $\Psiv_X$ for the Hankel-DMD algorithm applied to the cylinder wake with $N=100$ and $M=120$. In essence, we are building a Krylov subspace by measuring the horizontal component of the velocity field at a single point. Middle: Eigenvalues computed using Hankel-DMD. Right: Convergence to the first $11$ eigenvalues with increasing amount of data $M$.}
\label{hankel_fig}
\end{figure}

As an example, we revisit the cylinder wake discussed in \cref{sec:example_cylinder}. For the observable $g$, we choose the horizontal velocity at a \textit{single point} in the middle of the channel, situated $4D$ downstream from the center of the cylinder. Initially setting $N=100$ and $M=120$, we plot the singular values of the data matrix $\Psiv_X$ on the left side of \cref{hankel_fig}. It is crucial to recognize that although the spectrum is pure point in this example, $g$ does not generate a finite-dimensional invariant subspace since $g$ projects non-trivially onto each eigenspace. Guided by these singular values, we apply \cref{alg:DMD_vanilla} with $r=39$, using the transposes of $\Psiv_X$ and $\Psiv_Y$ as the snapshot matrices. The eigenvalues are illustrated in the middle of \cref{hankel_fig} and should be compared with a subset of the eigenvalues from \cref{wake1}. On the right side of \cref{hankel_fig}, we present the relative $\ell^2$ error for the first $11$ eigenvalues as a function of $M$. The convergence is remarkable. Nonetheless, we must stress that this example is rather straightforward. Systems such as the Lorenz system, as discussed in \cref{sec:intro_lorenz} and tackled in \cite[Section 4.1]{arbabi2017ergodic}, pose a substantially more significant challenge. This is further exemplified by a slow decay of singular values in the data matrices.

\subsubsection{Hankel Alternative View of Koopman (HAVOK)}

Consider the (truncated) SVD of the transpose of the matrix $\Psiv_X$ in \eqref{hankel_data_matrices},
$$
\Psiv_X^\top\approx \Uv \mathbf{\Sigma}\Vv^*,\quad \Uv\in\mathbb{C}^{N\times r},\quad\mathbf{\Sigma}\in\mathbb{R}^{r\times r},\quad\Vv\in\mathbb{C}^{M\times r}.
$$
We can view the columns of the matrix $\Vv$ as coordinates for a state $\vv=[v_1\hspace{2mm}v_2\hspace{2mm}\cdots\hspace{2mm}v_r]$. If our discrete-time dynamical system corresponds to sampling a continuous-time dynamical system, DMD/EDMD results in a linear regression model
$$
\frac{d \vv}{dt}=\hat{\Kv}\vv,\quad\text{for some matrix}\quad\hat{\Kv}\in\mathbb{C}^{r\times r}.
$$
This can be very effective for weakly nonlinear systems \citep{champion2019discovery} and if $r$ is sufficiently large to capture an almost invariant subspace \citep{arbabi2017ergodic}. However, it can be challenging to identify a small (approximately) closed linear model for chaotic systems.

\begin{figure}
\centering
\includegraphics[width=0.85\textwidth,trim={0mm 0mm 0mm 0mm},clip]{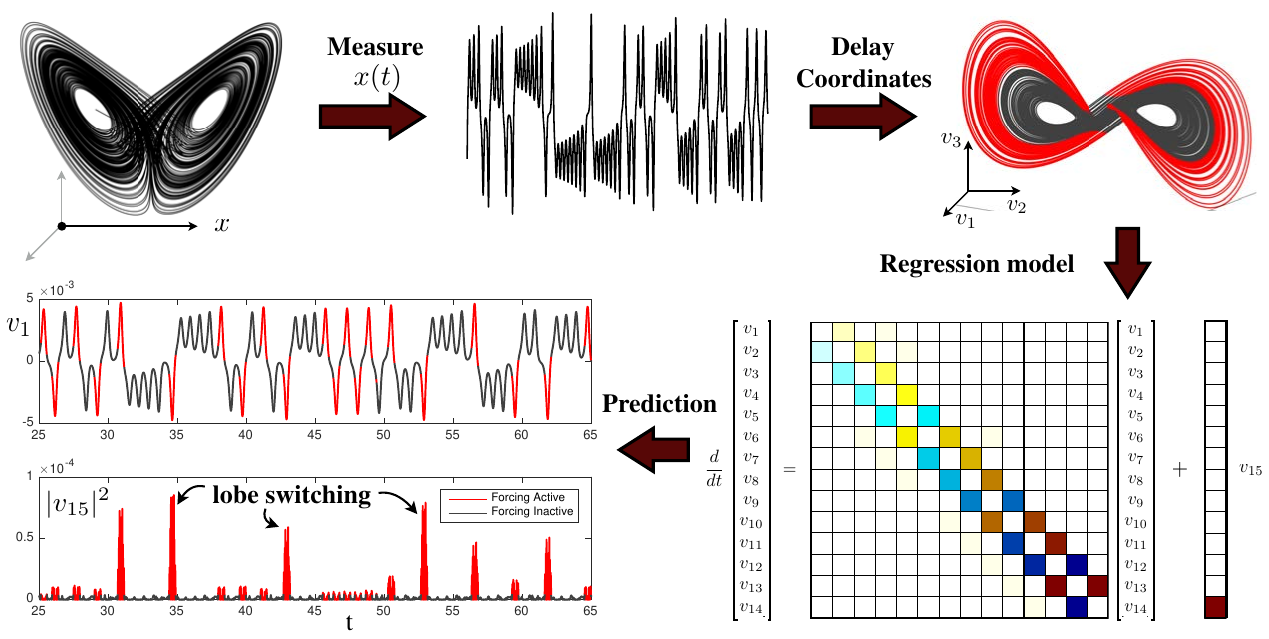}
\caption{Decomposition of chaos into a linear system with forcing. A time series $x(t)$ is stacked into a Hankel matrix, the SVD of which yields a hierarchy of  \emph{eigen} time series that produce a delay-embedded attractor. A best-fit linear regression model is obtained on the delay coordinates; the linear fit for the first $r-1$ variables is excellent, but the last coordinate $v_r$ is not well-modeled as linear. Instead, $v_r$ is an input that forces the first ${r-1}$ variables. Rare forcing events correspond to lobe switching in the chaotic dynamics. Reproduced with permission from \citep{brunton2021modern}, copyright © 2022 Society for Industrial and Applied Mathematics, all rights reserved.}
\label{HAVOK_fig}
\end{figure}

An alternative, known as the \textit{Hankel Alternative View of Koopman} (HAVOK) framework, proposed by \cite{brunton2017chaos}, is to build a linear model on the first $r - 1$ variables $\tilde{\vv}=[v_1\hspace{2mm}v_2\hspace{2mm}\cdots\hspace{2mm}v_{r-1}]$ variables and impose the last variable, $v_r$, as a forcing term:
$$
\frac{d \tilde{\vv}}{dt}=\tilde{\Kv}\tilde{\vv}+\Bv v_r,\quad\text{for some matrices}\quad\tilde{\Kv}\in\mathbb{C}^{r-1\times r-1},\Bv\in\mathbb{C}^{r-1\times 1}.
$$
Here, $v_r$ acts as an input forcing to the linear dynamics of the model, which approximates the nonlinear dynamics of the original system. Typically, the statistics of $v_r$ are non-Gaussian. For instance, in \cref{HAVOK_fig}, we summarize the results for the Lorenz system. The long tails in the statistics of $v_r$ correspond to rare-event forcing that drives lobe switching. For strategies on using HAVOK in systems with multiple time scales, see \citep{champion2019discovery}. \cite{hirsh2021structured} established connections between HAVOK and the Frenet--Serret frame from differential geometry, motivating a more accurate computational modeling approach.

\subsection{Controlling projection errors: Residual Dynamic Mode Decomposition (ResDMD)}
\label{sec:ResDMD}

We saw in \cref{sec:EDMD} that EDMD builds a finite matrix approximation of the Koopman operator. In particular, for a dictionary $\{\psi_1,\ldots,\psi_{N}\}$ forming a finite-dimensional subspace $V_N=\mathrm{span}\{\psi_1,\ldots,\psi_{N}\}$, the EDMD matrix corresponds to the projected operator $\mathcal{P}_{V_{N}}\mathcal{K}\mathcal{P}_{V_{N}}^*$. Care must be taken when discretizing or truncating an infinite-dimensional operator to a finite matrix to compute spectral properties. In general, several well-studied pitfalls include:
\begin{itemize}
\item \textbf{Spectral Pollution:} This term describes false eigenvalues that accumulate at points not in the spectrum as the discretization size increases.
\item \textbf{Spectral Invisibility:} Discretizing an operator can cause us to miss parts of its spectrum, even as the size of the discretization increases.
\item \textbf{Lack of Verification:} Even if a method converges as the discretization parameter grows, how much of the output can we trust for a finite discretization size?
\item \textbf{Continuous Spectra:} Discretizing to a finite matrix results in a discrete set of eigenvalues. How can we recover continuous spectra?
\end{itemize}
We have already encountered these effects in this review (e.g., spectral pollution as discussed in \cref{duffing2}), and they are well-known throughout the Koopman literature. In the following section, we will delve into strategies to mitigate these issues, focusing on controlling projection errors when transitioning from $\mathcal{K}$ to $\mathcal{P}_{V_{N}}\mathcal{K}\mathcal{P}_{V_{N}}^*$. The algorithm that does this is \textit{Residual DMD} (ResDMD), introduced by \cite{colbrook2021rigorous}.

\subsubsection{The algorithm}

The main idea behind ResDMD is to compute an infinite-dimensional residual. We follow the notation of \cref{sec:EDMD} that described EDMD. Consider an observable $g=\Psiv\gv\in V_N$, which we aim to be an approximate eigenfunction of $\mathcal{K}$ with an approximate eigenvalue $\lambda$. For now, the method of determining the pair $(\lambda,g)$ is left unspecified. In connection with pseudospectra and the approximate point spectrum discussed in \cref{sec:spectra_crash_course}, a way to measure the suitability of the candidate pair $(\lambda,g)$ is through the relative residual
\begin{equation}
\label{residual1}
\frac{\|(\mathcal{K}-\lambda I)g\|}{\|g\|}=\sqrt{\frac{\int_{\Omega}|[\mathcal{K}g](\xv)-\lambda g(\xv)|^2\ \mathrm{d} \omega(\xv)}{\int_{\Omega}|g(\xv)|^2\ \mathrm{d} \omega(\xv)}}=\sqrt{\frac{\langle \mathcal{K}g,\mathcal{K}g \rangle-\lambda\langle g,\mathcal{K}g \rangle-\overline{\lambda}\langle \mathcal{K}g,g \rangle+|\lambda|^2\langle g,g \rangle}{\langle g,g \rangle}}.
\end{equation}
For instance, if $\mathcal{K}$ is a normal operator (one that commutes with its adjoint), then
$$
\mathrm{dist}(\lambda,\mathrm{Sp}(\mathcal{K}))=\inf_{f}\frac{\|(\mathcal{K}-\lambda I)f\|}{\|f\|}\leq\frac{\|(\mathcal{K}-\lambda I)g\|}{\|g\|}.
$$
In the case of a non-normal $\mathcal{K}$, the residual in \eqref{residual1} is closely related to the concept of pseudospectra. Adopting the quadrature interpretation of EDMD, we can define a finite data approximation of the relative residual as:
$$
\mathrm{res}(\lambda,g)=\|(\Dv^{1/2}\Psiv_Y-\lambda\Dv^{1/2}\Psiv_X)\gv\|_{\ell^2}/\|\Dv^{1/2}\Psiv_X\gv\|_{\ell^2}.
$$
We then have
\begin{align}
[\mathrm{res}(\lambda,g)]^2&=\frac{\gv^*\left[\Psiv_Y^*\Dv\Psiv_Y-\lambda \Psiv_Y^*\Dv\Psiv_X -\overline{\lambda} \Psiv_X^*\Dv\Psiv_Y + |\lambda|^2\Psiv_X^*\Dv\Psiv_X \right]\gv}{\gv^*\Psiv_X^*\Dv\Psiv_X\gv}\notag\\
&=\frac{\gv\left[\Psiv_Y^*\Dv\Psiv_Y-\lambda\Av^*-\overline{\lambda}\Av+|\lambda|^2\Gv\right]\gv}{\gv^*\Gv\gv},\label{residual2}
\end{align}
where $\Gv$ and $\Av$ are the same matrices from \eqref{eq_EDMD_corr_matrices}. The right-hand side of \eqref{residual2} has an additional matrix $\Lv:=\Psiv_Y^*\Dv\Psiv_Y$. Under the assumption that the quadrature rule converges, this matrix approximates $\mathcal{K}^*\mathcal{K}$:
\begin{equation}
	\label{quad_convergence2}
	\lim_{M\rightarrow\infty}\Lv_{jk} = \langle \mathcal{K}\psi_k,\mathcal{K}\psi_j \rangle.
\end{equation}
Comparing \eqref{residual1} and the square-root of \eqref{residual2}, we observe that
$$
\lim_{M\rightarrow\infty}\mathrm{res}(\lambda,g)=\|(\mathcal{K}-\lambda I)g\|/\|g\|.
$$
Note that there is no approximation or projection on the right-hand side of this equation. Consequently, we can compute an \textit{infinite-dimensional residual} directly using finite matrices, achieving exactness in the limit of large data sets. ResDMD leverages this residual in a suite of algorithms to compute various spectral properties of $\mathcal{K}$, two of which are the focus of the subsequent discussion.

As a first approach, we can implement the EDMD algorithm (\cref{alg:EDMD}) to generate candidate eigenpairs $(\lambda,g)$, followed by the computation of residuals. This process is outlined in \cref{alg:ResDMD1}. Importantly, this approach is not more computationally demanding than EDMD itself. Additionally, it is worth noting that one is not restricted to using EDMD exclusively for selecting candidate eigenpairs; any suitable method can be employed. We can avoid spectral pollution by setting a threshold to discard residuals that exceed a certain tolerance. This also serves as a validation mechanism for the computations. However, it is crucial to recognize that \cref{alg:ResDMD1}, when relying on EDMD for computing candidate eigenpairs, does not inherently circumvent the issue of spectral invisibility. To address this, we need to consider approaches that approximate the pseudospectrum.

\begin{algorithm}[t]
\textbf{Input:} Snapshot data $\Xv\in\mathbb{C}^{d\times M}$ and $\Yv\in\mathbb{C}^{d\times M}$, quadrature weights $\{w_m\}_{m=1}^{M}$, and a dictionary of functions $\{\psi_j\}_{j=1}^{N}$.\\
\vspace{-4mm}
\begin{algorithmic}[1]
\State Compute the matrices $\Psiv_X$ and $\Psiv_Y$ defined in \eqref{psidef} and $\Dv=\mathrm{diag}(w_1,\ldots,w_{M})$.
\State Compute the EDMD matrix $\Kv= (\Dv^{1/2}\Psiv_X)^\dagger \Dv^{1/2}\Psiv_Y\in\mathbb{C}^{N\times N}$.
\State Compute the eigendecomposition $
\Kv\Vv=\Vv\mathbf{\Lambda}$. The columns of $\Vv=[\vv_1\cdots\vv_n]$ are eigenvector coefficients and $\mathbf{\Lambda}$ is a diagonal matrix of eigenvalues $\lambda_1,\ldots,\lambda_n$.
\State For each eigenpair $(\lambda_j,\vv_j)$ compute $
\mathrm{res}(\lambda_j,\Psiv\vv_j){=} \|(\Dv^{1/2}\Psiv_Y{-}\lambda_j\Dv^{1/2}\Psiv_X)\vv_j\|_{\ell^2}/\|\Dv^{1/2}\Psiv_X\vv_j\|_{\ell^2}.$
\end{algorithmic} \textbf{Output:} The eigenvalues $\mathbf{\Lambda}$, eigenvector coefficients $\Vv\in\mathbb{C}^{N\times N}$ and residuals $\{\mathrm{res}(\lambda_j,\Psiv\vv_j)\}$.
\caption{ResDMD for computing residuals \citep{colbrook2021rigorous}.}
\label{alg:ResDMD1}
\end{algorithm}

For computing pseudospectra, working in the standard $\ell^2$ norm is beneficial instead of the norm induced by the matrix $\Gv$. We compute an economy QR decomposition of the data matrix
$$
\Dv^{1/2}\mathbf{\Psi}_X=\Qv\Rv,\quad  \Qv\in\mathbb{C}^{M\times N},\Rv\in\mathbb{C}^{N\times N},
$$
where $\Qv$ has orthonormal columns and $\Rv$ is upper triangular with positive diagonals. Letting $\wv=\Rv\gv$, we have
$$
\|\Dv^{1/2}\Psiv_X\gv\|_{\ell^2}^2=\gv^*\Rv^*\Qv^*\Qv\Rv\gv=\gv^*\Rv^*\Rv\gv=\wv^*\wv=\|\wv\|_{\ell^2}^2.
$$
Consequently, the residual can be expressed as:
\begin{equation}
\label{final_residual}
\mathrm{res}(z,g)=\|(\Dv^{1/2}\Psiv_Y\Rv^{-1}-z\Qv)\wv\|_{\ell^2}/\|\wv\|_{\ell^2}.
\end{equation}
For a given $z\in\mathbb{C}$, our objective is to minimize this residual, which corresponds to finding the smallest singular value of the matrix $(\Dv^{1/2}\Psiv_Y\Rv^{-1}-z\Qv)\in\mathbb{C}^{M\times N}$. Denoting the smallest singular value by $\sigma_{\mathrm{\inf}}$, we must do this for various values of $z$. Given that $M>N$, a computational advantage is gained by considering the $N\times N$ matrix $(\Dv^{1/2}\Psiv_Y\Rv^{-1}-z\Qv)^*(\Dv^{1/2}\Psiv_Y\Rv^{-1}-z\Qv)$ and computing
$$
\sqrt{\sigma_{\mathrm{\inf}}((\Dv^{1/2}\Psiv_Y\Rv^{-1}-z\Qv)^*(\Dv^{1/2}\Psiv_Y\Rv^{-1}-z\Qv))}=\sigma_{\mathrm{\inf}}(\Dv^{1/2}\Psiv_Y\Rv^{-1}-z\Qv).
$$
Typically, computing singular values in this manner is not recommended due to the potential loss of precision owing to the square root. However, in most applications, the resulting error is significantly smaller than the errors inherent in the data matrices or the quadrature approximation of inner products. If precision becomes a concern, $\sigma_{\mathrm{\inf}}(\Dv^{1/2}\Psiv_Y\Rv^{-1}-\lambda\Qv)$ can be computed directly in subsequent algorithms. Since $\Qv^*\Qv=\Iv$, we have
\begin{align*}
&(\Dv^{1/2}\Psiv_Y\Rv^{-1}-z\Qv)^*(\Dv^{1/2}\Psiv_Y\Rv^{-1}-z\Qv)\\
&\quad\quad\quad\quad=(\Rv^*)^{-1}\Psiv_Y^*\Dv\Psiv_Y\Rv^{-1}-z (\Rv^*)^{-1}\Psiv_Y^*\Dv^{1/2}\Qv-\overline{z}\Qv^*\Dv^{1/2}\Psiv_Y\Rv^{-1}+|z|^2\Iv.
\end{align*}
The minimum singular values of this matrix are then computed across a grid of $z$ values. This procedure is detailed in \cref{alg:ResDMD2}, where the approximation of the $\epsilon$-pseudospectrum is defined as the set of grid points where the minimized residual falls below $\epsilon$. If required, the algorithm can also be extended to compute $\epsilon$-pseudoeigenfunctions (discussed in \cref{sec:spectra_crash_course}).

\begin{algorithm}[t]
\textbf{Input:} Snapshot data $\Xv\in\mathbb{C}^{d\times M}$ and $\Yv\in\mathbb{C}^{d\times M}$, quadrature weights $\{w_m\}_{m=1}^{M}$, dictionary of functions $\{\psi_j\}_{j=1}^{N}$, accuracy goal $\epsilon>0$, and grid of points $\{z_\ell\}_{\ell=1}^k\subset\mathbb{C}$.\\
\vspace{-4mm}
\begin{algorithmic}[1]
\State Compute the matrices $\Psiv_X$ and $\Psiv_Y$ defined in \eqref{psidef} and $\Dv=\mathrm{diag}(w_1,\ldots,w_{M})$.
\State Compute an economy QR decomposition $\Dv^{1/2}\mathbf{\Psi}_X=\Qv\Rv$, where $\Qv\in\mathbb{C}^{M\times N},\Rv\in\mathbb{C}^{N\times N}$.
\State Compute $\Cv_2=(\Rv^*)^{-1}\Psiv_Y^*\Dv\Psiv_Y\Rv^{-1}$ and $\Cv_1=\Qv^*\Dv^{1/2}\Psiv_Y\Rv^{-1}$.
\State Compute $\tau_\ell=\sigma_{\mathrm{\inf}}(\Cv_2-z_\ell \Cv_1^*-\overline{z_\ell} \Cv_1+|z_\ell|^2\Iv)$ for $\ell=1,\ldots, k$ ($\sigma_{\mathrm{\inf}}$ is smallest singular value).

\noindent(If wanted, compute the corresponding right-singular vectors $\wv_\ell$ and set $\vv_j=\Rv^{-1}\wv_j$.)
\end{algorithmic} \textbf{Output:} Estimate of the pseudospectrum $\{z_\ell:\tau_\ell<\epsilon\}$ (if wanted, corresponding pseudoeigenfunctions $\{\mathbf{\Psi}\vv_\ell:\tau_\ell<\epsilon\}$).
\caption{ResDMD for computing pseudospectra \citep{colbrook2021rigorous}. One can also compute the singular values directly (of a $\mathbb{C}^{M\times N}$ matrix) without the square root.}
\label{alg:ResDMD2}
\end{algorithm}

\subsubsection{Convergence theory}

\cite{colbrook2021rigorous} present several convergence results concerning ResDMD. We have already discussed that if the quadrature rule underlying EDMD converges,
$$
\lim_{M\rightarrow\infty}\mathrm{res}(\lambda,g)=\|(\mathcal{K}-\lambda I)g\|/\|g\|.
$$
Therefore, we can avoid spectral pollution in the large data limit by selecting eigenpairs with small residuals (as computed in \cref{alg:ResDMD1}). Let $\Gamma^{\epsilon}_{N,M}$ be the output $\{z_\ell:\tau_\ell<\epsilon\}$ of \cref{alg:ResDMD2}. With a minor modification for the boundary case where $\tau=\epsilon$,
$$
\lim_{M\rightarrow\infty}\Gamma^{\epsilon}_{N,M}=:\Gamma^{\epsilon}_{N}\subset \mathrm{Sp}_{\epsilon}(\mathcal{K}).
$$
In other words, ResDMD provides verified approximations of pseudospectra. Moreover, under mild conditions on the dictionary and an $N$-dependent grid $\{z_\ell\}_{\ell=1}^k$,
$$
\lim_{N\rightarrow\infty}\Gamma^{\epsilon}_{N}=\mathrm{Cl}\left(\left\{\lambda\in\mathbb{C}:\exists g\in L^2(\Omega,\omega)\text{ such that }\|g\|=1,\|(\mathcal{K}-\lambda I)g\|<\epsilon\right\}\right).
$$
As $\epsilon\downarrow 0$, the set on the right-hand side converges to the approximate point spectrum $\mathrm{Sp}_{\mathrm{ap}}(\mathcal{K})$. Thus, ResDMD allows us to compute $\mathrm{Sp}_{\mathrm{ap}}(\mathcal{K})$ via a convergent algorithm. \cite{colbrook2021rigorous} further discuss alterations that allow the computation of the full pseudospectrum $\mathrm{Sp}_{\epsilon}(\mathcal{K})$, and consequently the complete spectrum $\mathrm{Sp}(\mathcal{K})$. In summary, ResDMD addresses the challenges of spectral pollution and spectral invisibility, providing a method for verified spectral computations of general Koopman operators.

A careful reader will note that a few of these algorithms require us to take several parameters successively to infinity. This was also the case for EDMD, as discussed in \cref{sec:EDMD_convergence}. These limits do not generally commute, and it may be impossible to rewrite them with fewer limits or develop a different algorithm that uses fewer limits. This is a generic feature of infinite-dimensional spectral problems \citep{colbrookthesis} and has given rise to the \textit{Solvability Complexity Index} \citep{hansen2011solvability,SCI_ref,colbrook2022computation,colbrook2022foundations}. We do not go into the details, but there are many open questions on the foundations of computing spectral properties of Koopman operators. In particular, lower bounds on the number of successive limits needed to compute spectra of Koopman operators is an ongoing research problem (see \cref{sec:open_problems}).

We have yet to discuss continuous spectra, the final pitfall mentioned in the bullet point list at the start of this section. Using ResDMD, \cite{colbrook2021rigorous} also provide an algorithm that computes spectral measures of Koopman operators associated with generic measure-preserving systems. This approach and others for spectral measures are discussed in \cref{sec:methods_for_spec_meas}.

\subsubsection{Examples}

We have already seen an example of ResDMD in action in \cref{sec:duffing_spectral_pollution}. Here, we present some examples from \citep{colbrook2021rigorous,colbrook2023residualJFM}.

\begin{figure}
\centering
\includegraphics[width=1\textwidth,trim={0mm 0mm 0mm 5mm},clip]{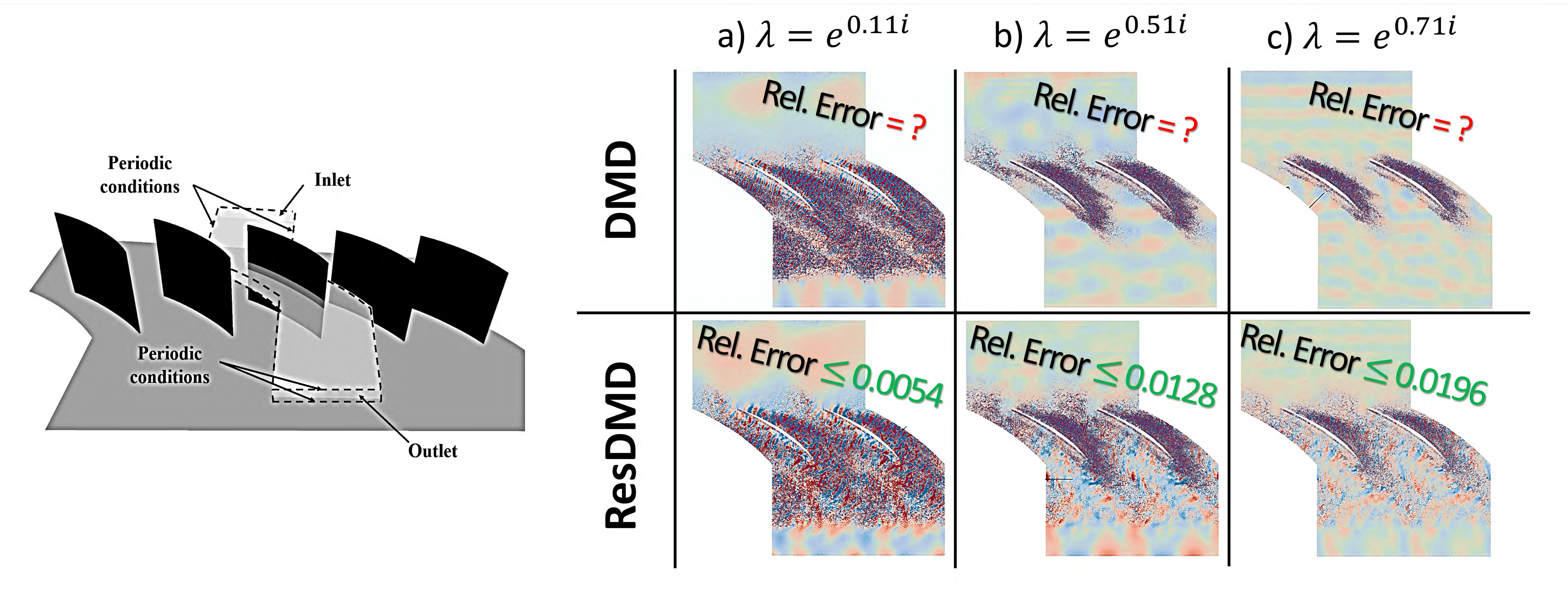}
\caption{Left: Large-scale wall-resolved turbulent flow past a periodic cascade of aerofoils. Right: Comparison of computed Koopman modes using ResDMD and DMD across various frequencies. ResDMD highlights stronger acoustic waves between cascades and larger-scale turbulent fluctuations past the trailing edge, which is crucial for understanding acoustic interactions with engine turbines and nearby structures. The residuals in ResDMD underscore its capability to capture nonlinear dynamics accurately and verifiably. Reproduced with permission from \citep{colbrook2021rigorous}.}
\label{fig:resdmd1}
\end{figure}

The first example we consider is a large-scale wall-resolved turbulent flow past a periodic cascade of aerofoils depicted on the left in \cref{fig:resdmd1}. This setup is motivated by ongoing efforts to mitigate noise sources from aerial vehicles. The data is collected from a high-fidelity simulation solving the fully nonlinear Navier–Stokes equations \citep{koch2021large}, with a Reynolds number of $3.88\times10^5$ and a Mach number of $0.07$. A two-dimensional slice of the pressure field is recorded at $295,122$ points across trajectories of length 798 and sampled every $2\times 10^{-5}$ seconds. ResDMD can be used with kernelized EDMD, and we use $N=250$ functions in our dictionary. \cref{fig:resdmd1} (right) shows the computed Koopman modes for a range of representative frequencies. We also show the corresponding Koopman modes computed using DMD. For the first column, ResDMD shows stronger acoustic waves between the cascades. Detecting these vibrations is essential as they can damage engine turbines \citep{parker1984acoustic}. ResDMD shows larger-scale turbulent fluctuations past the trailing edge for the second and third columns. This can be crucial for understanding acoustic interactions with nearby structures such as subsequent blade rows \citep{woodley1999resonant}. The residuals for ResDMD are small, particularly given the enormous state-space dimension. This example demonstrates two benefits of ResDMD compared with DMD: (1) ResDMD can capture the nonlinear dynamics (just like EDMD), and (2) it computes residuals, thus providing an accuracy certificate.

\begin{figure}
\centering
\includegraphics[width=1\textwidth,trim={0mm 0mm 0mm 0mm},clip]{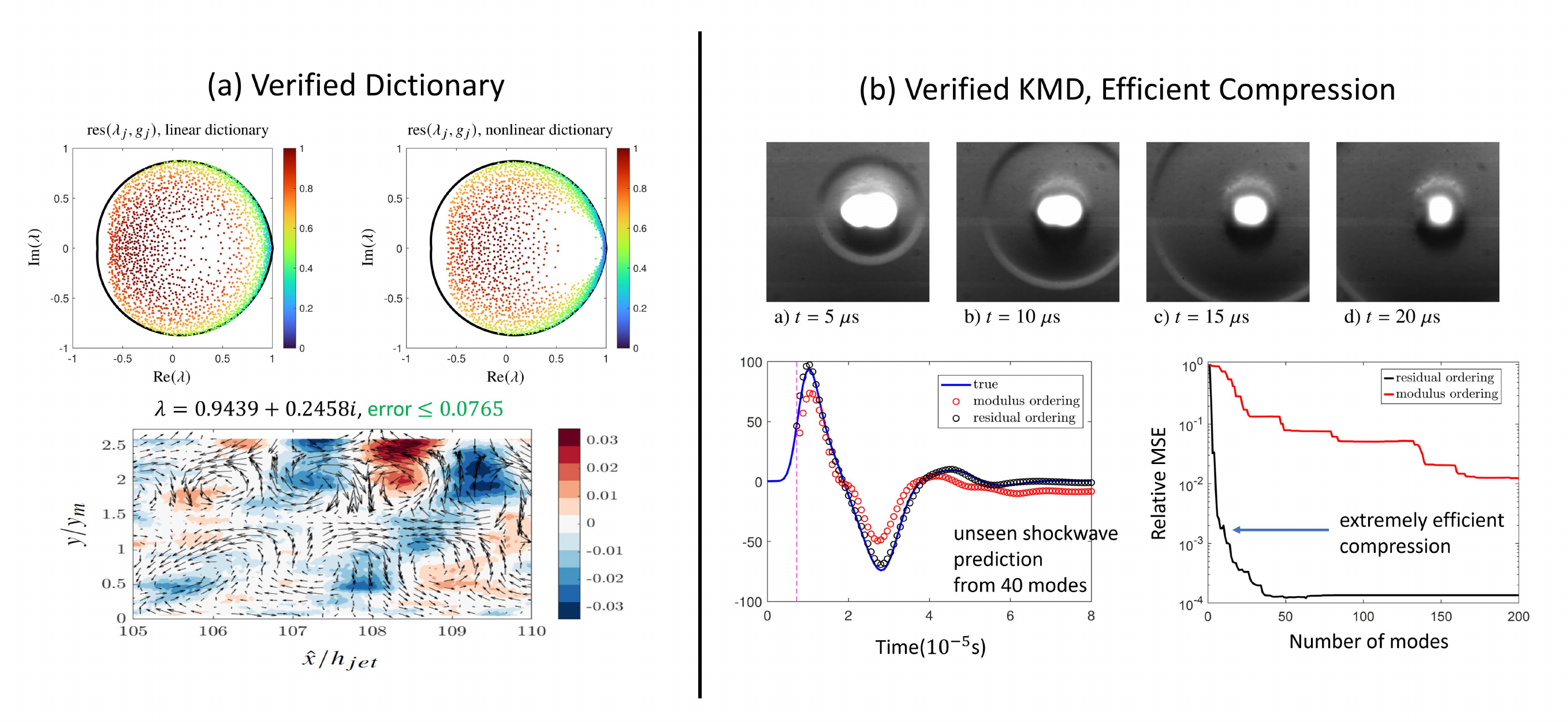}
\caption{ResDMD applications in validating EDMD and KMD. Left: Comparison of two dictionaries for turbulent boundary layer flow analysis, with the nonlinear dictionary showing smaller residuals and revealing verified transient modes (bottom of the figure). Right: Demonstration of KMD's ability to capture a highly nonlinear shockwave, where ordering modes by residual values enables efficient data compression and precise shockwave prediction. Reproduced with permission from \citep{colbrook2023residualJFM}.}
\label{fig:resdmd2}
\end{figure}

With its capability to verifiably compute spectra, ResDMD can be employed for validating dictionaries in methods like EDMD and verifying the efficacy of KMD itself. We present two illustrative examples in \cref{fig:resdmd2}, adapted from \citep{colbrook2023residualJFM}. The comparison of two dictionaries used for analyzing turbulent boundary layer flow is shown on the left. Here, the nonlinear dictionary demonstrates smaller residuals, leading to the identification of verified transient modes, as depicted at the bottom of the figure. On the right, the figure illustrates the proficiency of KMD in capturing a highly nonlinear shockwave. By ordering modes based on their residual values, we achieve efficient data compression and accurate prediction of the shockwave dynamics.

\section{Variants that Preserve Structure}
\label{sec:structure_preserving_methods}

One of the most exciting recent developments in DMD is the introduction of methods that preserve the structures of the underlying dynamical system in \eqref{eq:DynamicalSystem}. When studying a system from a data-driven perspective, it is often the case that one possesses partial knowledge of the system's underlying physics. Methods that leverage this structure typically exhibit a greater resistance to noise, better generalization, and demand less data for training. Structure-preserving algorithms have a deep-rooted history in geometric integration \citep{3-540-30663-3} and have recently gained traction in data-driven methods \citep{celledoni2021structure,greydanus2019hamiltonian,hernandez2021structure,hesthaven2022reduced,karniadakis2021physics,loiseau2018constrained,otto2023unified}. In the context of DMD, this area is burgeoning. We will concentrate on three methods:
\begin{itemize}
	\item \textbf{Physics-Informed DMD:} This provides a framework for incorporating symmetries into DMD through additional constraints in the least-squares problem \eqref{DMD_opt_vanilla}. The original paper focused on five fundamental physical principles: conservation,
self-adjointness,
localization,
causality,
and shift-equivariance. The idea is far more general and has ushered in a new wave of DMD methods.
	\item \textbf{Measure-Preserving EDMD:} This enforces measure-preserving EDMD truncations, leading to a Galerkin method whose eigendecomposition converges to the spectral quantities of Koopman operators (including spectral measures and continuous spectra) for general measure-preserving dynamical systems.	Like EDMD, it can be used with any dictionary of observables and with different data types. Preserving the measure is crucial for convergence, recovering the correct dynamical behavior, stability, robustness to noise, and improved qualitative and long-time behavior.
	\item \textbf{Compactification:} These methods for continuous-time measure-preserving systems are based on the compactification of the Koopman generator or its resolvent. They automatically lead to skew-adjoint approximations whose spectral properties converge to that of the Koopman generator. Additionally, approximations are expressed in a well-conditioned basis of kernel eigenvectors computed from trajectory data.
\end{itemize}
Subsequently, we will discuss additional DMD methods based on preserving structure. The methods we discuss open the door to future extensions to more general structure-preserving methods for Koopman operators and data-driven dynamical systems.

\subsection{Physics-Informed Dynamic Mode Decomposition (piDMD)}

\subsubsection{The framework}

\textit{Physics-Informed DMD} (piDMD), introduced by \cite{baddoo2023physics}, provides an overarching framework for integrating physical principles -- such as symmetries, invariances and conservation laws -- into DMD. The idea is to replace the optimization problem in \eqref{DMD_opt_vanilla} by a \textit{constrained} optimization problem
\begin{equation}
\label{DMD_opt_piDMD}
\min_{\Kv_{\mathrm{piDMD}}\in\mathcal{M}} \left\|\Yv-\Kv_{\mathrm{piDMD}}\Xv\right\|_{\mathrm{F}}.
\end{equation}
The matrix manifold $\mathcal{M}$ is dictated by the known physics of the system in \eqref{eq:DynamicalSystem}. One selects $\mathcal{M}$ so that its members satisfy certain symmetries of the system. The optimization problem in \eqref{DMD_opt_piDMD} is known as a Procrustes problem\footnote{In Greek mythology, Procrustes was a bandit who would stretch or amputate the limbs of his victims to force them to fit onto his bed. Herein, $\Xv$ plays the role of Procrustes' victim, $\Yv$ is the bed, and $\Kv_{\mathrm{piDMD}}$ is the `treatment' (stretching or amputation).}, which comprises of finding the optimal transformation between two matrices subject to certain constraints. Numerous exact solutions exist for Procrustes problems, including the notable cases of orthogonal matrices \citep{schonemann1966generalized}, and symmetric matrices \citep{higham1988symmetric}. When exact solutions are not possible, algorithmic solutions can be effective \citep{boumal2014manopt}. Procrustes analysis finds relevance in many fields, as detailed in the monograph of \cite{gower2004procrustes}. 

To apply piDMD, we first identify the system's known or suspected physical properties. Once the physical principles we wish to enforce are determined, these laws must be translated into the matrix manifold where the linear model will be constrained. With a defined target matrix manifold, we numerically solve the relevant Procrustes problem in \eqref{DMD_opt_piDMD}. The concluding step encompasses extracting physical information from the refined model. For instance, one might analyze the spectrum, DMD modes, and the related KMD.

\cite{baddoo2023physics} focus on five fundamental physical principles: conservation,
self-adjointness,
localization,
causality,
and shift-equivariance.
Several closed-form solutions and efficient algorithms for the corresponding piDMD optimizations are derived. With fewer degrees of freedom, piDMD models are typically less prone to overfitting, require less training data, and are often less computationally expensive to build than standard DMD models. This reduction in the size of required training data is connected with the problem of matrix recovery from matrix-vector products, whereby enforcing structures reduces the number of queries needed \citep{halikiasstructured}. A fundamental issue related to the DMD algorithm is the fact that low-rank matrices are not provably recoverable from snapshot pairs (without access to adjoints) until there are at least as many pairs as state dimensions \citep[Thm. 2.5]{otto2023model}.

\subsubsection{Examples}

\begin{figure}[t]
	\centering
	\includegraphics[width=\linewidth]{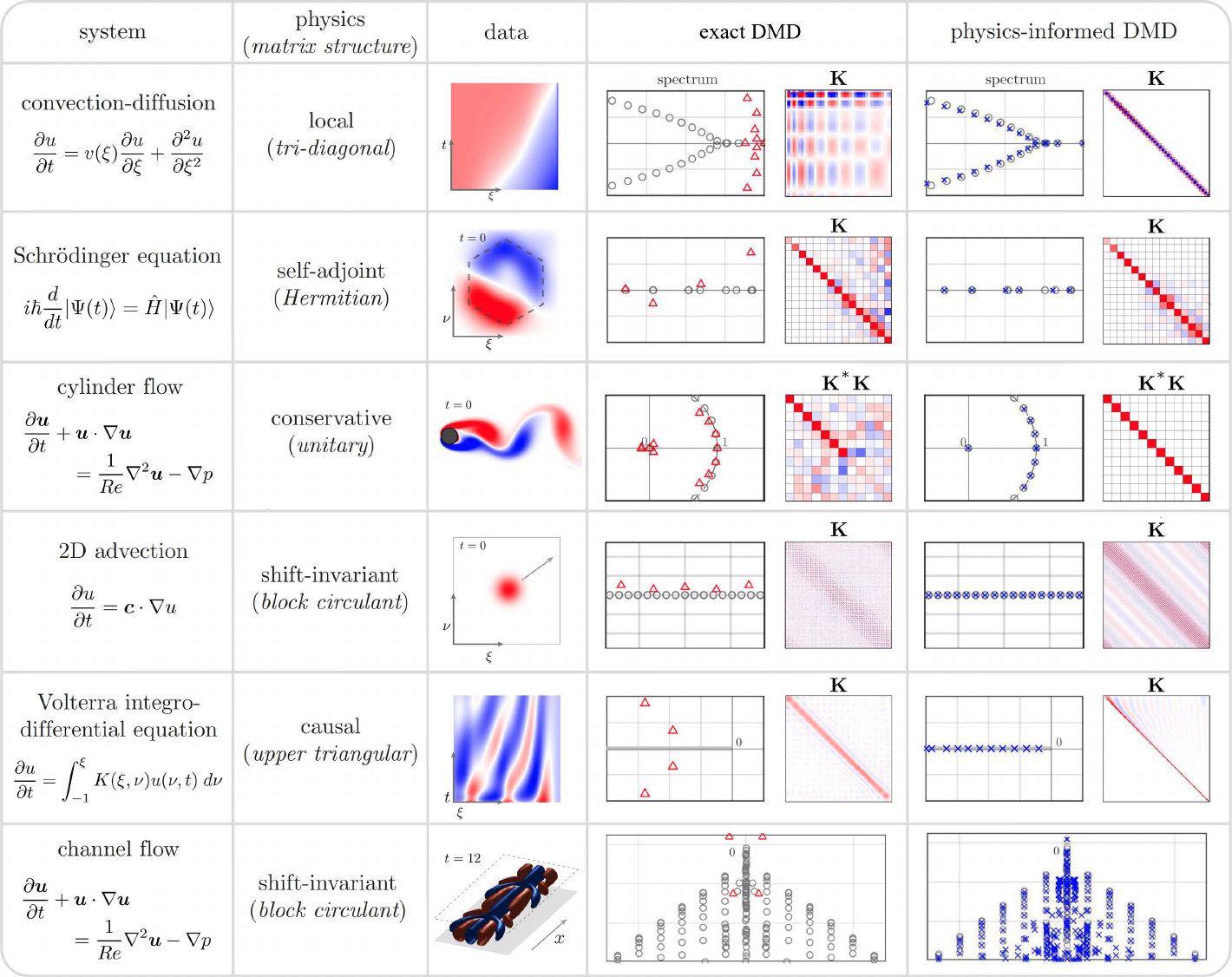}
	\caption{A comparison of the models learned by exact DMD (\cref{alg:DMD_vanilla}) and piDMD for a range of applications. The structure of the model matrices is also illustrated. In the spectrum subplots, the true eigenvalues are shown as 
{\protect\tikz \protect\draw[draw = gray,line width = 1pt] (.5,0) circle (3pt);},
the DMD eigenvalues as
{\protect\tikz \protect\draw[draw=red,thick] (0,0) --
	(0.2cm,0) -- (0.1cm,0.2cm) -- (0,0);}
, and the piDMD eigenvalues as 
{\protect\tikz \protect\draw[draw = blue, thick] (.5,0) node[blue, cross=3.5pt] {};}. In each case,
	the eigenvalues of piDMD are more accurate than exact DMD. Reproduced with permission from \citep{baddoo2023physics}.}
	\label{fig:piDMD}
\end{figure}

To showcase the breadth of piDMD, \cref{fig:piDMD} shows six physical examples. \cite{baddoo2023physics} provide full experimental details for each example. Each row corresponds to a different system, and the corresponding constraint is listed in the second column. Exact DMD (\cref{alg:DMD_vanilla}) is compared to piDMD in terms of the computed matrix $\Kv$ and the eigenvalues. In general, constraining the matrix $\Kv$ to lie on the appropriate manifold $\mathcal{M}$ leads to more accurate approximations of the eigenvalues. The advantage of preserving structure is striking!

\subsubsection{Future work}
\label{sec:piDMD_future}

We have only started to tap the potential of adding constraints in the optimization problem in \eqref{DMD_opt_vanilla}. This idea will likely be an active research area over the next few years. With that in mind, it is worth mentioning several challenges and directions of future work pointed out by \cite{baddoo2023physics}:
\begin{itemize}
	\item \textbf{Knowing the physics:} In some scenarios where the physics is poorly understood, determining suitable physical laws to impose on the model can be challenging. Is it possible to learn symmetries and then incorporate them as constraints?
\item \textbf{Complicated manifolds:} For problems with intricate geometries and multiple dimensions, interpreting the physical principle as a matrix manifold can be a roadblock, as the manifold can become exceedingly complicated.
\item \textbf{Regularizers:} In many applications, such as when the data are very noisy or the physical laws and constraints are only approximately understood, it may be more appropriate to merely encourage $\Kv$ towards $\mathcal{M}$, e.g., through a regularizer.
\item \textbf{Nonlinear observables:} It is not always clear how to extend the approach of piDMD to nonlinear observables and EDMD. Such an approach is crucial for strongly nonlinear systems to maintain the connection with Koopman operators. For example, if the dictionary consists of the state vector $\xv$, an upper triangular matrix $\Kv$ can have a clear meaning in terms of causality. But how should one incorporate causality into other choices of dictionaries? The manifold can depend on the chosen dictionary in a highly complex manner.
\item \textbf{Convergence:} In connection with the previous point in this list, proving the convergence of piDMD in the large data limit or large dictionary limit typically requires a Galerkin interpretation as in \cref{sec:Galerkin_perspective}. This connection is not always immediate.
\end{itemize}
We expect these last two points, in particular, to lead to many exciting future works.

\subsection{Measure-Preserving Extended Dynamic Mode Decomposition (mpEDMD)}
\label{sec:mpEDMD}

\textit{Measure-Preserving EDMD} (mpEDMD), introduced by \cite{colbrook2023mpedmd}, enforces that the EDMD approximation is measure-preserving. The system being measure-preserving is equivalent to the Koopman operator $\mathcal{K}$ being an isometry. Namely, $\|\mathcal{K}g\|=\|g\|$ for any observable $g\in L^2(\Omega,\omega)$. The mpEDMD algorithm is simple and robust, with no tuning parameters. We outline the method, discuss its convergence properties, and end with two examples. Note that we do not need to assume the system is ergodic or invertible.

\subsubsection{The algorithm}
\label{sec:mpedmd_alg}

We follow the notation of \cref{sec:EDMD} that described EDMD. Recall that we have a dictionary $\{\psi_1,\ldots,\psi_{N}\}$, i.e., a list of observables, in the space $L^2(\Omega,\omega)$. These observables form a finite-dimensional subspace $V_N=\mathrm{span}\{\psi_1,\ldots,\psi_{N}\}$. Our starting point is the observation that the Gram matrix $\Gv=\Psiv_X^*\Dv\Psiv_X$ in \eqref{eq_EDMD_corr_matrices} provides an approximation of the inner product $\langle \cdot,\cdot\rangle$ on $L^2(\Omega,\omega)$. Namely, we have the following inner product induced by $\Gv$:
\begin{equation}
\label{inner_product_G}
\hv^*\Gv\gv=\sum_{j,k=1}^N\overline{\hv_j}\gv_k\Gv_{j,k}\approx \sum_{j,k=1}^N  \overline{\hv_j}\gv_k\langle \psi_k,\psi_j \rangle=\langle \Psiv \gv,\Psiv \hv \rangle.
\end{equation}
If the convergence in \eqref{quad_convergence} holds, then the left-hand side of \eqref{inner_product_G} converges to the right-hand side as $M\rightarrow\infty$. Hence, if $g=\Psiv\gv\in V_N$ and we approximate the action of $\mathcal{K}$ on $V_N$ by a matrix $\Kv$,
$$
\|g\|^2\approx\gv^*\Gv\gv,\quad \|\mathcal{K}g\|^2\approx\|\Psiv \Kv\gv\|^2\approx \gv^*\Kv^*\Gv\Kv\gv.
$$
Since $\mathcal{K}$ is an isometry, $\|g\|^2=\|\mathcal{K}g\|^2.$ Therefore, it is natural to enforce
$$
	\gv^*\Gv\gv=\gv^*\Kv^*\Gv\Kv\gv\quad \forall \gv\in\mathbb{C}^N.
$$
This condition holds if and only if $\Kv^*\Gv\Kv=\Gv$. Returning to the optimization problem in \eqref{eq:ContinuousLeastSquaresProblem}, we now make two changes. First, we set $\Cv=\Gv^{1/2}$ so that $\|\Cv\gv\|_{\ell^2}=\sqrt{\gv^*\Gv\gv}\approx\|g\|$. Second, we enforce the additional constraint $\Kv^*\Gv\Kv=\Gv$. This leads us to the optimization problem
\begin{equation}
	\label{EDMD_opt_prob2bb}
	\underset{\underset{\Kv^*\Gv\Kv=\Gv}{\Kv\in\mathbb{C}^{N\times N}}}{\min} \int_\Omega \left\|\Psiv(\Fv(\xv))\Gv^{-1/2} - \Psiv(\xv)\Kv\Gv^{-1/2}\right\|^2_{\ell^2}\ \mathrm{d}\omega(\xv).
\end{equation}
In a nutshell, we enforce that our Galerkin approximation is an isometry with respect to the learned, data-driven inner product induced by $\Gv$. After applying the quadrature rule we used for EDMD, the discretized version of \eqref{EDMD_opt_prob2bb} is
\begin{equation}
\label{EDMD_opt_prob3}
\underset{\underset{\Kv^*\Gv\Kv=\Gv}{\Kv\in\mathbb{C}^{N\times N}}}{\min}\sum_{m=1}^{M} w_m\left\|\Psiv(\yv^{(m)})\Gv^{-1/2}-\Psiv(\xv^{(m)})\Kv\Gv^{-1/2}\right\|^2_{\ell^2}.
\end{equation}
Letting $\Kv=\Gv^{-1/2}\Bv\Gv^{1/2}$ for some matrix $\Bv$, the problem in \eqref{EDMD_opt_prob3} is equivalent to
\begin{equation}
\label{EDMD_opt_prob4}
\underset{\underset{\Bv^*\Bv=\Iv}{\Bv\in\mathbb{C}^{N\times N}}}{\min} \left\|\Dv^{1/2}\Psiv_X\Gv^{-1/2}\Bv-\Dv^{1/2}\Psiv_Y\Gv^{-1/2}\right\|^2_{\mathrm{F}},
\end{equation}
where $\Iv$ denotes the identity matrix. The optimization problem in \eqref{EDMD_opt_prob4} is known as the \textit{orthogonal Procrustes problem} \citep{schonemann1966generalized,arun1992unitarily}. The predominant method for computing a solution is via an SVD. First, we compute an SVD of
$$
\Gv^{-1/2}\Psiv_Y^*\Dv\Psiv_X\Gv^{-1/2}=\Gv^{-1/2}\Av^*\Gv^{-1/2}=\Uv_1\mathbf{\Sigma} \Uv_2^*,
$$
where $\Av=\Psiv_X^*\Dv\Psiv_Y$ is the matrix from \eqref{eq_EDMD_corr_matrices}. A solution of \eqref{EDMD_opt_prob4} is then $\Bv=\Uv_2\Uv_1^*$ and we take $\Kv=\Gv^{-1/2}\Uv_2\Uv_1^*\Gv^{1/2}$.

Since $\Kv$ is similar to a unitary matrix, its eigenvalues lie along the unit circle. For stability purposes, the best way to compute the eigendecomposition of $\Kv$ is to do so for the unitary matrix $\Uv_2\Uv_1^*$. To numerically ensure an orthonormal basis of eigenvectors, we use MATLAB's \texttt{schur} command in the examples of this paper. It is also beneficial to replace the square root $\Gv^{1/2}$ with a suitable upper triangular matrix $\Rv$ such that $\Gv=\Rv^*\Rv$. Such an upper triangular matrix can be computed using an economy QR decomposition of the data matrix as
$$
\Dv^{1/2}\mathbf{\Psi}_X=\Qv\Rv,\quad  \Qv\in\mathbb{C}^{M\times N},\Rv\in\mathbb{C}^{N\times N},
$$
where $\Qv$ has orthonormal columns and $\Rv$ is upper triangular with positive diagonals. This leads to a mathematically equivalent algorithm but is faster and more numerically robust in practice.\footnote{I am indebted to Zlatko Drmač for pointing this out.} The computation of $\Kv$ and its eigendecomposition is summarized in \cref{alg:mpEDMD}. Arguing as we did for EDMD, we obtain a KMD via
$$
g\approx \Psiv \Vv\left[\Vv^{-1}(\Dv^{1/2}\Psiv_X)^{\dagger}\Dv^{1/2}\left(g(\xv^{(1)}),\ldots,g(\xv^{(M)})\right)^\top\right], \quad g\in L^2(\Omega,\omega).
$$
Explicitly applied to the state vector $\xv$, we have (transposed) Koopman modes
$$
\mathbf{\Phi}^\top=\Vv^{-1}(\Dv^{1/2}\Psiv_X)^{\dagger}\Dv^{1/2}\left(\xv^{(1)},\ldots,\xv^{(M)}\right)^\top\in\mathbb{C}^{N\times d}.
$$
Note that mpEDMD can be used with generic choices of dictionary that generate $\Gv$ and $\Av$.

\begin{algorithm}[t]
\textbf{Input:} Snapshot data $\Xv\in\mathbb{C}^{d\times M}$ and $\Yv\in\mathbb{C}^{d\times M}$, quadrature weights $\{w_m\}_{m=1}^{M}$, and a dictionary of functions $\{\psi_j\}_{j=1}^{N}$.\\
\vspace{-4mm}
\begin{algorithmic}[1]
\State Compute the matrices $\Psiv_X$ and $\Psiv_Y$ defined in \eqref{psidef} and $\Dv=\mathrm{diag}(w_1,\ldots,w_{M})$.
\State Compute an economy QR decomposition $\Dv^{1/2}\mathbf{\Psi}_X=\Qv\Rv$, where $\Qv\in\mathbb{C}^{M\times N},\Rv\in\mathbb{C}^{N\times N}$.
\State Compute an SVD of $(\Rv^{-1})^{*}\Psiv_Y^*\Dv^{1/2}\Qv=\Uv_1\mathbf{\Sigma} \Uv_2^*$.
\State Compute the eigendecomposition $\Uv_2\Uv_1^*=\hat{\Vv}\mathbf{\Lambda} \hat{\Vv}^*$ (via a Schur decomposition).
\State Compute $\Kv=\Rv^{-1}\Uv_2\Uv_1^*\Rv$ and $\Vv=\Rv^{-1}\hat{\Vv}$.
\end{algorithmic} \textbf{Output:} Koopman matrix $\Kv$, with eigenvectors $\Vv$ and eigenvalues $\mathbf{\Lambda}$.
\caption{The mpEDMD algorithm \citep{colbrook2023mpedmd}.}
\label{alg:mpEDMD}
\end{algorithm}

Finally, the relationship between mpEDMD and piDMD is worth commenting on. For conservative systems, piDMD enforces the DMD matrix in \eqref{DMD_opt_vanilla} to be orthogonal and uses linear observables. This implicitly assumes that these linear observables (and the coordinates used) are orthonormal in $L^2(\Omega,\omega)$, an assumption that typically does not hold. In contrast, mpEDMD works in a data-driven inner product space induced by $\Gv$. The resulting Gram matrix of the observables must be included in a measure-preserving discretization; otherwise, the wrong measure may be preserved (see the example of turbulent flow in \citep{colbrook2023mpedmd} where mpEDMD and piDMD are contrasted). Thus, we can think of the relationship between mpEDMD and piDMD as akin to the relationship between EDMD and DMD (see the discussion in \cref{sec:galerkin_interp,sec:EDMD_algorithm}), with an additional difference arising from the use of the inner product arising from the Gram matrix $\Gv$.

\subsubsection{Convergence theory}
\label{sec:mpEDMD_convergence}

Several convergence results for mpEDMD are proven in \citep{colbrook2023mpedmd}. First, echoing \cref{sec:EDMD_convergence}, we can consider the two limits $M\rightarrow\infty$ and $N\rightarrow\infty$. Assuming that the quadrature rule underlying EDMD converges, i.e., \eqref{quad_convergence} holds, the EDMD matrix corresponds to $\mathcal{P}_{V_{N}}\mathcal{K}\mathcal{P}_{V_{N}}^*$. In contrast, the mpEDMD matrix corresponds to the unitary part of a \textit{polar decomposition} of $\mathcal{P}_{V_{N}}\mathcal{K}\mathcal{P}_{V_{N}}^*$. Call this matrix $\Kv_N$. Under a natural density assumption of $V_N$ as $N\rightarrow\infty$, $\Kv_N$ converges strongly to $\mathcal{K}$, meaning that \eqref{SOT_Koopman} holds.

We can consider the spectral measures from \cref{sec:spectra_crash_course} for a measure-preserving system. These spectral measures provide a diagonalization of the Koopman operator $\mathcal{K}$ and form the foundation of the KMD. As the dictionary $\{\psi_1,\ldots,\psi_{N}\}$ used in EDMD becomes richer, the spectral measures computed by EDMD do not typically converge in any sense to that of $\mathcal{K}$.\footnote{Even more fundamentally, the eigenvalues of EDMD typically lie within and accumulate within the unit disk. So, the measures are not even on the same space.} This contrasts the spectral measures of mpEDMD \citep{colbrook2023mpedmd}. The critical step in making this convergence work is that mpEDMD provides a unitary Galerkin approximation of $\mathcal{K}$.

The mpEDMD algorithm leads to the following approximations of spectral measures, where $\Kv$, $\Vv=[\vv_1 \,\cdots\,\vv_N]$ and $\mathbf{\Lambda}=\mathrm{diag}(\lambda_1,\ldots,\lambda_N)$ are the outputs of \cref{alg:mpEDMD}. To approximate the spectral measure $\mathcal{E}$, we consider the spectral measure, $\mathcal{E}_{N,M}$, of the matrix $\Kv$ on the Hilbert space $\mathbb{C}^N$ with the inner product in \eqref{inner_product_G} induced by $G$:
$$
\mathrm{d}\mathcal{E}_{N,M}(\lambda)=\sum_{j=1}^{N}\vv_j\vv_j^*\Gv\delta(\lambda-\lambda_j)\,\mathrm{d}\lambda.
$$
Let $g\in L^2(\Omega,\omega)$ with $\|g\|=1$. We approximate $\mu_g$ by $\smash{\mu_{\gv}^{(N,M)}}$, where
$$
\mathrm{d}\mu_{\gv}^{(N,M)}(\lambda)=\sum_{j=1}^N\delta(\lambda-\lambda_j)|\vv_j^*\Gv\gv|^2\mathrm{d}\lambda
$$
and $\gv$ is normalized so that $\gv^*\Gv\gv=1$. Since $\{\Gv^{1/2}\vv_j\}_{j=1}^N$ is an orthonormal basis for $\mathbb{C}^N$, $\mu_{\gv}^{(N,M)}$ is a probability measure on the unit circle $\mathbb{T}$.

The most natural way for measures to converge is in a \textit{weak sense}. We say that a sequence of measures $\mu_n$ converges weakly to a measure $\mu$ on $\mathbb{T}$ if for any continuous function $\phi:\mathbb{T}\rightarrow\mathbb{C}$,
$$
\lim_{n\rightarrow\infty}\int_{\mathbb{T}} \phi(\lambda)\ \mathrm{d}\mu_n(\lambda)=\int_{\mathbb{T}} \phi(\lambda)\ \mathrm{d}\mu(\lambda).
$$
This convergence is captured by the so-called Wasserstein 1 metric between probability measures:
$$
W_1(\mu,\nu)\coloneqq\sup\left\{\int_{\mathbb{T}}\phi(\lambda)\,d(\mu-\nu)(\lambda):\phi:\mathbb{T}\rightarrow\mathbb{R}\text{ Lip. cts., Lip. constant $\leq1$}\right\}.
$$
Under mild conditions on the dictionary as $N\rightarrow\infty$, mpEDMD has the following convergence properties. If $\lim_{N\rightarrow\infty}\Psiv \gv_N=g$ and $\phi:\mathbb{T}\rightarrow\mathbb{C}$ is continuous, then
$$
\lim_{N\rightarrow\infty}\limsup_{M\rightarrow\infty}\left\|\int_{\mathbb{T}} \phi(\lambda)\ \mathrm{d}\mathcal{E}(\lambda)g-\Psiv \int_{\mathbb{T}}\phi(\lambda)\ \mathrm{d}\mathcal{E}_{N,M}(\lambda) \gv_N\right\|=0.
$$
Moreover, for the scalar-valued spectral measures,
\begin{equation}\label{weak_weak_convergence2}
\lim_{N\rightarrow\infty}\limsup_{M\rightarrow\infty}W_1\left(\mu_g,\smash{\mu_{\gv}^{(N,M)}}\right)=0.
\end{equation}
Furthermore, if $\{g,\mathcal{K}g,\ldots,\mathcal{K}^{L_N-1}g\}\subset V_N$ (time-delay embedding), then
\begin{equation}\label{weak_weak_convergence3}
\limsup_{M\rightarrow\infty}W_1\left(\mu_g,\smash{\mu_{\gv}^{(N,M)}}\right)\lesssim{\log(L_N)}/{L_N}.
\end{equation}
The bound in \eqref{weak_weak_convergence3} provides an explicit convergence rate in the number of delays used.

Further properties of mpEDMD proven in \citep{colbrook2023mpedmd} include respecting invariance subspace properties of $\mathcal{K}$, well-conditioning of the matrix $\Kv$ and its eigendecomposition (these do not hold for EDMD in general), convergence of KMDs, and 
$$
\lim_{N\rightarrow\infty}\limsup_{M\rightarrow\infty} \sup_{\lambda\in \mathrm{Sp}_{\mathrm{ap}}(\mathcal{K})}\mathrm{dist}(\lambda,\{\lambda_1,\ldots,\lambda_N\})=0.
$$
In other words, we avoid spectral invisibility and do not miss parts of the spectrum. We can also combine with the techniques of \cref{sec:ResDMD} to avoid spectral pollution. Finally, in connection with \cref{sec:noise_DMD_robustness}, the solution to the orthogonal Procrustes problem \eqref{EDMD_opt_prob4} is also the solution to the corresponding constrained total least squares problem \citep{arun1992unitarily}. Hence, in a similar vein to tlsDMD in \cref{sec:tlsDMD}, mpEDMD is optimally robust when noise is present in both data matrices in \eqref{EDMD_opt_prob4} \citep{van1991total}.

\subsubsection{Examples}

We consider two examples of mpEDMD. The first shows the convergence to spectral measures for a system with continuous spectra. The second shows the conservation of energy and statistics for a turbulent boundary layer flow, where the snapshots are collected experimentally.

\paragraph{Convergence to continuous spectra}

We first revisit the Lorenz system from \cref{sec:intro_lorenz}, but now with a discrete time step of $\Delta t=0.1$. An arbitrary observable is chosen as
$$
g(\xv) =g(x,y,z)= c\tanh((xy-3z)/5),
$$
where $c$ is a normalization constant ensuring $\|g\|=1$. We employ delay-embedding to construct a Krylov subspace $V_N=\{g,\mathcal{K}g,\ldots,\mathcal{K}^{N-1}g\}$. The matrices $\Psiv_X$ and $\Psiv_Y$ are computed by evaluating $g$ pointwise at the snapshot matrices of $\xv$. A set of $M=10^4$ snapshots is collected along a single trajectory following an initial burn-in period. It is important to recall that the spectrum of the Koopman operator is continuous, featuring an embedded trivial eigenvalue at $\lambda=1$. Therefore, we demonstrate the convergence of the mpEDMD approximation of $\mu_g$. For visualization purposes, we transition from variables in $\mathbb{T}$ to complex-arguments in the interval $[-\pi,\pi)$. For a probability measure $\mu$ on $\mathbb{T}$, its cumulative distribution function (cdf) on $[-\pi,\pi)$ is defined as
$$
F_\mu(\theta)=\mu\left(\{\exp(it):\pi\leq t\leq \theta\}\right).
$$
One can express the metric $W_1$ in terms of these cdfs \citep{hundrieser2022statistics}.

\begin{figure}
\centering
\raisebox{-0.5\height}{\includegraphics[width=0.24\textwidth,trim={0mm 0mm 0mm 0mm},clip]{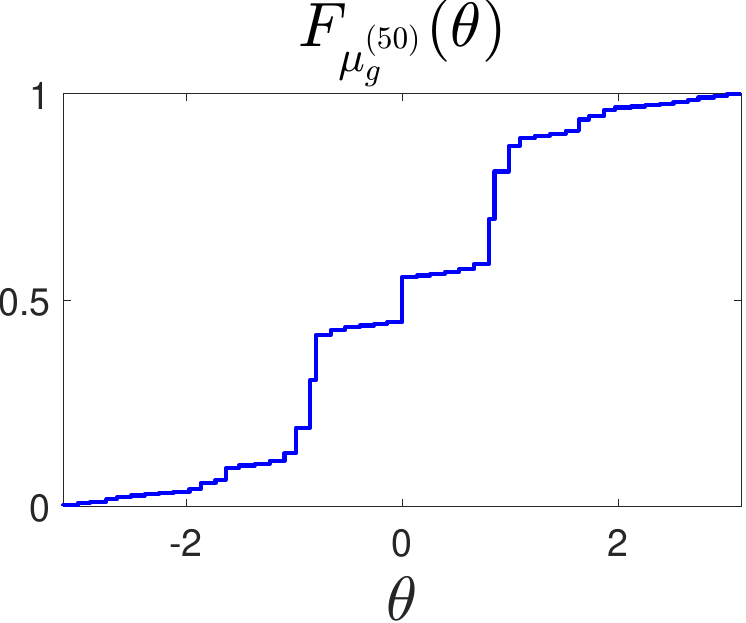}}
\raisebox{-0.5\height}{\includegraphics[width=0.24\textwidth,trim={0mm 0mm 0mm 0mm},clip]{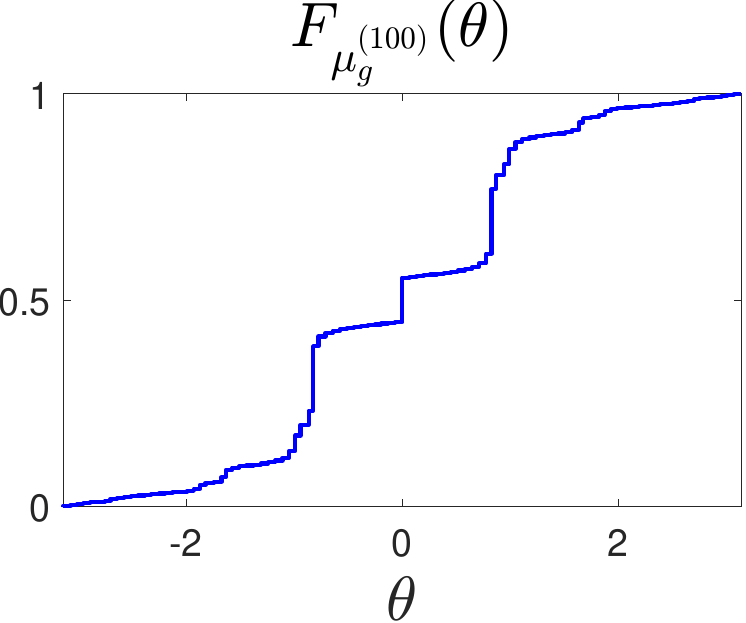}}
\raisebox{-0.5\height}{\includegraphics[width=0.24\textwidth,trim={0mm 0mm 0mm 0mm},clip]{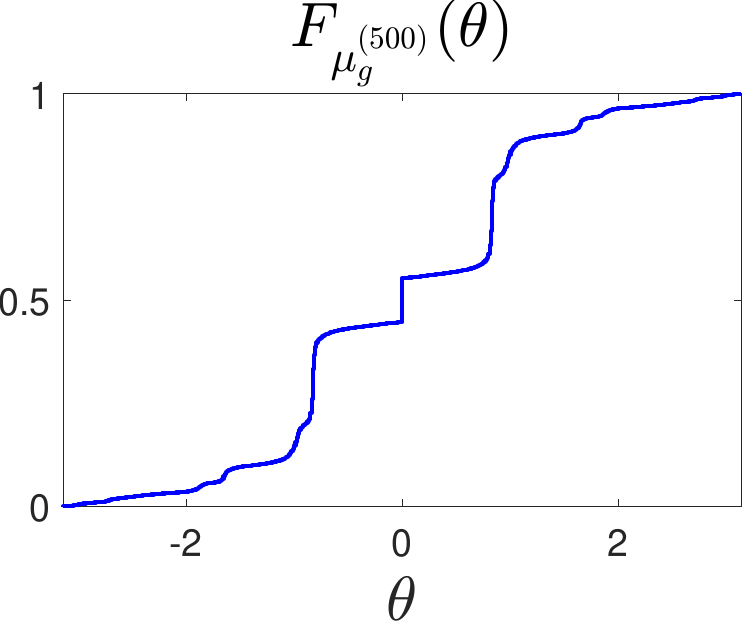}}
\raisebox{-0.5\height}{\includegraphics[width=0.24\textwidth,trim={0mm 0mm 0mm 0mm},clip]{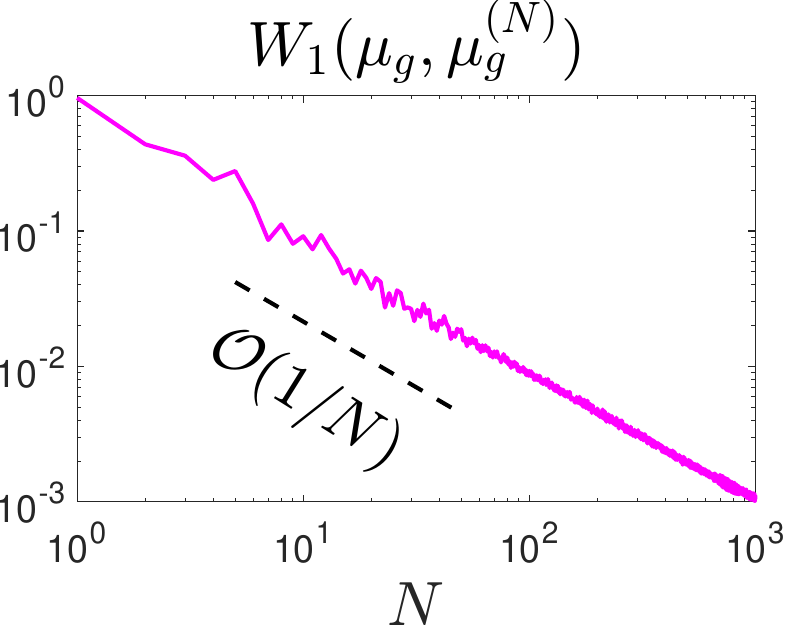}}
\caption{First three subplots: The cdfs computed by mpEDMD for various values of $N$. Far-right: The $W_1$ metric between the spectral measure computed by mpEDMD, $\mu_g^{(N)}$, and the spectral measure of the Koopman operator, $\mu_g$. The $W_1$ distance to $\mu_{g}$ is computed by comparing to an approximation with larger $N$ selected large enough to have a negligible effect on the shown errors.}
\label{fig_mpEDMD1}
\end{figure}

\cref{fig_mpEDMD1} displays the cdfs of $\mu_g^{(N)}=\mu_\gv^{(N,10^4)}$ for various choices of $N$, illustrating the convergence of spectral measures. Notably, there is a discontinuity in the cdfs at $\theta=0$, corresponding to the eigenvalue at $\lambda=1$. Away from this value, the cdfs show pointwise convergence. The error measured in the Wasserstein 1 metric is depicted on the far right of the plot. Consistent with \eqref{weak_weak_convergence3}, this error decreases as $\mathcal{O}(1/N)$. For similar analyses and examples regarding the projection-valued spectral measures, see \citep{colbrook2023mpedmd}.

\paragraph{Conservation of energy and statistics for turbulent boundary layer flow}

\begin{figure}
  \centering
     \includegraphics[width=1\textwidth,trim={0mm 0mm 0mm 0mm},clip]{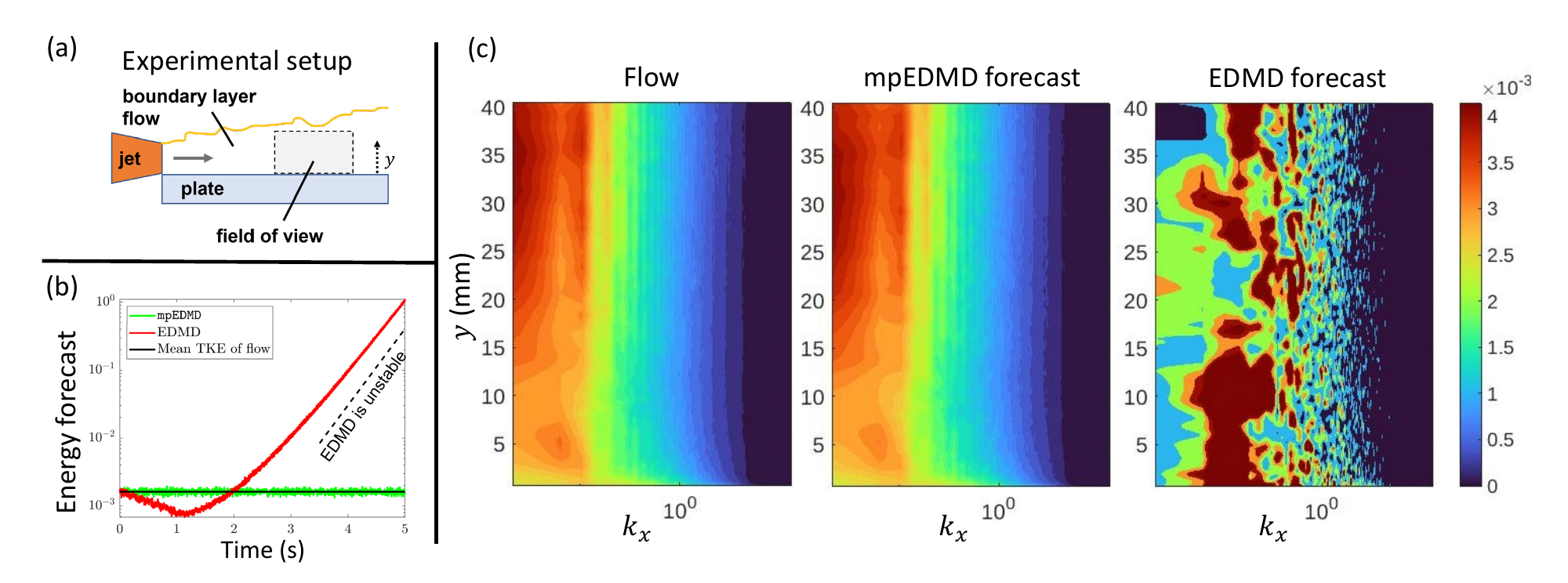}
\caption{(a) Experimental wall-jet boundary layer flow setup with Reynolds number $6.4\times10^4$. (b) Horizontal averages of the forecasts for turbulent kinetic energy, which show the stability of mpEDMD. (c) Wavenumber spectra measure the energy content of various turbulent structures as a function of their size, thus providing an efficient measure of a flow reconstruction method's performance over various spatial scales. This demonstrates the importance of structure-preserving discretizations (mpEDMD). Reproduced with permission from \citep{colbrook2023mpedmd}, copyright © 2023 Society for Industrial and Applied Mathematics, all rights reserved.}
\label{fig:mpEDMD2}
\end{figure}

We now examine the boundary layer formed by a thin jet injecting air onto a smooth, flat wall, as depicted in panel (a) of \cref{fig:mpEDMD2}. The experiments are conducted in the wind tunnel at Virginia Tech \citep{szoke2021flow}. We utilize a two-component, time-resolved particle image velocimetry system to capture $10^3$ snapshots of the two-dimensional velocity field of the wall-jet flow. These snapshots are taken over a spatial grid and a period of 1 second. The jet velocity is set at $U_{\mathrm{jet}}=50$m/s, corresponding to a jet Reynolds number of $6.4 \times 10^4$. The field-of-view spans approximately 75mm by 40mm, with a spatial resolution of approximately $\Delta x=\Delta y\approx0.24$mm. This setup leads to a dimension $d=102,300$ in \eqref{eq:DynamicalSystem}. A full SVD of the data matrix is employed to create a dictionary. The flow exhibits zero pressure gradient turbulent boundary layer characteristics within the region between the wall and the velocity profile peak at approximately $y=15.5$mm. Above this region, the flow is dominated by a two-dimensional shear layer characterized by large, energetic flow structures. This scenario presents a significant challenge for conventional DMD approaches due to the multiple turbulent scales within the boundary layer.

We investigate the conservation of energy and flow statistics using the KMD for future state predictions. We consider the velocity profiles predicted by mpEDMD and EDMD over 5 seconds, five times the observation window, starting from an initial state $\xv_0$ randomly selected from the trajectory data. The results are averaged over 100 such random initializations. Panel (b) of \cref{fig:mpEDMD2} shows the turbulent kinetic energy (TKE) of the predictions, averaged in the homogenous horizontal direction and normalized by $U_{\mathrm{jet}}^2$. The instability of EDMD is evident. In contrast, mpEDMD preserves the inner product associated with the TKE.

To examine the statistics of the predictions, panel (c) of \cref{fig:mpEDMD2} presents the wavenumber spectrum. This spectrum is computed by applying the Fourier transform to the spatial autocorrelations of the predictions in the horizontal direction, as detailed in \citep[Chapter 8]{glegg2017aeroacoustics}. The wavenumber spectrum provides insights into the energy content of various turbulent structures based on their size. It also serves as an efficient measure of a flow-reconstruction method's performance across different spatial scales. The wavenumber spectrum derived from mpEDMD shows remarkable alignment with the actual flow, demonstrating its efficacy. In contrast, EDMD completely fails to capture the correct turbulent statistics.

\subsection{Compactification methods for continuous-time systems}
\label{sec:compactification_methods}

For continuous-time invertible measure-preserving systems, the Koopman generator $\mathcal{L}$, as defined in \eqref{eq:gen_def}, is skew-adjoint. A sophisticated suite of methods exists aimed at approximating such generators through compactification. Working in continuous time presents at least two advantages. First, the generator $\mathcal{L}$ is skew-adjoint, while the Koopman operators $\mathcal{K}_{\Delta t}$ are unitary. Developing projection methods that preserve skew-adjointness is generally much more straightforward than preserving unitarity (although we have seen that mpEDMD leads to an appropriate unitary discretization). Second, by computing the spectral properties of the generator $\mathcal{L}$, we are no longer constrained by the need to select a specific discrete time step. We gain comprehensive spectral information for the entire family of Koopman operators $\{\mathcal{K}_{\Delta t}:\Delta t>0\}$.

\cite{das2021reproducing} developed an approach based on a one-parameter family of reproducing kernels, $\{p_\tau:\tau>0\}$, satisfying mild regularity assumptions. This method utilizes corresponding integral operators to perturb the Koopman generator $\mathcal{L}$ to a compact operator on the corresponding RKHS, $\mathcal{H}_\tau$. Assuming ergodic flow, \cite{das2021reproducing} constructed a one-parameter family of skew-adjoint compact operators, $W_\tau : \mathcal{H}_\tau \rightarrow \mathcal{H}_\tau$, where $W_\tau = P_\tau \mathcal{L} P_\tau^*$ and $P_\tau : L^2(\Omega,\omega) \to \mathcal H_\tau$ is the integral operator defined by
$$
    [P_\tau g](\xv') = \int_{\Omega} p_\tau(\xv',\xv) g(\xv)\,\mathrm{d}\omega(\xv).
$$
The operators $W_\tau$ are unitarily equivalent to $\mathcal{L}_\tau = G_\tau^{1/2} \mathcal{L} G_\tau^{1/2}$ acting on $L^2(\Omega,\omega)$, with $G_\tau = P_\tau^* P_\tau$. The operators $\mathcal{L}_\tau$ are compact, skew-adjoint, and converge in the strong resolvent sense to the generator $\mathcal{L}$ as $\tau \rightarrow 0$. Since each $\mathcal{L}_\tau$ is compact, its spectrum can be computed by projection onto finite-dimensional subspaces without spectral pollution and without missing parts of the spectrum in the limit of infinite discretization size \citep{SCI_ref}. This procedure yields approximate Koopman eigenvalues and eigenfunctions, which have been demonstrated to lie within the $\epsilon$-pseudospectrum of the Koopman operator, with the value of $\epsilon$ dependent on an RKHS-induced Dirichlet energy functional. In particular, approximate eigenfunctions with small Dirichlet energy as $\tau\rightarrow 0$ are approximately cyclical, slowly decorrelating observables under potentially mixing dynamics. It is important to note that two limits are implicitly involved here: the first concerns the parameter that controls the projection size used to approximate spectra of $\mathcal{L}_\tau$, and the second is as $\tau$ approaches zero. Another potential limitation of this method is its use of finite-difference schemes on time-ordered data. Although the error from these approximations can be controlled in the limit of a vanishing sampling interval via RKHS regularity, finite differencing generally reduces noise robustness.

Another approach involves the \textit{resolvent} of the generator, $(\mathcal{L}-zI)^{-1}$, where $z \in \mathbb C\backslash{i\mathbb{R}}$. By taking the Laplace transform of the Koopman semigroup, we can observe that \citep{susuki2021koopman}
\begin{equation}
\label{eq:resolventintegral}
 (\mathcal{L}-zI)^{-1} = - \int_0^\infty e^{-zt} \mathcal{K}_{\Delta t} \, \mathrm{d}t, \quad \mathrm{Re}(z) > 0.
\end{equation}
\cite{valva2023consistent} combine the compactification approach from \citep{das2021reproducing} and the integral representation of the resolvent used in \citep{susuki2021koopman} to construct a compact operator that acts as the resolvent of a skew-adjoint operator. The result is a family of skew-adjoint unbounded operators with compact resolvents, whose spectral measures converge weakly to those of $\mathcal{L}$. This method not only preserves skew-adjointness but also eliminates the need for finite-difference approximations of the generator by using a quadrature approximation for the integral in \eqref{eq:resolventintegral}. It offers a flexible framework that allows for the control of approximation accuracy by varying $z$ in relation to the sampling interval and the timespan of the training data. Additionally, the finite-rank operators are expressed in a well-conditioned basis of kernel eigenvectors, computed from trajectory data with convergence assurances in the large-data limit \citep{das2021reproducing}. These basis vectors are particularly well-suited to invariant measures supported on sets with complex geometries (e.g., fractal attractors) that are embedded in high-dimensional ambient spaces.

\subsubsection{Example}

\begin{figure}
\centering
\raisebox{-0.5\height}{\includegraphics[width=0.9\textwidth,trim={0mm 2mm 0mm 2mm},clip]{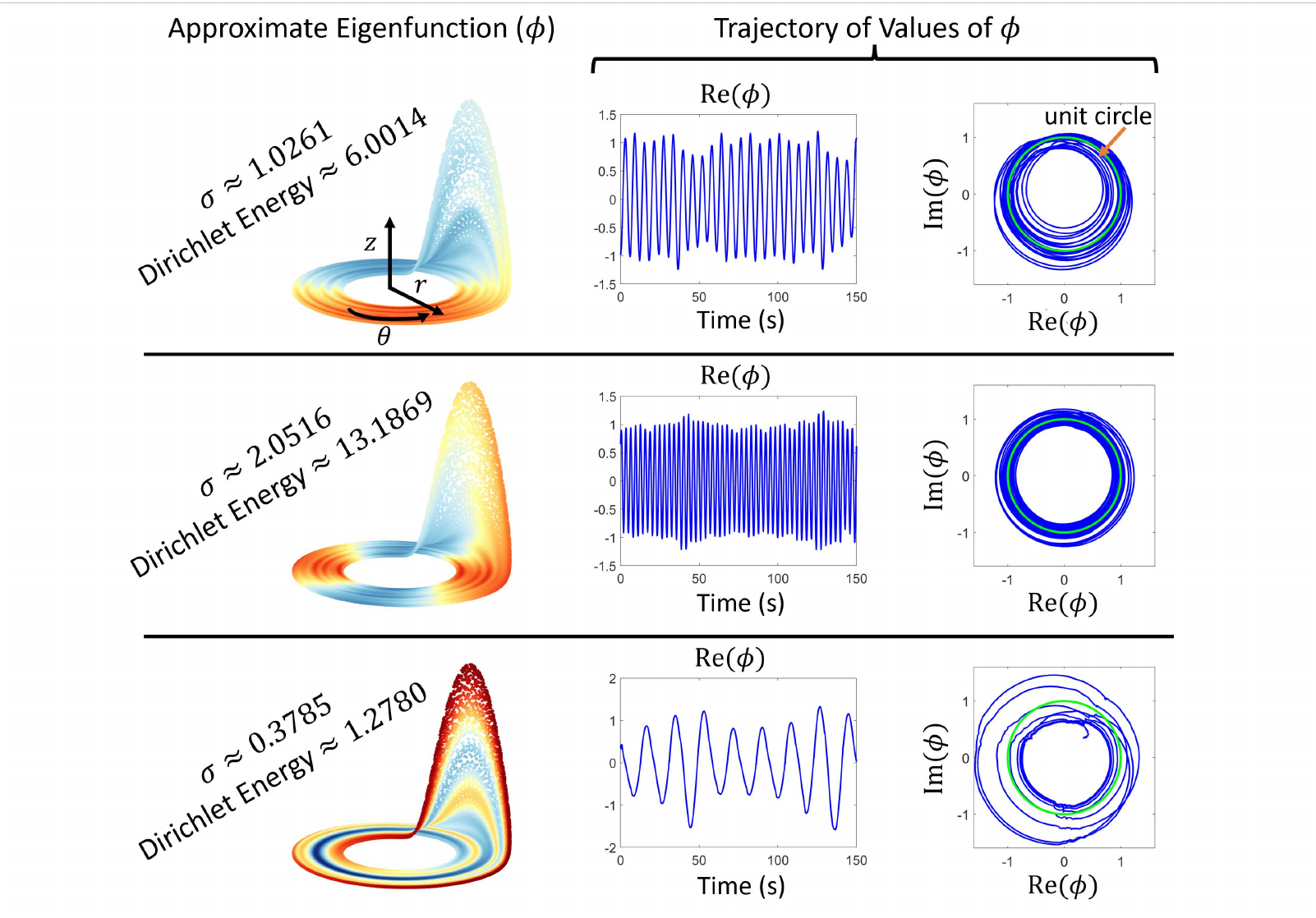}}
\caption{Approximate eigenfunctions of the R\"ossler system computed using the method of \citep{valva2023consistent}. On the right, we show the trajectories of these functions. A trajectory on the unit circle corresponds to coherent periodic behavior. The bottom row displays an approximate eigenfunction with radial variability and larger deviations from the unit circle.}
\label{fig:rossler}
\end{figure}

The R\"ossler system \citep{rossler1976equation} consists of the following three coupled ordinary differential equations:
$$
\dot{x}=-y-z,\quad\dot{y}=x+0.1y,\quad \dot{z}=0.1+z(x-14).
$$
We consider the dynamics of $\xv=(x,y,z)$ on the R\"ossler attractor. The R\"ossler system is often viewed as a simplified analog of the Lorenz (63) system. However, despite the simplicity of its governing equations, it exhibits complex dynamical characteristics. Theorems on the existence and measure-theoretic mixing properties of the R\"ossler system analogous to those for the Lorenz system have not been established. Nevertheless, the system has been studied extensively through analytical and numerical techniques, supporting the hypothesis that the R\"ossler system is mixing \citep{peifer2005mixing}. Assuming this, it follows that zero is the only eigenvalue of $\mathcal{L}$, corresponding to a constant eigenfunction and is simple. The integral in \eqref{eq:resolventintegral} is approximated by truncating the domain of integration (taking advantage of the exponential decay in the integrand) and Simpson's quadrature rule. Full algorithmic details of the method are given in \citep[Algorithm 1]{valva2023consistent}. Data is collected along a single trajectory of length $64,000$ with time-step $\Delta t=0.04$. We use MATLAB's \texttt{ode45} command to collect the data after an initial burn-in time to ensure that the initial point is (approximately) on the attractor. The dictionary consists of 2,001 data-driven kernel eigenfunctions, and the smoothing parameter is set as $\tau=2\times 10^{-6}$.

\cref{fig:rossler} shows three approximate eigenfunctions along with the corresponding values of $\sigma$, so that $i\sigma$ lies in the spectrum of $\mathcal{L}$. To the figure's right, we illustrate the trajectory of these approximate eigenfunctions. The chaotic behavior of the R\"ossler system predominantly occurs in the $(r=\sqrt{x^2+y^2}, z)$ coordinates, while the evolution of the azimuthal angle $\theta$ in the $z = 0$ plane proceeds at a near-constant angular frequency, approximately equal to one in natural time units. We suspect this distinction contributes to the challenge of capturing the approximate eigenfunctions in \cref{fig:rossler}. The first two approximate eigenfunctions are highly coherent, predominantly functions of the azimuthal phase angle, and evolve near-periodically over several Lyapunov times. The third approximate eigenfunction (bottom row) exhibits manifest radial variability in state space and amplitude-modulated time series. Additionally, resolving the radial direction may be more challenging in a data-driven basis, as most variability in the input data occurs in the azimuthal or vertical directions.

\begin{figure}
\centering
\raisebox{-0.5\height}{\includegraphics[width=0.45\textwidth,trim={0mm 0mm 0mm 0mm},clip]{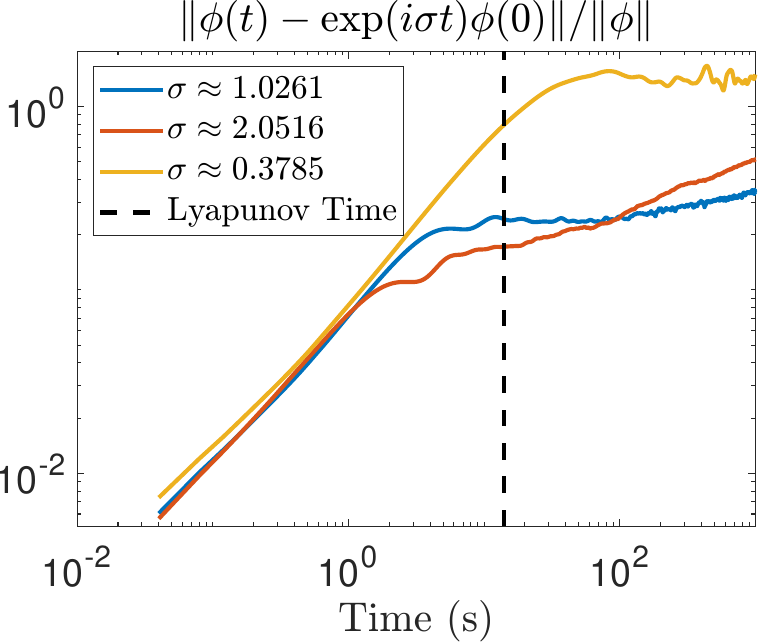}}
\caption{Relative residuals of the approximate eigenfunctions in \cref{fig:rossler} plotted as a function of time. The residuals increase steadily up to the characteristic Lyapunov timescale.}
\label{fig:rossler2}
\end{figure}

We have also presented the Dirichlet energies, indicative of the function's variability or roughness. The approximate eigenfunction corresponding to a frequency of $\sigma=0.3785$ demonstrates relatively low variability than the others. Furthermore, the increase in energy from the eigenfunction with frequency 1.0261 to its harmonic with frequency 2.0516 also mirrors this variability increase. Denoting each approximate eigenfunction's trajectory by $\phi(t)$, \cref{fig:rossler2} displays the relative residual $\|\phi(t)-\exp(i\sigma t)\phi(0)\|/\|\phi\|$. These residuals steadily increase up to the characteristic Lyapunov time of the system and exhibit a larger residual for the third approximate eigenfunction. To summarize, as an a posteriori metric, Dirichlet energy is generally independent from pseudospectral residuals. The information it provides can be useful in supervised learning tasks, e.g., when performing out-of-sample evaluation of the eigenfunctions in prediction problems. In practice, eigenfunctions with small Dirichlet energy also tend to have small pseudospectral residuals, though the precise ordering obtained from the two approaches may differ. 

\subsection{Further methods}
\label{sec:structure_pres_further_meth}

\cite{huang2018data} were among the first to enforce structure in DMD by introducing \textit{Naturally Structured DMD}. This variant ensures positivity and offers the added option of the Markov property. In another early paper, \cite{salova2019koopman} investigated dynamical systems with symmetries characterized by a finite group. Utilizing representation theory, the authors demonstrated that the Koopman operator and its EDMD approximations can be block diagonalized using a symmetry basis. This basis respects the isotypic component structure related to the underlying symmetry group and the actions of its elements, providing insights into suitable dictionaries. However, the data must align with the system's symmetries to achieve an exact block-diagonal approximation matrix. In an earlier work, \cite{sharma2016correspondence} connected the spatiotemporal symmetries of the Navier--Stokes equation with its spatial and temporal Koopman operators. \cite{kaiser2018discovering} presented a method to detect conservation laws using Koopman operator approximations, which can subsequently be employed to control Hamiltonian systems. 

We saw that mpEDMD and compactification methods are well-suited to measure-preserving systems. \cite{govindarajan2019approximation,govindarajan2021approximation} provide another approach that is similar to the Ulam approximation of the Perron--Frobenius operator \citep{ulam1960collection,li1976finite}. They proposed \textit{periodic approximations} for Koopman operators under conditions where $\Omega$ is compact, $\omega$ is absolutely continuous with respect to the Lebesgue measure, and the system is both measure-preserving and invertible \citep{govindarajan2019approximation}. This framework was developed into an algorithm for systems on tori and extended to continuous-time systems in \citep{govindarajan2021approximation}. The technique hinges on constructing a periodic approximation of the dynamics via a state-space partition, thus enabling the approximation of the Koopman operator's action through a permutation. The concept of periodic approximations has roots in the works of \cite{halmos1944approximation} and \cite{lax1971approximation}. This method yields measures that converge weakly to the spectral measures of the Koopman operator. Furthermore, periodic approximations are positive operators and uphold the multiplicative structure of the Koopman operator, i.e., $\mathcal{K}(fg) = (\mathcal{K}f)(\mathcal{K}g)$. A significant unresolved question is how these results can be generalized to handle systems that are not necessarily invariant with respect to a Lebesgue absolutely continuous measure, such as those defined on intricate domains like chaotic attractors, and how to develop efficient schemes in high dimension.

DMD fails in translational problems, such as wave-like phenomena, moving interfaces, and moving shocks \citep{kutz2016dynamic}. This limitation can be attributed to the dominant advection behavior propagating through the entire high-dimensional domain. This propagation makes establishing a global spatiotemporal basis challenging within a low-dimensional subspace. Drawing inspiration from Lagrangian POD \citep{mojgani2017lagrangian}, \cite{lu2020lagrangian} introduced \textit{Lagrangian DMD} that constructs a reduced-order model within the Lagrangian framework. Temporally evolving characteristic lines are selected as a central observable, and a low-dimensional structure in the Lagrangian framework is identified. \textit{Port-Hamiltonian DMD} \citep{morandin2023port} adapts the DMD within the port-Hamiltonian systems framework to ensure the system satisfies a dissipation inequality. \textit{Symmetric DMD} \citep{cohen2020mode} mandates the dynamics matrix to be symmetric. \textit{Constrained DMD} \citep{krake2022constrained} ensures the presence of specific frequencies by incorporating constraints into DMD.

On the transfer operator side, \cite{mehta2006symmetry} addresses symmetries of the Perron--Frobenius operator in relation to the admissible symmetry properties of attractors. \textit{Constrained Ulam DMD} \citep{goswami2018constrained} uses a minimization problem with constraints that guarantee a positive operator with a row sum equal to one. Beyond DMD, \cite{mardt2020deep} developed deep learning Markov and Koopman models with physical constraints. \cite{pan2020physics} learn continuous-time Koopman operators with deep neural networks and enforce stability by ensuring that eigenvalues have non-positive real parts. \cite{hirsh2021structured} presented a theoretical connection between time-delay embedding models and the Frenet--Serret frame (intrinsic coordinates formed by applying the Gram--Schmidt procedure to the derivatives of the trajectory) from differential geometry. This was used to develop structured HAVOK models.

\section{Further Topics and Open Problems}
\label{sec:conclusion}

We conclude this review of DMD with further topics that are connected with DMD and Koopman operators, followed by outlining some future challenges in the field.

\subsection{Transfer operators}
\label{sec:transfer_ops}

Perron--Frobenius operators, also known as Ruelle \citep{ruelle1968statistical} or transfer operators, act on measures through pullbacks. When considering appropriately chosen spaces of observables and measures, the Koopman and Perron--Frobenius operators emerge as dual pairs, thereby offering equivalent information. In the context of ergodic dynamical systems, natural spaces of observables are typically $L^2$ spaces of complex-valued scalar functions associated with invariant probability measures, while natural spaces of measures involve complex measures with $L^2$ densities. The operational distinction between Koopman and Perron--Frobenius operators has led to the development of two distinct families of approximation techniques. However, recent works, such as \citep{2158-2491_2016_1_51}, have started to bridge this gap. For a comprehensive introduction to Perron--Frobenius operators, see the textbook of \cite{bollt2013applied}.

Data-driven techniques employing Perron--Frobenius operators began with the seminal work of \cite{dellnitz1999approximation}. A prevalent approach in these methods involves approximating the Perron--Frobenius operator's spectrum using Ulam's method \citep{ulam1960collection,li1976finite}. This method involves partitioning the state space into a finite set of disjoint subsets. The transition probabilities between these subsets are then estimated by counting transitions observed in extensive simulations or experimental data. The derived transition probability matrix essentially serves as a Galerkin projection of a smoothed compact transfer operator, slightly perturbed by noise. The matrix's eigenvectors, corresponding to eigenvalues at or near the unit circle, are instrumental in identifying coherent sets. Subsequent research by Dellnitz, Froyland, and colleagues \citep{dellnitz2000isolated,froyland2003detecting,froyland2007ulam,froyland2008unwrapping,froyland2014detecting} focused on specific system classes with quasi-compact Perron--Frobenius operators. Their work rigorously demonstrated Ulam's method's efficacy in accurately approximating isolated Perron--Frobenius eigenvalues and their associated eigenfunctions.

The Perron--Frobenius operator has been used to analyze the global behavior of dynamical systems across various fields. Its applications span molecular dynamics \citep{schutte2013metastability,schutte2023overcoming}, fluid dynamics \citep{froyland2014well,froyland2016optimal}, meteorology and atmospheric sciences \citep{tantet2015early,tantet2018crisis,froyland2021spectral}, as well as engineering \citep{vaidya2010nonlinear,ober2015multiobjective}. Various toolboxes, such as GAIO \citep{dellnitz2001algorithms}, can compute almost invariant sets or metastable states. These toolboxes utilize adaptive box discretizations of the state space to approximate the system's behavior efficiently. However, it is important to note that this approach is generally more suited to low-dimensional problems.

\subsection{Continuous spectra and spectral measures}
\label{sec:methods_for_spec_meas}

In \cref{sec:spectra_crash_course}, we saw how Koopman operators associated with measure-preserving systems have spectral measures that provide a KMD. In the course of this review, we have met several methods that converge to spectral measures: mpEDMD in \cref{sec:mpEDMD} (see \cref{fig_mpEDMD1} for an example), methods based on compactification in \cref{sec:compactification_methods}, and partitioning of statespace to obtain periodic approximations in \cref{sec:structure_pres_further_meth}. There are other methods that are not based on the eigenvalues of a finite matrix. \cite{korda2020data} approximate the moments of spectral measures of ergodic systems using \eqref{eq:birkhoff}, and then use the Christoffel--Darboux kernel to analyze the atomic and absolutely continuous parts of the spectrum. They also compute the spectral projection on a given segment of the unit circle. See also \citep{arbabi2017study}, who use harmonic averaging and Welch's method \citep{welch1967use} to compute the discrete and continuous spectrum of the Koopman operator for post-transient flows. Using the resolvent operator and ResDMD, \cite{colbrook2021rigorous} compute smoothed approximations of spectral measures associated with general measure-preserving dynamical systems. They prove explicit high-order convergence theorems for the computation of spectral measures in various senses, including computing the density of the continuous spectrum, spectral projections of subsets of the unit circle, and the discrete spectrum. These smoothing techniques can also be used for self-adjoint operators \citep{colbrook2021computing}.

However, we end this discussion with the following warning to the reader about recovering atomic parts of spectral measures that should be kept in mind for all of the above methods. As soon as the spectral measure $\mu_g$ has atoms (i.e., $\mathcal{K}$ has eigenvalues, and $g$ is not orthogonal to all the eigenspaces), the map $\lambda\mapsto \mu_g(\{\lambda\})$ is discontinuous. One can prove that, in general, separating the point spectrum from the rest of the spectrum, either in terms of spectral measures or spectral sets, is impossible for \textit{any} algorithm. This holds even for simple classes of operators~\citep{colbrook2021computingCIMP}, unless we know apriori that the spectrum is discrete in a region of interest \citep[Section 7.3]{colbrook2021computing}. An excellent example and discussion of this point is provided on the second page of \citep{govindarajan2019approximation}. This is one reason the above methods can only compute spectral measures in a weak or setwise sense. It also helps explain why many of these techniques involve some form of smoothing.

\subsection{Stochastic dynamical systems}
\label{sec:stochastic_systems}

Stochastic dynamical systems are widely used to model and study systems that evolve under the influence of both deterministic and random effects. It is common to replace \eqref{eq:DynamicalSystem} with a discrete-time Markov process
\begin{equation}
\label{eq:stochDynamicalSystem} 
\xv_{n+1} = \Fv(\xv_n,\tau_{n}), \qquad n= 0,1,2,\ldots, 
\end{equation}
where $\{\tau_n\}\in\Omega_s$ are independent and identically distributed random variables with distribution $\rho$ supported on $\Omega_s$, and $\Fv : \Omega\times\Omega_s\rightarrow\Omega$ is a function. The stochastic Koopman operator (also called the Kolmogorov operator) is the expectation:
\begin{equation}\label{def:K_1}
[\mathcal{K}g](\xv) = \int_{\Omega_s} g(\Fv(\xv,\tau))\,\mathrm{d}\rho(\tau).  
\end{equation}
In contrast to the deterministic case, stochastic Koopman operators typically have discrete spectra due to diffusion. A primary focus has been the challenge of noisy observables in EDMD-type methods \citep{takeishi2017subspace,wanner2022robust}, and debiasing DMD~\citep{hemati2017biasing,dawson2016characterizing,takeishi2017subspace}. \cite{vcrnjaric2020koopman} developed a stochastic Hankel-DMD algorithm for numerical approximations of the stochastic Koopman operator. \cite{klus2020data} used gEDMD to derive models for SDEs with applications in control. \cite{sinha2020robust} provided an explicit optimization-based approximation of stochastic Koopman operators. \cite{wu2020variational} developed a variational approach for Markov processes that finds optimal feature mappings and optimal Markovian models of the dynamics from the top singular components of the Koopman operator. The definition in \eqref{def:K_1} involves an expectation. \cite{colbrook2023beyond} demonstrated the benefits and necessity of going beyond expectations of trajectories. They incorporated the concept of variance into the Koopman framework, establishing its relationship with batched Koopman operators. This led to an extension of ResDMD, resulting in convergence to the spectral properties of stochastic Koopman operators, a Koopman analog of a variance-bias decomposition, and the concept of variance-pseudospectra as a measure of statistical coherency. 

\subsection{Some open problems}
\label{sec:open_problems}

There has been substantial interest in Koopmanism over the last decade, and we expect this interest only to grow. This is an exciting time to be working in this field, which is at the crossroads of dynamical systems theory, data analysis, spectral theory, and computational analysis. We end with some open problems that the author believes will lead to important breakthroughs in the coming years:
\begin{itemize}	
\item \textbf{Banach spaces:} We have focused on Koopman operators defined on the Hilbert space $L^2(\Omega,\omega)$. In some cases, it is more appropriate to consider function spaces that are Banach spaces \citep{mohr2014construction,mezic2020spectrum}. Computational tools for infinite-dimensional spectral problems on Banach spaces are less developed than those for Hilbert spaces. An exception is the transfer operator community, which has developed methods for quasi-compact Perron--Frobenius operators (see \cref{sec:transfer_ops}). Another challenge is the development of the theory of the KMD in the absence of spectral theorems. These issues are expected to be particularly significant in applying the Koopman framework and DMD to transient and off-attractor dynamics.

For some chaotic systems under appropriate conditions, the eigenvalues in the large subspace limit of EDMD correspond to the eigenvalues of transfer operators in suitable function spaces, as discussed in \citep{slipantschuk2020dynamic, wormell2023orthogonal, bandtlow2023edmd}. Reconciling these generalized eigenfunctions, the $L^2$ flavor of DMD methods, and appropriate Banach spaces is a key open problem. These questions could be crucial for understanding and improving the effectiveness of EDMD, including guidance on the choice of dictionary.
		
		\item \textbf{Choice of dictionary:} One of the most significant open problems in EDMD is the selection of observables or dictionaries. At present, this process can be considered more of an art than a science. While well-conditioned bases can be constructed for some systems, as outlined in \cref{sec:compactification_methods}, this task often presents a substantial challenge. This is also true for methods based on delay-embedding, where the choice of delay itself is a classical problem with many available heuristics.

		\item \textbf{Foundations:} All of the convergence results for DMD and Koopman operators rely on algorithms that depend on several parameters, with successive limits of these parameters taken to achieve convergence. This is not accidental and occurs generically in infinite-dimensional spectral computations \citep{SCI_ref}. It is possible to classify the difficulty of computational problems, including data-driven ones. To date, the Koopman community has only provided upper bounds, i.e., algorithms that converge for specific classes of problems. A significant open problem is the development of lower bounds, i.e., universal impossibility results that indicate an intrinsic difficulty in a problem that cannot be overcome by any algorithm. Such results are beginning to emerge in the world of deep learning, particularly regarding the existence vs. trainability of neural networks \citep{colbrook2022difficulty}. We expect them to be equally fruitful in the Koopman context.
	
Lower bounds are essential for several reasons. First, they prevent futile searches for non-existent algorithms. Second, they often elucidate why certain algorithms cannot exist. When combined with upper bounds, this knowledge can lead to natural assumptions about the dynamical systems or the data required to achieve our computational goals.
	\item \textbf{Further structure-preserving methods:} We have already outlined several open problems stimulated by piDMD in \cref{sec:piDMD_future}. Understanding the relationships between structures or symmetries in dynamical systems and their manifestation in the Koopman spectrum lies at the forefront of current knowledge. A crucial challenge is extending constraints applicable to DMD with linear observables to EDMD with nonlinear observables, which will be instrumental in applying structure-preserving methods to nonlinear systems effectively. Another related open problem is establishing the convergence of constrained DMD methods to the spectral properties of Koopman operators and the convergence of KMDs.
	
	\item \textbf{Verified control:} An exciting development area in modern Koopman theory is its use for controlling nonlinear systems. In \cref{sec:extensions_for_control}, we discussed how the practice of Koopman control currently surpasses the theoretical understanding. Only a limited number of systems with a known Koopman-invariant subspace and verifiable eigenfunctions exist. Thus, developing methods to validate Koopman models for control purposes is a crucial problem. Successfully addressing this issue will likely lead to further insights and enhancements in the practice of Koopman control.
		
\end{itemize}



\small
\section*{Acknowledgements}

I am grateful to numerous friends for their support and stimulating discussions during the completion of this work. The following is in no particular order. Thank you Steven Brunton and Nathan Kutz, for your thoughts and perspectives on an early version of this manuscript, in particular, \cref{tab:summary}. Thank you Igor Mezi{\'c}, for our discussions regarding spectral theory and computation. Thank you Benjamin Erichson, for running the numerical experiments behind \cref{fig_SAT_modes,fig_SAT_times} and discussing rDMD with me. Thank you Claire Valva and Dimitrios Giannakis, for discussing compactification methods with me, and thank you Claire for running the numerical experiments behind \cref{fig:rossler}. Thank you Zlatko Drmač, for our discussions regarding the numerical linear algebra aspects of DMD. Thank you Julia Slipantschuk and Caroline Wormell for our discussions regarding numerical methods for transfer operators. Thank you
Zlatko Drmač,
Dimitrios Giannakis,
Igor Mezi{\'c}, and
Samuel Otto,
for your insightful feedbacks on an initial draft.
\normalsize

\begin{spacing}{.9}
\setlength{\bibsep}{6.25pt}
\bibliographystyle{agsm}
{\small\linespread{0.9}\selectfont{}
\bibliography{bib_file_DMD3}}

@phdthesis{colbrookthesis,
  title={The {Foundations of Infinite-Dimensional Spectral Computations}},
  author={Colbrook, M. J.},
  year={2020},
  school={University of Cambridge}
}

@Article{hansen2011solvability,
  author    = {Hansen, Anders},
  journal   = {Journal of the American Mathematical Society},
  title     = {On the solvability complexity index, the {$n$}-pseudospectrum and approximations of spectra of operators},
  year      = {2011},
  issn      = {1088-6834},
  month     = jul,
  number    = {1},
  pages     = {81--124},
  volume    = {24},
  doi       = {10.1090/s0894-0347-2010-00676-5},
  publisher = {American Mathematical Society (AMS)},
}

@Article{colbrook2022computation,
  author    = {Colbrook, Matthew J.},
  journal   = {Foundations of Computational Mathematics},
  title     = {On the computation of geometric features of spectra of linear operators on {H}ilbert spaces},
  year      = {2022},
  issn      = {1615-3383},
  month     = dec,
  pages     = {1--82},
  doi       = {10.1007/s10208-022-09598-0},
  publisher = {Springer Science and Business Media LLC},
}

@Article{colbrook2022foundations,
  author    = {Colbrook, Matthew J. and Hansen, Anders C.},
  journal   = {Journal of the European Mathematical Society},
  title     = {The foundations of spectral computations via the solvability complexity index hierarchy},
  year      = {2022},
  issn      = {1435-9855},
  month     = nov,
  number    = {12},
  pages     = {4639--4728},
  volume    = {25},
  doi       = {10.4171/jems/1289},
  publisher = {European Mathematical Society - EMS - Publishing House GmbH},
}

@Article{colbrook2021computing,
  author    = {Colbrook, Matthew and Horning, Andrew and Townsend, Alex},
  journal   = {SIAM Review},
  title     = {Computing spectral measures of self-adjoint operators},
  year      = {2021},
  issn      = {1095-7200},
  month     = jan,
  number    = {3},
  pages     = {489--524},
  volume    = {63},
  doi       = {10.1137/20m1330944},
  publisher = {Society for Industrial \& Applied Mathematics (SIAM)},
}

@Article{colbrook2021computingCIMP,
  author    = {Colbrook, Matthew J.},
  journal   = {Communications in Mathematical Physics},
  title     = {Computing spectral measures and spectral types},
  year      = {2021},
  issn      = {1432-0916},
  month     = apr,
  number    = {1},
  pages     = {433--501},
  volume    = {384},
  doi       = {10.1007/s00220-021-04072-4},
  publisher = {Springer Science and Business Media LLC},
}

@Article{colbrook2019compute,
  author    = {Colbrook, Matthew J. and Roman, Bogdan and Hansen, Anders C.},
  journal   = {Physical Review Letters},
  title     = {How to compute spectra with error control},
  year      = {2019},
  issn      = {1079-7114},
  month     = jun,
  number    = {25},
  pages     = {250201},
  volume    = {122},
  doi       = {10.1103/physrevlett.122.250201},
  publisher = {American Physical Society (APS)},
}

@Article{SCI_ref,
  author        = {Ben-Artzi, Jonathan and Colbrook, Matthew J. and Hansen, Anders C. and Nevanlinna, Olavi and Seidel, Markus},
  journal       = {arXiv preprint arXiv:1508.03280},
  title         = {Computing spectra - {O}n the solvability complexity index hierarchy and towers of algorithms},
  year          = {2020},
  month         = aug,
  archiveprefix = {arXiv},
  copyright     = {arXiv.org perpetual, non-exclusive license},
  doi           = {10.48550/ARXIV.1508.03280},
  eprint        = {1508.03280},
  primaryclass  = {cs.CC},
  publisher     = {arXiv},
}

@Article{mezic2004comparison,
  author    = {Mezi{\'c}, I. and Banaszuk, A.},
  journal   = {Physica D: Nonlinear Phenomena},
  title     = {Comparison of systems with complex behavior},
  year      = {2004},
  issn      = {0167-2789},
  month     = oct,
  number    = {1-2},
  pages     = {101--133},
  volume    = {197},
  doi       = {10.1016/j.physd.2004.06.015},
  publisher = {Elsevier BV},
}

@Article{dellnitz1999approximation,
  author    = {Dellnitz, Michael and Junge, Oliver},
  journal   = {SIAM Journal on Numerical Analysis},
  title     = {On the approximation of complicated dynamical behavior},
  year      = {1999},
  issn      = {1095-7170},
  month     = jan,
  number    = {2},
  pages     = {491--515},
  volume    = {36},
  doi       = {10.1137/s0036142996313002},
  publisher = {Society for Industrial \& Applied Mathematics (SIAM)},
}

@Article{froyland2007ulam,
  author    = {Froyland, Gary},
  journal   = {Discrete and Continuous Dynamical Systems},
  title     = {On {U}lam approximation of the isolated spectrum and eigenfunctions of hyperbolic maps},
  year      = {2007},
  issn      = {1553-5231},
  number    = {3},
  pages     = {671--689},
  volume    = {17},
  doi       = {10.3934/dcds.2007.17.671},
  publisher = {American Institute of Mathematical Sciences (AIMS)},
}

@Article{giannakis2021delay,
  author    = {Giannakis, Dimitrios},
  journal   = {Research in the Mathematical Sciences},
  title     = {Delay-coordinate maps, coherence, and approximate spectra of evolution operators},
  year      = {2021},
  issn      = {2197-9847},
  month     = jan,
  number    = {1},
  pages     = {1--33},
  volume    = {8},
  doi       = {10.1007/s40687-020-00239-y},
  publisher = {Springer Science and Business Media LLC},
}

@Article{mezicAMS,
  author    = {Mezi{\'c}, I.},
  journal   = {Notices of the American Mathematical Society},
  title     = {Koopman operator, geometry, and learning of dynamical systems},
  year      = {2021},
  issn      = {1088-9477},
  month     = aug,
  number    = {07},
  pages     = {1},
  volume    = {68},
  doi       = {10.1090/noti2306},
  publisher = {American Mathematical Society (AMS)},
}

@Article{susuki2021koopman,
  author    = {Susuki, Yoshihiko and Mauroy, Alexandre and Mezi{\'c}, Igor},
  journal   = {SIAM Journal on Applied Dynamical Systems},
  title     = {Koopman resolvent: {A} {L}aplace-domain analysis of nonlinear autonomous dynamical systems},
  year      = {2021},
  issn      = {1536-0040},
  month     = jan,
  number    = {4},
  pages     = {2013--2036},
  volume    = {20},
  doi       = {10.1137/20m1335935},
  publisher = {Society for Industrial \& Applied Mathematics (SIAM)},
}

@Article{mezic2020numerical,
  author    = {Mezi{\'c}, Igor},
  journal   = {Mathematics},
  title     = {On numerical approximations of the {K}oopman operator},
  year      = {2022},
  issn      = {2227-7390},
  month     = apr,
  number    = {7},
  pages     = {1180},
  volume    = {10},
  doi       = {10.3390/math10071180},
  publisher = {MDPI AG},
}

@Article{mezic2020spectrum,
  author    = {Mezi{\'c}, Igor},
  journal   = {Journal of Nonlinear Science},
  title     = {Spectrum of the {K}oopman operator, spectral expansions in functional spaces, and state-space geometry},
  year      = {2020},
  issn      = {1432-1467},
  month     = dec,
  number    = {5},
  pages     = {2091--2145},
  volume    = {30},
  doi       = {10.1007/s00332-019-09598-5},
  publisher = {Springer Science and Business Media LLC},
}

@Article{brunton2017chaos,
  author    = {Brunton, Steven L. and Brunton, Bingni W. and Proctor, Joshua L. and Kaiser, Eurika and Kutz, J. Nathan},
  journal   = {Nature Communications},
  title     = {Chaos as an intermittently forced linear system},
  year      = {2017},
  issn      = {2041-1723},
  month     = may,
  number    = {1},
  pages     = {1--9},
  volume    = {8},
  doi       = {10.1038/s41467-017-00030-8},
  publisher = {Springer Science and Business Media LLC},
}

@Book{mohri2018foundations,
  author    = {Mohri, Mehryar and Rostamizadeh, Afshin and Talwalkar, Ameet},
  publisher = {The MIT Press},
  title     = {Foundations of machine learning},
  year      = {2018},
  address   = {Cambridge, Massachusetts},
  edition   = {Second edition},
  isbn      = {9780262351362},
  series    = {Adaptive computation and machine learning},
  pagetotal = {486},
}

@Book{friedman2017elements,
  author    = {Hastie, Trevor and Tibshirani, Robert and Friedman, Jerome},
  publisher = {Springer New York},
  title     = {The elements of statistical learning},
  year      = {2009},
  isbn      = {9780387848587},
  doi       = {10.1007/978-0-387-84858-7},
  issn      = {2197-568X},
  journal   = {Springer Series in Statistics},
}

@Book{hey2009fourth,
  editor    = {Hey, Anthony J. G. and Tansley, Stewart and Tolle, Kristin Michele},
  publisher = {Microsoft Research},
  title     = {The fourth paradigm},
  year      = {2009},
  address   = {Redmond, Wash.},
  isbn      = {9780982544204},
  volume    = {1},
  ppn_gvk   = {1648419534},
  subtitle  = {Data intensive scientific discovery},
}

@Article{trefethen1993hydrodynamic,
  author    = {Trefethen, Lloyd N. and Trefethen, Anne E. and Reddy, Satish C. and Driscoll, Tobin A.},
  journal   = {Science},
  title     = {Hydrodynamic stability without eigenvalues},
  year      = {1993},
  issn      = {1095-9203},
  month     = jul,
  number    = {5121},
  pages     = {578--584},
  volume    = {261},
  doi       = {10.1126/science.261.5121.578},
  publisher = {American Association for the Advancement of Science (AAAS)},
}

@Article{caflisch1998monte,
  author    = {Caflisch, Russel E.},
  journal   = {Acta Numerica},
  title     = {Monte {C}arlo and quasi-{M}onte {C}arlo methods},
  year      = {1998},
  issn      = {1474-0508},
  month     = jan,
  pages     = {1--49},
  volume    = {7},
  doi       = {10.1017/s0962492900002804},
  publisher = {Cambridge University Press (CUP)},
}

@Article{woodley1999resonant,
  author    = {Woodley, B. M. and Peake, N.},
  journal   = {Journal of Fluid Mechanics},
  title     = {Resonant acoustic frequencies of a tandem cascade. {P}art 2. {R}otating blade rows},
  year      = {1999},
  issn      = {1469-7645},
  month     = aug,
  pages     = {241--256},
  volume    = {393},
  doi       = {10.1017/s0022112099005613},
  publisher = {Cambridge University Press (CUP)},
}

@Article{parker1984acoustic,
  author    = {Parker, R.},
  journal   = {Journal of Sound and Vibration},
  title     = {Acoustic resonances and blade vibration in axial flow compressors},
  year      = {1984},
  issn      = {0022-460X},
  month     = feb,
  number    = {4},
  pages     = {529--539},
  volume    = {92},
  doi       = {10.1016/0022-460x(84)90196-2},
  publisher = {Elsevier BV},
}

@Article{korda2020data,
  author    = {Korda, Milan and Putinar, Mihai and Mezi{\'c}, Igor},
  journal   = {Applied and Computational Harmonic Analysis},
  title     = {Data-driven spectral analysis of the {K}oopman operator},
  year      = {2020},
  issn      = {1063-5203},
  month     = mar,
  number    = {2},
  pages     = {599--629},
  volume    = {48},
  doi       = {10.1016/j.acha.2018.08.002},
  publisher = {Elsevier BV},
}

@Article{zaslavsky2002chaos,
  author    = {Zaslavsky, George M.},
  journal   = {Physics Reports},
  title     = {Chaos, fractional kinetics, and anomalous transport},
  year      = {2002},
  issn      = {0370-1573},
  month     = dec,
  number    = {6},
  pages     = {461--580},
  volume    = {371},
  doi       = {10.1016/s0370-1573(02)00331-9},
  publisher = {Elsevier BV},
}

@Article{brunton2016koopman,
  author    = {Brunton, Steven L. and Brunton, Bingni W. and Proctor, Joshua L. and Kutz, J. Nathan},
  journal   = {PLOS ONE},
  title     = {Koopman Invariant Subspaces and Finite Linear Representations of Nonlinear Dynamical Systems for Control},
  year      = {2016},
  issn      = {1932-6203},
  month     = feb,
  number    = {2},
  pages     = {e0150171},
  volume    = {11},
  doi       = {10.1371/journal.pone.0150171},
  publisher = {Public Library of Science (PLoS)},
}

@Book{shields1973theory,
  author    = {Shields, Paul C.},
  publisher = {University of Chicago Press, Chicago},
  title     = {The Theory of {B}ernoulli Shifts},
  year      = {1973},
}

@Book{walters2000introduction,
  author    = {Walters, Peter},
  publisher = {Springer},
  title     = {An introduction to ergodic theory},
  year      = {2000},
  address   = {New York},
  edition   = {1. softcover printing},
  isbn      = {0387951520},
  series    = {Graduate texts in mathematics},
  volume    = {79},
  pagetotal = {250},
  ppn_gvk   = {324378890},
}

@Book{dubrovin2012modern,
  author    = {Dubrovin, Boris A. and Fomenko, Anatolij Timofeevi{\v{c}} and Novikov, Serge{\u\i} Petrovich},
  publisher = {Springer New York},
  title     = {{Modern Geometry - Methods and Applications Part I. The Geometry of Surfaces, Transformation Groups, and Fields}},
  year      = {1984},
  isbn      = {9781468499469},
  volume    = {104},
  doi       = {10.1007/978-1-4684-9946-9},
  issn      = {0072-5285},
  journal   = {Graduate Texts in Mathematics},
}

@Article{brunton2021modern,
  author    = {Brunton, Steven L. and Budi{\v{s}}i{\'c}, Marko and Kaiser, Eurika and Kutz, J. Nathan},
  journal   = {SIAM Review},
  title     = {Modern {K}oopman Theory for Dynamical Systems},
  year      = {2022},
  issn      = {1095-7200},
  month     = may,
  number    = {2},
  pages     = {229--340},
  volume    = {64},
  doi       = {10.1137/21m1401243},
  publisher = {Society for Industrial \& Applied Mathematics (SIAM)},
}

@Article{lewin2010spectral,
  author    = {Lewin, Mathieu and S{\'e}r{\'e}, {\'E}ric},
  journal   = {Proceedings of the London Mathematical Society},
  title     = {Spectral pollution and how to avoid it},
  year      = {2010},
  issn      = {0024-6115},
  month     = dec,
  number    = {3},
  pages     = {864--900},
  volume    = {100},
  doi       = {10.1112/plms/pdp046},
  publisher = {Wiley},
}

@Article{budivsic2012applied,
  author    = {Budi{\v{s}}i{\'c}, Marko and Mohr, Ryan and Mezi{\'c}, Igor},
  journal   = {Chaos: An Interdisciplinary Journal of Nonlinear Science},
  title     = {Applied {K}oopmanism},
  year      = {2012},
  issn      = {1089-7682},
  month     = dec,
  number    = {4},
  pages     = {047510},
  volume    = {22},
  doi       = {10.1063/1.4772195},
  publisher = {AIP Publishing},
}

@Article{taira2017modal,
  author    = {Taira, Kunihiko and Brunton, Steven L. and Dawson, Scott T. M. and Rowley, Clarence W. and Colonius, Tim and McKeon, Beverley J. and Schmidt, Oliver T. and Gordeyev, Stanislav and Theofilis, Vassilios and Ukeiley, Lawrence S.},
  journal   = {AIAA Journal},
  title     = {Modal analysis of fluid flows: {A}n overview},
  year      = {2017},
  issn      = {1533-385X},
  month     = dec,
  number    = {12},
  pages     = {4013--4041},
  volume    = {55},
  doi       = {10.2514/1.j056060},
  publisher = {American Institute of Aeronautics and Astronautics (AIAA)},
}

@Article{klus2020data,
  author    = {Klus, Stefan and N{\"u}ske, Feliks and Peitz, Sebastian and Niemann, Jan-Hendrik and Clementi, Cecilia and Sch{\"u}tte, Christof},
  journal   = {Physica D: Nonlinear Phenomena},
  title     = {Data-driven approximation of the {K}oopman generator: {M}odel reduction, system identification, and control},
  year      = {2020},
  issn      = {0167-2789},
  month     = may,
  pages     = {132416},
  volume    = {406},
  doi       = {10.1016/j.physd.2020.132416},
  publisher = {Elsevier BV},
}

@Article{welch1967use,
  author    = {Welch, P.},
  journal   = {IEEE Transactions on Audio and Electroacoustics},
  title     = {The use of fast {F}ourier transform for the estimation of power spectra: {A} method based on time averaging over short, modified periodograms},
  year      = {1967},
  issn      = {0018-9278},
  month     = jun,
  number    = {2},
  pages     = {70--73},
  volume    = {15},
  doi       = {10.1109/tau.1967.1161901},
  publisher = {Institute of Electrical and Electronics Engineers (IEEE)},
}

@Article{arbabi2017study,
  author    = {Arbabi, Hassan and Mezi{\'c}, Igor},
  journal   = {Physical Review Fluids},
  title     = {Study of dynamics in post-transient flows using {K}oopman mode decomposition},
  year      = {2017},
  issn      = {2469-990X},
  month     = dec,
  number    = {12},
  pages     = {124402},
  volume    = {2},
  doi       = {10.1103/physrevfluids.2.124402},
  publisher = {American Physical Society (APS)},
}

@Book{beer1993topologies,
  author    = {Beer, Gerald},
  publisher = {Springer Netherlands},
  title     = {Topologies on Closed and Closed Convex Sets},
  year      = {1993},
  isbn      = {9789401581493},
  volume    = {268},
  doi       = {10.1007/978-94-015-8149-3},
}

@InProceedings{bruder2019modeling,
  author     = {Bruder, Daniel and Gillespie, Brent and David Remy, C. and Vasudevan, Ram},
  booktitle  = {Robotics: Science and Systems XV},
  title      = {Modeling and control of soft robots using the {K}oopman operator and model predictive control},
  year       = {2019},
  month      = jun,
  publisher  = {Robotics: Science and Systems Foundation},
  series     = {RSS2019},
  collection = {RSS2019},
  doi        = {10.15607/rss.2019.xv.060},
  journal    = {arXiv preprint arXiv:1902.02827},
}

@article{scholkopf2001kernel,
  title={The kernel trick for distances},
  author={Scholkopf, Bernhard},
  journal={Adv. Neur. Info. Proc. Syst.},
  pages={301--307},
  year={2001},
  publisher={MIT; 1998}
}

@Article{mezic2013analysis,
  author    = {Mezi{\'c}, Igor},
  journal   = {Annual Review of Fluid Mechanics},
  title     = {Analysis of fluid flows via spectral properties of the {K}oopman operator},
  year      = {2013},
  issn      = {1545-4479},
  month     = jan,
  number    = {1},
  pages     = {357--378},
  volume    = {45},
  doi       = {10.1146/annurev-fluid-011212-140652},
  publisher = {Annual Reviews},
}

@Article{chen2012variants,
  author    = {Chen, Kevin K. and Tu, Jonathan H. and Rowley, Clarence W.},
  journal   = {Journal of Nonlinear Science},
  title     = {Variants of dynamic mode decomposition: {B}oundary condition, {K}oopman, and {F}ourier analyses},
  year      = {2012},
  issn      = {1432-1467},
  month     = apr,
  number    = {6},
  pages     = {887--915},
  volume    = {22},
  doi       = {10.1007/s00332-012-9130-9},
  publisher = {Springer Science and Business Media LLC},
}

@Article{wynn2013optimal,
  author    = {Wynn, A. and Pearson, D. S. and Ganapathisubramani, B. and Goulart, P. J.},
  journal   = {Journal of Fluid Mechanics},
  title     = {Optimal mode decomposition for unsteady flows},
  year      = {2013},
  issn      = {1469-7645},
  month     = sep,
  pages     = {473--503},
  volume    = {733},
  doi       = {10.1017/jfm.2013.426},
  publisher = {Cambridge University Press (CUP)},
}

@Article{schmid2010dynamic,
  author    = {Schmid, Peter J.},
  journal   = {Journal of Fluid Mechanics},
  title     = {Dynamic mode decomposition of numerical and experimental data},
  year      = {2010},
  issn      = {1469-7645},
  month     = jul,
  pages     = {5--28},
  volume    = {656},
  doi       = {10.1017/s0022112010001217},
  publisher = {Cambridge University Press (CUP)},
}

@Article{li1976finite,
  author    = {Li, Tien-Yien},
  journal   = {Journal of Approximation Theory},
  title     = {Finite approximation for the {F}robenius-{P}erron operator. {A} solution to {U}lam's conjecture},
  year      = {1976},
  issn      = {0021-9045},
  month     = jun,
  number    = {2},
  pages     = {177--186},
  volume    = {17},
  doi       = {10.1016/0021-9045(76)90037-x},
  publisher = {Elsevier BV},
}

@Article{mauroy2013isostables,
  author    = {Mauroy, Alexandre and Mezi{\'c}, Igor and Moehlis, Jeff},
  journal   = {Physica D: Nonlinear Phenomena},
  title     = {Isostables, isochrons, and {K}oopman spectrum for the action--angle representation of stable fixed point dynamics},
  year      = {2013},
  issn      = {0167-2789},
  month     = oct,
  pages     = {19--30},
  volume    = {261},
  doi       = {10.1016/j.physd.2013.06.004},
  publisher = {Elsevier},
}

@Article{lorenz1963deterministic,
  author    = {Lorenz, Edward N.},
  journal   = {Journal of the Atmospheric Sciences},
  title     = {Deterministic Nonperiodic Flow},
  year      = {1963},
  issn      = {1520-0469},
  month     = mar,
  number    = {2},
  pages     = {130--141},
  volume    = {20},
  doi       = {10.1175/1520-0469(1963)020<0130:dnf>2.0.co;2},
  publisher = {American Meteorological Society},
}

@Article{klus2020eigendecompositions,
  author    = {Klus, Stefan and Schuster, Ingmar and Muandet, Krikamol},
  journal   = {Journal of Nonlinear Science},
  title     = {Eigendecompositions of transfer operators in reproducing kernel {H}ilbert spaces},
  year      = {2020},
  issn      = {1432-1467},
  month     = aug,
  number    = {1},
  pages     = {283--315},
  volume    = {30},
  doi       = {10.1007/s00332-019-09574-z},
  publisher = {Springer Science and Business Media LLC},
}

@Article{2158-2491_2016_1_51,
  author    = {Stefan Klus and Peter Koltai and Christof Sch\"{u}tte},
  journal   = {Journal of Computational Dynamics},
  title     = {On the numerical approximation of the {P}erron-{F}robenius and {K}oopman operator},
  year      = {2016},
  issn      = {2158-2491},
  month     = sep,
  number    = {1},
  pages     = {1--12},
  volume    = {3},
  doi       = {10.3934/jcd.2016003},
  publisher = {American Institute of Mathematical Sciences (AIMS)},
}

@Article{koopman1932dynamical,
  author    = {Koopman, Bernard O. and von Neumann, John},
  journal   = {Proceedings of the National Academy of Sciences},
  title     = {Dynamical Systems of Continuous Spectra},
  year      = {1932},
  issn      = {1091-6490},
  month     = mar,
  number    = {3},
  pages     = {255--263},
  volume    = {18},
  doi       = {10.1073/pnas.18.3.255},
  publisher = {National Academy of Sciences},
}

@Article{koopman1931hamiltonian,
  author    = {Koopman, Bernard O.},
  journal   = {Proceedings of the National Academy of Sciences},
  title     = {Hamiltonian systems and transformation in {H}ilbert space},
  year      = {1931},
  issn      = {1091-6490},
  month     = may,
  number    = {5},
  pages     = {315--318},
  volume    = {17},
  doi       = {10.1073/pnas.17.5.315},
  publisher = {National Academy of Sciences},
}

@InProceedings{koch2021large,
  author    = {Koch, R{\'e}gis and Sanjos{\'e}, Marl{\`e}ne and Moreau, Stephane},
  booktitle = {AIAA AVIATION 2021 FORUM},
  title     = {Large-Eddy Simulation of a Linear Compressor Cascade with Tip Gap: {A}erodynamic and Acoustic Analysis},
  year      = {2021},
  month     = jul,
  pages     = {2312},
  publisher = {American Institute of Aeronautics and Astronautics},
  doi       = {10.2514/6.2021-2312},
}

@Article{klus2018data,
  author    = {Klus, Stefan and N{\"u}ske, Feliks and Koltai, P{\'e}ter and Wu, Hao and Kevrekidis, Ioannis and Sch{\"u}tte, Christof and No{\'e}, Frank},
  journal   = {Journal of Nonlinear Science},
  title     = {Data-driven model reduction and transfer operator approximation},
  year      = {2018},
  issn      = {1432-1467},
  month     = jan,
  number    = {3},
  pages     = {985--1010},
  volume    = {28},
  doi       = {10.1007/s00332-017-9437-7},
  publisher = {Springer Science and Business Media LLC},
}

@Article{nuske2014variational,
  author    = {N{\"u}ske, Feliks and Keller, Bettina G. and P{\'e}rez-Hern{\'a}ndez, Guillermo and Mey, Antonia S. J. S. and No{\'e}, Frank},
  journal   = {Journal of Chemical Theory and Computation},
  title     = {Variational Approach to Molecular Kinetics},
  year      = {2014},
  issn      = {1549-9626},
  month     = mar,
  number    = {4},
  pages     = {1739--1752},
  volume    = {10},
  doi       = {10.1021/ct4009156},
  publisher = {American Chemical Society (ACS)},
}

@Article{klus2018kernel,
  author    = {Klus, Stefan and Bittracher, Andreas and Schuster, Ingmar and Sch{\"u}tte, Christof},
  journal   = {The Journal of Chemical Physics},
  title     = {A kernel-based approach to molecular conformation analysis},
  year      = {2018},
  issn      = {1089-7690},
  month     = dec,
  number    = {24},
  pages     = {244109},
  volume    = {149},
  doi       = {10.1063/1.5063533},
  publisher = {AIP Publishing},
}

@Article{mardt2018vampnets,
  author    = {Mardt, Andreas and Pasquali, Luca and Wu, Hao and No{\'e}, Frank},
  journal   = {Nature Communications},
  title     = {{VAMP}nets for deep learning of molecular kinetics},
  year      = {2018},
  issn      = {2041-1723},
  month     = jan,
  number    = {1},
  pages     = {1--11},
  volume    = {9},
  doi       = {10.1038/s41467-017-02388-1},
  publisher = {Springer Science and Business Media LLC},
}

@Book{ulam1960collection,
  author    = {Ulam, Stanislaw M.},
  publisher = {Science Editions},
  title     = {Problems in Modern Mathematics},
  year      = {1960},
}

@Article{otto2019linearly,
  author    = {Otto, Samuel E. and Rowley, Clarence W.},
  journal   = {SIAM Journal on Applied Dynamical Systems},
  title     = {Linearly Recurrent Autoencoder Networks for Learning Dynamics},
  year      = {2019},
  issn      = {1536-0040},
  month     = jan,
  number    = {1},
  pages     = {558--593},
  volume    = {18},
  doi       = {10.1137/18m1177846},
  publisher = {Society for Industrial & Applied Mathematics (SIAM)},
}

@Article{baddoo2021kernel,
  author    = {Baddoo, Peter J. and Herrmann, Benjamin and McKeon, Beverley J. and Brunton, Steven L.},
  journal   = {Proceedings of the Royal Society A: Mathematical, Physical and Engineering Sciences},
  title     = {Kernel learning for robust dynamic mode decomposition: {L}inear and nonlinear disambiguation optimization},
  year      = {2022},
  issn      = {1471-2946},
  month     = apr,
  number    = {2260},
  pages     = {20210830},
  volume    = {478},
  doi       = {10.1098/rspa.2021.0830},
  publisher = {The Royal Society},
}

@Article{lusch2018deep,
  author    = {Lusch, Bethany and Kutz, J. Nathan and Brunton, Steven L.},
  journal   = {Nature Communications},
  title     = {Deep learning for universal linear embeddings of nonlinear dynamics},
  year      = {2018},
  issn      = {2041-1723},
  month     = nov,
  number    = {1},
  pages     = {1--10},
  volume    = {9},
  doi       = {10.1038/s41467-018-07210-0},
  publisher = {Springer Science and Business Media LLC},
}

@InProceedings{yeung2019learning,
  author       = {Yeung, Enoch and Kundu, Soumya and Hodas, Nathan},
  booktitle    = {2019 American Control Conference (ACC)},
  title        = {Learning deep neural network representations for {K}oopman operators of nonlinear dynamical systems},
  year         = {2019},
  month        = jul,
  organization = {IEEE},
  pages        = {4832--4839},
  publisher    = {IEEE},
  doi          = {10.23919/acc.2019.8815339},
}

@Article{li2017extended,
  author    = {Li, Qianxiao and Dietrich, Felix and Bollt, Erik M. and Kevrekidis, Ioannis G.},
  journal   = {Chaos: An Interdisciplinary Journal of Nonlinear Science},
  title     = {Extended dynamic mode decomposition with dictionary learning: {A} data-driven adaptive spectral decomposition of the {K}oopman operator},
  year      = {2017},
  issn      = {1089-7682},
  month     = oct,
  number    = {10},
  pages     = {103111},
  volume    = {27},
  doi       = {10.1063/1.4993854},
  publisher = {AIP Publishing},
}

@Article{mezic2005spectral,
  author    = {Mezi{\'c}, Igor},
  journal   = {Nonlinear Dynamics},
  title     = {Spectral properties of dynamical systems, model reduction and decompositions},
  year      = {2005},
  issn      = {1573-269X},
  month     = aug,
  number    = {1},
  pages     = {309--325},
  volume    = {41},
  doi       = {10.1007/s11071-005-2824-x},
  publisher = {Springer Science and Business Media LLC},
}

@Article{pan2020structure,
  author    = {Pan, Shaowu and Duraisamy, Karthik},
  journal   = {Chaos: An Interdisciplinary Journal of Nonlinear Science},
  title     = {On the structure of time-delay embedding in linear models of non-linear dynamical systems},
  year      = {2020},
  issn      = {1089-7682},
  month     = jul,
  number    = {7},
  pages     = {073135},
  volume    = {30},
  doi       = {10.1063/5.0010886},
  publisher = {AIP Publishing},
}

@Article{kamb2020time,
  author    = {Kamb, Mason and Kaiser, Eurika and Brunton, Steven L. and Kutz, J. Nathan},
  journal   = {SIAM Journal on Applied Dynamical Systems},
  title     = {Time-delay observables for {K}oopman: {T}heory and applications},
  year      = {2020},
  issn      = {1536-0040},
  month     = jan,
  number    = {2},
  pages     = {886--917},
  volume    = {19},
  doi       = {10.1137/18m1216572},
  publisher = {Society for Industrial & Applied Mathematics (SIAM)},
}

@Article{eckmann1985ergodic,
  author    = {Eckmann, J.-P. and Ruelle, D.},
  journal   = {Reviews of Modern Physics},
  title     = {Ergodic theory of chaos and strange attractors},
  year      = {1985},
  issn      = {0034-6861},
  month     = jul,
  number    = {3},
  pages     = {617--656},
  volume    = {57},
  doi       = {10.1103/revmodphys.57.617},
  publisher = {American Physical Society (APS)},
}

@Article{tufillaro1992experimental,
  author    = {Tufillaro, Nicholas B. and Abbot, Tyler and Reilly, Jeremiah and Hickey, F. Roger},
  journal   = {American Journal of Physics},
  title     = {An Experimental Approach to Nonlinear Dynamics and Chaos.},
  year      = {1993},
  issn      = {1943-2909},
  month     = oct,
  number    = {10},
  pages     = {958--959},
  volume    = {61},
  doi       = {10.1119/1.17380},
  publisher = {American Association of Physics Teachers (AAPT)},
}

@Article{arbabi2017ergodic,
  author    = {Arbabi, Hassan and Mezi{\'c}, Igor},
  journal   = {SIAM Journal on Applied Dynamical Systems},
  title     = {Ergodic Theory, Dynamic Mode Decomposition, and Computation of Spectral Properties of the {K}oopman Operator},
  year      = {2017},
  issn      = {1536-0040},
  month     = jan,
  number    = {4},
  pages     = {2096--2126},
  volume    = {16},
  doi       = {10.1137/17m1125236},
  publisher = {Society for Industrial & Applied Mathematics (SIAM)},
}

@Article{kaiser2017data,
  author    = {Kaiser, Eurika and Kutz, J. Nathan and Brunton, Steven L.},
  journal   = {Machine Learning: Science and Technology},
  title     = {Data-driven discovery of {K}oopman eigenfunctions for control},
  year      = {2021},
  issn      = {2632-2153},
  month     = jun,
  number    = {3},
  pages     = {035023},
  volume    = {2},
  doi       = {10.1088/2632-2153/abf0f5},
  publisher = {IOP Publishing},
}

@Article{luzzatto2005lorenz,
  author    = {Luzzatto, Stefano and Melbourne, Ian and Paccaut, Frederic},
  journal   = {Communications in Mathematical Physics},
  title     = {The {L}orenz attractor is mixing},
  year      = {2005},
  issn      = {1432-0916},
  month     = aug,
  number    = {2},
  pages     = {393--401},
  volume    = {260},
  doi       = {10.1007/s00220-005-1411-9},
  publisher = {Springer Science and Business Media LLC},
}

@Book{halmos2017lectures,
  author    = {Halmos, Paul R.},
  publisher = {Dover Publications},
  title     = {Lectures on ergodic theory},
  year      = {2017},
  address   = {Mineola, New York},
  edition   = {Unabridged republication of the work originally published in 1956},
  isbn      = {9780486826844},
  series    = {Dover Books on Mathematics},
  pagetotal = {1113},
  ppn_gvk   = {1007362278},
}

@Article{giannakis2019data,
  author    = {Giannakis, Dimitrios},
  journal   = {Applied and Computational Harmonic Analysis},
  title     = {Data-driven spectral decomposition and forecasting of ergodic dynamical systems},
  year      = {2019},
  issn      = {1063-5203},
  month     = sep,
  number    = {2},
  pages     = {338--396},
  volume    = {47},
  doi       = {10.1016/j.acha.2017.09.001},
  publisher = {Elsevier BV},
}

@Article{brunton2016extracting,
  author    = {Brunton, Bingni W. and Johnson, Lise A. and Ojemann, Jeffrey G. and Kutz, J. Nathan},
  journal   = {Journal of Neuroscience Methods},
  title     = {Extracting spatial-temporal coherent patterns in large-scale neural recordings using dynamic mode decomposition},
  year      = {2016},
  issn      = {0165-0270},
  month     = jan,
  pages     = {1--15},
  volume    = {258},
  doi       = {10.1016/j.jneumeth.2015.10.010},
  publisher = {Elsevier BV},
}

@Article{das2021reproducing,
  author    = {Das, Suddhasattwa and Giannakis, Dimitrios and Slawinska, Joanna},
  journal   = {Applied and Computational Harmonic Analysis},
  title     = {Reproducing kernel {H}ilbert space compactification of unitary evolution groups},
  year      = {2021},
  issn      = {1063-5203},
  month     = sep,
  pages     = {75--136},
  volume    = {54},
  doi       = {10.1016/j.acha.2021.02.004},
  publisher = {Elsevier BV},
}

@Article{poincare1899methodes,
  author    = {Poincar{\'e}, Henri},
  journal   = {Il Nuovo Cimento},
  title     = {Les M{\'e}thodes Nouvelles de la M{\'e}canique C{\'e}leste},
  year      = {1899},
  issn      = {1827-6121},
  month     = jul,
  number    = {1},
  pages     = {128--130},
  volume    = {10},
  doi       = {10.1007/bf02742713},
  publisher = {Springer Science and Business Media LLC},
}

@Article{zhao2016analog,
  author    = {Zhao, Zhizhen and Giannakis, Dimitrios},
  journal   = {Nonlinearity},
  title     = {Analog forecasting with dynamics-adapted kernels},
  year      = {2016},
  issn      = {1361-6544},
  month     = aug,
  number    = {9},
  pages     = {2888--2939},
  volume    = {29},
  doi       = {10.1088/0951-7715/29/9/2888},
  publisher = {IOP Publishing},
}

@Article{das2019delay,
  author    = {Das, Suddhasattwa and Giannakis, Dimitrios},
  journal   = {Journal of Statistical Physics},
  title     = {Delay-coordinate maps and the spectra of {K}oopman operators},
  year      = {2019},
  issn      = {1572-9613},
  month     = apr,
  number    = {6},
  pages     = {1107--1145},
  volume    = {175},
  doi       = {10.1007/s10955-019-02272-w},
  publisher = {Springer Science and Business Media LLC},
}

@Article{giannakis2021learning,
  author    = {Giannakis, Dimitrios and Henriksen, Amelia and Tropp, Joel A. and Ward, Rachel},
  journal   = {SIAM Journal on Applied Dynamical Systems},
  title     = {Learning to Forecast Dynamical Systems from Streaming Data},
  year      = {2023},
  issn      = {1536-0040},
  month     = may,
  number    = {2},
  pages     = {527--558},
  volume    = {22},
  doi       = {10.1137/21m144983x},
  publisher = {Society for Industrial & Applied Mathematics (SIAM)},
}

@Article{burov2021kernel,
  author    = {Burov, Dmitry and Giannakis, Dimitrios and Manohar, Krithika and Stuart, Andrew},
  journal   = {Multiscale Modeling \& Simulation},
  title     = {Kernel analog forecasting: {M}ultiscale test problems},
  year      = {2021},
  issn      = {1540-3467},
  month     = jan,
  number    = {2},
  pages     = {1011--1040},
  volume    = {19},
  doi       = {10.1137/20m1338289},
  publisher = {Society for Industrial \& Applied Mathematics (SIAM)},
}

@Article{birkhoff1931proof,
  author    = {Birkhoff, George D.},
  journal   = {Proceedings of the National Academy of Sciences},
  title     = {Proof of the Ergodic Theorem},
  year      = {1931},
  issn      = {1091-6490},
  month     = dec,
  number    = {12},
  pages     = {656--660},
  volume    = {17},
  doi       = {10.1073/pnas.17.2.656},
  publisher = {National Acad Sciences},
}

@Article{govindarajan2019approximation,
  author    = {Govindarajan, Nithin and Mohr, Ryan and Chandrasekaran, Shivkumar and Mezi{\'c}, Igor},
  journal   = {SIAM Journal on Applied Dynamical Systems},
  title     = {On the Approximation of {K}oopman Spectra for Measure Preserving Transformations},
  year      = {2019},
  issn      = {1536-0040},
  month     = jan,
  number    = {3},
  pages     = {1454--1497},
  volume    = {18},
  doi       = {10.1137/18m1175094},
  publisher = {Society for Industrial & Applied Mathematics (SIAM)},
}

@Book{arnold1989mathematical,
  author    = {Arnold, V. I.},
  publisher = {Springer New York},
  title     = {Mathematical Methods of Classical Mechanics},
  year      = {1989},
  isbn      = {9781475720631},
  doi       = {10.1007/978-1-4757-2063-1},
  issn      = {0072-5285},
  journal   = {Graduate Texts in Mathematics},
  volum     = {60},
}

@Article{das2020koopman,
  author    = {Das, Suddhasattwa and Giannakis, Dimitrios},
  journal   = {Applied and Computational Harmonic Analysis},
  title     = {Koopman spectra in reproducing kernel {H}ilbert spaces},
  year      = {2020},
  issn      = {1063-5203},
  month     = sep,
  number    = {2},
  pages     = {573--607},
  volume    = {49},
  doi       = {10.1016/j.acha.2020.05.008},
  publisher = {Elsevier BV},
}

@Article{drmac2018data,
  author    = {Drma{\v{c}}, Zlatko and Mezi{\'c}, Igor and Mohr, Ryan},
  journal   = {SIAM Journal on Scientific Computing},
  title     = {Data Driven Modal Decompositions: {A}nalysis and Enhancements},
  year      = {2018},
  issn      = {1095-7197},
  month     = jan,
  number    = {4},
  pages     = {A2253--A2285},
  volume    = {40},
  doi       = {10.1137/17m1144155},
  publisher = {Society for Industrial & Applied Mathematics (SIAM)},
}

@Article{kachurovskii1996rate,
  author    = {Kachurovskii, Alexander Grigoryevich},
  journal   = {Russian Mathematical Surveys},
  title     = {The rate of convergence in ergodic theorems},
  year      = {1996},
  issn      = {1468-4829},
  month     = aug,
  number    = {4},
  pages     = {653--703},
  volume    = {51},
  doi       = {10.1070/rm1996v051n04abeh002964},
  publisher = {Steklov Mathematical Institute},
}

@Article{shargorodsky2008level,
  author    = {Shargorodsky, E.},
  journal   = {Bulletin of the London Mathematical Society},
  title     = {On the level sets of the resolvent norm of a linear operator},
  year      = {2008},
  issn      = {0024-6093},
  month     = may,
  number    = {3},
  pages     = {493--504},
  volume    = {40},
  doi       = {10.1112/blms/bdn038},
  publisher = {Wiley},
}

@Book{trefethen2005spectra,
  author    = {Trefethen, Lloyd N. and Embree, Mark},
  publisher = {Princeton University Press},
  title     = {Spectra and Pseudospectra: {T}he Behavior of Nonnormal Matrices and Operators},
  year      = {2005},
  isbn      = {9780691213101},
  month     = jan,
  doi       = {10.1515/9780691213101},
}

@article{roch1996c,
  title={{$C^*$}-algebra techniques in numerical analysis},
  author={Roch, Steffen and Silbermann, Bernd},
  journal={J. Oper. Theory},
  pages={241--280},
  year={1996},
  publisher={JSTOR}
}

@Article{williams2015data,
  author    = {Williams, Matthew O. and Kevrekidis, Ioannis G. and Rowley, Clarence W.},
  journal   = {Journal of Nonlinear Science},
  title     = {A data--driven approximation of the {K}oopman operator: {E}xtending dynamic mode decomposition},
  year      = {2015},
  issn      = {1432-1467},
  month     = jun,
  number    = {6},
  pages     = {1307--1346},
  volume    = {25},
  doi       = {10.1007/s00332-015-9258-5},
  publisher = {Springer Science and Business Media LLC},
}

@Article{korda2018convergence,
  author    = {Korda, Milan and Mezi{\'c}, Igor},
  journal   = {Journal of Nonlinear Science},
  title     = {On convergence of extended dynamic mode decomposition to the {K}oopman operator},
  year      = {2018},
  issn      = {1432-1467},
  month     = nov,
  number    = {2},
  pages     = {687--710},
  volume    = {28},
  doi       = {10.1007/s00332-017-9423-0},
  publisher = {Springer Science and Business Media LLC},
}

@Book{daubechies1992ten,
  author    = {Daubechies, Ingrid},
  publisher = {Society for Industrial and Applied Mathematics},
  title     = {Ten Lectures on Wavelets},
  year      = {1992},
  isbn      = {9781611970104},
  month     = jan,
  doi       = {10.1137/1.9781611970104},
}

@Book{katznelson2004introduction,
  author    = {Katznelson, Yitzhak},
  publisher = {Cambridge University Press},
  title     = {An Introduction to Harmonic Analysis},
  year      = {2004},
  doi       = {10.1017/cbo9781139165372},
}

@Article{williams2015kernel,
  author    = {O. Williams, Matthew and W. Rowley, Clarence and G. Kevrekidis, Ioannis},
  journal   = {Journal of Computational Dynamics},
  title     = {A kernel-based method for data-driven {K}oopman spectral analysis},
  year      = {2015},
  issn      = {2158-2505},
  number    = {2},
  pages     = {247--265},
  volume    = {2},
  doi       = {10.3934/jcd.2015005},
  publisher = {American Institute of Mathematical Sciences (AIMS)},
}

@Book{kutz2016dynamic,
  author    = {Kutz, J. Nathan and Brunton, Steven L. and Brunton, Bingni W. and Proctor, Joshua L.},
  publisher = {Society for Industrial and Applied Mathematics},
  title     = {Dynamic Mode Decomposition: {D}ata-Driven Modeling of Complex Systems},
  year      = {2016},
  isbn      = {9781611974508},
  month     = nov,
  doi       = {10.1137/1.9781611974508},
}

@Article{rowley2009spectral,
  author    = {Rowley, Clarence W. and Mezi{\'c}, Igor and Bagheri, Shervin and Schlatter, Philipp and Henningson, Dan S.},
  journal   = {Journal of Fluid Mechanics},
  title     = {Spectral analysis of nonlinear flows},
  year      = {2009},
  issn      = {1469-7645},
  month     = nov,
  pages     = {115--127},
  volume    = {641},
  doi       = {10.1017/s0022112009992059},
  publisher = {Cambridge University Press (CUP)},
}

@Article{tu2014dynamic,
  author    = {H. Tu, Jonathan and W. Rowley, Clarence and M. Luchtenburg, Dirk and L. Brunton, Steven and Nathan Kutz, J.},
  journal   = {Journal of Computational Dynamics},
  title     = {On dynamic mode decomposition: Theory and applications},
  year      = {2014},
  issn      = {2158-2505},
  number    = {2},
  pages     = {391--421},
  volume    = {1},
  doi       = {10.3934/jcd.2014.1.391},
  publisher = {American Institute of Mathematical Sciences (AIMS)},
}

@Article{halmos1963does,
  author    = {Halmos, P. R.},
  journal   = {The American Mathematical Monthly},
  title     = {What Does the Spectral Theorem Say?},
  year      = {1963},
  issn      = {1930-0972},
  month     = mar,
  number    = {3},
  pages     = {241--247},
  volume    = {70},
  doi       = {10.1080/00029890.1963.11990075},
  fjournal  = {Amer. Math. Monthly},
  publisher = {Informa UK Limited},
}

@Article{hemati2017biasing,
  author    = {Hemati, Maziar S. and Rowley, Clarence W. and Deem, Eric A. and Cattafesta, Louis N.},
  journal   = {Theoretical and Computational Fluid Dynamics},
  title     = {De-biasing the dynamic mode decomposition for applied {K}oopman spectral analysis of noisy datasets},
  year      = {2017},
  issn      = {1432-2250},
  month     = apr,
  number    = {4},
  pages     = {349--368},
  volume    = {31},
  doi       = {10.1007/s00162-017-0432-2},
  publisher = {Springer Science and Business Media LLC},
}

@Article{dawson2016characterizing,
  author    = {Dawson, Scott T. M. and Hemati, Maziar S. and Williams, Matthew O. and Rowley, Clarence W.},
  journal   = {Experiments in Fluids},
  title     = {Characterizing and correcting for the effect of sensor noise in the dynamic mode decomposition},
  year      = {2016},
  issn      = {1432-1114},
  month     = feb,
  number    = {3},
  pages     = {1--19},
  volume    = {57},
  doi       = {10.1007/s00348-016-2127-7},
  publisher = {Springer Science and Business Media LLC},
}

@Article{bottcher1983finite,
  author    = {B{\"o}ttcher, Albrecht and Silbermann, Bernd},
  journal   = {Mathematische Nachrichten},
  title     = {The finite section method for {T}oeplitz operators on the quarter-plane with piecewise continuous symbols},
  year      = {1983},
  issn      = {1522-2616},
  number    = {1},
  pages     = {279--291},
  volume    = {110},
  doi       = {10.1002/mana.19831100120},
  publisher = {Wiley},
}

@Book{pazy2012semigroups,
  author    = {Pazy, Amnon},
  publisher = {Springer},
  title     = {Semigroups of linear operators and applications to partial differential equations},
  year      = {2010},
  address   = {New York},
  edition   = {3},
  isbn      = {978-0-387-90845-8},
  series    = {Applied mathematical sciences},
  volume    = {44},
  pagetotal = {279},
  ppn_gvk   = {620391103},
}

@InProceedings{schmid2009dynamic,
  author    = {Schmid, Peter J.},
  booktitle = {8th International Symposium on Particle Image Velocimetry (PIV09)},
  title     = {Dynamic mode decomposition of experimental data},
  year      = {2009},
}

@Article{proctor2016dynamic,
  author    = {Proctor, Joshua L. and Brunton, Steven L. and Kutz, J. Nathan},
  journal   = {SIAM Journal on Applied Dynamical Systems},
  title     = {Dynamic Mode Decomposition with Control},
  year      = {2016},
  issn      = {1536-0040},
  month     = jan,
  number    = {1},
  pages     = {142--161},
  volume    = {15},
  doi       = {10.1137/15m1013857},
  publisher = {Society for Industrial & Applied Mathematics (SIAM)},
}

@Article{vcrnjaric2020koopman,
  author    = {{\v{C}}rnjari{\'c}-{\v{Z}}ic, Nelida and Ma{\'c}e{\v{s}}i{\'c}, Senka and Mezi{\'c}, Igor},
  journal   = {Journal of Nonlinear Science},
  title     = {Koopman Operator Spectrum for Random Dynamical Systems},
  year      = {2020},
  issn      = {1432-1467},
  month     = sep,
  number    = {5},
  pages     = {2007--2056},
  volume    = {30},
  doi       = {10.1007/s00332-019-09582-z},
  publisher = {Springer Science and Business Media LLC},
}

@Article{wanner2022robust,
  author    = {Wanner, Mathias and Mezi{\'c}, Igor},
  journal   = {SIAM Journal on Applied Dynamical Systems},
  title     = {Robust approximation of the stochastic {K}oopman operator},
  year      = {2022},
  issn      = {1536-0040},
  month     = jul,
  number    = {3},
  pages     = {1930--1951},
  volume    = {21},
  doi       = {10.1137/21m1414425},
  publisher = {Society for Industrial & Applied Mathematics (SIAM)},
}

@Article{coifman2006diffusion,
  author    = {Coifman, Ronald R. and Lafon, St{\'e}phane},
  journal   = {Applied and Computational Harmonic Analysis},
  title     = {Diffusion maps},
  year      = {2006},
  issn      = {1063-5203},
  month     = jul,
  number    = {1},
  pages     = {5--30},
  volume    = {21},
  doi       = {10.1016/j.acha.2006.04.006},
  publisher = {Elsevier BV},
}

@Article{takeishi2017subspace,
  author    = {Takeishi, Naoya and Kawahara, Yoshinobu and Yairi, Takehisa},
  journal   = {Physical Review E},
  title     = {Subspace dynamic mode decomposition for stochastic {K}oopman analysis},
  year      = {2017},
  issn      = {2470-0053},
  month     = sep,
  number    = {3},
  pages     = {033310},
  volume    = {96},
  doi       = {10.1103/physreve.96.033310},
  publisher = {American Physical Society (APS)},
}

@Article{sinha2020robust,
  author    = {Sinha, Subhrajit and Huang, Bowen and Vaidya, Umesh},
  journal   = {Journal of Nonlinear Science},
  title     = {On robust computation of {K}oopman operator and prediction in random dynamical systems},
  year      = {2020},
  issn      = {1432-1467},
  month     = nov,
  number    = {5},
  pages     = {2057--2090},
  volume    = {30},
  doi       = {10.1007/s00332-019-09597-6},
  publisher = {Springer},
}

@Article{hesthaven2022reduced,
  author    = {Hesthaven, Jan S. and Pagliantini, Cecilia and Rozza, Gianluigi},
  journal   = {Acta Numerica},
  title     = {Reduced basis methods for time-dependent problems},
  year      = {2022},
  issn      = {1474-0508},
  month     = may,
  pages     = {265--345},
  volume    = {31},
  doi       = {10.1017/s0962492922000058},
  publisher = {Cambridge University Press (CUP)},
}

@Article{nuske2023finite,
  author    = {N{\"u}ske, Feliks and Peitz, Sebastian and Philipp, Friedrich and Schaller, Manuel and Worthmann, Karl},
  journal   = {Journal of Nonlinear Science},
  title     = {Finite-data error bounds for {K}oopman-based prediction and control},
  year      = {2023},
  issn      = {1432-1467},
  month     = nov,
  number    = {1},
  pages     = {14},
  volume    = {33},
  doi       = {10.1007/s00332-022-09862-1},
  publisher = {Springer Science and Business Media LLC},
}

@Article{rowley2005model,
  author    = {Rowley, Clarence W.},
  journal   = {International Journal of Bifurcation and Chaos},
  title     = {Model reduction for fluids, using balanced proper orthogonal decomposition},
  year      = {2005},
  issn      = {1793-6551},
  month     = mar,
  number    = {03},
  pages     = {997--1013},
  volume    = {15},
  doi       = {10.1142/s0218127405012429},
  publisher = {World Scientific Pub Co Pte Lt},
}

@Article{berkooz1993proper,
  author    = {Berkooz, Gal and Holmes, Philip and Lumley, John L.},
  journal   = {Annual Review of Fluid Mechanics},
  title     = {The Proper Orthogonal Decomposition in the Analysis of Turbulent Flows},
  year      = {1993},
  issn      = {1545-4479},
  month     = jan,
  number    = {1},
  pages     = {539--575},
  volume    = {25},
  doi       = {10.1146/annurev.fl.25.010193.002543},
  publisher = {Annual Reviews},
}

@Article{wu2020variational,
  author    = {Wu, Hao and No{\'e}, Frank},
  journal   = {Journal of Nonlinear Science},
  title     = {Variational approach for learning {M}arkov processes from time series data},
  year      = {2020},
  issn      = {1432-1467},
  month     = aug,
  number    = {1},
  pages     = {23--66},
  volume    = {30},
  doi       = {10.1007/s00332-019-09567-y},
  publisher = {Springer Science and Business Media LLC},
}

@Article{klus2020kernel,
  author    = {Klus, Stefan and N{\"u}ske, Feliks and Hamzi, Boumediene},
  journal   = {Entropy},
  title     = {Kernel-based approximation of the {K}oopman generator and {S}chr{\"o}dinger operator},
  year      = {2020},
  issn      = {1099-4300},
  month     = jun,
  number    = {7},
  pages     = {722},
  volume    = {22},
  doi       = {10.3390/e22070722},
  publisher = {MDPI AG},
}

@Article{lu2020prediction,
  author    = {Lu, Hannah and Tartakovsky, Daniel M.},
  journal   = {SIAM Journal on Scientific Computing},
  title     = {Prediction Accuracy of Dynamic Mode Decomposition},
  year      = {2020},
  issn      = {1095-7197},
  month     = jan,
  number    = {3},
  pages     = {A1639--A1662},
  volume    = {42},
  doi       = {10.1137/19m1259948},
  publisher = {Society for Industrial & Applied Mathematics (SIAM)},
}

@Article{duke2012error,
  author    = {Duke, Daniel and Soria, Julio and Honnery, Damon},
  journal   = {Experiments in Fluids},
  title     = {An error analysis of the dynamic mode decomposition},
  year      = {2012},
  issn      = {1432-1114},
  month     = dec,
  number    = {2},
  pages     = {529--542},
  volume    = {52},
  doi       = {10.1007/s00348-011-1235-7},
  publisher = {Springer Science and Business Media LLC},
}

@Article{mollenhauer2022kernel,
  author  = {Mollenhauer, Mattes and Klus, Stefan and Sch{\"u}tte, Christof and Koltai, P{\'e}ter},
  journal = {Journal of Machine Learning Research},
  title   = {Kernel autocovariance operators of stationary processes: {E}stimation and convergence},
  year    = {2022},
  number  = {327},
  pages   = {1--34},
  volume  = {23},
}

@Article{carleman1932application,
  author    = {Carleman, Torsten},
  journal   = {Acta Mathematica},
  title     = {Application de la th{\'e}orie des {\'e}quations int{\'e}grales lin{\'e}aires aux syst{\`e}mes d'{\'e}quations diff{\'e}rentielles non lin{\'e}aires},
  year      = {1932},
  issn      = {0001-5962},
  number    = {0},
  pages     = {63--87},
  volume    = {59},
  doi       = {10.1007/bf02546499},
  publisher = {International Press of Boston},
}

@Article{gavish2014optimal,
  author    = {Gavish, Matan and Donoho, David L.},
  journal   = {IEEE Transactions on Information Theory},
  title     = {The optimal hard threshold for singular values is {$4/\sqrt{3}$}},
  year      = {2014},
  issn      = {1557-9654},
  month     = aug,
  number    = {8},
  pages     = {5040--5053},
  volume    = {60},
  doi       = {10.1109/tit.2014.2323359},
  publisher = {Institute of Electrical and Electronics Engineers (IEEE)},
}

@Article{udell2019big,
  author    = {Udell, Madeleine and Townsend, Alex},
  journal   = {SIAM Journal on Mathematics of Data Science},
  title     = {Why Are Big Data Matrices Approximately Low Rank?},
  year      = {2019},
  issn      = {2577-0187},
  month     = jan,
  number    = {1},
  pages     = {144--160},
  volume    = {1},
  doi       = {10.1137/18m1183480},
  publisher = {Society for Industrial & Applied Mathematics (SIAM)},
}

@InBook{drmavc2020dynamic,
  author    = {Drma{\v{c}}, Zlatko},
  pages     = {161--194},
  publisher = {Springer International Publishing},
  title     = {Dynamic Mode Decomposition - {A} Numerical Linear Algebra Perspective},
  year      = {2020},
  isbn      = {9783030357139},
  booktitle = {The Koopman Operator in Systems and Control},
  doi       = {10.1007/978-3-030-35713-9_7},
  issn      = {1610-7411},
  journal   = {The Koopman Operator in Systems and Control: Concepts, Methodologies, and Applications},
}

@Article{drmac2019data,
  author    = {Drma{\v{c}}, Zlatko and Mezi{\'c}, Igor and Mohr, Ryan},
  journal   = {SIAM Journal on Scientific Computing},
  title     = {Data driven {K}oopman spectral analysis in {V}andermonde--{C}auchy form via the {DFT}: {N}umerical method and theoretical insights},
  year      = {2019},
  number    = {5},
  pages     = {A3118--A3151},
  volume    = {41},
  doi       = {10.1137/18m1227688},
  publisher = {SIAM},
}

@Article{colbrook2023residualJFM,
  author    = {Colbrook, Matthew J. and Ayton, Lorna J. and Sz{\H{o}}ke, M{\'a}t{\'e}},
  journal   = {Journal of Fluid Mechanics},
  title     = {Residual dynamic mode decomposition: {R}obust and verified {K}oopmanism},
  year      = {2023},
  issn      = {1469-7645},
  month     = jan,
  pages     = {A21},
  volume    = {955},
  doi       = {10.1017/jfm.2022.1052},
  publisher = {Cambridge University Press (CUP)},
}

@Article{jackson1987finite,
  author    = {Jackson, C. P.},
  journal   = {Journal of Fluid Mechanics},
  title     = {A finite-element study of the onset of vortex shedding in flow past variously shaped bodies},
  year      = {1987},
  issn      = {1469-7645},
  month     = sep,
  number    = {1},
  pages     = {23--45},
  volume    = {182},
  doi       = {10.1017/s0022112087002234},
  publisher = {Cambridge University Press (CUP)},
}

@Article{zebib1987stability,
  author    = {Zebib, A.},
  journal   = {Journal of Engineering Mathematics},
  title     = {Stability of viscous flow past a circular cylinder},
  year      = {1987},
  issn      = {1573-2703},
  number    = {2},
  pages     = {155--165},
  volume    = {21},
  doi       = {10.1007/bf00127673},
  publisher = {Springer Science and Business Media LLC},
}

@Article{noack2003hierarchy,
  author    = {Noack, Bernd R. and Afanasiev, Konstantin and Morzynski, Marek and Tadmor, Gilead and Thiele, Frank},
  journal   = {Journal of Fluid Mechanics},
  title     = {A hierarchy of low-dimensional models for the transient and post-transient cylinder wake},
  year      = {2003},
  issn      = {1469-7645},
  month     = dec,
  pages     = {335--363},
  volume    = {497},
  doi       = {10.1017/s0022112003006694},
  publisher = {Cambridge University Press (CUP)},
}

@Article{noack1994global,
  author    = {Noack, Bernd R. and Eckelmann, Helmut},
  journal   = {Journal of Fluid Mechanics},
  title     = {A global stability analysis of the steady and periodic cylinder wake},
  year      = {1994},
  issn      = {1469-7645},
  month     = jul,
  pages     = {297--330},
  volume    = {270},
  doi       = {10.1017/s0022112094004283},
  publisher = {Cambridge University Press (CUP)},
}

@Article{rowley2017model,
  author    = {Rowley, Clarence W. and Dawson, Scott T. M.},
  journal   = {Annual Review of Fluid Mechanics},
  title     = {Model Reduction for Flow Analysis and Control},
  year      = {2017},
  issn      = {1545-4479},
  month     = jan,
  number    = {1},
  pages     = {387--417},
  volume    = {49},
  doi       = {10.1146/annurev-fluid-010816-060042},
  publisher = {Annual Reviews},
}

@Article{bagheri2013koopman,
  author    = {Bagheri, Shervin},
  journal   = {Journal of Fluid Mechanics},
  title     = {Koopman-mode decomposition of the cylinder wake},
  year      = {2013},
  issn      = {1469-7645},
  month     = jun,
  pages     = {596--623},
  volume    = {726},
  doi       = {10.1017/jfm.2013.249},
  publisher = {Cambridge University Press (CUP)},
}

@Article{taira2020modal,
  author    = {Taira, Kunihiko and Hemati, Maziar S. and Brunton, Steven L. and Sun, Yiyang and Duraisamy, Karthik and Bagheri, Shervin and Dawson, Scott T. M. and Yeh, Chi-An},
  journal   = {AIAA Journal},
  title     = {Modal analysis of fluid flows: {A}pplications and outlook},
  year      = {2020},
  issn      = {1533-385X},
  month     = mar,
  number    = {3},
  pages     = {998--1022},
  volume    = {58},
  doi       = {10.2514/1.j058462},
  publisher = {American Institute of Aeronautics and Astronautics (AIAA)},
}

@InProceedings{jozsa2016validation,
  author     = {J{\'o}zsa, Tam{\'a}s and Sz{\H{o}}ke, M{\'a}t{\'e} and Teschner, Tom-Robin and K{\"o}n{\"o}zsy, L{\'a}szl{\'o} Z. and Moulitsas, Irene},
  booktitle  = {Proceedings of the VII European Congress on Computational Methods in Applied Sciences and Engineering (ECCOMAS Congress 2016)},
  title      = {Validation and verification of a {2D} lattice {B}oltzmann solver for incompressible fluid flow},
  year       = {2016},
  publisher  = {Institute of Structural Analysis and Antiseismic Research School of Civil Engineering National Technical University of Athens (NTUA) Greece},
  series     = {ECCOMAS Congress 2016},
  doi        = {10.7712/100016.1869.10678},
}

@Article{szHoke2017performance,
  author    = {Sz{\H{o}}ke, M{\'a}t{\'e} and Jozsa, Tamas Istvan and Kolesz{\'a}r, {\'A}d{\'a}m and Moulitsas, Irene and K{\"o}n{\"o}zsy, L{\'a}szl{\'o}},
  journal   = {ACM Transactions on Mathematical Software},
  title     = {Performance evaluation of a two-dimensional lattice {B}oltzmann solver using {CUDA} and {PGAS} {UPC} based parallelisation},
  year      = {2017},
  issn      = {1557-7295},
  month     = jul,
  number    = {1},
  pages     = {1--22},
  volume    = {44},
  doi       = {10.1145/3085590},
  publisher = {Association for Computing Machinery (ACM)},
}

@inproceedings{takens2006detecting,
  title={Detecting strange attractors in turbulence},
  author={Takens, Floris},
  booktitle={Dynamical Systems and Turbulence, Warwick 1980: proceedings of a symposium held at the University of Warwick 1979/80},
  pages={366--381},
  year={2006},
  organization={Springer}
}

@Article{wu2021challenges,
  author    = {Wu, Ziyou and Brunton, Steven L. and Revzen, Shai},
  journal   = {Journal of The Royal Society Interface},
  title     = {Challenges in dynamic mode decomposition},
  year      = {2021},
  issn      = {1742-5662},
  month     = dec,
  number    = {185},
  pages     = {20210686},
  volume    = {18},
  doi       = {10.1098/rsif.2021.0686},
  publisher = {The Royal Society},
}

@Article{zhang2019online,
  author    = {Zhang, Hao and Rowley, Clarence W. and Deem, Eric A. and Cattafesta, Louis N.},
  journal   = {SIAM Journal on Applied Dynamical Systems},
  title     = {Online Dynamic Mode Decomposition for Time-Varying Systems},
  year      = {2019},
  issn      = {1536-0040},
  month     = jan,
  number    = {3},
  pages     = {1586--1609},
  volume    = {18},
  doi       = {10.1137/18m1192329},
  publisher = {Society for Industrial & Applied Mathematics (SIAM)},
}

@Article{macesic2018koopman,
  author    = {Ma{\'c}e{\v{s}}i{\'c}, Senka and {\v{C}}rnjari{\'c}-{\v{Z}}ic, Nelida and Mezi{\'c}, Igor},
  journal   = {SIAM Journal on Applied Dynamical Systems},
  title     = {Koopman operator family spectrum for nonautonomous systems},
  year      = {2018},
  issn      = {1536-0040},
  month     = jan,
  number    = {4},
  pages     = {2478--2515},
  volume    = {17},
  doi       = {10.1137/17m1133610},
  publisher = {Society for Industrial & Applied Mathematics (SIAM)},
}

@Article{mezic2016koopman,
  author    = {Mezi{\'c}, Igor and Surana, Amit},
  journal   = {IFAC-PapersOnLine},
  title     = {Koopman mode decomposition for periodic/quasi-periodic time dependence},
  year      = {2016},
  issn      = {2405-8963},
  number    = {18},
  pages     = {690--697},
  volume    = {49},
  doi       = {10.1016/j.ifacol.2016.10.246},
  publisher = {Elsevier BV},
}

@Article{giannakis2020extraction,
  author    = {Giannakis, Dimitrios and Das, Suddhasattwa},
  journal   = {Physica D: Nonlinear Phenomena},
  title     = {Extraction and prediction of coherent patterns in incompressible flows through space-time {K}oopman analysis},
  year      = {2020},
  issn      = {0167-2789},
  month     = jan,
  pages     = {132211},
  volume    = {402},
  doi       = {10.1016/j.physd.2019.132211},
  publisher = {Elsevier BV},
}

@Article{giannakis2019spatiotemporal,
  author    = {Giannakis, Dimitrios and Ourmazd, Abbas and Slawinska, Joanna and Zhao, Zhizhen},
  journal   = {Journal of Nonlinear Science},
  title     = {Spatiotemporal Pattern Extraction by Spectral Analysis of Vector-Valued Observables},
  year      = {2019},
  issn      = {1432-1467},
  month     = may,
  number    = {5},
  pages     = {2385--2445},
  volume    = {29},
  doi       = {10.1007/s00332-019-09548-1},
  publisher = {Springer Science and Business Media LLC},
}

@Article{kutz2016multiresolution,
  author    = {Kutz, J. Nathan and Fu, Xing and Brunton, Steven L.},
  journal   = {SIAM Journal on Applied Dynamical Systems},
  title     = {Multiresolution Dynamic Mode Decomposition},
  year      = {2016},
  issn      = {1536-0040},
  month     = jan,
  number    = {2},
  pages     = {713--735},
  volume    = {15},
  doi       = {10.1137/15m1023543},
  publisher = {Society for Industrial & Applied Mathematics (SIAM)},
}

@Article{nyquist1928certain,
  author    = {Nyquist, H.},
  journal   = {Transactions of the American Institute of Electrical Engineers},
  title     = {Certain topics in telegraph transmission theory},
  year      = {1928},
  issn      = {0096-3860},
  month     = apr,
  number    = {2},
  pages     = {617--644},
  volume    = {47},
  doi       = {10.1109/t-aiee.1928.5055024},
  publisher = {Institute of Electrical and Electronics Engineers (IEEE)},
}

@Article{shannon1948mathematical,
  author    = {Shannon, C. E.},
  journal   = {Bell System Technical Journal},
  title     = {A mathematical theory of communication},
  year      = {1948},
  issn      = {0005-8580},
  month     = jul,
  number    = {3},
  pages     = {379--423},
  volume    = {27},
  doi       = {10.1002/j.1538-7305.1948.tb01338.x},
  publisher = {Institute of Electrical and Electronics Engineers (IEEE)},
}

@Article{jovanovic2014sparsity,
  author    = {Jovanovi{\'c}, Mihailo R. and Schmid, Peter J. and Nichols, Joseph W.},
  journal   = {Physics of Fluids},
  title     = {Sparsity-promoting dynamic mode decomposition},
  year      = {2014},
  issn      = {1089-7666},
  month     = feb,
  number    = {2},
  volume    = {26},
  doi       = {10.1063/1.4863670},
  publisher = {AIP Publishing},
}

@Article{donoho2006compressed,
  author    = {Donoho, D. L.},
  journal   = {IEEE Transactions on Information Theory},
  title     = {Compressed sensing},
  year      = {2006},
  issn      = {0018-9448},
  month     = apr,
  number    = {4},
  pages     = {1289--1306},
  volume    = {52},
  doi       = {10.1109/tit.2006.871582},
  publisher = {Institute of Electrical and Electronics Engineers (IEEE)},
}

@Article{candes2006robust,
  author    = {Candes, E. J. and Romberg, J. and Tao, T.},
  journal   = {IEEE Transactions on Information Theory},
  title     = {Robust uncertainty principles: {E}xact signal reconstruction from highly incomplete frequency information},
  year      = {2006},
  issn      = {0018-9448},
  month     = feb,
  number    = {2},
  pages     = {489--509},
  volume    = {52},
  doi       = {10.1109/tit.2005.862083},
  publisher = {Institute of Electrical and Electronics Engineers (IEEE)},
}

@Book{foucart2013invitation,
  author    = {Foucart, Simon and Rauhut, Holger},
  publisher = {Springer New York},
  title     = {A Mathematical Introduction to Compressive Sensing},
  year      = {2013},
  isbn      = {9780817649487},
  doi       = {10.1007/978-0-8176-4948-7},
  issn      = {2296-5017},
  journal   = {Applied and Numerical Harmonic Analysis},
}

@Book{adcock2021compressive,
  author    = {Adcock, Ben and Hansen, Anders C.},
  publisher = {Cambridge University Press},
  title     = {Compressive Imaging: {S}tructure, Sampling, Learning},
  year      = {2021},
  isbn      = {9781108421614},
  month     = jul,
  doi       = {10.1017/9781108377447},
}

@Article{candes2008introduction,
  author    = {Candes, E. J. and Wakin, M. B.},
  journal   = {IEEE Signal Processing Magazine},
  title     = {An Introduction To Compressive Sampling},
  year      = {2008},
  issn      = {1053-5888},
  month     = mar,
  number    = {2},
  pages     = {21--30},
  volume    = {25},
  doi       = {10.1109/msp.2007.914731},
  publisher = {Institute of Electrical and Electronics Engineers (IEEE)},
}

@Article{tu2014spectral,
  author    = {Tu, Jonathan H. and Rowley, Clarence W. and Kutz, J. Nathan and Shang, Jessica K.},
  journal   = {Experiments in Fluids},
  title     = {Spectral analysis of fluid flows using sub-{N}yquist-rate {PIV} data},
  year      = {2014},
  issn      = {1432-1114},
  month     = sep,
  number    = {9},
  pages     = {1--13},
  volume    = {55},
  doi       = {10.1007/s00348-014-1805-6},
  publisher = {Springer Science and Business Media LLC},
}

@Article{brunton2016compressed,
  author    = {Brunton, Steven L. and Proctor, Joshua L. and Tu, Jonathan H. and Kutz, J. Nathan},
  journal   = {Journal of computational dynamics},
  title     = {Compressed sensing and dynamic mode decomposition},
  year      = {2016},
  number    = {2},
  pages     = {165--191},
  volume    = {2},
  publisher = {Journal of Computational Dynamics},
}

@Article{candes2006near,
  author    = {Candes, Emmanuel J. and Tao, Terence},
  journal   = {IEEE Transactions on Information Theory},
  title     = {Near-optimal signal recovery from random projections: {U}niversal encoding strategies?},
  year      = {2006},
  issn      = {0018-9448},
  month     = dec,
  number    = {12},
  pages     = {5406--5425},
  volume    = {52},
  doi       = {10.1109/tit.2006.885507},
  publisher = {Institute of Electrical and Electronics Engineers (IEEE)},
}

@Article{needell2009cosamp,
  author    = {Needell, D. and Tropp, J. A.},
  journal   = {Applied and Computational Harmonic Analysis},
  title     = {Co{S}a{MP}: {I}terative signal recovery from incomplete and inaccurate samples},
  year      = {2009},
  issn      = {1063-5203},
  month     = may,
  number    = {3},
  pages     = {301--321},
  volume    = {26},
  doi       = {10.1016/j.acha.2008.07.002},
  publisher = {Elsevier},
}

@Article{kou2017improved,
  author    = {Kou, Jiaqing and Zhang, Weiwei},
  journal   = {European Journal of Mechanics - B/Fluids},
  title     = {An improved criterion to select dominant modes from dynamic mode decomposition},
  year      = {2017},
  issn      = {0997-7546},
  month     = mar,
  pages     = {109--129},
  volume    = {62},
  doi       = {10.1016/j.euromechflu.2016.11.015},
  publisher = {Elsevier BV},
}

@Article{erichson2019compressed,
  author    = {Erichson, N. Benjamin and Brunton, Steven L. and Kutz, J. Nathan},
  journal   = {Journal of Real-Time Image Processing},
  title     = {Compressed dynamic mode decomposition for background modeling},
  year      = {2019},
  issn      = {1861-8219},
  month     = nov,
  number    = {5},
  pages     = {1479--1492},
  volume    = {16},
  doi       = {10.1007/s11554-016-0655-2},
  publisher = {Springer Science and Business Media LLC},
}

@Article{huhn2023parametric,
  author    = {Huhn, Quincy A. and Tano, Mauricio E. and Ragusa, Jean C. and Choi, Youngsoo},
  journal   = {Journal of Computational Physics},
  title     = {Parametric dynamic mode decomposition for reduced order modeling},
  year      = {2023},
  issn      = {0021-9991},
  month     = feb,
  pages     = {111852},
  volume    = {475},
  doi       = {10.1016/j.jcp.2022.111852},
  publisher = {Elsevier BV},
}

@InProceedings{sinha2019computation,
  author       = {Sinha, Subhrajit and Vaidya, Umesh and Yeung, Enoch},
  booktitle    = {2019 American Control Conference (ACC)},
  title        = {On computation of {K}oopman operator from sparse data},
  year         = {2019},
  month        = jul,
  organization = {IEEE},
  pages        = {5519--5524},
  publisher    = {IEEE},
  doi          = {10.23919/acc.2019.8814861},
}

@InBook{clainche2018spatio,
  author    = {Vega, Jos{\'e} M. and Le Clainche, Soledad},
  pages     = {121--157},
  publisher = {Elsevier},
  title     = {Spatio-temporal {K}oopman decomposition},
  year      = {2021},
  volume    = {28},
  booktitle = {Higher Order Dynamic Mode Decomposition and Its Applications},
  doi       = {10.1016/b978-0-12-819743-1.00011-2},
  journal   = {Journal of Nonlinear Science},
}

@Article{manohar2019optimized,
  author    = {Manohar, Krithika and Kaiser, Eurika and Brunton, Steven L. and Kutz, J. Nathan},
  journal   = {Multiscale Modeling \& Simulation},
  title     = {Optimized Sampling for Multiscale Dynamics},
  year      = {2019},
  issn      = {1540-3467},
  month     = jan,
  number    = {1},
  pages     = {117--136},
  volume    = {17},
  doi       = {10.1137/17m1162366},
  publisher = {Society for Industrial & Applied Mathematics (SIAM)},
}

@InProceedings{takeishi2017bayesian,
  author     = {Takeishi, Naoya and Kawahara, Yoshinobu and Tabei, Yasuo and Yairi, Takehisa},
  booktitle  = {Proceedings of the Twenty-Sixth International Joint Conference on Artificial Intelligence},
  title      = {Bayesian Dynamic Mode Decomposition},
  year       = {2017},
  month      = aug,
  pages      = {2814--2821},
  publisher  = {International Joint Conferences on Artificial Intelligence Organization},
  series     = {IJCAI-2017},
  collection = {IJCAI-2017},
  doi        = {10.24963/ijcai.2017/392},
}

@Article{loiseau2018constrained,
  author    = {Loiseau, Jean-Christophe and Brunton, Steven L.},
  journal   = {Journal of Fluid Mechanics},
  title     = {Constrained sparse {G}alerkin regression},
  year      = {2018},
  issn      = {1469-7645},
  month     = jan,
  pages     = {42--67},
  volume    = {838},
  doi       = {10.1017/jfm.2017.823},
  publisher = {Cambridge University Press (CUP)},
}

@Article{khosravi2023representer,
  author    = {Khosravi, Mohammad},
  journal   = {IEEE Transactions on Automatic Control},
  title     = {Representer Theorem for Learning {K}oopman Operators},
  year      = {2023},
  issn      = {2334-3303},
  month     = may,
  number    = {5},
  pages     = {2995--3010},
  volume    = {68},
  doi       = {10.1109/tac.2023.3242325},
  publisher = {Institute of Electrical and Electronics Engineers (IEEE)},
}

@Article{noack2016recursive,
  author    = {Noack, Bernd R. and Stankiewicz, Witold and Morzy{\'n}ski, Marek and Schmid, Peter J.},
  journal   = {Journal of Fluid Mechanics},
  title     = {Recursive dynamic mode decomposition of transient and post-transient wake flows},
  year      = {2016},
  issn      = {1469-7645},
  month     = nov,
  pages     = {843--872},
  volume    = {809},
  doi       = {10.1017/jfm.2016.678},
  publisher = {Cambridge University Press (CUP)},
}

@Article{klus2018tensor,
  author    = {Klus, Stefan and Gel{\ss}, Patrick and Peitz, Sebastian and Sch{\"u}tte, Christof},
  journal   = {Nonlinearity},
  title     = {Tensor-based dynamic mode decomposition},
  year      = {2018},
  issn      = {1361-6544},
  month     = jun,
  number    = {7},
  pages     = {3359--3380},
  volume    = {31},
  doi       = {10.1088/1361-6544/aabc8f},
  publisher = {IOP Publishing},
}

@Article{klus_towards,
  author    = {Stefan Klus and Christof Sch{\"u}tte},
  journal   = {Journal of Computational Dynamics},
  title     = {Towards tensor-based methods for the numerical approximation of the {P}erron--{F}robenius and {K}oopman operator},
  year      = {2016},
  issn      = {2158-2491},
  month     = nov,
  number    = {2},
  pages     = {139--161},
  volume    = {3},
  doi       = {10.3934/jcd.2016007},
  publisher = {American Institute of Mathematical Sciences (AIMS)},
}

@Article{tissot2014model,
  author    = {Tissot, Gilles and Cordier, Laurent and Benard, Nicolas and Noack, Bernd R.},
  journal   = {Comptes Rendus M{\'e}canique},
  title     = {Model reduction using dynamic mode decomposition},
  year      = {2014},
  issn      = {1631-0721},
  month     = jun,
  number    = {6-7},
  pages     = {410--416},
  volume    = {342},
  doi       = {10.1016/j.crme.2013.12.011},
  publisher = {Cellule MathDoc/CEDRAM},
}

@Article{andreuzzi2023dynamic,
  author    = {Andreuzzi, Francesco and Demo, Nicola and Rozza, Gianluigi},
  journal   = {SIAM Journal on Applied Dynamical Systems},
  title     = {A dynamic mode decomposition extension for the forecasting of parametric dynamical systems},
  year      = {2023},
  issn      = {1536-0040},
  month     = aug,
  number    = {3},
  pages     = {2432--2458},
  volume    = {22},
  doi       = {10.1137/22m1481658},
  publisher = {Society for Industrial & Applied Mathematics (SIAM)},
}

@Article{hirsh2020centering,
  author    = {Hirsh, Seth M. and Harris, Kameron Decker and Kutz, J. Nathan and Brunton, Bingni W.},
  journal   = {SIAM Journal on Applied Dynamical Systems},
  title     = {Centering data improves the dynamic mode decomposition},
  year      = {2020},
  issn      = {1536-0040},
  month     = jan,
  number    = {3},
  pages     = {1920--1955},
  volume    = {19},
  doi       = {10.1137/19m1289881},
  publisher = {Society for Industrial & Applied Mathematics (SIAM)},
}

@Article{bai2020dynamic,
  author    = {Bai, Zhe and Kaiser, Eurika and Proctor, Joshua L. and Kutz, J. Nathan and Brunton, Steven L.},
  journal   = {AIAA Journal},
  title     = {Dynamic mode decomposition for compressive system identification},
  year      = {2020},
  issn      = {1533-385X},
  month     = feb,
  number    = {2},
  pages     = {561--574},
  volume    = {58},
  doi       = {10.2514/1.j057870},
  publisher = {American Institute of Aeronautics and Astronautics (AIAA)},
}

@Article{williams2016extending,
  author    = {Williams, Matthew O. and Hemati, Maziar S. and Dawson, Scott T. M. and Kevrekidis, Ioannis G. and Rowley, Clarence W.},
  journal   = {IFAC-PapersOnLine},
  title     = {Extending data-driven {K}oopman analysis to actuated systems},
  year      = {2016},
  issn      = {2405-8963},
  number    = {18},
  pages     = {704--709},
  volume    = {49},
  doi       = {10.1016/j.ifacol.2016.10.248},
  publisher = {Elsevier BV},
}

@Article{champion2019discovery,
  author    = {Champion, Kathleen P. and Brunton, Steven L. and Kutz, J. Nathan},
  journal   = {SIAM Journal on Applied Dynamical Systems},
  title     = {Discovery of nonlinear multiscale systems: {S}ampling strategies and embeddings},
  year      = {2019},
  issn      = {1536-0040},
  month     = jan,
  number    = {1},
  pages     = {312--333},
  volume    = {18},
  doi       = {10.1137/18m1188227},
  publisher = {Society for Industrial & Applied Mathematics (SIAM)},
}

@Article{alla2017nonlinear,
  author    = {Alla, Alessandro and Kutz, J. Nathan},
  journal   = {SIAM Journal on Scientific Computing},
  title     = {Nonlinear model order reduction via dynamic mode decomposition},
  year      = {2017},
  issn      = {1095-7197},
  month     = jan,
  number    = {5},
  pages     = {B778--B796},
  volume    = {39},
  doi       = {10.1137/16m1059308},
  publisher = {Society for Industrial & Applied Mathematics (SIAM)},
}

@Article{drmac2020least,
  author    = {Drma{\v{c}}, Zlatko and Mezi{\'c}, Igor and Mohr, Ryan},
  journal   = {SIAM Journal on Scientific Computing},
  title     = {On least squares problems with certain {V}andermonde--{K}hatri--{R}ao structure with applications to {DMD}},
  year      = {2020},
  issn      = {1095-7197},
  month     = jan,
  number    = {5},
  pages     = {A3250--A3284},
  volume    = {42},
  doi       = {10.1137/19m1288474},
  publisher = {Society for Industrial & Applied Mathematics (SIAM)},
}

@Article{jiang2022correcting,
  author    = {Jiang, Lijian and Liu, Ningxin},
  journal   = {Journal of Computational Physics},
  title     = {Correcting noisy dynamic mode decomposition with {K}alman filters},
  year      = {2022},
  issn      = {0021-9991},
  month     = jul,
  pages     = {111175},
  volume    = {461},
  doi       = {10.1016/j.jcp.2022.111175},
  publisher = {Elsevier BV},
}

@Article{pan2015accuracy,
  author    = {Pan, Chong and Xue, Dong and Wang, Jinjun},
  journal   = {Experiments in Fluids},
  title     = {On the accuracy of dynamic mode decomposition in estimating instability of wave packet},
  year      = {2015},
  issn      = {1432-1114},
  month     = aug,
  number    = {8},
  pages     = {1--15},
  volume    = {56},
  doi       = {10.1007/s00348-015-2015-6},
  publisher = {Springer Science and Business Media LLC},
}

@Article{le2017higher,
  author    = {Le Clainche, Soledad and Vega, Jos{\'e} M.},
  journal   = {SIAM Journal on Applied Dynamical Systems},
  title     = {Higher order dynamic mode decomposition},
  year      = {2017},
  issn      = {1536-0040},
  month     = jan,
  number    = {2},
  pages     = {882--925},
  volume    = {16},
  doi       = {10.1137/15m1054924},
  publisher = {Society for Industrial & Applied Mathematics (SIAM)},
}

@Article{askham2018variable,
  author    = {Askham, Travis and Kutz, J. Nathan},
  journal   = {SIAM Journal on Applied Dynamical Systems},
  title     = {Variable projection methods for an optimized dynamic mode decomposition},
  year      = {2018},
  issn      = {1536-0040},
  month     = jan,
  number    = {1},
  pages     = {380--416},
  volume    = {17},
  doi       = {10.1137/m1124176},
  publisher = {Society for Industrial & Applied Mathematics (SIAM)},
}

@Article{azencot2019consistent,
  author    = {Azencot, Omri and Yin, Wotao and Bertozzi, Andrea},
  journal   = {SIAM Journal on Applied Dynamical Systems},
  title     = {Consistent dynamic mode decomposition},
  year      = {2019},
  issn      = {1536-0040},
  month     = jan,
  number    = {3},
  pages     = {1565--1585},
  volume    = {18},
  doi       = {10.1137/18m1233960},
  publisher = {Society for Industrial & Applied Mathematics (SIAM)},
}

@Article{sashidhar2022bagging,
  author    = {Sashidhar, Diya and Kutz, J. Nathan},
  journal   = {Philosophical Transactions of the Royal Society A: Mathematical, Physical and Engineering Sciences},
  title     = {Bagging, optimized dynamic mode decomposition for robust, stable forecasting with spatial and temporal uncertainty quantification},
  year      = {2022},
  issn      = {1471-2962},
  month     = jun,
  number    = {2229},
  pages     = {20210199},
  volume    = {380},
  doi       = {10.1098/rsta.2021.0199},
  publisher = {The Royal Society},
}

@Article{halko2011finding,
  author    = {Halko, N. and Martinsson, P. G. and Tropp, J. A.},
  journal   = {SIAM Review},
  title     = {Finding structure with randomness: {P}robabilistic algorithms for constructing approximate matrix decompositions},
  year      = {2011},
  issn      = {1095-7200},
  month     = jan,
  number    = {2},
  pages     = {217--288},
  volume    = {53},
  doi       = {10.1137/090771806},
  publisher = {Society for Industrial & Applied Mathematics (SIAM)},
}

@Article{rowley2006linear,
  author    = {Rowley, Clarence W. and Williams, David R. and Colonius, Tim and Murray, Richard M. and Macmynowski, Douglas G.},
  journal   = {Journal of Fluid Mechanics},
  title     = {Linear models for control of cavity flow oscillations},
  year      = {2006},
  issn      = {1469-7645},
  month     = jan,
  number    = {1},
  pages     = {317--330},
  volume    = {547},
  doi       = {10.1017/s0022112005007299},
  publisher = {Cambridge University Press (CUP)},
}

@Book{higham2008functions,
  author    = {Higham, Nicholas J.},
  publisher = {Society for Industrial and Applied Mathematics},
  title     = {Functions of matrices: {T}heory and computation},
  year      = {2008},
  isbn      = {9780898717778},
  month     = jan,
  doi       = {10.1137/1.9780898717778},
}

@Article{bagheri2014effects,
  author    = {Bagheri, Shervin},
  journal   = {Physics of Fluids},
  title     = {Effects of weak noise on oscillating flows: {L}inking quality factor, {F}loquet modes, and {K}oopman spectrum},
  year      = {2014},
  issn      = {1089-7666},
  month     = sep,
  number    = {9},
  volume    = {26},
  doi       = {10.1063/1.4895898},
  publisher = {AIP Publishing},
}

@Article{nonomura2018dynamic,
  author    = {Nonomura, Taku and Shibata, Hisaichi and Takaki, Ryoji},
  journal   = {AIP Advances},
  title     = {Dynamic mode decomposition using a {K}alman filter for parameter estimation},
  year      = {2018},
  issn      = {2158-3226},
  month     = oct,
  number    = {10},
  volume    = {8},
  doi       = {10.1063/1.5031816},
  publisher = {AIP Publishing},
}

@Article{nonomura2019extended,
  author    = {Nonomura, Taku and Shibata, Hisaichi and Takaki, Ryoji},
  journal   = {PLOS ONE},
  title     = {Extended-{K}alman-filter-based dynamic mode decomposition for simultaneous system identification and denoising},
  year      = {2019},
  issn      = {1932-6203},
  month     = feb,
  number    = {2},
  pages     = {e0209836},
  volume    = {14},
  doi       = {10.1371/journal.pone.0209836},
  editor    = {Xin, Baogui},
  publisher = {Public Library of Science (PLoS)},
}

@Article{hemati2014dynamic,
  author    = {Hemati, Maziar S. and Williams, Matthew O. and Rowley, Clarence W.},
  journal   = {Physics of Fluids},
  title     = {Dynamic mode decomposition for large and streaming datasets},
  year      = {2014},
  issn      = {1089-7666},
  month     = nov,
  number    = {11},
  volume    = {26},
  doi       = {10.1063/1.4901016},
  publisher = {AIP Publishing},
}

@InProceedings{hemati2016improving,
  author    = {Hemati, Maziar and Deem, Eric and Williams, Matthew and Rowley, Clarence W. and Cattafesta, Louis N.},
  booktitle = {54th AIAA Aerospace Sciences Meeting},
  title     = {Improving separation control with noise-robust variants of dynamic mode decomposition},
  year      = {2016},
  month     = jan,
  pages     = {1103},
  publisher = {American Institute of Aeronautics and Astronautics},
  doi       = {10.2514/6.2016-1103},
}

@Article{golub1973differentiation,
  author    = {Golub, G. H. and Pereyra, V.},
  journal   = {SIAM Journal on Numerical Analysis},
  title     = {The differentiation of pseudo-inverses and nonlinear least squares problems whose variables separate},
  year      = {1973},
  issn      = {1095-7170},
  month     = apr,
  number    = {2},
  pages     = {413--432},
  volume    = {10},
  doi       = {10.1137/0710036},
  publisher = {Society for Industrial & Applied Mathematics (SIAM)},
}

@Book{pereyra2010exponential,
  author    = {Pereyra, Victor and Scherer, Godela},
  publisher = {Bentham eBooks},
  title     = {Exponential data fitting and its applications},
  year      = {2010},
  address   = {Sharjah, UAE},
  isbn      = {9781608050482},
  pagetotal = {195},
  ppn_gvk   = {1749037033},
}

@Article{gueniat2015dynamic,
  author    = {Gu{\'e}niat, Florimond and Mathelin, Lionel and Pastur, Luc R.},
  journal   = {Physics of Fluids},
  title     = {A dynamic mode decomposition approach for large and arbitrarily sampled systems},
  year      = {2015},
  issn      = {1089-7666},
  month     = feb,
  number    = {2},
  volume    = {27},
  doi       = {10.1063/1.4908073},
  publisher = {AIP Publishing},
}

@Article{leroux2016dynamic,
  author    = {Leroux, Romain and Cordier, Laurent},
  journal   = {Experiments in Fluids},
  title     = {Dynamic mode decomposition for non-uniformly sampled data},
  year      = {2016},
  issn      = {1432-1114},
  month     = may,
  number    = {5},
  pages     = {94},
  volume    = {57},
  doi       = {10.1007/s00348-016-2165-1},
  publisher = {Springer Science and Business Media LLC},
}

@Book{breiman2017classification,
  author    = {Breiman, Leo and Friedman, Jerome H. and Olshen, Richard A. and Stone, Charles J.},
  publisher = {Routledge},
  title     = {Classification and regression trees},
  year      = {2017},
  isbn      = {9781315139470},
  month     = oct,
  doi       = {10.1201/9781315139470},
}

@Article{nuske2021tensor,
  author    = {N{\"u}ske, Feliks and Gel{\ss}, Patrick and Klus, Stefan and Clementi, Cecilia},
  journal   = {Physica D: Nonlinear Phenomena},
  title     = {Tensor-based computation of metastable and coherent sets},
  year      = {2021},
  issn      = {0167-2789},
  month     = dec,
  pages     = {133018},
  volume    = {427},
  doi       = {10.1016/j.physd.2021.133018},
  publisher = {Elsevier BV},
}

@Article{fujii2019dynamic,
  author    = {Fujii, Keisuke and Kawahara, Yoshinobu},
  journal   = {Neural Networks},
  title     = {Dynamic mode decomposition in vector-valued reproducing kernel {H}ilbert spaces for extracting dynamical structure among observables},
  year      = {2019},
  issn      = {0893-6080},
  month     = sep,
  pages     = {94--103},
  volume    = {117},
  doi       = {10.1016/j.neunet.2019.04.020},
  publisher = {Elsevier BV},
}

@article{kawahara2016dynamic,
  title={Dynamic mode decomposition with reproducing kernels for {K}oopman spectral analysis},
  author={Kawahara, Yoshinobu},
  journal={Advances in neural information processing systems},
  volume={29},
  year={2016}
}

@Article{redman2021koopman,
  author    = {Redman, William T.},
  journal   = {Chaos: An Interdisciplinary Journal of Nonlinear Science},
  title     = {On {K}oopman mode decomposition and tensor component analysis},
  year      = {2021},
  issn      = {1089-7682},
  month     = may,
  number    = {5},
  volume    = {31},
  doi       = {10.1063/5.0046325},
  publisher = {AIP Publishing},
}

@Article{sharma2016correspondence,
  author    = {Sharma, Ati S. and Mezi{\'c}, Igor and McKeon, Beverley J.},
  journal   = {Physical Review Fluids},
  title     = {Correspondence between {K}oopman mode decomposition, resolvent mode decomposition, and invariant solutions of the {N}avier-{S}tokes equations},
  year      = {2016},
  issn      = {2469-990X},
  month     = jul,
  number    = {3},
  pages     = {032402},
  volume    = {1},
  doi       = {10.1103/physrevfluids.1.032402},
  publisher = {American Physical Society (APS)},
}

@Article{herrmann2021data,
  author    = {Herrmann, Benjamin and Baddoo, Peter J. and Semaan, Richard and Brunton, Steven L. and McKeon, Beverley J.},
  journal   = {Journal of Fluid Mechanics},
  title     = {Data-driven resolvent analysis},
  year      = {2021},
  issn      = {1469-7645},
  month     = may,
  pages     = {A10},
  volume    = {918},
  doi       = {10.1017/jfm.2021.337},
  publisher = {Cambridge University Press (CUP)},
}

@Article{towne2018spectral,
  author    = {Towne, Aaron and Schmidt, Oliver T. and Colonius, Tim},
  journal   = {Journal of Fluid Mechanics},
  title     = {Spectral proper orthogonal decomposition and its relationship to dynamic mode decomposition and resolvent analysis},
  year      = {2018},
  issn      = {1469-7645},
  month     = may,
  pages     = {821--867},
  volume    = {847},
  doi       = {10.1017/jfm.2018.283},
  publisher = {Cambridge University Press (CUP)},
}

@Article{scherl2020robust,
  author    = {Scherl, Isabel and Strom, Benjamin and Shang, Jessica K. and Williams, Owen and Polagye, Brian L. and Brunton, Steven L.},
  journal   = {Physical Review Fluids},
  title     = {Robust principal component analysis for modal decomposition of corrupt fluid flows},
  year      = {2020},
  issn      = {2469-990X},
  month     = may,
  number    = {5},
  pages     = {054401},
  volume    = {5},
  doi       = {10.1103/physrevfluids.5.054401},
  publisher = {American Physical Society (APS)},
}

@Article{martinsson2020randomized,
  author    = {Martinsson, Per-Gunnar and Tropp, Joel A.},
  journal   = {Acta Numerica},
  title     = {Randomized numerical linear algebra: {F}oundations and algorithms},
  year      = {2020},
  issn      = {1474-0508},
  month     = may,
  pages     = {403--572},
  volume    = {29},
  doi       = {10.1017/s0962492920000021},
  publisher = {Cambridge University Press (CUP)},
}

@Article{erichson2016randomized,
  author    = {Erichson, N. Benjamin and Donovan, Carl},
  journal   = {Computer Vision and Image Understanding},
  title     = {Randomized low-rank dynamic mode decomposition for motion detection},
  year      = {2016},
  issn      = {1077-3142},
  month     = may,
  pages     = {40--50},
  volume    = {146},
  doi       = {10.1016/j.cviu.2016.02.005},
  publisher = {Elsevier BV},
}

@Article{bistrian2017randomized,
  author    = {Bistrian, Diana Alina and Navon, Ionel Michael},
  journal   = {International Journal for Numerical Methods in Engineering},
  title     = {Randomized dynamic mode decomposition for nonintrusive reduced order modelling},
  year      = {2017},
  issn      = {1097-0207},
  month     = feb,
  number    = {1},
  pages     = {3--25},
  volume    = {112},
  doi       = {10.1002/nme.5499},
  publisher = {Wiley},
}

@Book{van1991total,
  author    = {Van Huffel, Sabine and Vandewalle, Joos},
  publisher = {Society for Industrial and Applied Mathematics},
  title     = {The total least squares problem: {C}omputational aspects and analysis},
  year      = {1991},
  isbn      = {9781611971002},
  month     = jan,
  doi       = {10.1137/1.9781611971002},
}

@Article{erichson2019randomized,
  author    = {Erichson, N. Benjamin and Mathelin, Lionel and Kutz, J. Nathan and Brunton, Steven L.},
  journal   = {SIAM Journal on Applied Dynamical Systems},
  title     = {Randomized dynamic mode decomposition},
  year      = {2019},
  issn      = {1536-0040},
  month     = jan,
  number    = {4},
  pages     = {1867--1891},
  volume    = {18},
  doi       = {10.1137/18m1215013},
  publisher = {Society for Industrial & Applied Mathematics (SIAM)},
}

@Book{hansen2006deblurring,
  author    = {Hansen, Per Christian and Nagy, James G. and O'Leary, Dianne P.},
  publisher = {Society for Industrial and Applied Mathematics},
  title     = {Deblurring images: {M}atrices, spectra, and filtering},
  year      = {2006},
  isbn      = {9780898718874},
  month     = jan,
  doi       = {10.1137/1.9780898718874},
}

@Article{rokhlin2010randomized,
  author    = {Rokhlin, Vladimir and Szlam, Arthur and Tygert, Mark},
  journal   = {SIAM Journal on Matrix Analysis and Applications},
  title     = {A randomized algorithm for principal component analysis},
  year      = {2010},
  issn      = {1095-7162},
  month     = jan,
  number    = {3},
  pages     = {1100--1124},
  volume    = {31},
  doi       = {10.1137/080736417},
  publisher = {Society for Industrial & Applied Mathematics (SIAM)},
}

@Article{gu2015subspace,
  author    = {Gu, M.},
  journal   = {SIAM Journal on Scientific Computing},
  title     = {Subspace iteration randomization and singular value problems},
  year      = {2015},
  issn      = {1095-7197},
  month     = jan,
  number    = {3},
  pages     = {A1139--A1173},
  volume    = {37},
  doi       = {10.1137/130938700},
  publisher = {Society for Industrial & Applied Mathematics (SIAM)},
}

@Article{ahmed2022dynamic,
  author    = {Ahmed, Shady E. and Dabaghian, Pedram H. and San, Omer and Bistrian, Diana A. and Navon, Ionel M.},
  journal   = {Physics of Fluids},
  title     = {Dynamic mode decomposition with core sketch},
  year      = {2022},
  issn      = {1089-7666},
  month     = jun,
  number    = {6},
  volume    = {34},
  doi       = {10.1063/5.0095163},
  publisher = {AIP Publishing},
}

@Article{rosenfeld2022dynamic,
  author    = {Rosenfeld, Joel A. and Kamalapurkar, Rushikesh and Gruss, L. Forest and Johnson, Taylor T.},
  journal   = {Journal of Nonlinear Science},
  title     = {Dynamic mode decomposition for continuous time systems with the {L}iouville operator},
  year      = {2022},
  issn      = {1432-1467},
  month     = dec,
  number    = {1},
  pages     = {1--30},
  volume    = {32},
  doi       = {10.1007/s00332-021-09746-w},
  publisher = {Springer Science and Business Media LLC},
}

@Article{philipp2023error,
  author  = {Philipp, Friedrich and Schaller, Manuel and Worthmann, Karl and Peitz, Sebastian and N{\"u}ske, Feliks},
  journal = {arXiv preprint arXiv:2301.08637},
  title   = {Error bounds for kernel-based approximations of the {K}oopman operator},
  year    = {2023},
}

@Article{degennaro2019scalable,
  author    = {DeGennaro, Anthony M. and Urban, Nathan M.},
  journal   = {SIAM Journal on Scientific Computing},
  title     = {Scalable extended dynamic mode decomposition using random kernel approximation},
  year      = {2019},
  issn      = {1095-7197},
  month     = jan,
  number    = {3},
  pages     = {A1482--A1499},
  volume    = {41},
  doi       = {10.1137/17m115414x},
  publisher = {Society for Industrial & Applied Mathematics (SIAM)},
}

@Article{alford2022deep,
  author    = {Alford-Lago, D. J. and Curtis, C. W. and Ihler, A. T. and Issan, O.},
  journal   = {Chaos: An Interdisciplinary Journal of Nonlinear Science},
  title     = {Deep learning enhanced dynamic mode decomposition},
  year      = {2022},
  issn      = {1089-7682},
  month     = mar,
  number    = {3},
  volume    = {32},
  doi       = {10.1063/5.0073893},
  publisher = {AIP Publishing},
}

@Article{li2021deep,
  author    = {Li, Mengnan and Jiang, Lijian},
  journal   = {Journal of Computational Physics},
  title     = {Deep learning nonlinear multiscale dynamic problems using {K}oopman operator},
  year      = {2021},
  issn      = {0021-9991},
  month     = dec,
  pages     = {110660},
  volume    = {446},
  doi       = {10.1016/j.jcp.2021.110660},
  publisher = {Elsevier BV},
}

@Article{takeishi2017learning,
  author  = {Takeishi, Naoya and Kawahara, Yoshinobu and Yairi, Takehisa},
  journal = {Advances in neural information processing systems},
  title   = {Learning {K}oopman invariant subspaces for dynamic mode decomposition},
  year    = {2017},
  volume  = {30},
}

@Article{wehmeyer2018time,
  author    = {Wehmeyer, Christoph and No{\'e}, Frank},
  journal   = {The Journal of Chemical Physics},
  title     = {Time-lagged autoencoders: {D}eep learning of slow collective variables for molecular kinetics},
  year      = {2018},
  issn      = {1089-7690},
  month     = mar,
  number    = {24},
  volume    = {148},
  doi       = {10.1063/1.5011399},
  publisher = {AIP Publishing},
}

@InProceedings{azencot2020forecasting,
  author       = {Azencot, Omri and Erichson, N. Benjamin and Lin, Vanessa and Mahoney, Michael},
  booktitle    = {International Conference on Machine Learning},
  title        = {Forecasting sequential data using consistent {K}oopman autoencoders},
  year         = {2020},
  organization = {PMLR},
  pages        = {475--485},
}

@Article{eivazi2021recurrent,
  author    = {Eivazi, Hamidreza and Guastoni, Luca and Schlatter, Philipp and Azizpour, Hossein and Vinuesa, Ricardo},
  journal   = {International Journal of Heat and Fluid Flow},
  title     = {Recurrent neural networks and {K}oopman-based frameworks for temporal predictions in a low-order model of turbulence},
  year      = {2021},
  issn      = {0142-727X},
  month     = aug,
  pages     = {108816},
  volume    = {90},
  doi       = {10.1016/j.ijheatfluidflow.2021.108816},
  publisher = {Elsevier BV},
}

@Article{xie2019graph,
  author    = {Xie, Tian and France-Lanord, Arthur and Wang, Yanming and Shao-Horn, Yang and Grossman, Jeffrey C.},
  journal   = {Nature Communications},
  title     = {Graph dynamical networks for unsupervised learning of atomic scale dynamics in materials},
  year      = {2019},
  issn      = {2041-1723},
  month     = jun,
  number    = {1},
  pages     = {2667},
  volume    = {10},
  doi       = {10.1038/s41467-019-10663-6},
  publisher = {Springer Science and Business Media LLC},
}

@Article{manojlovic2020applications,
  author  = {Manojlovi{\'c}, Iva and Fonoberova, Maria and Mohr, Ryan and Andrej{\v{c}}uk, Aleksandr and Drma{\v{c}}, Zlatko and Kevrekidis, Yannis and Mezi{\'c}, Igor},
  journal = {arXiv preprint arXiv:2006.11765},
  title   = {Applications of {K}oopman mode analysis to neural networks},
  year    = {2020},
}

@Article{dogra2020optimizing,
  author  = {Dogra, Akshunna S. and Redman, William},
  journal = {Advances in Neural Information Processing Systems},
  title   = {Optimizing neural networks via {K}oopman operator theory},
  year    = {2020},
  pages   = {2087--2097},
  volume  = {33},
}

@Article{dietrich2020koopman,
  author    = {Dietrich, Felix and Thiem, Thomas N. and Kevrekidis, Ioannis G.},
  journal   = {SIAM Journal on Applied Dynamical Systems},
  title     = {On the {K}oopman operator of algorithms},
  year      = {2020},
  issn      = {1536-0040},
  month     = jan,
  number    = {2},
  pages     = {860--885},
  volume    = {19},
  doi       = {10.1137/19m1277059},
  publisher = {Society for Industrial & Applied Mathematics (SIAM)},
}

@Article{colbrook2023beyond,
  author        = {Colbrook, Matthew J. and Li, Qin and Raut, Ryan V. and Townsend, Alex},
  journal       = {arXiv preprint arXiv:2308.10697},
  title         = {Beyond expectations: {R}esidual dynamic mode decomposition and variance for stochastic dynamical systems},
  year          = {2023},
  month         = aug,
  archiveprefix = {arXiv},
  copyright     = {arXiv.org perpetual, non-exclusive license},
  doi           = {10.48550/ARXIV.2308.10697},
  eprint        = {2308.10697},
  primaryclass  = {math.DS},
  publisher     = {arXiv},
}

@Article{colbrook2021rigorous,
  author    = {Colbrook, Matthew J. and Townsend, Alex},
  journal   = {Communications on Pure and Applied Mathematics},
  title     = {Rigorous data-driven computation of spectral properties of {K}oopman operators for dynamical systems},
  year      = {2023},
  issn      = {1097-0312},
  month     = jul,
  number    = {1},
  pages     = {221--283},
  volume    = {77},
  doi       = {10.1002/cpa.22125},
  publisher = {Wiley},
}

@Article{colbrook2023mpedmd,
  author    = {Colbrook, Matthew J.},
  journal   = {SIAM Journal on Numerical Analysis},
  title     = {The mp{EDMD} algorithm for data-driven computations of measure-preserving dynamical systems},
  year      = {2023},
  issn      = {1095-7170},
  month     = jun,
  number    = {3},
  pages     = {1585--1608},
  volume    = {61},
  doi       = {10.1137/22m1521407},
  publisher = {Society for Industrial & Applied Mathematics (SIAM)},
}

@Article{sechi2021estimation,
  author    = {Sechi, Renata and Sikorski, Alexander and Weber, Marcus},
  journal   = {Multiscale Modeling \& Simulation},
  title     = {Estimation of the {K}oopman generator by {N}ewton's extrapolation},
  year      = {2021},
  issn      = {1540-3467},
  month     = jan,
  number    = {2},
  pages     = {758--774},
  volume    = {19},
  doi       = {10.1137/20m1333006},
  publisher = {Society for Industrial & Applied Mathematics (SIAM)},
}

@Article{mauroy2019koopman,
  author    = {Mauroy, Alexandre and Goncalves, Jorge},
  journal   = {IEEE Transactions on Automatic Control},
  title     = {Koopman-based lifting techniques for nonlinear systems identification},
  year      = {2020},
  issn      = {2334-3303},
  month     = jun,
  number    = {6},
  pages     = {2550--2565},
  volume    = {65},
  doi       = {10.1109/tac.2019.2941433},
  publisher = {Institute of Electrical and Electronics Engineers (IEEE)},
}

@Article{mezic2002ergodic,
  author    = {Mezi{\'c}, Igor and Sotiropoulos, Fotis},
  journal   = {Physics of Fluids},
  title     = {Ergodic theory and experimental visualization of invariant sets in chaotically advected flows},
  year      = {2002},
  issn      = {1089-7666},
  month     = may,
  number    = {7},
  pages     = {2235--2243},
  volume    = {14},
  doi       = {10.1063/1.1480266},
  publisher = {AIP Publishing},
}

@Article{juang1985eigensystem,
  author    = {Juang, Jer-Nan and Pappa, Richard S.},
  journal   = {Journal of Guidance, Control, and Dynamics},
  title     = {An eigensystem realization algorithm for modal parameter identification and model reduction},
  year      = {1985},
  issn      = {1533-3884},
  month     = sep,
  number    = {5},
  pages     = {620--627},
  volume    = {8},
  doi       = {10.2514/3.20031},
  publisher = {American Institute of Aeronautics and Astronautics (AIAA)},
}

@Article{broomhead1989time,
  author    = {Broomhead, David S. and Jones, Roger},
  journal   = {Proceedings of the Royal Society of London. A. Mathematical and Physical Sciences},
  title     = {Time-series analysis},
  year      = {1989},
  number    = {1864},
  pages     = {103--121},
  volume    = {423},
  publisher = {The Royal Society London},
}

@Article{hirsh2021structured,
  author    = {Hirsh, Seth M. and Ichinaga, Sara M. and Brunton, Steven L. and Nathan Kutz, J. and Brunton, Bingni W.},
  journal   = {Proceedings of the Royal Society A: Mathematical, Physical and Engineering Sciences},
  title     = {Structured time-delay models for dynamical systems with connections to {F}renet--{S}erret frame},
  year      = {2021},
  issn      = {1471-2946},
  month     = oct,
  number    = {2254},
  pages     = {20210097},
  volume    = {477},
  doi       = {10.1098/rspa.2021.0097},
  publisher = {The Royal Society},
}

@Article{alexander2020operator,
  author    = {Alexander, Romeo and Giannakis, Dimitrios},
  journal   = {Physica D: Nonlinear Phenomena},
  title     = {Operator-theoretic framework for forecasting nonlinear time series with kernel analog techniques},
  year      = {2020},
  issn      = {0167-2789},
  month     = aug,
  pages     = {132520},
  volume    = {409},
  doi       = {10.1016/j.physd.2020.132520},
  publisher = {Elsevier BV},
}

@Article{kostic2022learning,
  author  = {Kostic, Vladimir and Novelli, Pietro and Maurer, Andreas and Ciliberto, Carlo and Rosasco, Lorenzo and Pontil, Massimiliano},
  journal = {Advances in Neural Information Processing Systems},
  title   = {Learning dynamical systems via {K}oopman operator regression in reproducing kernel {H}ilbert spaces},
  year    = {2022},
  pages   = {4017--4031},
  volume  = {35},
}

@Article{ikeda2022koopman,
  author    = {Ikeda, Masahiro and Ishikawa, Isao and Schlosser, Corbinian},
  journal   = {Chaos: An Interdisciplinary Journal of Nonlinear Science},
  title     = {Koopman and {P}erron--{F}robenius operators on reproducing kernel {B}anach spaces},
  year      = {2022},
  issn      = {1089-7682},
  month     = dec,
  number    = {12},
  volume    = {32},
  doi       = {10.1063/5.0094889},
  publisher = {AIP Publishing},
}

@InProceedings{mezic2015applications,
  author    = {Mezi{\'c}, Igor},
  booktitle = {2015 54th IEEE Conference on Decision and Control (CDC)},
  title     = {On applications of the spectral theory of the {K}oopman operator in dynamical systems and control theory},
  year      = {2015},
  month     = dec,
  pages     = {7034--7041},
  publisher = {IEEE},
  doi       = {10.1109/cdc.2015.7403328},
}

@Article{mauroy2016global,
  author    = {Mauroy, Alexandre and Mezi{\'c}, Igor},
  journal   = {IEEE Transactions on Automatic Control},
  title     = {Global stability analysis using the eigenfunctions of the {K}oopman operator},
  year      = {2016},
  issn      = {1558-2523},
  month     = nov,
  number    = {11},
  pages     = {3356--3369},
  volume    = {61},
  doi       = {10.1109/tac.2016.2518918},
  publisher = {Institute of Electrical and Electronics Engineers (IEEE)},
}

@Book{hill1986introduction,
  title={An introduction to statistical thermodynamics},
  author={Hill, Terrell L.},
  year={1986},
  publisher={Courier Corporation}
}

@Article{kryloff1937theorie,
  author    = {Kryloff, Nicolas and Bogoliouboff, Nicolas},
  journal   = {The Annals of Mathematics},
  title     = {La th{\'e}orie g{\'e}n{\'e}rale de la mesure dans son application {\`a} l'{\'e}tude des syst{\`e}mes dynamiques de la m{\'e}canique non lin{\'e}aire},
  year      = {1937},
  issn      = {0003-486X},
  month     = jan,
  number    = {1},
  pages     = {65--113},
  volume    = {38},
  doi       = {10.2307/1968511},
  publisher = {JSTOR},
}

@Book{conway2019course,
  author    = {Conway, John B.},
  publisher = {Springer New York},
  title     = {A course in functional analysis},
  year      = {2007},
  edition   = {2},
  isbn      = {9781475743838},
  volume    = {96},
  doi       = {10.1007/978-1-4757-4383-8},
  issn      = {0072-5285},
  journal   = {Graduate Texts in Mathematics},
}

@Book{kantz2004nonlinear,
  author    = {Kantz, Holger and Schreiber, Thomas},
  publisher = {Cambridge Univ. Press},
  title     = {Nonlinear time series analysis},
  year      = {2006},
  address   = {Cambridge},
  edition   = {Second Edition},
  isbn      = {9780521821506},
  series    = {Cambridge nonlinear science series},
  pagetotal = {369},
  ppn_gvk   = {639454070},
}

@Book{mauroy2020koopman,
  author    = {Mauroy, Alexandre and Susuki, Y. and Mezi{\'c}, I.},
  publisher = {Springer International Publishing AG},
  title     = {Koopman operator in systems and control},
  year      = {2020},
  address   = {Cham},
  isbn      = {9783030357139},
  number    = {v.484},
  series    = {Lecture Notes in Control and Information Sciences Ser.},
  pagetotal = {1568},
  ppn_gvk   = {1751910016},
  subtitle  = {Concepts, Methodologies, and Applications},
}

@Article{otto2021koopman,
  author    = {Otto, Samuel E. and Rowley, Clarence W.},
  journal   = {Annual Review of Control, Robotics, and Autonomous Systems},
  title     = {Koopman operators for estimation and control of dynamical systems},
  year      = {2021},
  issn      = {2573-5144},
  month     = may,
  number    = {1},
  pages     = {59--87},
  volume    = {4},
  doi       = {10.1146/annurev-control-071020-010108},
  publisher = {Annual Reviews},
}

@InProceedings{arbabi2018data,
  author       = {Arbabi, Hassan and Korda, Milan and Mezi{\'c}, Igor},
  booktitle    = {2018 IEEE Conference on Decision and Control (CDC)},
  title        = {A data-driven {K}oopman model predictive control framework for nonlinear partial differential equations},
  year         = {2018},
  month        = dec,
  organization = {IEEE},
  pages        = {6409--6414},
  publisher    = {IEEE},
  doi          = {10.1109/cdc.2018.8619720},
}

@InBook{peitz2020feedback,
  author    = {Peitz, Sebastian and Klus, Stefan},
  pages     = {257--282},
  publisher = {Springer},
  title     = {Feedback control of nonlinear {PDE}s using data-efficient reduced order models based on the {K}oopman operator},
  year      = {2020},
  isbn      = {9783030357139},
  booktitle = {The Koopman Operator in Systems and Control},
  doi       = {10.1007/978-3-030-35713-9_10},
  issn      = {1610-7411},
  journal   = {The Koopman Operator in Systems and Control: Concepts, Methodologies, and Applications},
}

@InProceedings{abraham2017model,
  author       = {Abraham, Ian and de la Torre, Gerardo and Murphey, Todd},
  booktitle    = {Robotics: Science and Systems XIII},
  title        = {Model-based control using {K}oopman operators},
  year         = {2017},
  month        = jul,
  organization = {MIT Press Journals},
  publisher    = {Robotics: Science and Systems Foundation},
  series       = {RSS2017},
  collection   = {RSS2017},
  doi          = {10.15607/rss.2017.xiii.052},
}

@Article{abraham2019active,
  author    = {Abraham, Ian and Murphey, Todd D.},
  journal   = {IEEE Transactions on Robotics},
  title     = {Active learning of dynamics for data-driven control using {K}oopman operators},
  year      = {2019},
  issn      = {1941-0468},
  month     = oct,
  number    = {5},
  pages     = {1071--1083},
  volume    = {35},
  doi       = {10.1109/tro.2019.2923880},
  publisher = {Institute of Electrical and Electronics Engineers (IEEE)},
}

@InProceedings{mamakoukas2019local,
  author     = {Mamakoukas, Giorgos and Castano, Maria and Tan, Xiaobo and Murphey, Todd},
  booktitle  = {Robotics: Science and Systems XV},
  title      = {Local {K}oopman operators for data-driven control of robotic systems},
  year       = {2019},
  month      = jun,
  publisher  = {Robotics: Science and Systems Foundation},
  series     = {RSS2019},
  collection = {RSS2019},
  doi        = {10.15607/rss.2019.xv.054},
}

@Article{haggerty2023control,
  author    = {Haggerty, David A. and Banks, Michael J. and Kamenar, Ervin and Cao, Alan B. and Curtis, Patrick C. and Mezi{\'c}, Igor and Hawkes, Elliot W.},
  journal   = {Science Robotics},
  title     = {Control of soft robots with inertial dynamics},
  year      = {2023},
  issn      = {2470-9476},
  month     = aug,
  number    = {81},
  pages     = {eadd6864},
  volume    = {8},
  doi       = {10.1126/scirobotics.add6864},
  publisher = {American Association for the Advancement of Science (AAAS)},
}

@Article{korda2018power,
  author    = {Korda, Milan and Susuki, Yoshihiko and Mezi{\'c}, Igor},
  journal   = {IFAC-PapersOnLine},
  title     = {Power grid transient stabilization using {K}oopman model predictive control},
  year      = {2018},
  issn      = {2405-8963},
  number    = {28},
  pages     = {297--302},
  volume    = {51},
  doi       = {10.1016/j.ifacol.2018.11.718},
  publisher = {Elsevier BV},
}

@Article{netto2018robust,
  author    = {Netto, Marcos and Mili, Lamine},
  journal   = {IEEE Transactions on Power Systems},
  title     = {A robust data-driven {K}oopman {K}alman filter for power systems dynamic state estimation},
  year      = {2018},
  issn      = {1558-0679},
  month     = nov,
  number    = {6},
  pages     = {7228--7237},
  volume    = {33},
  doi       = {10.1109/tpwrs.2018.2846744},
  publisher = {Institute of Electrical and Electronics Engineers (IEEE)},
}

@InProceedings{hasnain2020steady,
  author       = {Hasnain, Aqib and Boddupalli, Nibodh and Balakrishnan, Shara and Yeung, Enoch},
  booktitle    = {2020 American Control Conference (ACC)},
  title        = {Steady state programming of controlled nonlinear systems via deep dynamic mode decomposition},
  year         = {2020},
  month        = jul,
  organization = {IEEE},
  pages        = {4245--4251},
  publisher    = {IEEE},
  doi          = {10.23919/acc45564.2020.9147218},
}

@Article{narasingam2019koopman,
  author    = {Narasingam, Abhinav and Kwon, Joseph Sang-Il},
  journal   = {AIChE Journal},
  title     = {Koopman {L}yapunov-based model predictive control of nonlinear chemical process systems},
  year      = {2019},
  issn      = {1547-5905},
  month     = aug,
  number    = {11},
  pages     = {e16743},
  volume    = {65},
  doi       = {10.1002/aic.16743},
  publisher = {Wiley},
}

@Article{korda2018linear,
  author    = {Korda, Milan and Mezi{\'c}, Igor},
  journal   = {Automatica},
  title     = {Linear predictors for nonlinear dynamical systems: {K}oopman operator meets model predictive control},
  year      = {2018},
  issn      = {0005-1098},
  month     = jul,
  pages     = {149--160},
  volume    = {93},
  doi       = {10.1016/j.automatica.2018.03.046},
  publisher = {Elsevier BV},
}

@Article{proctor2018generalizing,
  author    = {Proctor, Joshua L. and Brunton, Steven L. and Kutz, J. Nathan},
  journal   = {SIAM Journal on Applied Dynamical Systems},
  title     = {Generalizing {K}oopman theory to allow for inputs and control},
  year      = {2018},
  issn      = {1536-0040},
  month     = jan,
  number    = {1},
  pages     = {909--930},
  volume    = {17},
  doi       = {10.1137/16m1062296},
  publisher = {Society for Industrial & Applied Mathematics (SIAM)},
}

@InProceedings{surana2016koopman,
  author       = {Surana, Amit},
  booktitle    = {2016 IEEE 55th Conference on Decision and Control (CDC)},
  title        = {Koopman operator based observer synthesis for control-affine nonlinear systems},
  year         = {2016},
  month        = dec,
  organization = {IEEE},
  pages        = {6492--6499},
  publisher    = {IEEE},
  doi          = {10.1109/cdc.2016.7799268},
}

@InProceedings{kaiser2018discovering,
  author       = {Kaiser, Eurika and Kutz, J. Nathan and Brunton, Steven L.},
  booktitle    = {2018 IEEE Conference on Decision and Control (CDC)},
  title        = {Discovering conservation laws from data for control},
  year         = {2018},
  month        = dec,
  organization = {IEEE},
  pages        = {6415--6421},
  publisher    = {IEEE},
  doi          = {10.1109/cdc.2018.8618963},
}

@Article{peitz2019koopman,
  author    = {Peitz, Sebastian and Klus, Stefan},
  journal   = {Automatica},
  title     = {Koopman operator-based model reduction for switched-system control of {PDE}s},
  year      = {2019},
  issn      = {0005-1098},
  month     = aug,
  pages     = {184--191},
  volume    = {106},
  doi       = {10.1016/j.automatica.2019.05.016},
  publisher = {Elsevier BV},
}

@InProceedings{sootla2016properties,
  author       = {Sootla, Aivar and Mauroy, Alexandre},
  booktitle    = {2016 American Control Conference (ACC)},
  title        = {Properties of isostables and basins of attraction of monotone systems},
  year         = {2016},
  month        = jul,
  organization = {IEEE},
  pages        = {7365--7370},
  publisher    = {IEEE},
  doi          = {10.1109/acc.2016.7526835},
}

@InProceedings{vaidya2007observability,
  author       = {Vaidya, Umesh},
  booktitle    = {2007 46th IEEE Conference on Decision and Control},
  title        = {Observability {G}ramian for nonlinear systems},
  year         = {2007},
  organization = {IEEE},
  pages        = {3357--3362},
  publisher    = {IEEE},
  doi          = {10.1109/cdc.2007.4434828},
}

@InProceedings{goswami2017global,
  author       = {Goswami, Debdipta and Paley, Derek A.},
  booktitle    = {2017 IEEE 56th Annual Conference on Decision and Control (CDC)},
  title        = {Global bilinearization and controllability of control-affine nonlinear systems: {A} {K}oopman spectral approach},
  year         = {2017},
  month        = dec,
  organization = {IEEE},
  pages        = {6107--6112},
  publisher    = {IEEE},
  doi          = {10.1109/cdc.2017.8264582},
}

@InProceedings{yeung2018koopman,
  author       = {Yeung, Enoch and Liu, Zhiyuan and Hodas, Nathan O.},
  booktitle    = {2018 Annual American Control Conference (ACC)},
  title        = {A {K}oopman operator approach for computing and balancing {G}ramians for discrete time nonlinear systems},
  year         = {2018},
  month        = jun,
  organization = {IEEE},
  pages        = {337--344},
  publisher    = {IEEE},
  doi          = {10.23919/acc.2018.8431738},
}

@Article{sinha2016operator,
  author    = {Sinha, S. and Vaidya, U. and Rajaram, R.},
  journal   = {Journal of Mathematical Analysis and Applications},
  title     = {Operator theoretic framework for optimal placement of sensors and actuators for control of nonequilibrium dynamics},
  year      = {2016},
  issn      = {0022-247X},
  month     = aug,
  number    = {2},
  pages     = {750--772},
  volume    = {440},
  doi       = {10.1016/j.jmaa.2016.03.058},
  publisher = {Elsevier BV},
}

@Article{sharma2019transfer,
  author    = {Sharma, Himanshu and Vaidya, Umesh and Ganapathysubramanian, Baskar},
  journal   = {Building and Environment},
  title     = {A transfer operator methodology for optimal sensor placement accounting for uncertainty},
  year      = {2019},
  issn      = {0360-1323},
  month     = may,
  pages     = {334--349},
  volume    = {155},
  doi       = {10.1016/j.buildenv.2019.03.054},
  publisher = {Elsevier BV},
}

@Article{mamakoukas2021derivative,
  author    = {Mamakoukas, Giorgos and Castano, Maria L. and Tan, Xiaobo and Murphey, Todd D.},
  journal   = {IEEE Transactions on Robotics},
  title     = {Derivative-based {K}oopman operators for real-time control of robotic systems},
  year      = {2021},
  issn      = {1941-0468},
  month     = dec,
  number    = {6},
  pages     = {2173--2192},
  volume    = {37},
  doi       = {10.1109/tro.2021.3076581},
  publisher = {Institute of Electrical and Electronics Engineers (IEEE)},
}

@Article{kaiser2018sparse,
  author    = {Kaiser, E. and Kutz, J. N. and Brunton, S. L.},
  journal   = {Proceedings of the Royal Society A: Mathematical, Physical and Engineering Sciences},
  title     = {Sparse identification of nonlinear dynamics for model predictive control in the low-data limit},
  year      = {2018},
  issn      = {1471-2946},
  month     = nov,
  number    = {2219},
  pages     = {20180335},
  volume    = {474},
  doi       = {10.1098/rspa.2018.0335},
  publisher = {The Royal Society},
}

@InProceedings{huang2018feedback,
  author       = {Huang, Bowen and Ma, Xu and Vaidya, Umesh},
  booktitle    = {2018 IEEE Conference on Decision and Control (CDC)},
  title        = {Feedback stabilization using {K}oopman operator},
  year         = {2018},
  month        = dec,
  organization = {IEEE},
  pages        = {6434--6439},
  publisher    = {IEEE},
  doi          = {10.1109/cdc.2018.8619727},
}

@InBook{huang2020data,
  author    = {Huang, Bowen and Ma, Xu and Vaidya, Umesh},
  pages     = {313--334},
  publisher = {Springer International Publishing},
  title     = {Data-driven nonlinear stabilization using {K}oopman operator},
  year      = {2020},
  isbn      = {9783030357139},
  booktitle = {The Koopman Operator in Systems and Control},
  doi       = {10.1007/978-3-030-35713-9_12},
  issn      = {1610-7411},
  journal   = {The Koopman Operator in Systems and Control: Concepts, Methodologies, and Applications},
}

@InProceedings{hemati2017dynamic,
  author    = {Hemati, Maziar and Yao, Huaijin},
  booktitle = {8th AIAA Theoretical Fluid Mechanics Conference},
  title     = {Dynamic mode shaping for fluid flow control: {N}ew strategies for transient growth suppression},
  year      = {2017},
  month     = jun,
  pages     = {3160},
  publisher = {American Institute of Aeronautics and Astronautics},
  doi       = {10.2514/6.2017-3160},
}

@Article{aloisio2022spectral,
  author        = {Aloisio, M. and Carvalho, S. L. and de Oliveira, C. R. and Souza, E.},
  journal       = {arXiv preprint arXiv:2209.05290},
  title         = {On spectral measures and convergence rates in von {N}eumann's ergodic theorem},
  year          = {2022},
  month         = sep,
  archiveprefix = {arXiv},
  copyright     = {Creative Commons Attribution 4.0 International},
  doi           = {10.48550/ARXIV.2209.05290},
  eprint        = {2209.05290},
  primaryclass  = {math.SP},
  publisher     = {arXiv},
}

@Article{young2002srb,
  author    = {Young, Lai-Sang},
  journal   = {Journal of Statistical Physics},
  title     = {What are {SRB} measures, and which dynamical systems have them?},
  year      = {2002},
  issn      = {0022-4715},
  number    = {5/6},
  pages     = {733--754},
  volume    = {108},
  doi       = {10.1023/a:1019762724717},
  publisher = {Springer Science and Business Media LLC},
}

@Article{berry2015nonparametric,
  author    = {Berry, Tyrus and Giannakis, Dimitrios and Harlim, John},
  journal   = {Physical Review E},
  title     = {Nonparametric forecasting of low-dimensional dynamical systems},
  year      = {2015},
  issn      = {1550-2376},
  month     = mar,
  number    = {3},
  pages     = {032915},
  volume    = {91},
  doi       = {10.1103/physreve.91.032915},
  publisher = {American Physical Society (APS)},
}

@Article{garcia2020error,
  author    = {Garc{\'\i}a Trillos, Nicol{\'a}s and Gerlach, Moritz and Hein, Matthias and Slep{\v{c}}ev, Dejan},
  journal   = {Foundations of Computational Mathematics},
  title     = {Error estimates for spectral convergence of the graph {L}aplacian on random geometric graphs toward the {L}aplace--{B}eltrami operator},
  year      = {2020},
  issn      = {1615-3383},
  month     = sep,
  number    = {4},
  pages     = {827--887},
  volume    = {20},
  doi       = {10.1007/s10208-019-09436-w},
  publisher = {Springer Science and Business Media LLC},
}

@Article{von2008consistency,
  author    = {von Luxburg, Ulrike and Belkin, Mikhail and Bousquet, Olivier},
  journal   = {The Annals of Statistics},
  title     = {Consistency of spectral clustering},
  year      = {2008},
  issn      = {0090-5364},
  month     = apr,
  number    = {2},
  pages     = {555--586},
  volume    = {36},
  doi       = {10.1214/009053607000000640},
  publisher = {Institute of Mathematical Statistics},
}

@InProceedings{giannakis2015spatiotemporal,
  author       = {Giannakis, Dimitrios and Slawinska, Joanna and Zhao, Zhizhen},
  booktitle    = {Feature Extraction: Modern Questions and Challenges},
  title        = {Spatiotemporal feature extraction with data-driven {K}oopman operators},
  year         = {2015},
  organization = {PMLR},
  pages        = {103--115},
}

@Article{goswami2018constrained,
  author    = {Goswami, Debdipta and Thackray, Emma and Paley, Derek A.},
  journal   = {IEEE Control Systems Letters},
  title     = {Constrained {U}lam dynamic mode decomposition: {A}pproximation of the {P}erron--{F}robenius operator for deterministic and stochastic systems},
  year      = {2018},
  issn      = {2475-1456},
  month     = oct,
  number    = {4},
  pages     = {809--814},
  volume    = {2},
  doi       = {10.1109/lcsys.2018.2849552},
  publisher = {IEEE},
}

@article{cohen2020mode,
  title={Mode decomposition for homogeneous symmetric operators},
  author={Cohen, Ido and Azencot, Omri and Lifshits, Pavel and Gilboa, Guy},
  journal={Preprint arXiv},
  year={2020}
}

@Article{krake2022constrained,
  author    = {Krake, Tim and Kl{\"o}tzl, Daniel and Eberhardt, Bernhard and Weiskopf, Daniel},
  journal   = {IEEE Transactions on Visualization and Computer Graphics},
  title     = {Constrained dynamic mode decomposition},
  year      = {2022},
  issn      = {2160-9306},
  number    = {1},
  pages     = {1--11},
  volume    = {29},
  doi       = {10.1109/tvcg.2022.3209437},
  publisher = {Institute of Electrical and Electronics Engineers (IEEE)},
}

@InProceedings{huang2018data,
  author       = {Huang, Bowen and Vaidya, Umesh},
  booktitle    = {2018 Annual American Control Conference (ACC)},
  title        = {Data-driven approximation of transfer operators: {N}aturally structured dynamic mode decomposition},
  year         = {2018},
  month        = jun,
  organization = {IEEE},
  pages        = {5659--5664},
  publisher    = {IEEE},
  doi          = {10.23919/acc.2018.8431409},
}

@InProceedings{mardt2020deep,
  author       = {Mardt, Andreas and Pasquali, Luca and No{\'e}, Frank and Wu, Hao},
  booktitle    = {Mathematical and Scientific Machine Learning},
  title        = {Deep learning {M}arkov and {K}oopman models with physical constraints},
  year         = {2020},
  organization = {PMLR},
  pages        = {451--475},
}

@Article{salova2019koopman,
  author    = {Salova, Anastasiya and Emenheiser, Jeffrey and Rupe, Adam and Crutchfield, James P. and D'Souza, Raissa M.},
  journal   = {Chaos: An Interdisciplinary Journal of Nonlinear Science},
  title     = {Koopman operator and its approximations for systems with symmetries},
  year      = {2019},
  issn      = {1089-7682},
  month     = sep,
  number    = {9},
  volume    = {29},
  doi       = {10.1063/1.5099091},
  publisher = {AIP Publishing},
}

@Article{mehta2006symmetry,
  author    = {Mehta, Prashant G. and Hessel-von Molo, Mirko and Dellnitz, Michael},
  journal   = {Journal of Difference Equations and Applications},
  title     = {Symmetry of attractors and the {P}erron--{F}robenius operator},
  year      = {2006},
  issn      = {1563-5120},
  month     = nov,
  number    = {11},
  pages     = {1147--1178},
  volume    = {12},
  doi       = {10.1080/10236190601045788},
  publisher = {Informa UK Limited},
}

@Article{lu2020lagrangian,
  author    = {Lu, Hannah and Tartakovsky, Daniel M.},
  journal   = {Journal of Computational Physics},
  title     = {Lagrangian dynamic mode decomposition for construction of reduced-order models of advection-dominated phenomena},
  year      = {2020},
  issn      = {0021-9991},
  month     = apr,
  pages     = {109229},
  volume    = {407},
  doi       = {10.1016/j.jcp.2020.109229},
  publisher = {Elsevier BV},
}

@Article{mojgani2017lagrangian,
  author        = {Mojgani, Rambod and Balajewicz, Maciej},
  journal       = {arXiv preprint arXiv:1701.04343},
  title         = {Lagrangian basis method for dimensionality reduction of convection dominated nonlinear flows},
  year          = {2017},
  month         = jan,
  archiveprefix = {arXiv},
  copyright     = {arXiv.org perpetual, non-exclusive license},
  doi           = {10.48550/ARXIV.1701.04343},
  eprint        = {1701.04343},
  primaryclass  = {physics.flu-dyn},
  publisher     = {arXiv},
}

@Article{pan2020physics,
  author    = {Pan, Shaowu and Duraisamy, Karthik},
  journal   = {SIAM Journal on Applied Dynamical Systems},
  title     = {Physics-informed probabilistic learning of linear embeddings of nonlinear dynamics with guaranteed stability},
  year      = {2020},
  issn      = {1536-0040},
  month     = jan,
  number    = {1},
  pages     = {480--509},
  volume    = {19},
  doi       = {10.1137/19m1267246},
  publisher = {Society for Industrial & Applied Mathematics (SIAM)},
}

@Article{morandin2023port,
  author    = {Morandin, Riccardo and Nicodemus, Jonas and Unger, Benjamin},
  journal   = {SIAM Journal on Scientific Computing},
  title     = {Port-{H}amiltonian dynamic mode decomposition},
  year      = {2023},
  issn      = {1095-7197},
  month     = jul,
  number    = {4},
  pages     = {A1690--A1710},
  volume    = {45},
  doi       = {10.1137/22m149329x},
  publisher = {Society for Industrial & Applied Mathematics (SIAM)},
}

@Article{baddoo2023physics,
  author    = {Baddoo, Peter J. and Herrmann, Benjamin and McKeon, Beverley J. and Nathan Kutz, J. and Brunton, Steven L.},
  journal   = {Proceedings of the Royal Society A: Mathematical, Physical and Engineering Sciences},
  title     = {Physics-informed dynamic mode decomposition},
  year      = {2023},
  issn      = {1471-2946},
  month     = mar,
  number    = {2271},
  pages     = {20220576},
  volume    = {479},
  doi       = {10.1098/rspa.2022.0576},
  publisher = {The Royal Society},
}

@Article{celledoni2021structure,
  author    = {Celledoni, E. and Ehrhardt, M. J. and Etmann, C. and McLachlan, R. I. and Owren, B. and Sch{\"o}nlieb, C.-B. and Sherry, F.},
  journal   = {European Journal of Applied Mathematics},
  title     = {Structure-preserving deep learning},
  year      = {2021},
  issn      = {1469-4425},
  month     = may,
  number    = {5},
  pages     = {888--936},
  volume    = {32},
  doi       = {10.1017/s0956792521000139},
  publisher = {Cambridge University Press (CUP)},
}

@Article{karniadakis2021physics,
  author    = {Karniadakis, George Em and Kevrekidis, Ioannis G. and Lu, Lu and Perdikaris, Paris and Wang, Sifan and Yang, Liu},
  journal   = {Nature Reviews Physics},
  title     = {Physics-informed machine learning},
  year      = {2021},
  issn      = {2522-5820},
  month     = may,
  number    = {6},
  pages     = {422--440},
  volume    = {3},
  doi       = {10.1038/s42254-021-00314-5},
  publisher = {Springer Science and Business Media LLC},
}

@article{greydanus2019hamiltonian,
  title={Hamiltonian neural networks},
  author={Greydanus, S. and Dzamba, M. and Yosinski, J.},
  journal={Advances in neural information processing systems},
  volume={32},
  year={2019}
}

@Article{hernandez2021structure,
  author    = {Hern{\'a}ndez, Q. and Bad{\'\i}as, A. and Gonz{\'a}lez, D. and Chinesta, F. and Cueto, E.},
  journal   = {Journal of Computational Physics},
  title     = {Structure-preserving neural networks},
  year      = {2021},
  issn      = {0021-9991},
  month     = feb,
  pages     = {109950},
  volume    = {426},
  doi       = {10.1016/j.jcp.2020.109950},
  publisher = {Elsevier BV},
}

@Book{3-540-30663-3,
  author    = {Hairer, E. and Lubich, C. and Wanner, G.},
  publisher = {Springer},
  title     = {Geometric numerical integration},
  year      = {2010},
  address   = {Berlin},
  edition   = {Second edition, first softcover print},
  isbn      = {9783642051579},
  number    = {31},
  series    = {Springer series in computational mathematics},
  pagetotal = {644},
  ppn_gvk   = {613185870},
  subtitle  = {Structure-preserving algorithms for ordinary differential equations},
}

@Article{govindarajan2021approximation,
  author    = {Govindarajan, Nithin and Mohr, Ryan and Chandrasekaran, Shivkumar and Mezic, Igor},
  journal   = {SIAM Journal on Applied Dynamical Systems},
  title     = {On the approximation of {K}oopman spectra of measure-preserving flows},
  year      = {2021},
  issn      = {1536-0040},
  month     = jan,
  number    = {1},
  pages     = {232--261},
  volume    = {20},
  doi       = {10.1137/19m1282908},
  publisher = {Society for Industrial & Applied Mathematics (SIAM)},
}

@Article{halmos1944approximation,
  author    = {Halmos, Paul R.},
  journal   = {Transactions of the American Mathematical Society},
  title     = {Approximation theories for measure preserving transformations},
  year      = {1944},
  issn      = {0002-9947},
  month     = jan,
  number    = {1},
  pages     = {1--18},
  volume    = {55},
  doi       = {10.2307/1990137},
  publisher = {JSTOR},
}

@Article{lax1971approximation,
  author    = {Lax, Peter D.},
  journal   = {Communications on Pure and Applied Mathematics},
  title     = {Approximation of measure preserving transformations},
  year      = {1971},
  issn      = {1097-0312},
  month     = mar,
  number    = {2},
  pages     = {133--135},
  volume    = {24},
  doi       = {10.1002/cpa.3160240204},
  publisher = {Wiley},
}

@Article{schonemann1966generalized,
  author    = {Sch{\"o}nemann, Peter H.},
  journal   = {Psychometrika},
  title     = {A generalized solution of the orthogonal procrustes problem},
  year      = {1966},
  issn      = {1860-0980},
  month     = mar,
  number    = {1},
  pages     = {1--10},
  volume    = {31},
  doi       = {10.1007/bf02289451},
  publisher = {Springer Science and Business Media LLC},
}

@Article{higham1988symmetric,
  author    = {Higham, Nicholas J.},
  journal   = {BIT},
  title     = {The symmetric {P}rocrustes problem},
  year      = {1988},
  issn      = {1572-9125},
  month     = mar,
  number    = {1},
  pages     = {133--143},
  volume    = {28},
  doi       = {10.1007/bf01934701},
  publisher = {Springer Science and Business Media LLC},
}

@Article{boumal2014manopt,
  author    = {Boumal, Nicolas and Mishra, Bamdev and Absil, P.-A. and Sepulchre, Rodolphe},
  journal   = {The Journal of Machine Learning Research},
  title     = {Manopt, a {MATLAB} toolbox for optimization on manifolds},
  year      = {2014},
  number    = {1},
  pages     = {1455--1459},
  volume    = {15},
  publisher = {JMLR. org},
}

@Book{gower2004procrustes,
  author    = {Gower, John C. and Dijksterhuis, Garmt B.},
  publisher = {Oxford University Press},
  title     = {Procrustes problems},
  year      = {2004},
  isbn      = {9780198510581},
  month     = jan,
  volume    = {30},
  doi       = {10.1093/acprof:oso/9780198510581.001.0001},
}

@Article{rossler1976equation,
  author    = {R{\"o}ssler, Otto E.},
  journal   = {Physics Letters A},
  title     = {An equation for continuous chaos},
  year      = {1976},
  issn      = {0375-9601},
  month     = jul,
  number    = {5},
  pages     = {397--398},
  volume    = {57},
  doi       = {10.1016/0375-9601(76)90101-8},
  publisher = {Elsevier BV},
}

@Article{peifer2005mixing,
  author    = {Peifer, M. and Schelter, B. and Winterhalder, M. and Timmer, J.},
  journal   = {Physical Review E},
  title     = {Mixing properties of the {R}{\"o}ssler system and consequences for coherence and synchronization analysis},
  year      = {2005},
  issn      = {1550-2376},
  month     = aug,
  number    = {2},
  pages     = {026213},
  volume    = {72},
  doi       = {10.1103/physreve.72.026213},
  publisher = {American Physical Society (APS)},
}

@Article{valva2023consistent,
  author        = {Valva, Claire and Giannakis, Dimitrios},
  journal       = {arXiv preprint arXiv:2309.00732},
  title         = {Consistent spectral approximation of {K}oopman operators using resolvent compactification},
  year          = {2023},
  month         = sep,
  archiveprefix = {arXiv},
  copyright     = {arXiv.org perpetual, non-exclusive license},
  doi           = {10.48550/ARXIV.2309.00732},
  eprint        = {2309.00732},
  primaryclass  = {math.DS},
  publisher     = {arXiv},
}
 \end{spacing}
\end{document}